%% file: ntrivodd.tex
\documentstyle[txmac,a4,%leqno,%
amssymb,%
bibenc,%
case,%
twoside,%
nocaphead2,%
epsf,%
alltt,%
myrot,%
mypic,times,mathptm]{article}

\advance\oddsidemargin by -1.9cm
\advance\evensidemargin by -1.9cm
\advance\textwidth by 3.8cm
\def\mynewtheo#1#2{%
\newtheorem{@#1}{#2}[section]%
\newenvironment{#1}{\begin{@#1}\rm}{\end{@#1}}}

\mynewtheo{lemma}{Lemma}
\mynewtheo{sublemma}{Sublemma}
\mynewtheo{exer}{Exercise}
\mynewtheo{theo}{Theorem}
\mynewtheo{rem}{Remark}
\mynewtheo{clc}{Calculation}
\mynewtheo{defi}{Definition}
\mynewtheo{conj}{Conjecture}
\mynewtheo{corr}{Corollary}
\mynewtheo{prop}{Proposition}
\mynewtheo{question}{Question}
\mynewtheo{exam}{Example}
\mynewtheo{conv}{Convention}

\makeatletter

\newenvironment{theorem}{\begin{theo}}{\end{theo}}

\begin{document}

\input{myeqn.tex}

\newcommand{\mybin}[2]{\text{$\Bigl(\begin{array}{@{}c@{}}#1\\#2%
\end{array}\Bigr)$}}
\newcommand{\mybinn}[2]{\text{$\biggl(\begin{array}{@{}c@{}}%
#1\\#2\end{array}\biggr)$}}

\def\overtwo#1{\mbox{\small$\mybin{#1}{2}$}}
\newcommand{\mybr}[2]{\text{$\Bigl\lfloor\mbox{%
\small$\displaystyle\frac{#1}{#2}$}\Bigr\rfloor$}}
\def\mybrtwo#1{\mbox{\mybr{#1}{2}}}

\def\myfrac#1#2{\raisebox{0.2em}{\small$#1$}\!/\!\raisebox{-0.2em}{\small$#2$}}
\def\ffrac#1#2{\mbox{\small$\ds\frac{#1}{#2}$}}

\def\Pos#1#2{{\diag{#1}{1}{1}{#2
\picmultiline{-8 1 -1 0}{0 1}{1 0}
\picmultiline{-8 1 -1 0}{0 0}{1 1}
}}}
\def\Neg#1#2{{\diag{#1}{1}{1}{#2
\picmultiline{-8 1 -1 0}{0 0}{1 1}
\picmultiline{-8 1 -1 0}{0 1}{1 0}
}}}
\def\Nul#1#2{{\diag{#1}{1}{1}{#2
\piccirclearc{0.5 1.4}{0.7}{-135 -45}
\piccirclearc{0.5 -0.4}{0.7}{45 135}
}}}
\def\Inf#1#2{{\diag{#1}{1}{1}{#2
\piccirclearc{0.5 1.4 x}{0.7}{135 -135}
\piccirclearc{0.5 -0.4 x}{0.7}{-45 45}
}}}
\def\pos{\Pos{0.5em}{\piclinewidth{10}}}
\def\neg{\Neg{0.5em}{\piclinewidth{10}}}
\def\nul{\Nul{0.5em}{\piclinewidth{10}}}
\def\inf{\Inf{0.5em}{\piclinewidth{10}}}

\def\noloop{{\diag{0.5cm}{0.5}{1}{\picline{0.25 0}{0.25 1}}}}

\def\llpoint#1{{\picfillgraycol{0}\picfilledcircle{#1}{0.15}{}}}
\def\lpoint#1{{\picfillgraycol{0}\picfilledcircle{#1}{0.08}{}}}
\def\point#1{{\picfillgraycol{0}\picfilledcircle{#1}{0.04}{}}}

\def\ReidI#1#2{
  \diag{0.5cm}{0.9}{1}{
    \pictranslate{0.4 0.5}{
      \picscale{#11 #21}{
        \picmultigraphics[S]{2}{1 -1}{
	  \picmulticurve{-6 1 -1 0}{0.5 -0.5}{0.5 0}{0.1 0.3}{-0.2 0.3}
	} %picmultigraphics
	\piccirclearc{-0.2 0}{0.3}{90 270}
      }%picscale
    }%pictranslate
  }%diag
}

\def\pt#1{{\picfillgraycol{0}\picfilledcircle{#1}{0.1}{}}}
\def\ppt#1{{\picfillgraycol{0}\picfilledcircle{#1}{0.06}{}}}

\def\cycl#1#2#3#4{\vrt{#1}\vrt{#2}\vrt{#3}\vrt{#4}%
\picline{#1}{#2}\picline{#2}{#3}\picline{#3}{#4}\picline{#4}{#1}%
}
\def\vrt#1{{\picfillgraycol{0}\picfilledcircle{#1}{0.09}{}}}

%usage ODD # of coords, then {}
\def\@curvepath#1#2#3{%
  \@ifempty{#2}{\piccurveto{#1 }{@stc}{@std}#3}%
    {\piccurveto{#1 }{#2 }{#2 #3 0.5 conv}
    \@curvepath{#3}}%
}
\def\curvepath#1#2#3{%
  \piccurve{#1 }{#2 }{#2 }{#2 #3 0.5 conv}%
  \picPSgraphics{/@stc [ #1 #2 -1 conv ] $ D /@std [ #1 ] $ D }%
  \@curvepath{#3}%
}

%usage EVEN # of coords, then {}
\def\@opencurvepath#1#2#3{%
  \@ifempty{#3}{\piccurveto{#1 }{#1 }{#2 }}%
    {\piccurveto{#1 }{#2 }{#2 #3 0.5 conv}\@opencurvepath{#3}}%
}
\def\opencurvepath#1#2#3{%
  \piccurve{#1 }{#2 }{#2 }{#2 #3 0.5 conv}%
  \@opencurvepath{#3}%
}
\def\epsfs#1#2{{\catcode`\_=11\relax\ifautoepsf\unitxsize#1\relax\else
\epsfxsize#1\relax\fi\epsffile{#2.eps}}}
\def\epsfsv#1#2{{\vcbox{\epsfs{#1}{#2}}}}
\def\vcbox#1{\setbox\@tempboxa=\hbox{#1}\parbox{\wd\@tempboxa}{\box
  \@tempboxa}}
\def\p{\epsfsv{2cm}}

\def\@test#1#2#3#4{%
  \let\@tempa\go@
  \@tempdima#1\relax\@tempdimb#3\@tempdima\relax\@tempdima#4\unitxsize\relax
  \ifdim \@tempdimb>\z@\relax
    \ifdim \@tempdimb<#2%
      \def\@tempa{\@test{#1}{#2}}%
    \fi
  \fi
  \@tempa
}

\def\go@#1\@end{}
\newdimen\unitxsize
\newif\ifautoepsf\autoepsftrue

\unitxsize4cm\relax
\def\epsfsize#1#2{\epsfxsize\relax\ifautoepsf
  {\@test{#1}{#2}{0.1 }{4   }
		{0.2 }{3   }
		{0.3 }{2   }
		{0.4 }{1.7 }
		{0.5 }{1.5 }
		{0.6 }{1.4 }
		{0.7 }{1.3 }
		{0.8 }{1.2 }
		{0.9 }{1.1 }
		{1.1 }{1.  }
		{1.2 }{0.9 }
		{1.4 }{0.8 }
		{1.6 }{0.75}
		{2.  }{0.7 }
		{2.25}{0.6 }
		{3   }{0.55}
		{5   }{0.5 }
		{10  }{0.33}
		{-1  }{0.25}\@end
		\ea}\ea\epsfxsize\the\@tempdima\relax
		\fi
		}

\def\rr#1{\unitxsize4.0cm\relax\epsffile{#1.eps}}
\def\vis#1#2{\hbox{\begin{tabular}{c}\rr{#2} \\ \ry{1.4em}$#1$\end{tabular} }}
\def\vis@#1#2{\hbox{\begin{tabular}{c}\rr{k-#1-#2} \end{tabular} }}

\author{A. Stoimenow\footnotemark[1]\\[2mm]
\small Research Institute for Mathematical Sciences, \\
\small Kyoto University, Kyoto 606-8502, Japan,\\
\small e-mail: {\tt stoimeno@kurims.kyoto-u.ac.jp}\\
\small WWW: {\hbox{\web|http://www.kurims.kyoto-u.ac.jp/~stoimeno|}}
}

{\def\thefootnote{\fnsymbol{footnote}}
\footnotetext[1]{Supported by 21st Century COE Program.}
}

\title{\large\bf \uppercase{Non-triviality of the Jones polynomial}
\\[1mm]
\uppercase{and the crossing numbers of amphicheiral knots}\\[4mm]
\phantom{\small\it This is
a preprint. I would be grateful for any comments and corrections.}}

\date{\phantom{\large Current version: \curv\ \ \ First version:
\makedate{1}{8}{2004}}}

\maketitle

\makeatletter

\let\vn\varnothing
\def\lpoint#1{{\picfillgraycol{0}\picfilledcircle{#1}{0.08}{}}}
\def\point#1{{\picfillgraycol{0}\picfilledcircle{#1}{0.04}{}}}
\let\ay\asymp
\let\pa\partial
\let\ap\alpha
\let\bt\beta
\let\Dl\Delta
\let\Gm\Gamma
\let\gm\gamma
\let\de\delta
\let\dl\delta
\let\eps\epsilon
\let\lm\lambda
\let\Lm\Lambda
\let\sg\sigma
\let\vp\varphi
\let\zt\zeta
\let\om\omega
\let\diagram\diag
\let\nb\nabla

\let\sm\setminus
\let\tl\tilde
\def\tD{\tl D}
\def\sgn{\mathop {\operator@font sgn}}
\def\spn{\mathop {\operator@font span}}
\def\Mc{\max\cf}
\def\Md{\max\deg}
\def\md{\min\deg}
\def\mc{\max\cf}
\def\mnc{\min\cf}
\def\vol{\text{\rm vol}\,}
\def\Ra{\Rightarrow}
\def\Lra{\Longrightarrow}
\def\lra{\longrightarrow}
\def\so{\Rightarrow}
\def\So{\Longrightarrow}
\def\nin{\not\in}
\let\tg\bigtriangleup 
\let\ds\displaystyle
\let\llra\longleftrightarrow
\let\reference\ref
\let\ul\underline
\let\es\enspace
\def\lfra{\leftrightarrow}
\def\cf{\text{\rm cf}\,}

\def\lb{\linebreak[0]}
\def\lz{\lb\verb}
\def\ssim{\stackrel{\ds \sim}{\vbox{\vskip-0.2em\hbox{$\scriptstyle
*$}}}}

\long\def\@makecaption#1#2{%
   % \tm
   \vskip \abovecaptionskip 
   {\let\label\@gobble
   \let\ignorespaces\@empty
   \xdef\@tempt{#2}%
   }%
   \ea\@ifempty\ea{\@tempt}{%
   \sbox\@tempboxa{%
      \fignr#1#2}%
      }{%
   \sbox\@tempboxa{%
      {\fignr#1:}\capt\ #2}%
      }%
   \ifdim \wd\@tempboxa >\captionwidth {%
      \centerline{\parbox{\captionwidth}{\unhbox \@tempboxa}}%
      %\rightskip=\@captionmargin\leftskip=\@captionmargin
      %\unhbox\@tempboxa\par
     }%
   \else
      \centerline{\box \@tempboxa}%
      % \hbox to\captionwidth{\hfil\box\@tempboxa\hfil}%
   \fi
   \vskip \belowcaptionskip
   }%
\def\fignr{\small\sffamily\bfseries}%
\def\capt{\small\sffamily}%

% \long\def\@makecaption#1#2{%
%    % \tm
%    \vskip 10pt
%    {\let\label\@gobble
%    \let\ignorespaces\@empty
%    \xdef\@tempt{#2}%
%    %\typeout{`#2'}%
%    }%
%    \ea\@ifempty\ea{\@tempt}{%
%    \setbox\@tempboxa\hbox{%
%       \fignr#1#2}%
%       }{%
%    \setbox\@tempboxa\hbox{%
%       {\fignr#1:}\capt\ #2}%
%       }%
%    \ifdim \wd\@tempboxa >\captionwidth {%
%       \rightskip=\@captionmargin\leftskip=\@captionmargin
%       \unhbox\@tempboxa\par}%
%    \else
%       \hbox to\captionwidth{\hfil\box\@tempboxa\hfil}%
%    \fi}%
% %
% \def\fignr{\small\sffamily\bfseries}%
% \def\capt{\small\sffamily}%

\newdimen\@captionmargin\@captionmargin2cm\relax
\newdimen\captionwidth\captionwidth0.8\hsize\relax

\def\eqref#1{(\protect\ref{#1})}

\def\proof{\@ifnextchar[{\@proof}{\@proof[\unskip]}}
\def\@proof[#1]{\noindent{\bf Proof #1.}\enspace}

\def\hint{\noindent Hint: }
\def\problem{\noindent{\bf Problem.} }

\def\@mt#1{\ifmmode#1\else$#1$\fi}
\def\qed{\hfill\@mt{\Box}}
\def\qqed{\hfill\@mt{\Box\enspace\Box}}

\def\hD{{\hat D}}
\def\bV{{\bar V}}
\def\bD{{\bar \Dl}}
\def\cU{{\cal U}}
\def\cC{{\cal C}}
\def\cP{{\cal P}}
\def\fg{{\frak g}}
\def\tr{\text{tr}}
\def\cZ{{\cal Z}}
\def\cD{{\cal D}}
\def\bR{{\Bbb R}}
\def\cE{{\cal E}}
\def\bZ{{\Bbb Z}}
\def\bN{{\Bbb N}}

\def\bysame{\same[\kern2cm]\,}

\def\br#1{\left\lfloor#1\right\rfloor}
\def\ag#1{\left\langle#1\right\rangle}
\def\BR#1{\left\lceil#1\right\rceil}

\def\abstractname{}

\@addtoreset {footnote}{page}

\renewcommand{\section}{%
   \@startsection
         {section}{1}{\z@}{-2.5ex \@plus -1ex \@minus -.2ex}%
               {3ex \@plus.2ex}{\Large\bf}%
}

\renewcommand{\subsubsection}{%
   \@startsection
         {subsubsection}{1}{\z@}{-1.5ex \@plus -1ex \@minus -.2ex}%
               {1ex \@plus.2ex}{\large\bf}%
}
\renewcommand{\@seccntformat}[1]{\csname the#1\endcsname .
\quad}

\def\bC{{\Bbb C}}
\def\bP{{\Bbb P}}

\makeatletter

\let\old@tl\~\def\~{\raisebox{-0.8ex}{\tt\old@tl{}}}
\let\lra\longrightarrow
\let\sm\setminus
\let\eps\varepsilon
\let\ex\exists
\let\fa\forall
\let\ps\supset

\def\rs#1{\raisebox{-0.4em}{$\big|_{#1}$}}

{\let\@noitemerr\relax
\vskip-2.7em\kern0pt\begin{abstract}
\noindent{\bf Abstract.}\enspace
Using an involved study of the Jones polynomial, we determine, as
our main result, the crossing numbers of (prime) amphicheiral knots.
As further applications, we show that several classes of links,
including semiadequate links and Whitehead doubles of semiadequate
knots, have non-trivial Jones polynomial. We also prove that there
are infinitely many positive knots with no positive minimal crossing
diagrams. Some relations to the twist number of a link, Mahler
measure and the hyperbolic volume are given, for example explicit
upper bounds on the volume for Montesinos and 3-braid links in terms
of their Jones polynomial. \\[2mm]
\end{abstract}
}\vspace{7mm}

{\parskip0.2mm\tableofcontents}
\vspace{7mm}

\section{Introduction}

The ultimate goal of this paper is to prove

\begin{theorem}\label{thm}
For each natural number $n\ge 15$, there exists a prime amphicheiral
knot of crossing number $n$.
\end{theorem}

This result thus determines the crossing numbers of prime amphicheiral
(or achiral) knots, i.e. those isotopic to their mirror images.
It can be considered as concluding a topic with an illustrious
history. The problem (given also as question 1.66 in \cite{Kirby})
goes back to the origins of knot theory some 120 years ago, when
Tait started compiling knot tables and sought the amphicheiral knots
therein (see for example \cite{HTW}). He actually might have believed
that only even $n$ occur, based on evidence from alternating knots
(which are highly over-represented in small crossing numbers).
The proof that alternating amphicheiral knots have even crossing
number was a major step that came with the discovery of the Jones
polynomial $V$ \cite{Jones} about 20 years ago. It was obtained as a
consequence of the crossing number result for alternating knots of
Kauffman, Murasugi and Thistlethwaite \cite{Kauffman,Murasugi,Thistle2}.
Their work also easily allows one to realize all even $n\ge 4$
(solving part (B) or the op.\ cit.\ question in Kirby's book).
The odd $n$ (part (A) of the question) are far more difficult to
deal with. A knot for $n=15$ was found by Hoste and Thistlethwaite
in a purely computational manner, during the compilation of knot
tables. However, exhaustive enumeration quickly becomes impossible
with increasing crossing numbers, and no generally applicable
methods are known to handle such examples. The approach
leading to our result will emerge from a much more advanced
treatment of the Jones polynomial, which has further noteworthy
implications.

The Jones polynomial, which has led to a significant change-of-face
of low-dimensional topology in recent years, remains despite its
popularity far from understood. Yet we do not know much about the
direct appearance of this invariant. For example, we can still not
answer Jones' question if it distinguishes all knots from the trivial
one (see \cite{Bi,Rolfsen,DH}). In \cite{Kauffman,Murasugi,Thistle2}
it was proved that the span (difference between minimal and
maximal degree) gives a lower bound for the crossing number,
which is exact on alternating diagrams, thus confirming a
century-old conjecture of Tait. For such links the coefficients
alternate in sign. In \cite{LickThis} (semi)adequate links were
introduced, motivated by further extensions of non-triviality and
crossing number results. Adequate links are a generalization of
alternating links, but in practice relatively few non-alternating
links are adequate. Semiadequate links are a much larger class.
It contains, beside alternating links, also other important
classes like positive and Montesinos links. In \cite{Thistle},
Thistlethwaite obtained a description of certain ``critical line''
coefficients of the Kauffman polynomial of
semiadequate links in terms of graph-theoretic invariants. For the
Jones polynomial it was observed in \cite{LickThis} that (almost by
definition) for semiadequate links the leading coefficient becomes
$\pm 1$. Apart from this insight, a while there was no further
information we had about coefficients of the polynomial.

This paper aims to provide some understanding of coefficients 2 and 3
of the Jones polynomial. Our starting point will be formulas for these
coefficients in semiadequate diagrams, shown in \S\reference{SSS}
(propositions \reference{pp21} and \reference{pv2}).
The formula for coefficient 2 was motivated by, and is a common
generalization of the one found by Dasbach-Lin \cite{DL} in
alternating diagrams and myself in \cite{posqp} for positive
diagrams. With these formulas we thus offer a unifying concept
of the two approaches. The formula for the third coefficient
was obtained independently and simultaneously by Dasbach-Lin in
\cite{DL2}. However, rather than just focusing on the formulas
themselves, our main contribution here will be to apply them to
several, incl. long-standing problems in knot theory.
These applications occupy the later sections, and can be
briefly grouped as follows.

% \begin{enumerate}
% \item
\em{Non-triviality of the Jones polynomial.}
Using the formula for coefficient 2 and some knowledge of properties
of positive links (some of which are prepared in \S\reference{SSS}),
we can conclude this non-triviality in many cases, including
semiadequate links (theorem \ref{tsqn}), Whitehead doubles
of semiadequate knots (proposition \reference{pWD}), 3-braid and
Montesinos links (corollary \reference{cntm}), and some
strongly $n$-trivial knots (proposition \reference{psnt}).
These results are treated in \S\ref{NT}
(except for 3-braids, which occupy \S\ref{S3B}). They extend
all previously known results 
(of \cite{Kauffman,Murasugi,Thistle2} for alternating,
\cite{LickThis} for adequate links, and \cite{Thistle2} for the
Kauffman polynomial).

\em{Estimates of hyperbolic volume.} The Jones polynomial coefficients
are related, as Dasbach-Lin \cite{DL} observe, to a quantity
called twist number. This number occurs in recent work of
Lackenby-Agol-Thurston \cite{Lackenby} as a way of estimating
hyperbolic volume. We can thus amplify Dasbach-Lin's
volume inequalities for alternating links, by similar (though
slightly more involved) upper estimates of the volume for adequate,
3-braid and Montesinos links in \S\ref{vest} (propositions
\reference{VOL3} and \reference{VOLM}).  Such relations aim
at providing a less exact, but more practical alternative to
the so-called Volume conjecture \cite{MM} (which states for a
general link an exact, but very involved, formula for the volume
in terms of the \em{colored} Jones polynomial). They also relate
to the open problem whether one can augment hyperbolic volume but
preserve the Jones polynomial. 

\em{3-braid links}. It turns out that 3-braid links are semiadequate
(corollary \reference{C3SQ}),
and we obtain several applications to them. We readily have
non-triviality of the polynomial (which,
% simplifies a cumbersome previous calculation for 3-braid knots (and
\em{inter alia}, also
gives a new proof of the faithfulness of the Burau representation).
We solve the problem of minimal conjugacy length in the 3-braid
group (theorem \reference{tq}), which was asked for by Birman-%
Menasco \cite{BirMen}, and extend their classification of
amphicheiral 3-braid links to those which are unorientedly
amphicheiral (theorem \ref{TU3}). We can also describe the 3-braid
links of unsharp Morton-Williams-Franks braid index inequality
(corollary \reference{CMWF}). %; a result claimed by El-Rifai

\em{Positive and amphicheiral knots}. The applications that require
the most substantial effort are the proof that infinitely many
positive knots have no minimal crossing diagrams in \S\ref{mpd}
(propositions \reference{ppp} and \reference{ppf}), and its
extension to show theorem \ref{thm}. In latter's proof, there
is a certain dichotomy between crossing numbers $n=15+4k$ and
$n=17+4k$, resulting from the somewhat different nature of the
examples we focus on. Accordingly, the proof of theorem \ref{thm}
is divided into the last two sections.

% that there are infinitely many amphicheiral
% knots of odd crossing number in \S\ref{mpe}. Latter result, theorem
% \reference{th15}, can be considered as solving (at least ``one half''
% of) a problem that goes back to some 120 years ago, when Tait started
% compiling knot tables and sought the amphicheiral knots therein.
% We extend the 15 crossing knot of Hoste and Thistlethwaite to an
% infinite family, settling crossing numbers $15+4k$, $k>0$.

The formulas of \S\reference{SSS} are the central and unifying theme
in all these results. The only part of the paper they do not directly
appear in is \S\ref{SMM}, which studies the relation between the
leading and trailing coefficients of the Jones polynomial, twist
numbers and Mahler measure. From this relation, Dasbach-Lin obtained,
and then further empirically speculated about, certain relations
between coefficients of the Jones polynomial and hyperbolic volume.
We give a qualitative improvement of the Dasbach-Lin result \cite{DL},
showing that \em{every} coefficient of the Jones polynomial gives
rise to a lower bound for the volume of alternating knots (corollary
\reference{ZZ}). This supplies with some more understanding their
experimental questions.

Besides polynomial non-triviality, and apart from the above particular
series of (achiral) examples, one should hope
to say something on crossing
numbers of semiadequate links in general. The only previous result,
in \cite{Thistle}, is that these links are non-trivial (i.e. their
crossing number is not 0). We develop a method to obtain estimates
on crossing numbers of semiadequate links in a separate paper
\cite{adeq}. Among others, we prove that a semiadequate link
has only finitely many reduced semiadequate diagrams, and there are
only finitely many semiadequate links with the same Kauffman polynomial.

\section{Basic preliminaries\label{BP}}

\begin{figure}[htb]
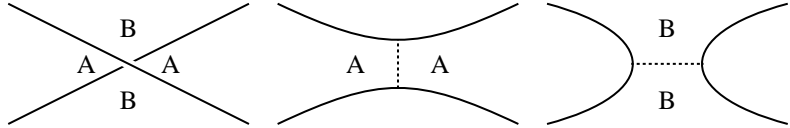

\[
\diag{8mm}{4}{2}{
   \picline{0 0}{4 2}
   \picmultiline{-5.0 1 -1.0 0}{0 2}{4 0}
   \picputtext{2.7 1}{A}
   \picputtext{1.3 1}{A}
   \picputtext{2 1.6}{B}
   \picputtext{2 0.4}{B}
} \quad
\diag{8mm}{4}{2}{
   \pictranslate{2 1}{
       \picmultigraphics[S]{2}{1 -1}{
           \piccurve{-2 1}{-0.3 0.2}{0.3 0.2}{2 1}
       }
       {\piclinedash{0.05 0.05}{0.01}
        \picline{0 -0.4}{0 0.4}
       }
   }
   \picputtext{2.7 1}{A}
   \picputtext{1.3 1}{A}
} \quad
\diag{8mm}{4}{2}{
   \pictranslate{2 1}{
       \picmultigraphics[S]{2}{-1 1}{
           \piccurve{2 -1}{0.1 -0.5}{0.1 0.5}{2 1}
       }
       {\piclinedash{0.05 0.05}{0.01}
        \picline{0 -0.6 x}{0 0.6 x}
       }
   }
   \picputtext{2 1.6}{B}
   \picputtext{2 0.4}{B}
}
\]
\caption{\label{figsplit}The A- and B-corners of a
crossing, and its both splittings. The corner A (resp. B)
is the one passed by the overcrossing strand when rotated 
counterclockwise (resp. clockwise) towards the undercrossing 
strand. A type A (resp.\ B) splitting is obtained by connecting 
the A (resp.\ B) corners of the crossing. It is useful to
put a ``trace'' of each splitted crossing as an arc connecting
the loops at the splitted spot.}
\end{figure}

The Jones polynomial is useful to define here via Kauffman's state
model \cite{Kauffman}. Recall, that the Kauffman bracket $\ag{D}$
of a link diagram $D$ is a Laurent polynomial in a variable
$A$, obtained by summing over all states $S$ the terms
\begin{eqn}\label{eq_12}
A^{\#A(S)-\#B(S)}\,\left(-A^2-A^{-2}\right)^{|S|-1}\,,
\end{eqn}
where a \em{state} is a choice of \em{splicings} (or
\em{splittings}) of type $A$ or $B$ for any
single crossing (see figure \ref{figsplit}), $\#A(S)$ and $\#B(S)$
denote the number of type A (resp. type B) splittings and $|S|$ the
number of (disjoint) circles obtained after all splittings in $S$.

We call the \em{$A$-state} the state in which all crossings are
$A$-spliced, and \em{$B$-state} is defined analogously. We call a trace
$a$ in the $A$-state \em{dual} to a trace $b$ in the $B$-state, if
$a$ and $b$ correspond to the same crossing as in figure \ref{figsplit}.

The Jones polynomial of a link $L$ can be specified from the
Kauffman bracket of some diagram $D$ of $L$ by
\begin{eqn}\label{conv}
V_L(t)\,=\,\left(-t^{-3/4}\right)^{-w(D)}\,\ag{D}
\raisebox{-0.6em}{$\Big |_{A=t^{-1/4}}$}\,,
\end{eqn}
with $w(D)$ being the writhe of $D$.

Let $S$ be the $A$-state of a diagram $D$ and $S'$ a state of $D$
with exactly one $B$-splicing. If $|S|>|S'|$ for all such $S'$,
we say that $D$ is \em{$A$-adequate}. Similarly one defines a
$B$-adequate diagram $D$. See \cite{LickThis,Thistle}. Then we set
a diagram to be
\begin{eqnarray*}
\mbox{adequate} & = & \mbox{A-semiadequate}\es
  \mbox{\em and}\es\mbox{B-semiadequate}\,, \\
\mbox{semiadequate} & = & %\mbox{\em either}\es
\mbox{A-semiadequate}\es\mbox{\em or}\es\mbox{B-semiadequate}\,, \\
\mbox{inadequate} & = & \mbox{\em neither}\es
  \mbox{A-semiadequate}\es\mbox{\em nor}\es\mbox{B-semiadequate}\,.
\end{eqnarray*}
(Note that inadequate is a stronger condition than not to be adequate.)

A link is called ($A$ or $B$-)adequate, if it has an ($A$ or
$B$-)adequate diagram. A link is \em{semiadequate} if it is $A$- or
$B$-adequate. A link is \em{inadequate}, if it is neither
$A$- nor $B$-adequate. 
% We call a link adequate resp.\ (A/B)-semiadequate if it has an
% adequate resp.\ (A/B)-semiadequate diagram.

As noted, semiadequate links are a much wider extension of the class
of alternating links than adequate links. For example, only 3
non-alternating knots in Rolfsen's tables \cite[appendix]{Rolfsen2}
are adequate, while all 55 are semiadequate.

Although the other polynomials will make only a secondary
appearance in the paper, we include a description at least to
clarify conventions.

The \em{skein (HOMFLY) polynomial} $P$ is a
Laurent polynomial in two variables $l$
and $m$ of oriented knots and links and can be defined
by being $1$ on the unknot and the (skein) relation
\begin{eqn}\label{1}
l^{-1}\,P\bigl(
\diag{5mm}{1}{1}{
\picmultivecline{0.18 1 -1.0 0}{1 0}{0 1}
\picmultivecline{0.18 1 -1.0 0}{0 0}{1 1}
}
\bigr)\,+\,
l \,P\bigl(
\diag{5mm}{1}{1}{
\picmultivecline{0.18 1 -1 0}{0 0}{1 1}
\picmultivecline{0.18 1 -1 0}{1 0}{0 1}
}
\bigr)\,=\,
-m\,P\bigl(
\diag{5mm}{1}{1}{
\piccirclevecarc{1.35 0.5}{0.7}{-230 -130}
\piccirclevecarc{-0.35 0.5}{0.7}{310 50}
}
\bigr)\,.
\end{eqn}
This convention uses the variables of \cite{LickMil}, but
differs from theirs by the interchange of $l$ and $l^{-1}$.
We call the three diagram fragments in \eqref{1} from left to
right a \em{positive} crossing, a \em{negative} crossing and a
\em{smoothed out} crossing (in the skein sense).

A diagram is called \em{positive}, if all its crossings are
positive. A(n oriented) link is positive, if it admits a
positive diagram (see for example \cite{pos}).

The \em{Kauffman polynomial} \cite{Kauffman2} $F$ is usually defined
via a regular isotopy invariant $\Lm(a,z)$ of unoriented links.

We use here a slightly different convention for the variables
in $F$, differing from \cite{Kauffman2,Thistle} by the interchange
of $a$ and $a^{-1}$. Thus in particular we have the relation
$F(D)(a,z)=a^{w(D)}\Lm(D)(a,z)$, where $w(D)$ is the writhe
of a link diagram $D$, and $\Lm(D)$ is the writhe-unnormalized
version of the polynomial. $\Lm$ is given in our convention by
the properties
\[
\begin{array}{c}
\Lambda\bigl(\Pos{0.5cm}{}\bigr)\ +\ \Lambda\bigl(\Neg{0.5cm}{}\bigr)\ =\ z\ \bigl(\ 
\Lambda\bigl(\Nul{0.5cm}{}\bigr)\ +\ \Lambda\bigl(\Inf{0.5cm}{}\bigr)\ \bigr)\,,\\[2mm]
\Lambda\bigl(\ \ReidI{-}{-}\bigr) = a^{-1}\ \Lambda\bigl(\noloop\bigr);\quad
\Lambda\bigl(\ \ReidI{-}{ }\bigr) = a\ \Lambda\bigl(\noloop\bigr)\,,\\[2mm]
\Lambda\bigl(\,\mbox{\Large $\bigcirc$}\,\bigr) = 1\,.
\end{array}
\]

Note that for $P$ and $F$ there are several other
variable conventions, differing from each other by possible inversion
and/or multiplication of some variable by some fourth root of unity.

The Jones polynomial $V$ can be obtained from $P$ and $F$ (in
our conventions) by the substitutions (with $i=\sqrt{-1}$; see 
\cite{LickMil} or \cite[\S III]{Kauffman2})
\begin{eqn}\label{xx}
V(t) \quad = \quad P(-it,i(t^{-1/2}-t^{1/2}))\quad =
\quad F(-t^{3/4},t^{1/4}+t^{-1/4})\,.
\end{eqn}

By $[P]_M$ we denote the coefficient of the monomial $M$ in the
polynomial $P$. If $P$ has a single variable, then we use the
exponent rather than the whole monomial for $M$. (For example,
$[V]_3=[V]_{t^3}$ for $V\in\bZ[t^{\pm 1}]$.)

\begin{defi}
Two crossings of a link diagram $D$ are twist equivalent if up
to flypes they form a clasp. (As in \cite{gen1}, this means that
they are either $\sim$-equivalent or $\ssim$-equivalent.)

Let $t(D)$ be the \em{twist number} of a diagram $D$,
which is the number of its twist equivalence classes. 
We call such equivalence classes also simply \em{twists}.
For a knot $K$ we define its twist number by
\[
t(K)\,:=\,\min\,\{\,t(D)\,:\,\mbox{$D$ is a diagram of $K$}\,\}\,.
\]
\end{defi}

(In order to distinguish the twist number from the variable of the
Jones polynomial, we will always write an argument in parentheses
behind it.)

Let $c(D)$ be the crossing number of a diagram $D$, and
$c(K)$ the crossing number of a knot or link $K$,
\[
c(K)\,:=\,\min\,\{\,c(D)\,:\,\mbox{$D$ is a diagram of $K$}\,\}\,.
\]

Thistlethwaite proved in \cite{Thistle} that (with our convention
for $\Lm$) for a link diagram $D$ of $c(D)$ crossings we have
$[\Lm(D)]_{z^la^m}\ne 0$ only if $l+|m|\le c(D)$, and that $D$ is
$A$ resp. $B$-adequate iff such a coefficient does not vanish for
some $l$ and $m$ with $l-m=c(D)$ resp. $l+m=c(D)$. These properties
imply most of his results, incl. the main one, the crossing number
minimality of adequate diagrams.

The coefficients of $\Lm(D)$ for which $l\pm m=c(D)$ form the
``critical line'' polynomials $\phi_{\mp}(D)$. Thistlethwaite
expresses these polynomials in terms of some graph invariants,
so they clearly encode combinatorial information of the diagram.
Unfortunately, we do not know how to interpret most of this
information, i.e., to say what tangible features of the diagram
it measures.

Still the minimal and maximal degree of $\phi_{\pm}(D)$ do have
a ``visual'' meaning. In \S 5 of \cite{adeq} we translated this
meaning to our present context. One degree can be expressed also
in terms of the writhe and Jones polynomial, thus giving a new
obstruction to semiadequacy (settling in particular the undecided
12 crossing knot in \cite{Thistle}). However, for semiadequate
links this obstruction (consequently vanishes and) gives no new
information. In contrast, the other degree (see theorem \ref{ozp}
below) gives a new invariant for semiadequacy. This will be used
crucially in the last sections, together with our work of ``decoding''
in such a visual way the Jones polynomial coefficients.

We also require occasionally the \em{Alexander polynomial} $\Dl$. We
define it here by being $1$ on the unknot and a relation involving
the diagrams occurring in \eqref{1} of triples of links $L_+$, $L_-$
and $L_0$, differing just at one crossing,
\begin{eqn}\label{srel}
\Dl(L_+)-\Dl(L_-)\,=\,\bigl(t^{1/2}-t^{-1/2}\bigr)\Dl(L_0)\,.
\end{eqn}
This relation is clearly a special case of \eqref{1}. Consequently,
there is the substitution formula
(see \cite{LickMil}; $i$ is again the complex unit),
\[
\Dl(t)\,=\,P(i,i(t^{1/2}-t^{-1/2}))\,,
\]
expressing $\Dl$ as a special case of $P$.

% The Jones polynomial $V(t)$ and Alexander polynomial $\Dl$ are
% used with the standard conventions, for example in \cite{coeff}.

\begin{defi}
When
\begin{eqn}\label{Vc}
V_K=a_0t^k+V_1t^{k+1}+\dots+a_dt^{k+d}
\end{eqn}
with $a_0\ne 0\ne a_d$ is the Jones polynomial of a knot or link $K$,
we will write for $d$ the \em{span} $\spn V_K$ of $V$, for $k$ the
\em{minimal degree} $\md V_K$ and for $k+d$ the \em{maximal degree}
$\Md V_K$. We will use throughout the paper the notation $V_i=V_i(K)=
a_i$ and $\bV_i=\bV_i(K)=a_{k-i}$ for the the $i+1$-st or $(i+1)$-last
coefficient of $V$ (since these terms will occur often, and to
abbreviate the clumsier alternative $[V]_{\md V+i}$ resp.
$[V]_{\Md V-i}$).
\end{defi}

Knots of $\le 10$ crossings will be denoted according to Rolfsen's
tables \cite[appendix]{Rolfsen2}, and for $\ge 11$ crossings according
to Hoste and Thistlethwaite's program KnotScape \cite{KnotScape},
appending for given crossing number the non-alternating
knots after the alternating ones.
The obverse (mirror image) of $K$ is denoted by $!K$. The mirroring
convention we use for knots in the tables, including for $\le 10$
crossings, is this of \cite{KnotScape}.

Dasbach-Lin (among other authors) noticed that $t(D)$ is an
invariant of alternating \em{diagrams} $D$ of an alternating knot $K$.
This follows, \em{inter alia}, from \cite{MenThis}. So we may define

\begin{defi}
The \em{alternating twist number} $t_a(K)$ is the twist
number of an alternating diagram $D$ of $K$.
\end{defi}

\begin{rem}
In general $t_a(K)>t(K)$, for example $7_4$ has $t_a(7_4)=3$ and
$t(7_4)=2$.
\end{rem}

In \cite{DL}, Dasbach-Lin consider $T_i(K):=|V_i|+|\bV_{i}|$ and prove

\begin{lemma}\label{lDL}(\cite{DL})
For an alternating knot $K$, we have $t_a(K)=T_1(K)$.
\end{lemma}

They apply this to obtain bounds on the volume of alternating knots
in terms of $|V_1|$ and $|\bV_{1}|$. This led us to consider these
values closer. In what follows we will explain the outcome. In
particular, one can extend some of their results in a weaker form
to adequate knots (\cite{Thistle}). On the other hand, $T_1$ seems
to be useful also in other situations, as we will see.

\section{Jones polynomial of (semi)adequate links\label{SSS}}

\subsection{The second coefficient}

We consider the bracket \cite{Kauffman} (rather than Tutte)
polynomial. The $A$-state of $D$, the state with all splicings
$A$, is denoted by $A(D)$. (In many cases we omit the argument
$D$ in this notation.). For us a state is always understood as
a planar picture of loops (solid lines) and traces connecting
these loops (dashed lines). Then it is clear that and how
to reconstruct $D$ from $A(D)$.

\begin{defi}
One fundamental object exploited in this paper is the \em{$A$-graph}
$G(A)=G(A(D))$ of $D$. It is defined as the planar graph with
vertices given by loops in the $A$-state of $D$, and edges given
by crossings of $D$. (The trace of each crossing connects two
loops.) The analogous terminology is set up also for the $B$-state.
% 
% Occasionally, by abuse of notation, we write $A(D)$ for $G(A(D))$.
% Still their distinction should be evident from
% the context. (For example a loop in $A(D)$ refers to the
% state, while a vertex in $A(D)$ to the graph.) We will
% point out when the separation between state and graph is
% essential (for example, in proposition \reference{pv2}).
\end{defi}

Clearly the $A$-state determines the $A$-graph, but not
conversely. Their distinction is relevant in some situations.
However, $G(A(D))$ (including its planar embedding) determines
$A(D)$ if $D$ is alternating; then sometimes $G(A(D))$ is called
the Tait graph of $D$. Note also that, for alternating $D$, the
duality of crossing traces between $A(D)$ and $B(D)$ corresponds 
to the duality (in the usual graph-theoretic sense) of edges in
the planar graphs $G(A(D))$ and $G(B(D))$.

Let $v(G)$ and $e(G)$ be the number
of vertices and edges of a graph $G$. Let $G'$ be $G$ with 
multiple edges removed (so that a simple edge remains). We 
call $G'$ the \em{reduction} of $G$.

We will write sometimes
\[
s_+(D)\,=\,v(G(A(D)))\,=\,v(G(A(D))')\,,\qquad
s_-(D)\,=\,v(G(B(D)))\,=\,v(G(B(D))')\,.
\]

The definition of $A$-adequate can be restated saying that
$G(A(D))$ has no edges connecting the same vertex. For $B$-splicings
the graph $G(B(D))$ and the property $B$-adequate are similarly
defined (and what is stated below proved).

In the following, we shall explain the second and third coefficient
of the Jones polynomial in semiadequate diagrams. Bae and Morton
\cite{Bmo} and Manchon \cite{Manchon} have done work in a different
direction, and studied the leading coefficients of the bracket (which
are $\pm 1$ in $A$-adequate diagrams) in more general situations.

Then we have

\begin{prop}\label{pp21}
If $D$ is $A$-adequate connected diagram, then in the representation
\eqref{Vc} of $V_D$ we have $V_0=\pm 1$, $V_1V_0\le 0$, and
\begin{eqn}\label{a1}
|V_{1}|=e(G(A(D))')-v(G(A(D))')+1=b_1(G(A(D))')
\end{eqn}
is the first Betti number of the reduced $A$-graph.
\end{prop}

\proof This will follow from the proof of proposition \ref{pv2}.
Alternatively, see \cite{posqp} where the formula is proved
for $D$ positive (then $A(D)$ is the Seifert graph).
The adaptation to $A$-adequate diagrams is easy. \qed

\begin{lemma}\label{lm1}
If $G$ is a planar simple graph (no multiple edges), then
$b_1(G)\le\br{\ffrac{2}{3}e(G)}-1$. This inequality is sharp
for proper $G$ when $e(G)>2$.
\end{lemma}

(Here and below $\br{x}$ is the largest integer not greater
than $x$.)

\proof The plane complement of $G$ has $b_1(G)+1$ cells
(including the one at infinity). 
Each cell has at least 3 edges, since $G$ has no
multiple edges, and each edge bounds at most two cells.
Now to make the inequality sharp, start with a triangle,
repeatedly connect (with 3 edges) a new inner vertex in one
of the faces of the graph. If $e\equiv 2\bmod 3$, finally
add a vertex on an edge $h$ between $v_{1,2}$ and connect it to 
$v$, where $v\ne v_{1,2}$ is a vertex of one of the faces that
$h$ bounds. \qed

In the following positivity arguments will be essential.
It is well-known that if a diagram $D$ is positive, then it
is $A$-adequate: the $A$-state of $D$ is just the Seifert
picture of $D$ (and the $A$-state loops are the Seifert circles).
Since $A$-adequacy is an unoriented condition, $D$ would remain
$A$-adequate even if we alter orientation of some components.
We say that an unoriented diagram $D$ \em{admits a positive
orientation}, or is \em{positively orientable}, if it arises
from such a diagram by forgetting orientation. The following lemma,
which will be repeatedly used below, specifies which $A$-adequate
diagrams are positively orientable.

\begin{lemma}\label{lmps}
Let $D$ be $A$-adequate. Then $D$ is positively orientable
iff $G(A(D))$ is a bipartite graph. 
\end{lemma}

\proof If $D$ is positively orientable, its graph $G(A(D))$ is the
Seifert graph of a (positive) diagram, and hence bipartite.
Conversely, if $G(A(D))$ is bipartite, it is possible to orient
the loops in $A(D)$ so that each trace looks locally like
\ \diag{5mm}{1}{1.3}{
  \picvecline{0 0}{0 1.3}
  \picvecline{1 0}{1 1.3}
  {\piclinedash{0.09 0.09}{0.01}
   \picline{0 0.6}{1 0.6}
  }
}\,.
Then it is clear that it is possible to extend this loop
orientation to an orientation of $D$, and with that orientation
$D$ becomes positive. \qed

\begin{corr}\label{crone}
Let $L$ be an $A$-adequate link, with an $A$-adequate connected
diagram $D$, then
\[
1-\br{\ffrac{2}{3}c(D)}\,\le\,V_0V_1\le\,\,0\,.
\]
If $V_1=0$, then $D$ admits a positive orientation, and (with this
orientation) $L$ is fibered.
\end{corr}

\proof We have $V_1=-V_0b_1(G(A(D))')$. Now $G(A(D))'$ is planar, so
to its $b_1$ we apply lemma \reference{lm1}.

Now if $V_1=0$, then $G(A(D))'$ is a tree, so in particular bipartite.
So the loops of the $A$-state of $D$ can be oriented alternatingly,
and then they become Seifert circles, and with the inherited orientation
the crossings become positive. That $L$ is then fibered follows from
\cite{posqp}. \qed

This gives a new semiadequacy test. Beside that Thistlethwaite's
condition on the positive critical line coefficients \cite{Thistle}
involves the Kauffman polynomial, which is considerably slower to
calculate, the Jones polynomial features sometimes prove essential,
as show the following examples.

\begin{exam}\label{x946}
The knots $!9_{46}$, $!9_{47}$ and $9_{48}$ (mirrored as in KnotScape)
have positive critical line coefficients of the Kauffman polynomial,
but $|V_0|=2$.
\end{exam}

\begin{exam}
The knot $14_{22068}$ has positive critical line coefficients,
but $V_0V_1>0$.
\end{exam}

% Here and below prime knots of $\ge 11$ crossings will be
% denoted according to Hoste and Thistlethwaite \cite{KnotScape},

The minimal number $e(n)=e_n$ of edges needed for a planar
simple graph to have given $b_1=n$ is by lemma \reference{lm1}
\[
e_n\,=\,\left\{\,\begin{array}{c@{\quad}c}
\ffrac{3}{2}(n+1) & \mbox{$n$ odd}\\[3mm]
\ffrac{3}{2}n+2 & \mbox{$n$ even}\end{array}\right.\,.
\]

For the next example, and for later discussion, it is useful to
define some properties of graphs (see also \cite{MS}).

\begin{defi}\label{SDE}
The \em{join} (or \em{block sum}, as called in \cite{Murasugi3})
`$*$' of two graphs is defined by
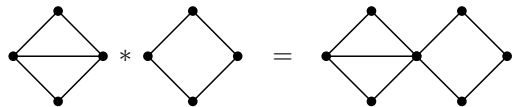
\[
\diag{6mm}{5}{2}{
  \picmultigraphics{2}{3 0}{\cycl{0 1}{1 0}{2 1}{1 2}}
  \picline{0 1}{2 1}
  \picputtext{2.5 1}{$*$}
}
\quad=\quad
\diag{6mm}{4}{2}{
  \picmultigraphics{2}{2 0}{\cycl{0 1}{1 0}{2 1}{1 2}}
  \picline{0 1}{2 1}
}
\]
This operation depends on the choice of a vertex in each one of
the graphs. %We call this vertex the \em{join vertex}. 

We call $v$ a \em{cut vertex} of a graph $G$, if $G$ gets
disconnected when deleting all edges incident to $v$ and
\em{additionally} $v$ itself. (When we delete an edge, we
understand that a vertex it is incident to is \em{not} to
be deleted too.)

Every connected non-trivial (i.e. with at least one edge) graph
$G$ can be written as a join $G_1*\dots*G_n$ for some non-trivial
connected graphs $G_i$, such that no $G_i$ has a cut vertex. We
call $G_i$ the \em{join factors} of the graph $G$. The number
$a(G)=n$ of join factors of $G$ is called \em{atom number} of $G$.
\end{defi}

\begin{exam}\label{x3.3}
The knot $14_{46350}$ has positive critical line coefficients, but
$V_0V_1=-8$. Similarly, $16_{484942}$, $16_{487600}$ and $16_{564314}$
have $V_0V_1=-10$. They are non-alternating. Since the $A$-state of a
non-alternating $A$-adequate diagram has a separating loop, we
see that the $A$-graph $G(A(D))$ of an $A$-adequate diagram $D$
must be the join of two graphs, i.e. have a cut vertex. Since 
\[
\min\,\bigl(\min\ \{\ e_a+e_b\,:\,a+b=8,\ a,b>0\ \},\ 
e_8+1\ \bigr)=15\,,
\]
we see that an $A$-adequate diagram of $14_{46350}$ must have
at least 15 crossings. Similarly an $A$-adequate diagram of the
16 crossing knots has at least 18 crossings. This is still a
non-trivial information in comparison to the Kauffman polynomial.
\end{exam}

In fact even without the assumption that the 16 crossing knots
in the previous example are non-alternating, we can deduce this.

\begin{corr}\label{cr3}
If $L$ is non-split and $e(|V_1|)>\spn V(L)$, then $L$ is
non-alternating.
\end{corr}

\proof If $D$ is an alternating diagram of $L$, then $c(D)=\spn V(L)$,
but we have $c(D)=e(A(D))\ge e(G(A(D))')\ge e(b_1(G(A(D))'))$, and since
$D$ is $A$-adequate, $e(b_1(G(A(D))'))=e(|V_1|)$. \qed

For the 16 crossing knots in example \reference{x3.3}
we have $e_{10}=17$, while $\spn V=14$.
So we see that these knots have no $A$-adequate minimal crossing
diagram, and so can not be alternating, not even adequate. Contrarily
their Jones polynomials are alternating, and monic on both sides.
Thus we have new information even compared to previous conditions
on the Jones polynomial.

\subsection{Conditions for positivity}

Note that in corollary \reference{crone} the case $V_1=0$
poses strong additional restrictions to $A$-semiadequacy. The
change of component orientation alters $V$ only by a positive unit
\cite{LickMil}, so that several positivity criteria still apply.
Among others, with $n(L)$ the number of components of $L$, we have
$V_0=(-1)^{n(L)-1}$, and for any $t\in(0,1]$ we have $(-1)^{n(L)-1}
V_L(t)\ge 0$ (see \cite{restr}). For knots positively orientable is
the same as positive, and we have a series of further properties.
For example $\md V_L=g(L)$, the genus of $L$, which is in particular
always positive. 

\begin{exam}
The knot $!11_{405}$ has $V_0=1$ and $V_1=0$, but $\md V=-2<0$.
The knot $!12_{1531}$ has $|V_0|=1$ and $V_1=0$, but $V_0=-1$, and
so is not $A$-semiadequate. (It has also $\md V<0$, but taking its
iterated connected sums with positive trefoils we can find an example
with $\md V>0$. Alternatively consider $!16_{1133621}$, where $\md
V=2$.) The knots have positive critical line polynomials.
\end{exam}

Similarly one can examine the inequalities for the Vassiliev
invariants
\begin{eqn}\label{v23}
v_2(K)\,=\,-\frac{1}{6}V''(1)\,\ge\,\frac{c(K)}{4}\,,
\quad\mbox{and}\quad\,
v_3(K)\,=\,-\frac{1}{12}V''(1)\,-\frac{1}{36}V'''(1)\,\ge\,
\frac{3}{8}c(K)-\frac{3}{4}\,,
\end{eqn}
proved in \cite{pos} (up to the different normalization of $v_3$
we use here). Again all these,
and the previous, conditions eventually dominate in the one or
other direction.

\begin{exam}
The knot $!15_{120559}$ has $V_0=1$ and $V_1=0$, and $v_2=5$ and
$v_3=16$, which satisfy the conditions of \cite{pos}, but $\md
V=0$. The knot $14_{26659}$ has $\md V=3$, $V_0=1$ and $V_1=0$,
and $v_2=4>14/4$, but $v_3=4<\myfrac{3}{8}\cdot 14-\myfrac{3}{4}$.
The knot $15_{136877}$ has $\md V=2$, $V_0=1$ and $V_1=0$,
and $v_3=5>\myfrac{3}{8}\cdot 15-\myfrac{3}{4}$, but $v_2=3<15/4$.
\end{exam}

Here is another application of the positive case. 

\begin{defi}
We call a link \em{weakly adequate}, if it is both $A$-adequate
and $B$-adequate.
\end{defi}

\begin{exam}
The knot of \cite{posex_bcr},
$11_{550}$, has a $B$-adequate 11 crossing diagram, which is not
$A$-adequate, and an $A$-adequate (because positive) 12 crossing
diagram. So it is weakly adequate but not adequate. Another such
example is Perko's knot $10_{161}$.
\end{exam}

\begin{prop}
A non-trivial weakly adequate \em{knot} has $T_1>0$.
\end{prop}

\proof $T_1(K)=0$ implies that both $K$ and its mirror image are
positive, which is impossible. \qed

\begin{exam}\label{x2.7}
Examples of adequate \em{links} with $T_1=0$ include the split union
of a figure-8-knot with an unknot, the Hopf link, and iterated
connected sums and disconnected twisted blackboard framed cables
of it. More such (incl. hyperbolic) examples will be given below.
\end{exam}

For knots a simple source of $T_1=0$ examples is

\begin{corr}
Let $K$ be a weakly adequate knot that has the Jones polynomial of a
$(m,n)$-torus knot. Then $m=2$. In particular, the only weakly adequate
torus knots are the $(2,n)$-torus knots.
\end{corr}

The second claim implies the claim for adequate torus knots, that
can be alternatively proved by combining work of \cite{Thistle} with
\cite{Murasugi2}, and in particular the property for alternating
torus knots, for which now numerous further proofs are known
(see \cite{restr} for some discussion).

\proof Since a torus knot is positive and fibered we see directly $V_1
=0$ for the Jones polynomial of any torus knot. Then (after possible
mirroring) $K$ is positive. Now we would like to prove that $\bV_{1}=0$
when $m,n>2$. We use Jones' formula for the polynomial of a torus knot
(see \cite[prop. 11.19]{Jones}).
\begin{equation}\label{Vtor}
V(T_{m,n})\,=\,-\frac{t^{m+n}-t^{m+1}-t^{n+1}+1}{t^2-1}
\,t^{(m-1)(n-1)/2}\,.
\end{equation}
Then, by conjugating and adjusting the degree, we see that $\bV_{1}=0$
is equivalent to the vanishing of the $t$-derivative of
$\ffrac{t^{m+n}-t^{m-1}-t^{n-1}+1}{1-t^2}$ at $t=0$. This is directly
to verify when $m,n>2$. Since one can recover
$(m,n)$ from $V(T_{m,n})$, the Jones polynomial distinguishes all
torus knots, and the second claim follows. \qed

\begin{exam}
It is not clear what non-torus knots can have the Jones polynomial
of a torus knot. For $m=2$ infinitely many such knots are known
from Birman \cite{Birman} (see \S\reference{S2q}),
but the case $m\ge 3$ is open. Such
examples do not occur among prime knots up to 16 crossings.
The simplest non-torus knot examples of $T_1=0$ are
$16_{1095701}$ and $16_{1318182}$ (which are hyperbolic).
\end{exam}

The next corollary determines asymptotically the minimal value
of $V_0V_1$ for positive link diagrams.

\begin{corr}
Let $L$ be a positive $n(L)$-component link, with a positive
connected diagram $D$. Then
\[
-\ffrac{1}{2}c(D)\,\le\,(-1)^{n(L)-1}V_1\,\le\,0\,,
\]
and the left inequality is asymptotically sharp (for large
$c(D)$).
\end{corr}

\proof Now $G(A(D))'$ is bipartite, so the smallest cell is a 4-gon. The
rest of the argument is as for lemma \reference{lm1}. To see sharpness
consider $A(D)$ to be the graph made by subdividing a square into
$n\times n$ smaller squares and let $n\to\infty$. \qed

\subsection{The third coefficient\label{S3c}}

The third coefficient in semiadequate diagrams can still be
identified in a relatively self-contained form. It depends, however,
on more than just $G(A(D))'$ or $G(A(D))$. Originally, the motivation
for the below formula was to estimate the volume of Montesinos links
(see \S\ref{vest}). However, the final sections show much more
important applications of it.

\parbox[t]{12.5cm}{
\begin{defi}
We call two edges $e_{1,2}$ in $G(A(D))'$ \em{intertwined},
if the following 3 conditions hold:
\begin{enumerate}
\item $e_{1,2}$ have a common vertex $v$.
\item The loop $l$ of $v$ in the $A$-state of $D$ separates the
loops $l_{1,2}$ of the other vertices $v_{1,2}$ of $e_{1,2}$.
\item $e_{1,2}$ correspond to at least a double edge in $A(D)$,
and there are traces of four crossings along $l$ that are connected
in cyclic order to $l_1,l_2,l_1,l_2$.
\end{enumerate}
\end{defi}
}
\enspace\parbox[t]{4.4cm}{
\[
\diag{9mm}{4}{4}{
  \picline{2 0}{2 4}
  {\piclinedash{0.1 0.1}{0.01}
   \picline{1 0.8}{2 0.8}
   \picline{3 1.6}{2 1.6}
   \picline{1 2.4}{2 2.4}
   \picline{3 3.2}{2 3.2}
  }
  \picfilledellipse{1.6 1 x}{0.3 1.2}{}
  \picfilledellipse{2.4 3 x}{0.3 1.2}{}
}
\]
}

% Here is an important point where the state $A(D)$ and graph $G(A(D))$
% must be separated.
It should be made clear that intertwinedness of edges in $G(A(D))'$
still depends on the \em{state} $A(D)$, rather than just the graph
$G(A(D))'$ or $G(A(D))$: in the graphs the intertwining information
becomes lost.

\begin{defi}
A \em{connection} in $A(D)$ is the set of traces between the
same two loops, i.e. an edge in $G(A(D))'$. A connection $e$
in $A(D)$ is said \em{multiple}, if it consists of at least
two crossing traces. More generally we can define the
\em{multiplicity} of a connection as the number of its traces.
\end{defi}

Then we can speak also of intertwined (multiple) connections. We will
later occasionally relax terminology even more and speak just of
intertwined loops $M$ and $N$, when their connections to a third
loop $L$ are intertwined. This is legitimate, because one can
determine $L$ from $(M,N)$ uniquely.

Call below a loop in $A(D)$ \em{separating} if it is connected
by crossing traces from either side. So if connections $(M,L)$
and $(N,L)$ are intertwined, then $L$ would be separating, and
connected by $M$ and $N$ from opposite sides. Clearly there can
be only one such $L$ for given $M$ and $N$. (Moreover, almost
throughout where we will apply this terminology below, there will
be in fact only one separating loop in $A(D)$.)

\begin{defi}
Define the \em{intertwining graph} $IG(A(D))$ to consist of vertices
given by multiple connections in $A(D)$ (or multiple edges in
$G(A(D))$), and edges connecting pairs of intertwined connections. 
\end{defi}

This is the second graph associated to $D$, which plays a
fundamental role in the whole paper.
Note that this is a simple graph (no multiple edges), but
not necessarily planar or connected. It may also be empty.
For example, for an alternating diagram this graph has no
edges, so is a (possibly empty) set of isolated vertices.
With the preceding remark, $IG(A(D))$ is determined by the
state $A(D)$, but not (in general) by the graph $G(A(D))$.

\begin{prop}\label{pv2}
If $D$ is $A$-adequate, then 
\[
V_0V_2\,=\,{|V_1|+1\choose 2}-\tg A(D)'+\chi(IG(A(D)))\,=\,
{|V_1|+1\choose 2}+e_{++}(A(D))-\tg A(D)'-\dl A(D)'\,,
\]
where $\tg A(D)'$ is the number of triangles (cycles of
length 3) in $G(A(D))'$, $\dl A(D)'=e(IG(A(D)))$ the number
of intertwined edge pairs in $G(A(D))'$, and $e_{++}(A(D))=
v(IG(A(D)))$ the number of multiple connections in $A(D)$.
\end{prop}

This formula was obtained independently and simultaneously by
Dasbach-Lin \cite{DL2}. Their proof is longer, but it unravels
the underlying combinatorics completely, while here we will help
ourselves with some (skein theoretic) work from \cite{posqp}.

\begin{rem}
Note that in alternating diagrams $\dl A(D)'=0$, since
there are no separating loops, while in positive diagrams
$\tg A(D)'=0$ because $G(A(D))'$ is bipartite. Note also that
a pair of intertwined edges does not occur in a triangle.
\end{rem}

\proof[of proposition \reference{pv2}]
The states that contribute to $V_2$ have $B$-splicings only in at
most two edges in $G(A(D))'$, or in 3 edges that form a triangle.

The zero edge case is the
$A$-state that defines $V_0$, and this gives a contribution of
$V_0{s_+-1\choose 2}$ to $V_2$ and $V_0(s_+-1)$ to $V_1$,
where $s_+(D)$ is the number of loops in the $A$-state.

The one edge case consists of the states that determine $V_1$.
If we have a connection $e$ in $A(D)$ of multiplicity $k$, then the
contributions of these states to $V_1$ are 
\begin{eqn}\label{opi}
V_0\left(-k+{k\choose 2}-{k\choose 3}+\dots\right)\,=-V_0\,.
\end{eqn}
With the term from the zero edge case this proves proposition
\reference{pp21}. Now the contributions to $V_0V_2$ are 
\[
-k(s_+-2)+{k\choose 2}(s_+-1)-{k\choose 3}s_++
{k\choose 4}(s_++1)-\dots\,=(s_+-2)\cdot
(-1)+\left[{k\choose 2}\cdot 1-
{k\choose 3}\cdot 2+{k\choose 4}\cdot 3-\dots\right]\,.
\]
The bracket evaluates to $1$ for multiple connections ($k\ge 2$)
and to $0$ for single ones ($k=1$).

Now we must consider the states that have $B$-splitting in two
edges in $G(A(D))'$. For most pairs of edges one has the product of
two evaluations as in \eqref{opi} (but with opposite sign, since
we have one loop less), which leads to ${e_+\choose 2}$,
with $e_+$ the number of edges of $G(A(D))'$.
So far we arrived, in collecting terms in $V_0V_2$, at
\begin{eqn}\label{imd}
{s_+-1\choose 2}-\big[e_+(s_+-2)-e_{++}\big]+{e_+\choose 2}=
{e_+-s_++2\choose 2}+e_{++}={|V_1|+1\choose 2}+e_{++}\,.
\end{eqn}
Herein all pairs of edges in $G(A(D))'$ were counted as if they
would contribute $1$. Now, there are two types of pairs of edges
where this may not be so, the intertwined edge pairs, and those
contained in triangles.

So now we need to take care of these ``degenerate pairs of edges''
in $G(A(D))'$, and of states with a triple of $B$-splicings in
triangles of $G(A(D))'$. Below we show how to find the corrections
to \eqref{imd} coming from these.

Consider triangles, and edge pairs contained in triangles. (Any
edge pair is contained in at most one triangle.) Fix an ordering
of edges in $G(A(D))'$, and then of traces in $A(D)$, so that if
$e_i<e_j$ in $A(D)$ then $e_i'<e_j'$ in $G(A(D))'$. Let $g_{1,2,3}
\in G(A(D))'$ be the edges of a fixed triangle with $g_1<g_2<g_3$.
We consider all states that have a $B$-splitting in preimages of
at least two of $g_{1,2,3}$. Set $H_{i}$ to be the set of preimage
traces of $g_{i}$ in $A(D)$. 

Let $i,j\in \{1,2,3\}$ with $i\ne j$. We sort the states by the
minimal (in $A(D)$ trace order) preimages $h_{i}\in H_i$,
$h_j\in H_j$ of $g_{i},g_j$ that are $B$-split. The case
that in all of $H_{1,2,3}$ a $B$-splicing is done is
incorporated into $i=1,j=2$. If one of $h_{i},h_j$ is not
maximal in $H_{i}$ resp.\ $H_j$, then the contributions
of the states that $B$-split any (possibly empty) set of
crossings in $H_{i,j}$ bigger than $h_{i,j}$, and possibly
in $H_3$ if $i=1,j=2$, is the product of two or three alternating
binomial coefficient sums, at least one of which evaluates
to $0$. So we need to count only contributions of states where
$h_{i,j}$ are both the maximal elements in $H_{i,j}$. However, if
$\{i,j\}=\{1,2\}$, we still have the non-trivial (and therefore
evaluating to 0) alternating binomial coefficient sum coming
from splicing edges in $H_3$. So the non-zero contributions
come from $\{i,j\}\ne \{1,2\}$, and they are two, rather than
3, as counted in ${e_+\choose 2}$. Thus the number of triangles
must be subtracted from \eqref{imd}.

Finally we need to take care of intertwined edge pairs. Here
the combinatorics is even messier, but we can help ourselves
with our previous work. Clearly, the contribution will not
depend on the remaining loops, and edges connecting them.
So we may evaluate it on a diagram what consists only of the
three relevant loops. Such are exactly the diagrams of prime
closed positive 3-braids. From \cite{posqp} we know that
$V_0V_2=1$ for such links. Since we can evaluate all other
terms easily, we obtain that intertwined edge pairs give no
contribution, and so must be subtracted from the ${e_+\choose 2}$
in \eqref{imd}. With this the desired formula is proved. \qed

Here is an application to positive links.

\begin{prop}\label{ppo}
If $K$ is a positive link, then with $\gm(K)=2\md V_K-|V_1|=
1-\chi(K)-|V_1|$, we have
\[
\min\left(\,0,\,\gm(K)-\left(\frac{\gm(K)}{2}\right)^2\,\right)\,
\le\,V_0V_2-{|V_1|+1\choose 2}\,\le\,\gm(K)\,.
\]
Moreover, these inequalities are sharp (i.e. equalities)
for links of arbitrarily small $\chi$.
\end{prop}

\proof If $D$ is a positive diagram of $K$, then $1-\chi=
b_1(G(A(D)))$ and $|V_1|=b_1(G(A(D))')$, so
\begin{eqn}\label{gme}
\gm(K)\,\ge\,e_{++}(A(D))\,.
\end{eqn}
We know also that $\tg A(D)'=0$. So it remains to see
\begin{eqn}\label{gmf}
\dl A(D)'\le \left(\frac{e_{++}}{2}\right)^2\,.
\end{eqn}
Now one can color the crossing traces (and then also
the edges in $G(A(D))'$) black and white, so that
traces of different (resp. same) colors connect to any
loop from different (resp. the same) side. Then only edge
pairs of opposite color can be intertwined in $G(A(D))'$.
For given number of multiple traces, the maximal number of
pairings occurs when all black traces are intertwined with
all white traces. In this situation the maximal number of
intertwinings is when exactly one half of the multiple traces
have either color.

To make the inequalities sharp, note that \eqref{gme} is sharp
whenever $A(D)$ has no $\ge 3$-ple traces. Many such diagrams
have $\dl=0$, for example the $(2,2,\dots,2)$-pretzel diagrams,
which thus realize the right inequality sharply.

To make the left inequality sharp we need to make
\eqref{gmf} sharp. Consider a positive diagram with one
separating Seifert circle $a$ and $n$ Seifert circles
$c_1,\dots,c_n$, $d_1,\dots,d_n$ on each side of $a$.
Let them be attached to $a$ by crossings in cyclic
order $c_1,\dots,c_n,d_1,\dots,d_n,c_n,\dots,c_1,d_n,
\dots,d_1$.  \qed

\begin{rem}\label{rNak}
Note also that $e_{++}(D)\le \myfrac{1}{2}c(D)$. This gives
some estimate on the crossing number of a positive diagram.
It seems, however, unlikely strong enough to prove
that some positive knots have no minimal (crossing number)
positive diagram. To find infinitely many such knots is a
problem was proposed by Nakamura \cite{Nakamura}.
We will later give a solution (see proposition \ref{ppp}),
although using more elaborate arguments.
\end{rem}

\begin{corr}\label{corf}
If $K$ is a fibered positive link, then
\begin{eqn}\label{8}
2\md V-\md V^2\,\le\,V_0V_2\,\le\,2\md V\,,
\end{eqn}
and when $K$ is prime, then $V_0V_2\,\le\,1$.
\end{corr}

\proof Now $|V_1|=0$. The inequalities \eqref{8} follow except
the cases of $(1-\chi)/2=\md V<2$ in the left inequality.
If $1-\chi\le 3$, then the left inequality is still valid when
\eqref{gme} is sharp. Since for a fibered positive link, $G(A(D))'$
is a tree, all connections in $A(D)$ can be assumed multiple
(otherwise $D$ can be reduced). So if \eqref{gme} is unsharp,
then at least a triple connection exists. Then for $1-\chi\le 3$
we have at most $3$ Seifert circles. Such diagrams are
of closed positive 3-braids. Then we know $V_0V_2=1$,
which satisfies the inequalities.

If now $K$ is prime, then $D$ is prime (see \cite{Ozawa}).
Now let $\star$ be an equivalence relation on
edges in $G(A(D))'$ defined by the transitive expansion
of the relation `intertwined'. That is, $k\star l$
if they can be related by a sequence of consecutively
intertwined edges. Now it remains to argue
that, since $G(A(D))'$ is a tree and $D$ is prime,
$\star$ must have a single equivalence class, and
therefore $\dl A(D)'\ge e_{++}-1$.

If $\star$ has $\ge 2$ equivalence classes, elements
$x_{1,2}$ in such classes connect some loop $l$ in
the $A$-state (from which side is not specified). For
a loop $x$ connecting $l$ by a connection of $n$ crossing
traces, inscribe an $n$-gon, consisting of $n$ chords,
into $l$ with vertices the endpoints of the traces.
Then if $x_{1,2}$ are in different $\star$-equivalence
classes, the chord diagram\footnote{Here we have chords
with coinciding basepoints, so this is not the sort
of chord diagram from Vassiliev invariant theory.}
in $l$ will be disconnected, with the chords coming
from $x_{1,2}$ in different components. So one can put
a interval $\gm$ that intersects $l$ twice and contains
parts of the chord diagram on both sides. Because
$G(A(D))'$ is a tree, none of the loops inside or
outside $l$ can be connected except through $l$.
So $\gm$ can be extended to a closed curve that
gives a non-trivial factor decomposition of $D$,
a contradiction.
\qed

\begin{rem}
The examples realizing the left inequality sharply
in proposition \reference{ppo} remain valid also for
corollary \reference{corf}. The right inequalities in the
corollary are sharp for positive braids. On the other hand,
positive braids allow an estimate also of $V_3$ \cite{posqp},
which seems hard in general. See also remark \reference{remV3}.
\end{rem}

\subsection{Cabling\label{Cab}}

In the following we briefly explain the effect of the
formulas for $V_{1,2}$ under cabling.

\begin{lemma}\label{lcab}
Let $D$ be an $A$-adequate link diagram (not necessarily reduced)
and $D_2$ the (blackboard framed) 2-parallel of $D$. Then
$V_0V_1(D_2)=V_0V_1(D)$ and
\[
V_0V_2(D_2)\,=\,\mybin{|V_1(D)|}{2}+1-\tg(A(D))\,.
\]
\end{lemma}

\proof To come from $D$ to $D_2$, one applies the
following modifications\\[6mm]
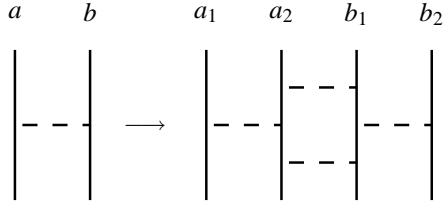
\[
\raisebox{0cm}{
\diag{1cm}{1}{2}{
  \picline{0 0}{0 2}
  \picline{1 0}{1 2}
  {\piclinedash{0.2}{0.1}
   \picline{1 1}{0 1}
  }
  \picputtext{0 2.5}{$a$}
  \picputtext{1 2.5}{$b$} 
}}
\quad\lra\quad
\raisebox{0cm}{
\diag{1cm}{3}{2}{
  \picmultigraphics{4}{1 0}{\picline{0 0}{0 2}}
  {\piclinedash{0.2}{0.1}
   \picline{1 1}{0 1}
   \picline{3 1}{2 1}
   \picline{2 1.5}{1 1.5}
   \picline{2 0.5}{1 0.5}
  }
  \picputtext{0 2.5}{$a_1$}
  \picputtext{1 2.5}{$a_2$} 
  \picputtext{2 2.5}{$b_1$}
  \picputtext{3 2.5}{$b_2$} 
}}
\]

Let $k(A(D))$ be the number of loops in $A(D)$ connected by
only one trace to some other loop. We call such loops isolated.

\begin{itemize}
\item All edges in $A(D_2)$ are multiple except $k(A(D))$
because all $a$-$b$ connections are doubled in $a_2-b_1$, and all 
loops in $A(D)$ except $k(A(D))$ are connected by at least two
traces to some other loop(s). So
\begin{eqn}\label{_1}
e_{++}(A(D_2))=e(G(A(D_2))')=v(A(D))+e(G(A(D))')-k(A(D))\,.
\end{eqn}

\item Triangles in $A(D_2)$ correspond bijectively to triangles
in $D$,
\begin{eqn}\label{_2}
\tg(A(D_2))=\tg(A(D))\,.
\end{eqn}

\item Now consider intertwined edges in $G(A(D_2))'$.
First, intertwined connections (or edge pairs in
$G(A(D_2))'$) correspond to 
\[
\bigl(\,\mbox{$(a_2,b_1)$-connection,
$(a_1,a_2)$-connection}\,\bigr)\,
\bigl(\,\mbox{$(a_2,b_1)$-connection,
$(b_1,b_2)$-connection}\,\bigr)\,,
\]
such that the $(a_1,a_2)$-connection resp. $(b_1,b_2)$-connection
are multiple, i.e. $a$ resp. $b$ is not an isolated loop in
$A(D)$. So for each edge\ \ $\diag{1em}{2}{1.3}{
\footnotesize\llpoint{0 0.5}\llpoint{2 0.5}
\picline{0 0.5}{2 0.5}\picputtext[d]{0 0.8}{$a$}
\picputtext[d]{2 0.8}{$b$}}$%
\ \ in $G(A(D))'$ there exist exactly two intertwined edge pairs
in $G(A(D_2))'$, but from these $2e(G(A(D))')$ pairs $k(A(D))$
must be excluded. (Here we use also that $a\ne b$.) So 
\begin{eqn}\label{_3}
\dl(G(A(D_2))')=2e(G(A(D))')-k(A(D))\,.
\end{eqn}
\end{itemize}

Then using \eqref{_1}, \eqref{_2}, and \eqref{_3},
% \tm
\begin{eqnarray}
\nonumber V_0V_2(D_2) & \,=\, & \mybin{|V_1(D_2)|+1}{2}-\tg(A(D_2))+
                      e_{++}(A(D_2))-\dl(A(D_2)) \\
\nonumber	 & \,=\, & \mybin{|V_1(D)|+1}{2} -\tg(A(D)) +
                      v(A(D))-e(G(A(D))') \\
\nonumber     & \,=\, & \mybin{|V_1(D)|+1}{2} -\tg(A(D)) +1 -|V_1(D)| \\
		 & \,=\, & \mybin{|V_1(D)|}{2} -\tg(A(D)) +1\,.
% this hack produces the qed-box instead of the eqn-label
\xdef\@tempk{\arabic{equation}}
\glet\eq@num\@eqnnum
\gdef\@eqnnum{\bgroup\let\r@fn\normalcolor
  \def\normalcolor##1(##2){\r@fn##1\qed}\eq@num\egroup}%
\edef\@currentlabel{\qed}
\end{eqnarray}%
\setcounter{equation}{\@tempk}%
\addtocounter{equation}{-1}%
\glet\@eqnnum\eq@num
% \qed
% \ntm

\begin{corr}\label{ccab}
The quantities $\tg(A(D))$ and $\chi(IG(A))=e_{++}(A(D))-\dl(A(D))$
are invariants of $A$-adequate diagrams $D$ of a link $L$. We call
them together with $\chi(G(A))$ the \em{A-state (semiadequacy)
invariants} of $L$.
\end{corr}

\proof Let $D_{1,2}$ be $A$-adequate diagrams of a link. By adding
kinks (which does not change the quantities) we may assume that 
$D_{1,2}$ have the same writhe. Then their 2-parallels are isotopic, 
and we see that $\tg(A(D))$ is invariant. The invariance of 
$\chi(IG(A(D)))$ follows from looking at $V_0V_2(D_i)$.
\qed

\begin{rem}
It is easily observed that the formulas remain correct
if we replace $2$-cables by $n$-cables for $n>2$.
(This means also that we obtain no new invariant information
from higher cables.)
\end{rem}

The advantage of lemma \reference{lcab} and corollary
\reference{ccab} is to extract invariant information out
of $V$ and its cables, which is much more direct to compute
than the whole (NP-hard; see \cite{JVW}) polynomial. This can,
for example, be applied to verify for a given knot a large
number of diagrams. % Such a verification was essential in
% the proof of the following theorem. (Details will appear in
% a separate paper with a contribution of Stefan Friedl.)
% 
% \begin{theorem}
% There exist mutants of different crossing number.
% \end{theorem}
% 
Note also that the meaning of the various quantities
encountered is rather visual from a semiadequate diagram~--
in contrast, for example, to the coefficients of
the Kauffman polynomial studied in \cite{Thistle2} in
this context. The work of sections \ref{mpd}, \ref{mpe}
and \reference{S10} exploits heavily this big advantage.
% We will use this convenient feature to prove 
% (with a much longer, but more down-to-earth argument)
% that there exist infinitely
% many odd crossing number achiral knots.

Let us note the following, which will be useful below,

\begin{corr}\label{cormut}
Semiadequacy invariants are mutation invariant.
\end{corr}

\proof It is well-known \cite{MorTr} that the Jones polynomial
and all its cables are mutation invariant. \qed

% \section{Applications}

\section{Non-triviality of the Jones polynomial\label{NT}}

\subsection{Semiadequate links}

An important application of corollary \reference{crone} is

\begin{theo}\label{tsqn}
There is no non-trivial semiadequate link $L$ with trivial Jones
polynomial (i.e., polynomial of the same component number unlink),
even up to units.
\end{theo}

This is an extension of the results of \cite{Kauffman,Murasugi,%
Thistle2} for alternating links, \cite{LickThis} for adequate
links, \cite{pos} for positive \em{knots}, and \cite{Thistle2}
for the Kauffman polynomial of semiadequate links. While the
existence of non-trivial links with trivial polynomial is now
settled for $n(L)\ge 2$ \cite{EKT}, the (most interesting)
case $n(L)=1$ remains open. In fact, except for the above
cited (meanwhile classical) results, and despite considerable
(including electronic) efforts \cite{Bi,Rolfsen,DH,pv4}, even
nicely defined general classes of knots on which one can exclude
trivial polynomial are scarce. Our theorem seems thus to
subsume all previous self-contained results in this direction.

\proof We use induction on $n(L)$. Let first $n(L)=1$. 
If the Jones polynomial of a \em{knot} is determined up
to units, it is uniquely determined, because $V(1)=1$ and $V'(1)=0$.
Thus $V=1$, and we have $V_1=0$. Now by corollary \reference{crone},
$L$ must be a positive knot, but there is no non-trivial positive knot
with trivial Jones polynomial \cite{pos}. Now let $n(L)>1$. By
\cite[corollary 3.2]{Thistle}, if $L$ is split, so is any semiadequate
diagram $D$ of $L$, and in this case we can argue by induction on the
number of components of the split parts of $D$. Else $D$ is connected
(non-split), but since  up to units $V=(-t^{1/2}-t^{-1/2})^{n(L)-1}$,
we have $V_0V_1=n(L)-1>0$ and a contradiction to corollary
\reference{crone}. \qed

\subsection{Montesinos links\label{Ml}}

\setbox\@tempboxa=\hbox{$M(\myfrac{3}{11},-\myfrac{1}{4},
\myfrac{2}{5},4)$}

\begin{figure}[htb]
\[
\begin{array}{c}
\epsfsv{5cm}{t1-dioph3} \\
\end{array}
\]
\caption{The Montesinos knot \unhbox\@tempboxa\ with
Conway notation $(213,-4,22,40)$.\label{figM}}
\end{figure}

To avoid confusion, let us fix some terminology relating to Montesinos
links. The details on Conway's notation can be found in his
original paper \cite{Conway}, or for example in \cite{Adams}.

Let the \em{continued} (or \em{iterated}) \em{fraction}
$[[s_1,\dots,s_r]]$ for integers $s_i$ be defined inductively by
$[[s]]=s$ and
\[
[[s_1,\dots,s_{r-1},s_r]]=s_r+\frac{1}{[[s_1,\dots,s_{r-1}]]}\,.
\]
The rational tangle $T(p/q)$ is the one with Conway notation
$c_{1}\ c_{2}\ \dots c_{n}$, when the $c_i$ are chosen so that
\begin{eqn}\label{ci}
[[c_{1},c_{2},c_{3},\dots,c_n]]=\frac{p}{q}\,.
\end{eqn}
One can assume without loss of generality that $(p,q)=1$, and $0<|p|<q$.
A \em{rational (or 2-bridge) link} $S(q,p)$ is the closure of $T(p/q)$.

Montesinos links (see e.g. \cite{BurZie}) are generalizations of
pretzel and rational links and special types of arborescent links. They
are denoted in the form $M(\frac{p_1}{q_1},\dots,\frac{p_n}{q_n},e)$,
where $e,q_i,p_i$ are integers, $(q_i,p_i)=1$ and $0<|p_i|<q_i$.
Sometimes $e$ is called the \em{integer part}, and the $\frac{p_i}
{q_i}$ are called \em{fractional parts}. They both together form the
\em{entries}. If $e=0$, it is omitted in the notation. Our convention
follows \cite{Oertel} and \em{may differ} from other authors' by the
sign of $e$ and/or multiplicative inversion of the fractional parts.
For example $M(\frac{p_1}{q_1},\dots,\frac{p_n}{q_n},e)$ is denoted 
as ${\frak m}(e;\frac{q_1}{p_1},\dots,$ $\frac{q_n}{p_n})$ in
\cite[definition 12.28]{BurZie} and as $M(-e;(q_1,p_1),\dots,$ $
(q_n,p_n))$ and the tables of \cite{Kawauchi}.

If all $|p_i|=1$, then the Montesinos link $M(\pm \frac{1}{q_1},\dots,
\pm \frac{1}{q_n},e)$ is called a \em{pretzel link} $P(\pm q_1,
\dots,\pm q_n,\eps,\dots,\eps)$, where $\eps=\sgn(e)$, and there
are $|e|$ copies of it. We also say it is a $(\pm q_1,\dots)$-pretzel
link.

To visualize the Montesinos link from a notation,
let $q_i/p_i$ be continued fractions of rational tangles
$c_{1,i}\dots c_{n_i,i}$ with $[[c_{1,i},c_{2,i},c_{3,i},\dots,
c_{l_i,i}]]=\frac{q_i}{p_i}$. Then $M(\frac{p_1}{q_1},
\dots,\frac{p_n}{q_n},e)$ is the link that corresponds to the
Conway notation
\begin{eqn}\label{cZ}
(c_{1,1}\dots c_{l_1,1}), (c_{1,2}\dots c_{l_2,2}), \dots,
(c_{1,n}\dots c_{l_n,n}), e\,0\,.
\end{eqn}
An example is shown in figure \reference{figM}.

An easy exercise shows that if $p_i>0$ resp. $p_i<0$, then 
\begin{eqn}\label{Mi}
M(\dots,p_i/q_i,\dots,e)\,=\,
M(\dots,(p_i\mp q_i)/q_i,\dots,e\pm 1)\,,
\end{eqn}
i.e. both forms represent the same link.
In our convention the identity \eqref{Mi} can be read naturally
to preserve the sum of all entries, and an integer entry
can be formally regarded as a fractional part. Theorem 12.29 in
\cite{BurZie} asserts that the entry sum, together with the vector
of the fractional parts, modulo $\bZ$ and up to cyclic permutations
and reversal, determine the isotopy class of a Montesinos link $L$.
So the number $n$ of fractional parts is an invariant of $L$; we call
it the \em{length} of $L$. If the length $n<3$, an easy observation
shows that the Montesinos link is in fact a rational link.

We can now conclude some work on Montesinos links that could only be
partially done in \cite{LickThis} using adequacy. Adequacy is clearly
too restrictive, for instance the knot $9_{46}$ is a pretzel knot that
is not adequate, but semiadequate; see example \reference{x946}.
(As the paper \cite{HTY} shows, the problem to determine the crossing
number of a general Montesinos link can unlikely be solved using
the Jones polynomial only; in \cite{LickThis} the Kauffman polynomial
was used.) 
% Thus the methods of
% \cite{LickThis} cannot apply to all Montesinos links, even although
% the word "all" appears explicitly (and overstatedly) in the abstract
% of \cite{LickThis}. (This point is important to notice, since
% otherwise some of the subsequent work here would not make sense.)

\begin{corr}\label{cntm}
No Montesinos link has trivial Jones polynomial up to units.
\end{corr}

\proof Let 
\begin{eqn}\label{Mont}
M(p_1/q_1,\dots,p_n/q_n,e)
\end{eqn}
be a representation chosen
up to mirroring so that $n\ge 2$, $0<p_i/q_i<1$. If $e\ge 0$, the
canonical diagram to this representation is alternating, if $e<-1$
it is semiadequate by splicing so that the $-e$ crossings create
$1-e$ loops in between them (and using that no $p_i/q_i$ represents
a single crossing tangle). Let this be, say, the $A$-state. Now
consider $e=-1$. Then write $M(p_1/q_1,\dots,p_{n-1}/q_{n-1},p_n'/q_n)$
such that $-1<p_n'/q_n<0$. If now $n=2$, we have a rational link,
so let $n>2$. Then since the tangle representing $p_n'/q_n$ does not
consist of a single crossing, the $B$-state shows semiadequacy. \qed

We can say something more on
the Jones polynomial, also in cases not covered in \cite{LickThis}.

For the rational link we have a new handy condition, which is
a simple straightforward observation.

\begin{prop}\label{prat}
If $L$ is a rational link, then $||V_1|-|\bV_{1}||\le 1$. \qed
\end{prop}

Let $D=M(p_1/q_1,\dots,p_n/q_n,e)$ be a Montesinos link diagram of
a representation chosen so that $n\ge 3$, $0<p_i/q_i<1$. (For
$n\le 2$ we obtain 2-bridge links.) Then if $e\ge 0$ or $e\le -n$,
the link is alternating. If $e<-1$ and $e>1-n$ then the link is
adequate and considered in \cite{LickThis}. We write $S(q,p)$
in Schubert form for the rational link with fraction $p/q$.

\begin{prop}\label{pmnt}
(1) If $e<-1$ in \eqref{Mont} then 
\[
|\bV_{1}|=|\bV_{1}(S(q_1,p_1)\#\dots\#S(q_n,p_n))|\,+\,
\left\{\,\begin{array}{c@{\quad}l}
1 & \mbox{if $e<-2$} \\
0 & \mbox{if $e=-2$} \end{array}\,\right.
\]
(2) If $e=-1$ and $\dl(L)\,:=\#\,\{\,i\,:p_i/q_i\ge 1/2\,\}-1$ is
non-zero, then $|\bV_0|=\dl$, and $\spn V(D)=c(D)-3$.
\end{prop}

\proof The diagram $D$ is $B$-adequate, and the $B$-state has
one separating loop that contains the $e$ twists. They contribute
to $b_1(A(D))$ depending of whether $e=-2$ or $e<-2$. This
explains part (1).

In part (2), the diagram $D$ is almost alternating. So it is neither
$A$- nor $B$-adequate. If we resolve the dealternator by the
bracket relation into diagrams $D_A$ and $D_B$, we observe that
$D_A$ and $D_B$ are both alternating (and in particular adequate).
The lowest terms in $A$ in $A\ag{D_A}$ and $A^{-1}\ag{D_B}$ will
cancel (since $D$ is not $B$-adequate), but $\dl$ is exactly the
difference of the second lowest terms. It is easy to see that
the highest terms in $A$ in $A\ag{D_A}$ and $A^{-1}\ag{D_B}$
again cancel (since $D$ is not $A$-adequate), but the second highest
terms differ by $\pm 1$ and do not cancel. So if $\dl\ne 0$,
the span of $V(D)$ is by two less than the difference between the
highest and lowest $t=A^{-1/4}$ terms, which is $(c(D)+s_+(D)+s_-(D)-
2)/2$. Since by direct check $s_+(D)+s_-(D)=c(D)$, the claim
follows. \qed

\begin{rem}
Note in particular that part (2) allows to construct Montesinos links
whose leading or trailing Jones polynomial coefficient is arbitrary
up to sign. (See also \cite{Manchon,HTY} for similar families of such
examples.)
\end{rem}

\begin{figure}[htb]
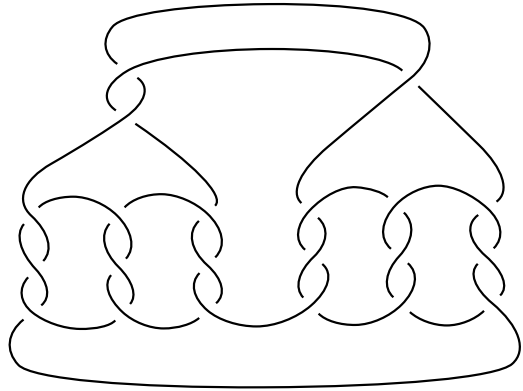

\begin{center}
\epsfsv{6cm}{t1-inaqarb}
\end{center}
\caption{\label{figinaqarb}An arborescent inadequate knot
(-3,-3,-3)2,(3,3,-3)1.}
\end{figure}

\begin{exam}
It was at first suggestive that one could extend corollary \ref{cntm}
also to
arborescent links, but a proof met difficulties that led to examples
like the one in figure \reference{figinaqarb}. This arborescent knot,
with Conway notation $(-3,-3,-3)2,(3,3,\lb -3)1$,
has Jones polynomial with leading and trailing coefficient $\pm 3$,
and so is neither $A$- nor $B$-adequate (though the Kauffman
polynomial is positive on both critical lines). Another example,
the arborescent knot $(3,3,3)-2,(3,3,3)-2,(-3,-3,-3)2$ has $V_0V_1=2$
and non-positive $B$-critical line of the Kauffman polynomial.
\end{exam}

% appears more technical and,
% if obtained, will be given at a later stage.
% 
% CHECK IF CEX CANDS ARE VALID!!!
% 
% 
% 
\subsection{Whitehead doubles}

Untwisted Whitehead doubles have trivial Alexander polynomial, 
and are one suggestive class of knots to look for trivial
Jones polynomial. (Practical calculations have shown that
the coefficients of the Jones polynomial of Whitehead doubles
are absolutely very small compared to their crossing number.)

\begin{prop}\label{pWD}
Let $K$ be a semiadequate non-trivial knot. Then the untwisted 
Whitehead doubles $Wh_{\pm}(K)$ of $K$ (with either clasp) have 
non-trivial Jones polynomial.
\end{prop}

\begin{rem} 
Because $V$ determines $v_2$ (from \eqref{v23}), and also
$v_2=\myfrac{1}{2}\Dl''(1)$, among Whitehead doubles only
untwisted ones may have trivial Jones polynomial.
\end{rem}

\proof Let $D$ be an $A$-adequate diagram of $K$. ($B$-adequate 
diagrams are dealt with similarly.) Without loss of generality
assume that $w(D)>0$ (possibly add positive kinks). If
one of $Wh_{\pm}(K)$ has $V=1$, then by a simple skein/bracket
relation argument, the blackboard framed, disconnected
2-cable diagram $\tl D$ of $D$ must have the
Jones polynomial $V_w=V(\tl D)$ of the $(2,2w)$-torus link
up to units. Since $D$ is $A$-adequate, so is $\tl D$,
and since $w>0$, the polynomial $V_w$ has $V_1=0$.
Thus $\tl D$ must admit a positive orientation. In particular,
all its components' self-crossings must be positive, that is,
$D$ is positive. But we know from \cite{pos} that the
untwisted Whitehead double of a positive knot has non-trivial 
Jones polynomial. \qed

This generalizes a result for adequate knots in
\cite{LickThis} and positive knots in \cite{pos} and  
considerably simplifies the quest for trivial
polynomial knots among Whitehead doubles. % after the
% conditions of \cite{pos,pv4}.
In particular, one can 
extend the verification of \cite{pv4} to establish that
no non-trivial knot of $\le 16$ crossings has untwisted 
Whitehead doubles with trivial Jones polynomial. 

\begin{rem}
We do not claim that $Wh_{\pm}(K)$ are themselves semiadequate.
For example, consider the alternating 11 crossing knots $11_{62}$
and $11_{247}$. The Jones polynomial of $Wh_{-}(11_{62})$ 
has $V_0=2$, while $Wh_{-}(11_{247})$ has Jones polynomial with
$\bV_{1}=0$ (but is clearly not negative). Thus neither $A$- nor
$B$-semiadequacy need to be preserved under taking \em{untwisted}
Whitehead doubles. (Note also that the Kauffman polynomial
semiadequacy test for such examples would have been
\em{considerably} less pleasant.)
\end{rem}

\subsection{Strongly $n$-trivial knots}

The next application concerns strongly $n$-trivial knots.
They were considered first around 1990 by Ohyama, and studied
more closely recently \cite{Torisu,HL,AK}. While one can
easily verify by calculating the Jones polynomial that
a given example is non-trivial, the proof of non-triviality
for a family of knots with arbitrarily large $n$ remained
open for a while. (For $n>2$ the Alexander polynomial is
trivial \cite{AK}.) A proof that partially features the Jones 
polynomial value $V(e^{\pi i/3})$ was given in \cite{tl},
but nothing about the Jones polynomial of the examples
directly could be said. How to evaluate the Jones polynomial
was also asked by Kalfagianni \cite{Kalfagianni}.
Now we can deal with a different class of examples, proving
the polynomial non-trivial.

\begin{prop}\label{psnt}
There exist for any $n$ strongly $n$-trivial knots with
non-trivial Jones polynomial.
\end{prop}

\proof Apply the construction of \cite{AK}
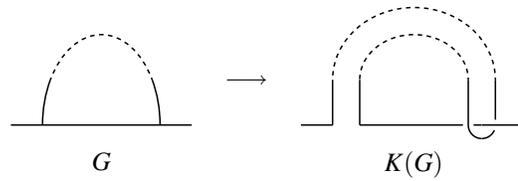
\[
\diag{6mm}{4}{3}{
  \picline{0 0.5}{4 0.5}
  \picellipsearc{2 0.5}{1.3 2}{0 30}
  \picellipsearc{2 0.5}{1.3 2}{150 180}
  \piclinedash{0.1 0.1}{0.05}
  \picellipsearc{2 0.5}{1.3 2}{30 150}
  \picputtext{2 -0.3}{$G$}
}
\quad\lra\quad
\diag{6mm}{5}{3}{
  \picPSgraphics{0 setlinecap}
  \picline{0 0.5}{0.7 0.5}
  \picline{1.3 0.5}{4 0.5}
  \picline{0.7 0.5}{0.7 1.5}
  \picline{1.3 0.5}{1.3 1.5}
  \picmultiarcangle{-7 1 -1 0}{3.7 1.5}{3.7 0.2}{4.0 0.2}{0.3}
  \picmultiarcangle{-7 1 -1 0}{4.0 0.2}{4.3 0.2}{4.3 1.5}{0.3}
  \picmultiline{-7 1 -1 0}{4 0.5}{5 0.5}
  \piclinedash{0.1 0.1}{0.05}
  \picellipsearc{2.5 1.5}{1.2 1.0}{0 180}
  \picellipsearc{2.5 1.5}{1.8 1.6}{0 180}
  \picputtext{2.5 -0.3}{$K(G)$}
}
\]
on the Suzuki graph $G=G_n$ obtained by gluing into
a circle the ends of the following family of
tangle diagrams (shown for $n=5$):
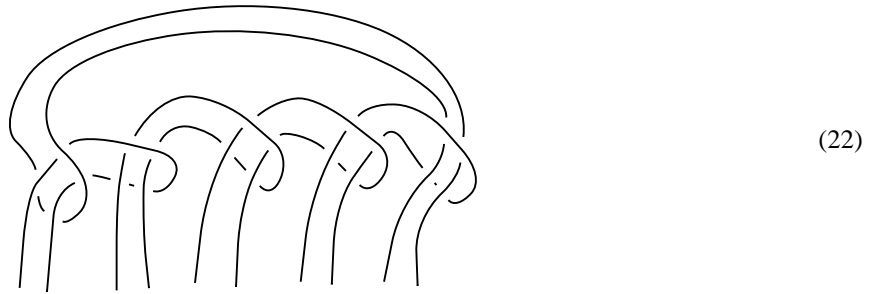
\begin{eqn}\label{dfg}
\epsfsv{5cm}{t1-str5triv}
\end{eqn}
(The signs of the clasps at the arrow hooks are
irrelevant.) Assume now some of the resulting knots
$K_n=K(G_n)$ for $n\ge 2$ has trivial
Jones polynomial. Then smoothing out the crossings
of the strong $n$-trivializer (the hook clasps),
we obtain a $(1+n)$-component link diagram $D=D_n$ with trivial
polynomial. This link diagram is regularly isotopic
to a diagram $D'=D'_n$ that is obtained from $G$ in
\eqref{dfg} by replacing $\diag{6mm}{2}{1}{
  \picline{0 0.3}{2 0.3}
  \picline{0.7 0.3}{0.7 1}
  \picline{1.3 0.3}{1.3 1}
}$ by $
\diag{6mm}{2}{1}{
  \picline{0 0.3}{0.7 0.3}
  \picline{1.3  0.3}{2 0.3}
  \picline{0.7 0.3}{0.7 1}
  \picline{1.3 0.3}{1.3 1}
}
$ and $\epsfsv{2cm}{t1-str5triv2}$ by $
\epsfsv{2cm}{t1-str5triv3}$. But by
a direct check, $D'$ is $A$-adequate. \qed

\begin{rem}
By using proposition \ref{pp21}, one can also show that
the polynomials of $K_n$ are pairwise distinct for different
$n$. While non-triviality of $K_n$ seems provable also using
\cite{HL,AK}, distinctness seems not so easy to see. For
example, I do not know if the genera of $K_n$ are all distinct
(although by \cite{HL} they have an increasing lower estimate).
\end{rem}

The next applications require a longer discussion and are put
into separate sections. In relation to non-triviality, let us
conclude (for now) by saying that from our present result
theorem \reference{tsqn}, and with some help
of Bae and Morton \cite{Bmo}, we proved in \cite{gener}:

\begin{theorem}(\cite{gener})
Non-trivial $k$-almost positive knots have non-trivial Jones
polynomial for $k\le 3$ and non-trivial skein polynomial for $k\le 4$.
\end{theorem}

\section{Minimal positive diagrams\label{mpd}}

% PUT 3BR, WH-DOUBLE, STR N-TRIV KNOT STUFF HERE
A further application relates to remark \reference{rNak}
and Nakamura's problem \cite{Nakamura}. We will exhibit
two types of infinite families of positive knots with
no minimal (crossing number) positive diagram. The first
example of such a knot was given in \cite{posex_bcr}. There
we noticed also that such examples give instances of a
negative solution to the problem of Cromwell \cite{Cromwell}
whether homogeneous knots have homogeneous minimal diagrams
(as do alternating knots).

The approach we use will be set forth in the last sections
for the exhibition of the families of odd crossing number achiral
knots. So this section can be considered a preparation for latter's
proof, though much more substantial argumentation remains.

\subsection{Genus three knots\label{G3}}

\begin{prop}\label{ppp}
There are infinitely many positive knots of genus 3 that
have no minimal positive diagram.
\end{prop}

\proof Consider the following family of positive knot diagrams of genus
3 (where the 6 half-twists in the group are indefinitely augmented).

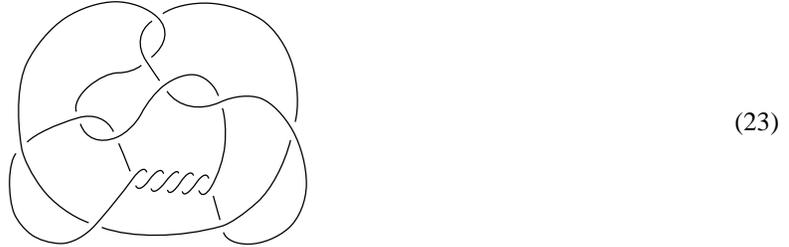
\begin{eqn}\label{dgg}
\epsfsv{3.6cm}{t1-takuji_family}
\end{eqn}

Let $D_n$ be the diagram in this family that has $2n$ crossings,
and $K_n$ the knot it represents.
As the original example $D_6$ in \cite{posex_bcr}, these diagrams
reduce to diagrams $\tl D_n$ on $2n-1$ crossings. The diagrams
$\tl D_n$ are almost positive, and can be checked to be
$B$-adequate. Now it follows from \cite{Thistle} that
then either (a) $c(K_n)=2n-2$, and $K_n$ has an adequate positive
diagram $D_n'$ of $2n-2$ crossings, or (b) $c(K_n)=2n-1$. In
case (b) we will show that no $2n-1$-crossing diagram
$D_n'$ of $K_n$ is positive.

Assume first (a). We use some computation based on the generator
theory explained in \cite{gen1,gen2}; for length reasons we repeat
only a minimum on details concerning the related notions. Two
crossings are \em{$\sim$-equivalent} if they form a reverse clasp up
to flypes. A \em{$\bar t_2'$-move} (or \em{twist}) adds two new
elements/crossings to some $\sim$-equivalence class of a diagram:
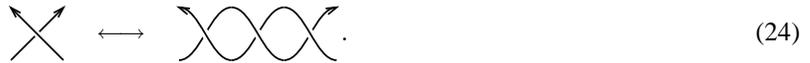
\begin{eqn}\label{move}
\diag{7mm}{1}{1}{
    \picmultivecline{0.12 1 -1.0 0}{0 0}{1 1}
    \picmultivecline{0.12 1 -1.0 0}{1 0}{0 1}
}\quad\llra\quad
\diag{7mm}{3}{2}{
  \picPSgraphics{0 setlinecap}
  \pictranslate{1 1}{
    \picrotate{-90}{
      \lbraid{0 -0.5}{1 1}
      \lbraid{0 0.5}{1 1}
      \lbraid{0 1.5}{1 1}
      \pictranslate{-0.5 0}{
      \picvecline{0.021 1.95}{0 2}
      \picvecline{0.021 -.95}{0 -1}
    }
    }
    }
}
\,.
\end{eqn}
(Such a move may also replace a positive crossing by 3 such.)
A \em{generator} is a diagram all whose $\sim$-equivalence classes
have at most two crossings; its \em{series} is the set of diagrams
obtained from the generator by crossing changes, and successive
flypes and $\bar t_2'$-moves. The importance of these concepts lies in
the fact that knot diagrams of given (canonical) genus decompose
into finitely many series \cite{gen1}. For genus 3, there are
4017 prime alternating generator knots \cite{gen2}.

A minimal crossing 
(adequate positive) diagram $D=D_n'$ lies in the series of some
genus $3$ generator $D'$. Since $K_n$ is prime by \cite{Ozawa},
so is $D$, and $D'$ has even crossing number. Prime genus $3$
generators are determined as explained in \cite{gen2}.

One calculates that $V_1(K_n)=-1$. Now we claim that
$D$ has at most one $\sim$-equivalence class of more than
2 crossings. Assume it has $2$ such classes. Then $G(A(D))'$
would contain two paths \ $\diag{0.6cm}{3}{1}{
  \picline{0 0.5}{3 0.5}
  \lpoint{0 0.5}
  \lpoint{1 0.5}
  \lpoint{2 0.5}
  \lpoint{3 0.5}
}$\ \ of length $\ge 3$. If the paths form parts of a single cycle,
and this is the only cycle in $G(A(D))'$, then the crossings
would be $\sim$-equivalent. So there exist two distinct
cycles, and then $V_2\le -2$, so the diagram cannot depict
any of the $K_n$. (Keep in mind that $V_0=+1$ in positive knot
diagrams.)

So $D$ is obtained from $D'$ by twisting (possibly repeatedly)
at a single crossing. Now adequacy is invariant under $\bar t_2'$
twists, after one twist is performed.

So we select prime even crossing genus 3 generators $D'$ with $V_1=0$
or $V_1=-1$, apply one $\bar t_2'$ twist at any crossing (\em{not}
simultaneously, that is, generating
a new diagram separately for each crossing), and check
whether $V_1=-1$ and the (twisted) diagram is adequate. Such
diagrams happen to occur only for one generator, $10_{154}$.
Now $\bar t_2'$ twists (even before a $\bar t_2'$ twist is applied
previously) are also easily seen to preserve $\bV_1$. But $\bV_1(K_n)
=-3$, while on all twisted diagrams $\bV_1$ equals $\bV_1(10_{154})=
-2$. So case (a) is ruled out.

Now turn to case (b). Assume $D$ is a (minimal crossing number)
positive diagram of $2n-1$ crossings. We consider odd crossing
number generators. Adequacy tests are no longer valid, but still we
can use $V_1$. So we check which positive generators have $V_1=-1$ or
$0$, apply a $\bar t_2'$ twist at any crossing, and check $V_1=-1$. We
can then also check whether $V_2=1$, as it is for $K_n$. Again $V_2$
does not change under further $\bar t_2'$ twists, after one twist is
performed. (Keep in mind that our diagrams $D$ are positive, and so
$\tg A(D)'=0$.)

Only 4 positive generators $D'$ produce diagrams $D$ that satisfy
these conditions. The positive generators all have a fragment that
admits a reducing move called second reduction move in \cite{pos}.
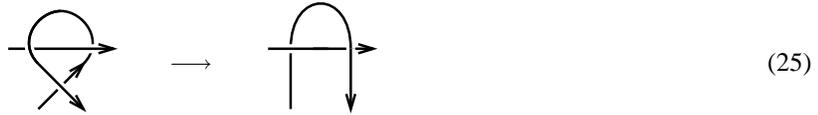
\begin{eqn}\label{second}
\diag{1cm}{2}{1.6}{
  \picline{0.3 1}{1 1}
  \picmulticurve{0.12 1 -1.0 0}{1.0 1.5}{0.7 1.5}{0.4 1.1}{0.7 0.8}
  \picmulticurve{0.12 1 -1.0 0}{1.3 0.8}{1.6 1.1}{1.3 1.5}{1.0 1.5}
  \picmultivecline{0.12 1 -1.0 0}{1 1}{1.7 1}
  \picvecline{0.7 0.2}{1.3 0.8}
  \picmultivecline{0.12 1 -1.0 0}{0.7 0.8}{1.3 0.2}
  \piccurve{1.0 1.5}{0.7 1.5}{0.4 1.1}{0.7 0.8}
}\quad\lra\quad
\diag{1cm}{2}{1.6}{
  \picvecline{1 1}{1.7 1}
  \picmulticurve{0.12 1 -1.0 0}{1.4 1}{1.4 1.8}{0.6 1.8}{0.6 1}
  \picmultivecline{0.12 1 -1.0 0}{1.4 1}{1.4 0.2}
  \picline{0.6 0.2}{0.6 1}
  \picmultiline{0.12 1 -1.0 0}{0.3 1}{1 1}
  \picline{0.9 1}{1.1 1}
}
\end{eqn}
Thus, even after arbitrary twists, such diagrams are not of minimal
crossing number. This contradiction completes the proof. \qed

In \cite{posqp}, we observed that for a positive link the condition
$V_1=0$ is equivalent to being fibered. The observation
made in the proof shows 

\begin{corr}
If for a positive link $V_0V_1=-1$, then the minimal genus
Seifert surface is unique and is obtained from a (positively and
at least twice full-)twisted annulus by iterated Hopf plumbing.
\end{corr}

\proof The canonical surface from a positive diagram is a minimal
genus surface. Then $V_0V_1=-1$ and the above argument show that
such a surface is obtained from a twisted annulus by Murasugi sum
of $(2,\,.\,)$-torus link fiber surfaces. The work of Kobayashi
\cite{Kobayashi} shows then that the surface is a unique minimal
genus surface. Then the argument in the proof of proposition 5.1
of \cite{GHY} shows that one can reduce Murasugi sum with the
$(2,\,.\,)$-torus link surfaces to Hopf plumbing. \qed

\subsection{Fibered knots}

Using refinements and extensions of the preceding arguments
we can settle the problem also for Nakamura's series of fibered
positive knots. This proof is longer, but uses less computation,
and initiates the technique needed later for our main result.

\begin{prop}\label{ppf}
There are infinitely many fibered positive knots that
have no minimal positive diagram.
\end{prop}

\proof Now we modify the initial example $K_6=11_{550}$ to diagrams
$D_{n}$ as in \cite{Nakamura} like
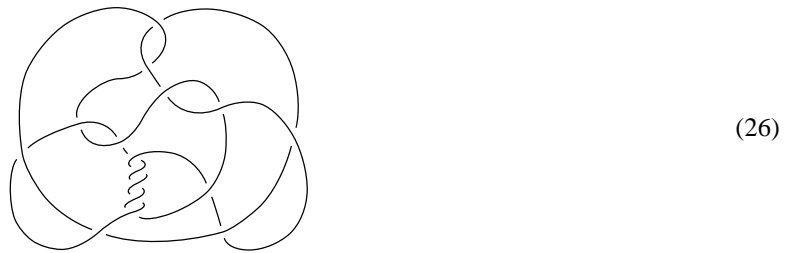
\begin{eqn}\label{dg2}
\epsfsv{3.6cm}{t1-takuji_family_2_}
\end{eqn}
Again let $K_n$ be the knot represented by $D_n$. The diagrams $D_n$
reduce to almost positive $B$-adequate diagrams of $2n-1$ crossings.
Now $K_n$ are fibered, of increasing genus. We have the previous
2 cases.

In case (a) we must have a positive adequate diagram with
5 Seifert circles ($=$loops in the $A$-state). Since $G(A(D))'$
is a tree, it has $4$ edges. All edges are multiple (a single
edges would give a nugatory crossing). Now $V_2=-1$, which
implies that we have 5 intertwined edge pairs. We recall the
intertwining graph $IG(A(D))$ defined by vertices being the (multiple)
edges in $G(A(D))'$, and edges given by intertwined edge pairs.
If we color the regions in the $A$ state of $D$ even-odd then
each edge in $G(A)'$ receives a color. Since only oppositely
colored edges can be intertwined, $IG(A(D))$ is a bipartite graph.
But no bipartite graph on $4$ vertices has $5$ edges. This
rules out case (a).

Now consider case (b) of $c(K_n)=2n-1$. Assume $K=K_n$ for some
$n$ has a positive diagram $D$ of $2n-1$ crossings. $D$ has 6
Seifert circles. We know also that $D$ is not (B-)adequate. By
the previous argument we see that $IG(A(D))$ has 5 vertices and
6 edges, so $IG(A(D))=K_{2,3}$ (the complete bipartite graph
on $2,3$ vertices) is the only option. This means
that $D$ has a single separating Seifert circle, with two
resp. three other Seifert circles, attached (by multiple
crossings) from either side, so that each pair of Seifert
circles from opposite sides is intertwined.

Let the \em{intertwining index} of an edge pair be half the
number of interchanged connections from eithers side of $l$.
For example the intertwining index of
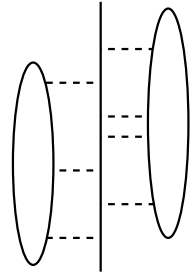
\[
\diag{9mm}{4}{4}{
  \picline{2 0}{2 4}
  {\piclinedash{0.1 0.1}{0.01}
   \picline{1 0.5}{2 0.5}
   \picline{3 1.0}{2 1.0}
   \picline{1 1.5}{2 1.5}
   \picline{3 2.0}{2 2.0}
   \picline{3 2.3}{2 2.3}
   \picline{1 2.8}{2 2.8}
   \picline{3 3.3}{2 3.3}
  }
  \picfilledellipse{1.6 1 x}{0.3 1.5}{}
  \picfilledellipse{2.2 3 x}{0.3 1.7}{}
}
\]
is $3$, and edges are intertwined iff their intertwining index is
$\ge 2$. (We will from now on, to save space, draw in diagrams only
a part of $l$ that contains its basepoints. So the straight line,
that represents $l$, is understood to be closed up.)

We assume in the following that an \em{edge} in $A(D)$ stands for a
possible (but non-empty) collection of parallel traces:
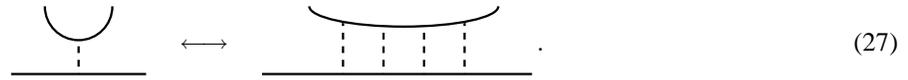
\begin{eqn}\label{parallel}
\diag{9mm}{2}{2.5}{
  \pictranslate{0 0.8}{
    \picrotate{90}{
      \pictranslate{-2 -2}{
        \picline{2 0}{2 2}
        {\piclinedash{0.1 0.1}{0.01}
         \picline{3 1.0}{2 1.0}
        }
        {\picgraycol{1}
          \picfill{
            \picellipsearc{1.0 3 x}{0.5 d}{90 270}
            \picPSgraphics{closepath}
          }
        }
        \picellipsearc{1.0 3 x}{0.5 d}{90 270}
      }
    }
  }
}\quad\llra\quad
\diag{9mm}{4}{2.5}{
  \pictranslate{2 0.8}{
    \picrotate{90}{
      \pictranslate{-2 -2}{
        \picline{2 0}{2 4}
        {\piclinedash{0.1 0.1}{0.01}
         \picline{3 1.0}{2 1.0}
         \picline{3 1.6}{2 1.6}
         \picline{3 2.2}{2 2.2}
         \picline{3 2.8}{2 2.8}
        }
        {\picgraycol{1}
        \picfill{
          \picellipsearc{1.9 3 x}{0.3 1.4}{90 270}
          \picPSgraphics{closepath}
        }
        }
        \picellipsearc{1.9 3 x}{0.3 1.4}{90 270}
      }
    }
  }
}\,.
\end{eqn}
(The term ``edge'' thus assumes a clear separation between
the state $A(D)$ and the graph $G(A(D))$.)
Traces between loops $a$ and $b$ are \em{parallel} if between
their basepoints on both $a$ and $b$ no traces connecting $a$
or $b$ to other loops occur. With this convention we identify,
and do not display several, parallel traces in diagrams. A
connection in $A(D)$, which is the set of all traces or edges
that connect the same two loops, in general decomposes into
several edges.

Now one can obtain $A(D)$ by starting with some single loop,
and then attaching new loops with all their traces.
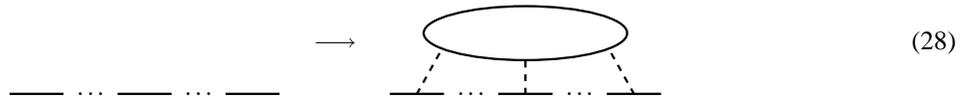
\begin{eqn}\label{*}\label{attach}
\diag{9mm}{4}{2.5}{
  \picline{0.5 0 x}{0.5 0.8 x}
  \picputtext{0.5 1.2 x}{$\dots$}
  \picline{0.5 1.6 x}{0.5 2.4 x}
  \picputtext{0.5 2.8 x}{$\dots$}
  \picline{0.5 3.2 x}{0.5 4.0 x}
}
\quad\lra\quad
\diag{9mm}{4}{2.5}{
  \picline{0.5 0 x}{0.5 0.8 x}
  \picputtext{0.5 1.2 x}{$\dots$}
  \picline{0.5 1.6 x}{0.5 2.4 x}
  \picputtext{0.5 2.8 x}{$\dots$}
  \picline{0.5 3.2 x}{0.5 4.0 x}
  {\piclinedash{0.1 0.1}{0.01}
   \picline{0.5 0.4 x}{1.4 0.9 x}
   \picline{0.5 2.0 x}{1.4 2.0 x}
   \picline{0.5 3.6 x}{1.4 3.1 x}
  }
  \picfilledellipse{2.0 1.4}{0.4 1.5 x}{}
}
\end{eqn}
Since $G(A(D))'$ is a tree, we can assume that we attach the
traces of the new loop to fragments of the same previous loop.

Let 
\begin{eqn}\label{gm}
\gm(D)\,:=\,c(D)+2-v(A(D))-v(B(D))\,.
\end{eqn}
(This number is always even, and non-negative.)
It is easy to see that \eqref{*} never reduces $\gm$. Our diagram
$D$ must, by calculation of $V(K_n)$, have the property that
$\spn V(D)=c(D)-2$, and therefore $\gm(D)\le 4$.

This means that if after some move \eqref{*}, $\gm$ becomes
$\ge 6$, then we can rule out any completion of this loop
insertion out. In particular, it is easy to see that no two
loops have intertwining index $\ge 4$.

Now we distinguish two cases.

\begin{caselist}
\case \label{cs1}
First assume there is a pair of loops of intertwining index 3.
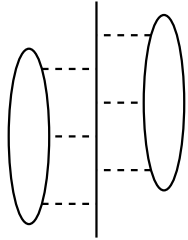
\[
\diag{9mm}{4}{3.5}{
  \picline{2 0}{2 3.5}
  {\piclinedash{0.1 0.1}{0.01}
   \picline{1 0.5}{2 0.5}
   \picline{3 1.0}{2 1.0}
   \picline{1 1.5}{2 1.5}
   \picline{3 2.0}{2 2.0}
   \picline{1 2.5}{2 2.5}
   \picline{3 3.0}{2 3.0}
  }
  \picfilledellipse{1.5 1 x}{0.3 1.3}{}
  \picfilledellipse{2.0 3 x}{0.3 1.3}{}
}
\]
Then $\gm=4$ already after inserting these two loops. So we are
not allowed to increase it when putting in the other two loops.
If we put a loop so as to create a new
region $A$, not connected by a trace from the other side of $l$,
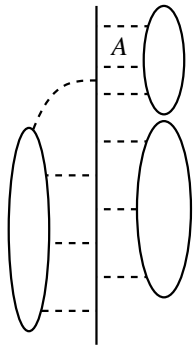
\[
\diag{9mm}{4}{5}{
  \picline{2 0}{2 5}
  {\piclinedash{0.1 0.1}{0.01}
   \picline{1 0.5}{2 0.5}
   \picline{3 1.0}{2 1.0}
   \picline{1 1.5}{2 1.5}
   \picline{3 2.0}{2 2.0}
   \picline{1 2.5}{2 2.5}
   \picline{3 3.0}{2 3.0}
   \picline{3 3.7}{2 3.7}
   \picline{3 4.1}{2 4.1}
   \picline{3 4.7}{2 4.7}
   \piccurve{1 3.0}{1.3 3.9}{1.6 3.9}{2. 3.9}
  }
  \picfilledellipse{1.7 1 x}{0.3 1.5}{}
  \picfilledellipse{2.0 3 x}{0.4 1.3}{}
  \picfilledellipse{4.2 3 x}{0.3 0.8}{}
  \picputtext{2.35 4.4}{$A$}
}
\]
then $A$ must be joined by such a trace after putting one of the
subsequent loops (otherwise the diagram will not be prime, or the
edges in the boundary would remain parallel, in opposition to our
convention). But one observes easily that this loop insertion would
augment $\gm$. 

Thus all $3$ remaining loops can be attached only between
connections from the other side of $l$. They must also be
intertwined among each other. One then observes that, up to
symmetries, the only option is
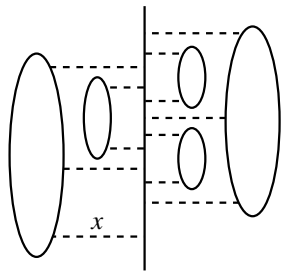
\[
\diag{9mm}{4}{3.9}{
  \picline{2 0}{2 3.9}
  {\piclinedash{0.1 0.1}{0.01}
   \picline{0.6 0.5}{2 0.5}
   \picline{3.4 1.0}{2 1.0}
   \picline{0.6 1.5}{2 1.5}
   \picline{3.4 2.25}{2 2.25}
   \picline{0.6 3.0}{2 3.0}
   \picline{3.4 3.5}{2 3.5}
   \picline{2.7 1.3}{2 1.3}
   \picline{1.3 1.8}{2 1.8}
   \picline{1.3 2.7}{2 2.7}
   \picline{2.7 2.0}{2 2.0}
   \picline{2.7 2.5}{2 2.5}
   \picline{2.7 3.2}{2 3.2}
  }
  \picfilledellipse{1.7 0.4 x}{0.4 1.5}{}
  \picfilledellipse{2.2 3.6 x}{0.4 1.4}{}
  \picfilledellipse{1.3 2.25}{0.2 0.6}{}
  \picfilledellipse{2.7 2.85}{0.2 0.45}{}
  \picfilledellipse{2.7 1.65}{0.2 0.45}{}
  \picputtext{1.3 0.7}{$x$}
}\quad.
\]
A direct check shows that this diagram is adequate if
and only if $x$ is a multiple edge. If $x$ is single,
then in $B(D)$ we have only one trace that connects the
same loop. By \cite{Bmo}, we conclude then that the extreme
$B$-term in the bracket vanishes. Then, since $\gm=4$,
\[
\spn V(D)\,<\,\frac{c(D)+v(A(D))+v(B(D))-2}{2}\,=\,c(D)-2\,,
\]
with a contradiction. So the argument in this case is complete.

\case Now assume all loops have intertwining index 2. Then our
$A$-state is obtained from 
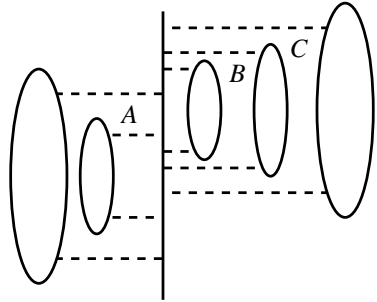
\[
\diag{11mm}{4}{3.5}{
  \picline{2 0}{2 3.5}
  {\piclinedash{0.1 0.1}{0.01}
   \picline{0.5 0.5}{2 0.5}
   \picline{0.5 2.5}{2 2.5}
   \picline{1.2 1.0}{2 1.0}
   \picline{1.2 2.0}{2 2.0}
   \picline{4.2 1.3}{2 1.3}
   \picline{4.2 3.3}{2 3.3}
   \picline{3.3 1.6}{2 1.6}
   \picline{3.3 3.0}{2 3.0}
   \picline{2.5 1.8}{2 1.8}
   \picline{2.5 2.8}{2 2.8}
  }
  \picfilledellipse{0.5 1.5}{0.34 1.3}{}
  \picfilledellipse{1.2 1.5}{0.2 0.7}{}
  \picfilledellipse{4.2 2.3}{0.34 1.3}{}
  \picfilledellipse{3.3 2.3}{0.2 0.8}{}
  \picfilledellipse{2.5 2.3}{0.2 0.6}{}
  \picputtext{1.6 2.25}{$A$}
  \picputtext{2.9 2.75}{$B$}
  \picputtext{3.65 3.05}{$C$}
}
\]
by adding edges. If we do not create equivalent edges, then we can
attach edges that connect a loop to one of regions $A$, $B$ or $C$.
(To connect to two, or all three of $A,B,C$ would violate \eqref{gm}.)
In all three cases, which we denote after the region concerned,
we can add at most two edges $x$ and/or $y$.
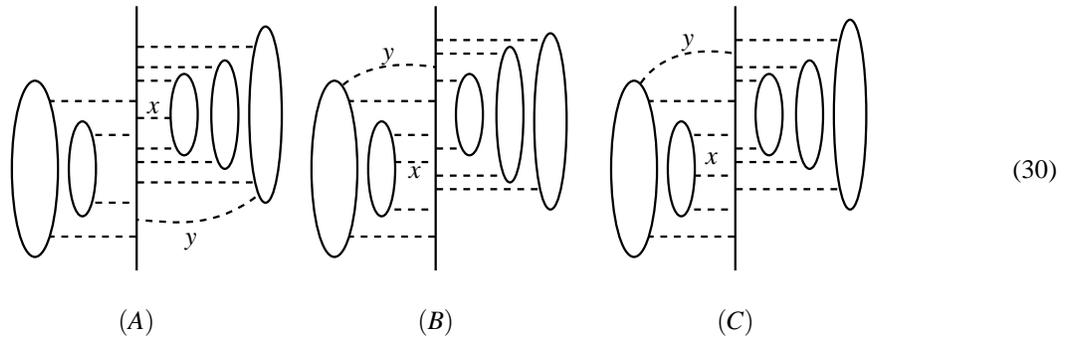
\begin{eqn}\label{ABC}
\begin{array}{c@{\quad}c@{\quad}c}
\diag{9mm}{4}{3.9}{
  \picline{2 0}{2 3.9}
  {\piclinedash{0.1 0.1}{0.01}
   \picline{0.5 0.5}{2 0.5}
   \picline{0.5 2.5}{2 2.5}
   \picline{1.2 1.0}{2 1.0}
   \picline{1.2 2.0}{2 2.0}
   \picline{3.9 1.3}{2 1.3}
   \picline{3.9 3.3}{2 3.3}
   \picline{3.3 1.6}{2 1.6}
   \picline{3.3 3.0}{2 3.0}
   \picline{2.7 1.8}{2 1.8}
   \picline{2.7 2.25}{2 2.25}
   \picputtext{2.25 2.4}{$x$}
   \picline{2.7 2.8}{2 2.8}
   \piccurve{3.9 1.2}{3.2 0.6}{2.6 0.7}{2 0.75}
   \picputtext{2.8 0.47}{$y$}
  }
  \picfilledellipse{0.5 1.5}{0.34 1.3}{}
  \picfilledellipse{1.2 1.5}{0.2 0.7}{}
  \picfilledellipse{3.9 2.3}{0.24 1.3}{}
  \picfilledellipse{3.3 2.3}{0.2 0.8}{}
  \picfilledellipse{2.7 2.3}{0.2 0.6}{}
}
 &
\diag{9mm}{4}{3.9}{
  \picline{2 0}{2 3.9}
  {\piclinedash{0.1 0.1}{0.01}
   \picline{0.5 0.5}{2 0.5}
   \picline{0.5 2.5}{2 2.5}
   \picline{1.2 0.9}{2 0.9}
   \picline{1.2 2.0}{2 2.0}
   \picline{3.7 1.2}{2 1.2}
   \picline{3.7 3.4}{2 3.4}
   \picline{3.1 1.4}{2 1.4}
   \picline{3.1 3.2}{2 3.2}
   \picline{2.5 1.8}{2 1.8}
   \picline{2.5 2.8}{2 2.8}
   \picputtext{1.69 1.45}{$x$}
   \picline{1.2 1.6}{2 1.6}
   \picputtext{1.3 3.2}{$y$}
   \piccurve{0.5 2.6}{0.8 3.0}{1.4 3.1}{2 3.0}
  }
  \picfilledellipse{0.5 1.5}{0.34 1.3}{}
  \picfilledellipse{1.2 1.5}{0.2 0.7}{}
  \picfilledellipse{3.7 2.2}{0.24 1.3}{}
  \picfilledellipse{3.1 2.3}{0.2 1.0}{}
  \picfilledellipse{2.5 2.3}{0.2 0.6}{}
} &
\diag{9mm}{4}{3.9}{
  \picline{2 0}{2 3.9}
  {\piclinedash{0.1 0.1}{0.01}
   \picline{0.5 0.5}{2 0.5}
   \picline{0.5 2.5}{2 2.5}
   \picline{1.2 0.9}{2 0.9}
   \picline{1.2 2.0}{2 2.0}
   \picline{3.7 1.2}{2 1.2}
   \picline{3.7 3.4}{2 3.4}
   \picline{3.1 1.6}{2 1.6}
   \picline{3.1 3.0}{2 3.0}
   \picline{2.5 1.8}{2 1.8}
   \picline{2.5 2.8}{2 2.8}
   \picputtext{1.65 1.65}{$x$}
   \picline{1.2 1.4}{2 1.4}
   \picputtext{1.3 3.4}{$y$}
   \piccurve{0.5 2.6}{0.8 3.2}{1.4 3.3}{2 3.2}
  }
  \picfilledellipse{0.5 1.5}{0.34 1.3}{}
  \picfilledellipse{1.2 1.5}{0.2 0.7}{}
  \picfilledellipse{3.7 2.3}{0.24 1.4}{}
  \picfilledellipse{3.1 2.3}{0.2 0.8}{}
  \picfilledellipse{2.5 2.3}{0.2 0.6}{}
}\\[16mm]
\\
(A) & (B) & (C) 
\end{array}
\end{eqn}
It is easy to see that adding only $x$ or only $y$ is equivalent,
and that the cases $B$ and $C$ are equivalent up to mutations.

\begin{caselist}
\case Neither $x$ nor $y$ is added. By direct verification, $D$
is adequate.

\case Both $x$ and $y$ are added. As in case \reference{cs1},
we observe the $B$-state of $D$. There are now two traces
that connect the same loops, which we call \em{inadequate},
but one has no \em{linked pair}, i.e. a pair of inadequate
traces at the same loop in the mutual position
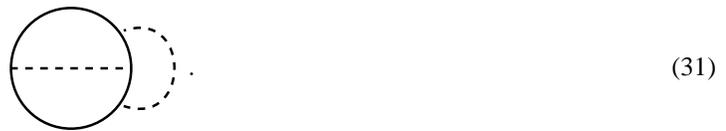
\begin{eqn}\label{lpt}
\diag{1cm}{2.5}{2.}{
  {\piclinedash{0.1 0.1}{0.01}
   \piccurve{1.7 0.5}{2.6 0.2}{2.6 1.8}{1.7 1.5}
  }
  \picfilledcircle{1 1}{0.8}{}
  {\piclinedash{0.1 0.1}{0.01}
   \picline{0.2 1}{1.8 1}
  }
}\,.  
\end{eqn}
Then again by replacing edges by (possibly multiple) traces, either
we have an adequate diagram, or if not, by \cite{Bmo} the extreme
$B$-state term in the bracket vanishes. Since again $\gm=4$,
we obtain $\spn V(D)<c(D)-2$.

\case The cases when, say, $x$ (but not $y$) is added are
more problematic. 

\begin{caselist}
\case\label{Acs} Consider first a trace connecting $A$. Name the other
crossings as indicated on the left of \eqref{four}. The resulting
link diagram $D'$ (if all edges are not yet replaced by
parallel multiple traces) is shown on the right. 
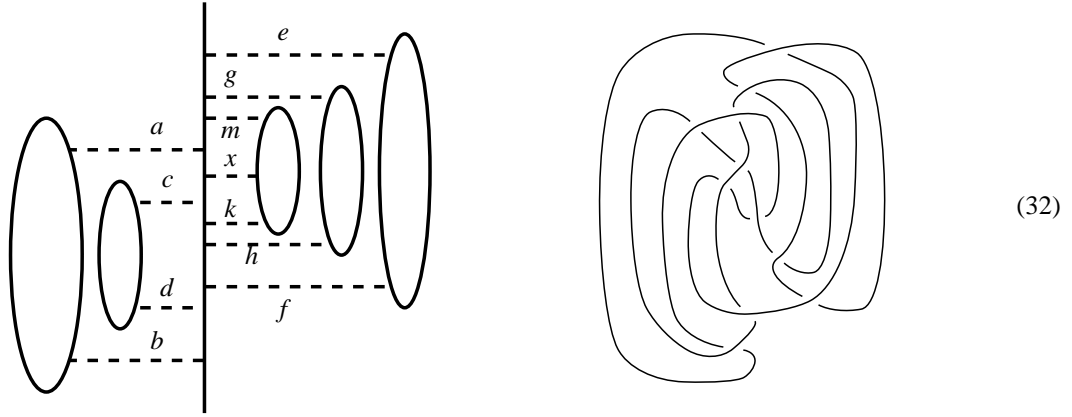
\begin{eqn}\label{four}
\diag{14mm}{4.3}{3.9}{
  \picline{2 0}{2 3.9}
  {\piclinedash{0.1 0.1}{0.01}
   \picline{0.5 0.5}{2 0.5}
   \picline{0.5 2.5}{2 2.5}
   \picline{1.2 1.0}{2 1.0}
   \picline{1.2 2.0}{2 2.0}
   \picline{3.9 1.2}{2 1.2}
   \picline{3.9 3.4}{2 3.4}
   \picline{3.3 1.6}{2 1.6}
   \picline{3.3 3.0}{2 3.0}
   \picline{2.7 1.8}{2 1.8}
   \picline{2.7 2.25}{2 2.25}
   \picline{2.7 2.8}{2 2.8}
   \picputtext{2.25 2.4}{$x$}
   \picputtext{1.65 2.2}{$c$}
   \picputtext{1.65 1.2}{$d$}
   \picputtext{1.55 0.7}{$b$}
   \picputtext{1.55 2.7}{$a$}
   \picputtext{2.75 3.6}{$e$}
   \picputtext{2.75 1.0}{$f$}
   \picputtext{2.45 1.5}{$h$}
   \picputtext{2.25 3.2}{$g$}
   \picputtext{2.25 1.95}{$k$}
   \picputtext{2.25 2.67}{$m$}
  }
  \picfilledellipse{0.5 1.5}{0.34 1.3}{}
  \picfilledellipse{1.2 1.5}{0.2 0.7}{}
  \picfilledellipse{3.9 2.3}{0.24 1.3}{}
  \picfilledellipse{3.3 2.3}{0.2 0.8}{}
  \picfilledellipse{2.7 2.3}{0.2 0.6}{}
}\rx{2cm}
\epsfsv{4.8cm}{t1-takuji_family_2}
\end{eqn}
Then (for single multiplicity edges)
in the $B$ state the trace of $x$ (that connects the
same loop), is linked with $a,b,c,d$. Latter four are
not linked among each other. So we must deal with the
situation where $x$ is a single crossing (and one of
$a,b,c,d$ is too). 

Now we require some computation.
Let $x=(x_1,\dots,x_{10})\in\bN_+^{10}$ with $x_i\ge 1$ and assign
to $x$ a diagram $D_x$ obtained from \eqref{four} by replacing
$a$ by $x_1$, $b$ by $x_2$, \dots, $m$ by $x_{10}$ parallel edges,
as in \eqref{parallel}. Define $\tl x\in\{1,2\}^{10}$ by
$(\tl x)_i=\min(x_i,2)$.

We generate the 1024 diagrams $D_{\tl x}$ obtained by taking $\tl
x\in\{1,2\}^{10}$, and calculate their Alexander polynomial. (In fact,
288 diagrams are enough, as we can assume up to flypes that $(2,1)
\not\in\{(a,b),(c,d),(e,f),(g,h)\}$, and $1\in\{a,b,c,d\}$.) 
Recall our normalization of $\Dl$ by the skein relation \eqref{srel},
or equivalently by the conditions $\Dl(1)=1$ and $\Dl(t)=\Dl(1/t)$.
Let $\bD_i=[\Dl]_{\Md\Dl-i}$. Then for all $D_{\tl x}$ we find
\begin{eqn}\label{strr}
\bD_0=\bD_1=1\,.
\end{eqn}
(That $\bD_1=1$ and $2\Md\Dl=1-\chi(D_{\tl x})$
is clear by theory \em{a priori}.) A direct skein argument
using \eqref{srel} shows then that for all $x\in\bN_+^{10}$ we
have $\bD_2(D_x)=\bD_2(D_{\tl x})$. Now again it follows by
computation that for our knots $\bD_2(K_n)\in \{-6,-7\}$. Contrarily
\begin{eqn}\label{str2}
\bD_2(D_x)<-7
\end{eqn}
unless all entries $x_i$ of $x$ are $1$ except at most one $i=
5,\dots,10$, that is, for at least nine $i\in\{1,\dots,10\}$
we have $x_i=1$. Because of \eqref{strr} it is enough to check
\eqref{str2} for the finitely many $x=\tl x$.

So we need to consider only diagrams $D_x$ where $x$ has at most
one entry $x_i>1$. But since $D'=D_{(1,\dots,1)}$ in \eqref{four}
is a 3-component link diagram, such $D_x$ will have at least two
components.

\case Finally, consider a trace connecting $B$. We name the edges
in the diagram \eqref{ABC} (C) as indicated below. (For a moment
think of $c'$ as not being there.)
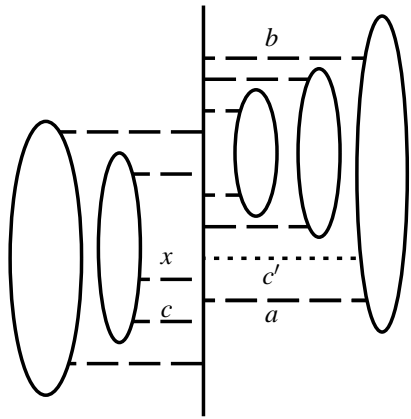
\[
\diag{14mm}{4}{3.9}{
  \picline{2 0}{2 3.9}
  {\piclinedash{0.3 0.1}{0.01}
   \picline{0.5 0.5}{2 0.5}
   \picline{0.5 2.7}{2 2.7}
   \picline{1.2 0.9}{2 0.9}
   \picline{1.2 2.3}{2 2.3}
   \picline{3.7 1.1}{2 1.1}
   {\piclinedash{0.05 0.08}{0.01}
    \picline{3.7 1.5}{2 1.5}
   }
   \picline{3.7 1.1}{2 1.1}
   \picline{3.7 3.4}{2 3.4}
   \picline{3.1 1.8}{2 1.8}
   \picline{3.1 3.2}{2 3.2}
   \picline{2.5 2.1}{2 2.1}
   \picline{2.5 2.9}{2 2.9}
   \picputtext{1.65 1.5}{$x$}
   \picputtext{1.65 1.0}{$c$}
   \picputtext{2.65 0.95}{$a$}
   \picputtext{2.65 1.35}{$c'$}
   \picputtext{2.65 3.6}{$b$}
   \picline{1.2 1.3}{2 1.3}
   % \picputtext{1.3 3.4}{$y$}
   % \piccurve{0.5 2.6}{0.8 3.2}{1.4 3.3}{2 3.2}
  }
  \picfilledellipse{0.5 1.5}{0.34 1.3}{}
  \picfilledellipse{1.2 1.6}{0.2 0.9}{}
  \picfilledellipse{3.7 2.3}{0.24 1.5}{}
  \picfilledellipse{3.1 2.5}{0.2 0.8}{}
  \picfilledellipse{2.5 2.5}{0.2 0.6}{}
}
\]
It suffices to consider only inadequate diagrams $D$ where $B(D)$ has
a linked pair \eqref{lpt} of traces. Now for a simple edge $x$ in $B
(D)$, i.e. one consisting of only one crossing trace, this trace is
linked with (the trace of) the edge $a$ or/and $b$, if these are
simple too,
and no further (potential) linked pairs occur. Thus we need to deal
with the cases when $x$ is single, and one of $a$ and $b$ is. Up
to a flype, we may assume that $a$ is single. Then using
Reidemeister III moves, one can move the crossings in $c$
into $c'$, and thus returns to case \reference{Acs}.
So our proof is now complete.
\end{caselist}
\end{caselist}
\qed
\end{caselist}

\section{3-Braids\label{S3B}}

\subsection{Semiadequacy of 3-braids}

Let $\bt$ be a braid word. We write
\[
\bt\,=\,\prod_{j=1}^k\sg_{p_j}^{q_j}
\]
with $\sg_i$ the Artin generators, $p_j\ne p_{j+1}$ and 
$q_j\ne 0$. We will deal mostly with 3-braids, since here
the most interesting applications are possible. For $3$-braids
$p_j$ are interchangingly $1$ and $2$. Thus $\bt$ is (at least
up to conjugacy) determined by the vector $(q_j)$ (up to cyclic
permutations), which is called the \em{Schreier vector}. It has
even length (unless the length is 1). \em{We will assume for the
rest of this section that the Schreier vector has length $>2$.}
(The braids $\sg_{1}^k\sg_2^l$ are easy to deal with directly.)
Any braid word $\bt$ gives a link diagram $\hat\bt$ under closure,
as a braid gives a link. Let $\dl=\sg_2\sg_1\sg_2=\sg_1\sg_2\sg_1$
be the square root of the center generator of $B_3$.

The following is shown by a careful observation:

\begin{lemma}\label{lemsaq}
A diagram of closed $3$-braid word $\bt$ (with Schreier vector
of length $>2$) is $A$-adequate if and only if it satisfies one of
the following two conditions:
\def\theenumi{\arabic{enumi}}
\begin{enumerate}
\item it is positive, or
\item 
\def\theenumii{\arabic{enumii}}
\begin{enumerate}
\item
it does not contain $\dl^{-1}=\sg_1^{-1}\sg_2^{-1}\sg_1^{-1}=
\sg_2^{-1}\sg_1^{-1}\sg_2^{-1}$ as subword and 
\item positive
entries in the Schreier vector are isolated (i.e., entries
cyclically before and after positive entries are negative).
\qed
\end{enumerate}
\end{enumerate}
\end{lemma}

\begin{corr}\label{caq}
A diagram of closed $3$-braid word $\bt$ is adequate if and only if
\def\theenumi{\arabic{enumi}}
\begin{enumerate}
\item it is positive and does not contain $\dl$ as subword, or
\item it is negative and does not contain $\dl^{-1}$ as subword, or
\item it is alternating.
\qed
\end{enumerate}
\end{corr}

\begin{rem} 
For positive/negative braids to contain $\dl^{\pm 1}$ is
equivalent to the Schreier vector containing no entries $\pm 1$.
\end{rem}

\begin{corr}
If a positive word $\bt$ contains $\dl$ as subword, then it is not
conjugate to a positive word that does not contain $\dl$ as subword.
\end{corr}

\proof The diagram $\hat\bt$ is adequate, and by \cite{Thistle} any
other same crossing number diagram must be adequate too. But if
$\bt'$ contains no $\dl$, then $\hat\bt'$ is not $B$-adequate.
\qed

\begin{theo}\label{tq}
Let $\bt$ be a $3$-braid word. Then the following are equivalent:
\def\theenumi{\alph{enumi}}
\def\labelenumi{(\theenumi)}
\begin{enumerate}
\item $\hat\bt$ is semiadequate.
\item $\bt$ has minimal length in its conjugacy class.
\item Some $\gm\in\{\bt,\bar\bt\}$ satisfies one of the conditions
enumerated in lemma \ref{lemsaq}.
\end{enumerate}
\end{theo}

%\ begin{rem}
Before the proof we make some remarks on consequences and related
problems.
The implication (c) $\Ra$ (b) says in particular that if a word $\bt$
has minimal length up to cyclic permutations, it has minimal
length up to conjugacy. We say then it has \em{minimal
conjugacy length}.
% \end{rem}

% \begin{rem}
There is a problem due to Stallings (see \cite[problem 1.8]{Kirby})
whether minimal length braid words are closed under end
extension (replacing a final $\sg_i^{\pm 1}$ by $\sg_i^{\pm k}$
for $k>1$). A natural modification of this problem is 

\begin{question}\label{qmn}
Are minimal conjugacy length words closed under (general)
extension (replacing any $\sg_i^{\pm 1}$ by $\sg_i^{\pm k}$)?
\end{question}

We can confirm this for 3-braids. To describe minimal length
or minimal conjugacy length words in general braid groups seems
extremely difficult, though.
% \end{rem}

% \begin{rem}
A further version of question \reference{qmn} is to replace
minimal conjugacy length by minimal length among Markov
equivalent braids (a) of the same or (b) of arbitrary strand
number. Using theorem \reference{tq} and \cite{BirMen} we
find that version (a) holds for 3-braids. A computer
quest yielded a negative answer to part (b) for 5-braids.

\begin{exam}
The knots $14_{42676}$ and $14_{42683}$ have
14 crossing 5-braid representations, which after twofold
extension give rise to 16 crossing braids that correspond
to 8 different 15 crossing knots.
% 15   102467 
% 15   153829 
% 15   156957 
% 15   204698 
% 15   218067 
% 15   219651 
% 15   224001 
% 15   224179 
Since these knots have 4-braid 15 crossing representations,
and the braid index of $14_{42676}$ and $14_{42683}$ is 5, we
see that neither minimal length, nor strand number minimality
are generally preserved under extension. The examples also
yield a negative answer to a similar problem: is the property
of a braid word to realize the crossing number of its closure
knot or link (i.e. the closed braid diagram to be of minimal
crossing number) preserved under extension? (Here somewhat
simpler examples were found previously.)
% The knot $13_{8850}$ has a $4$-braid word
% ${1,1,2,-3,2,  -3,-2,   -1,2,-1,2,-3,2}$. Applying
% extension twice, we can obtain the length-15-word 
% ${1,1,2,-3,2,2,-3,-2,-2,-1,2,-1,2,-3,2}$ that represents 
% the knot $14_{43087}$. 
\end{exam}

The proof of theorem \reference{tq} is mainly contained in
the following two lemmas.

\begin{lemma}\label{lem1}
If $\bt$ has minimal length up to cyclic permutations and
contains $\dl$ (resp. $\dl^{-1}$), then $\bt$ is positive (resp.
negative).
\end{lemma}

\proof It suffices to prove in the positive case.
The proof is given constructively by direct braid 
word transformation. We assume $\bt$ is not positive and
find a word reduction. To accomplish this, apart from braid 
relations $\sg_{1}^{\pm_1 1}\sg_{2}^{\pm_2 1}\sg_{1}^{\pm_2 1}
\lfra\sg_{2}^{\pm_2 1}\sg_{1}^{\pm_2 1}\sg_{2}^{\pm 1}$
(with the $\pm$ signs chosen consistently with their indices
and possible interchange of $\sg_{1,2}$) we use the
``wave'' moves $\sg_{1}^{\mp 1}\sg_2^k\sg_1^{\pm 1}
\lfra\sg_{2}^{\pm 1}\sg_1^k\sg_2^{\mp 1}$.

We assume that a negative letter is on the right of $\dl$ in $\bt$.
The case when it is on the left is obtained by reversing the
order of the letters in the word. The sequence of words below
gives a sample word reduction (an entry $\pm i$ for $i>0$ means
$\sg_{i}^{\pm 1}$). More formally one uses that one can `pull'
letters through $\dl$, interchanging indices $1\lfra 2$,
and can bring $\dl$ close to the negative letter, where it
cancels. Here is a simple example:\\
\def\QED{\hfill{\normalsize\qed}}
\hbox to 0.9\textwidth{\vbox{\small
\begin{alltt}
1 2 1 1 1 2 2 2 1 1 1 2 2 2 -1
1 2 1 1 1 2 2 2 1 1 -2 1 1 1 2
1 2 1 1 1 2 2 -1 2 2 1 1 1 1 2
1 2 1 1 -2 1 1 2 2 2 1 1 1 1 2
1 -1 2 2 1 1 1 2 2 2 1 1 1 1 2
     2 2 1 1 1 2 2 2 1 1 1 1 2\QED
\end{alltt}%
}}%\hfill$\Box$

\begin{lemma}\label{lem2}
If $\bt$ has Schreier vector that has length-2-subsequences
of both positive and negative entries, then $\bt$ is word reducible.
\end{lemma}

\proof Take a positive and negative length-2-subsequence $p,q$ of
the Schreier vector that have minimal distance to each other. Assume
w.l.o.g. $p$ comes before $q$ and is of the form $\sg_1^k\sg_2^l$.
Then the entries in the Schreier vector between (and including)
the second entry for $p$ and the first entry for $q$ must alternate
in sign. Then $q$ is of the form $\sg_1^{-k'}\sg_2^{-l'}$.

One brings then the subword `1 2' from the pair $p$ of consecutive
positive entries close to $q$
using $\sg_1\sg_2^k\sg_1^{-1}=\sg_2^{-1}\sg_1^k
\sg_2$ and $\sg_1\sg_2\sg_1^{-k}=\sg_2^{-k}\sg_1\sg_2$, and then
cancels.  Again we give a sample word reduction.\\
\hbox to 0.9\textwidth{\vbox{
{\small
\begin{alltt}
1 1 1 2 2 2 -1 -1 -1     2 2 2 -1 -1 -1      2 2 2 -1 -1 -1 -2 -2 -2
1 1 -2 1 1 1 2 -1 -1     2 2 2 -1 -1 -1      2 2 2 -1 -1 -1 -2 -2 -2
1 1 -2 1 1     -2 -2 1 2 2 2 2 -1 -1 -1      2 2 2 -1 -1 -1 -2 -2 -2
1 1 -2 1 1     -2 -2 -2 1 1 1 1 2 -1 -1      2 2 2 -1 -1 -1 -2 -2 -2
1 1 -2 1 1     -2 -2 -2 1 1 1     -2 -2  1 2 2 2 2 -1 -1 -1 -2 -2 -2
1 1 -2 1 1     -2 -2 -2 1 1 1     -2 -2 -2 1 1 1 1  2 -1 -1 -2 -2 -2
1 1 -2 1 1     -2 -2 -2 1 1 1     -2 -2 -2 1 1 1 -2 -2 1  2 -2 -2 -2
1 1 -2 1 1     -2 -2 -2 1 1 1     -2 -2 -2 1 1 1 -2 -2 1       -2 -2\QED
\end{alltt}%
}}}%\hfill \qed

\proof[of theorem \reference{tq}]
(a) $\Ra$ (b). If $\bt$ is not of minimal length in its conjugacy
class, $\hat\bt$ is regularly isotopic to a diagram of smaller
crossing number. This contradicts the result of Thistlethwaite
\cite[corollary 3.1]{Thistle}.

(b) $\Ra$ (c).
This is accomplished by lemmas \reference{lem1} and \reference{lem2}.

(c) $\Ra$ (a). Use lemma \ref{lemsaq}. \qed

% \subsection{Jones polynomial}
\subsection{Some applications}

Now we put together some consequences of the preceding treatment of
semiadequacy of $3$-braids. The most obvious one is

\begin{corr}\label{C3SQ}
$3$-braid links are semiadequate, and so have
non-trivial Jones polynomials up to units. \qed
\end{corr}

The material in this section was originally motivated by the
question whether there is a non-trivial $3$-braid knot with
trivial Jones polynomial. I was initially aware of Birman's
claim of a negative answer (see \cite{Birman,3br}). However,
Birman did not give any reference or further comment, and I
knew of no proof. Only after the present proof was obtained,
I was pointed to the paper \cite{Takahashi}. There explicit
(but very complicated) formulas for polynomials of closed
$3$-braids are given, and the Jones polynomial non-triviality
result (for knots) is claimed after an 18-page long calculation.

Apart from dealing with links too, our proof should be considered
more conceptual. We thereby also extend the result on the
non-triviality of the skein polynomial of \cite{3br}, although
our approach here is quite different from \cite{Birman} or
\cite{3br}. % This does not mean that it should be kept isolated.
% In fact, by using the work in \cite{3br}, we will later prove
% that there are only finitely many closed 3-braids with the same
% Jones polynomial.

The feature of excluding trivial polynomials up to units should
be contrasted to the existence of a (2-component) 4-braid link
with trivial polynomial up to units. Such a link, the closure of
$(\sg_2\sg_1\sg_3\sg_2)^2\sg_1^3\sg_3^{-3}$, is given in \cite{EKT}.
It also relates to the problem of the (possible) faithfulness of the
Burau representation. It is known that the Burau matrix determines
the Jones polynomial for $3$- or $4$-braids. So a non-faithful Burau
representation on $4$-braids will (in practice; see \cite{Bi}) imply
a $4$-braid knot or link $L$ with trivial polynomial, but the above
example cautions about a possible reverse conclusion. To put our work
into this context, it seems more likely (though not easy) to develop
a combinatorial proof that no $L$ exists and deduce that Burau is
faithful (as we can now for 3-braids, providing a numerous proof)
rather than the other way around.

\begin{rem}\label{remx}
There seems no easy way to decide what of $A$- or $B$-adequate
3-braid representations exist for a given link, starting from
another semiadequate diagram of it. For example $5_2$ and $10_{161}$
are $A$-adequate, but have only $B$-adequate 3-braid representations.
Latter knot thus shows that a weakly adequate 3-braid link may not be
adequate, while former knot shows that an adequate 3-braid knot may
not have an adequate 3-braid (word) representation.
\end{rem}

Contrarily, if
a link has 3-braid representations, one of which is $A$-adequate
and one of which is $B$-adequate, then it also has an adequate
representation. This follows, because $A$- and $B$-adequacy always
appear in minimal length 3-braid words. Using this remark, we can 
extend a result of Birman-Menasco \cite{BirMen}.

\begin{defi}
We call an oriented link \em{positively/negatively (orientedly)
amphicheiral, or achiral}, if it is isotopic to its mirror image
with the orientation of no/all components reversed, and
components (possibly and arbitrarily) permuted.
We call the link \em{unorientedly achiral}, if it is isotopic
to its mirror image with the orientation of some (possibly no,
but also not necessarily all) components reversed (and components
possibly permuted).
\end{defi}

It is clear that oriented achirality (of some sign) implies
unoriented one, and for \em{knots} both notions coincide.
The Hopf link is an example of a
link which is unorientedly achiral, but not orientedly so.
% Birman-Menasco described the orientedly achiral 3-braid links
% in \cite{BirMen}. We can settle now also the unorientedly achiral
% case.

\begin{theorem}\label{TU3}
Let $L$ be an unorientedly achiral 3-braid link. Then $L=\hat\bt$,
where $\bt\in B_3$ is either (a) an alternating braid with Schreier
vector admitting a dihedral (anti)symmetry, or (b) $\bt=\sg_1^{k}
\sg_2^{l}$ and $|k|,|l|\le 2$.
\end{theorem}

\proof We assume $\bt$ is a minimal length word up to conjugacy.
The cases that $\bt$ has a Schreier vector of length two or less
are easy to deal with, and lead to case (b). So we assume the
Schreier vector of $\bt$ is of length at least $4$.

By the above remark,
semiadequacy appears in minimal length braid words.
We know that any minimal length word is $A$-adequate or
$B$-adequate, but since $L$ is achiral and semiadequacy is
independent on component orientation, $\bt$ must be both,
i.e. it must be adequate.

In case $\bt$ is alternating, we can apply \cite{MenThis}.
An easy observation shows that the diagram is determined
by the Schreier vector up to dihedral moves (cyclic permutations
and reversal of order), and the only possible flypes in
a diagram of a closed 3-braid occur at Schreier vector of length
$4$, with an entry $\pm 1$ (see \cite{BirMen}). In that case, however,
the flypes just reverse the orientation of (all components
of) the link. So we obtain no new symmetries. We arrive at case (a).

It remains to exclude the case that $\bt$ is not alternating.
In that event, we know that $\bt$ is positive or negative,
and has no $\Dl$ as subword. Positive (or negative) links
are not orientedly achiral, so the possibility that $L$ is a knot
or that the achirality is oriented is easily ruled out.

It is clear that all components of an unorientedly achiral link
are achiral knots, or mirror images in pairs. Latter option does
not lead to anything new for closed 3-braids, so we ignore it.
Then if $L=\hat\bt$ is a 2-component link, two of the strands
$x,y$ of $\bt$ form (under closure) a positive achiral knot. This
knot must be trivial, and so $x,y$ have only one common crossing.
Then an easy observation shows that $\bt$ contains $\Dl$ as subword,
a contradiction.

So assume that $L=\hat\bt$ is a 3-component link. Then
a linking number argument shows that two of the components
of $L$ must have no common crossing. This means that in
the Schreier vector all entries are even and (up to
mirroring) positive. Then reversal of component
orientation makes the (say) positive diagram $D=\hat\bt$
into a negative diagram $D'$. Now $D'$ must be orientedly
isotopic to the mirror image of $D$. We know that in
positive and negative diagrams the canonical Seifert
surface has minimal genus. So to find a contradiction,
we just need to count the Seifert circles in $D$ in $D'$
and see that they are not equally many. \qed

In Birman-Menasco \cite{BirMen} the oriented achirality result
was obtained (which leads to case (a) above), so our result is a
slight extension of theirs. However, the proof in \cite{BirMen}
requires their
extremely involved proof of the classification of closed 3-braids.
It is therefore very useful to obtain simpler proofs of at least
consequences of their work. So more important than the drop of
orientations here is the complete bypassing of the method in
\cite{BirMen}, even though the insight in \cite{BirMen} motivated
the present proof, and with \cite{MenThis} we make use
of another quite substantial result.

With a similar argument, we easily see

\begin{corr}\label{c5.4}
A non-split $3$-braid \em{knot} with $T_1=0$ is a closed
positive or negative 3-braid.
\end{corr}

\proof Take a minimal word length 3-braid representation
$\bt$. The closure of $\bt$ is a semiadequate diagram, 
since $V_1=\bV_{1}=0$, it must be positive. \qed

The case of \em{links} is quite different, and we will describe
them now, since
in the following it will be somewhat important to understand
how one can make 3-braids more complicated without altering $T_1$.

\begin{prop}\label{pexc}
Assume $\bt$ is an $A$-adequate $3$-braid word and $V_1=0$. Then,
up to interchange $\sg_1\lfra\sg_2$,
\def\theenumi{\alph{enumi}}
\def\labelenumi{(\theenumi)}
\begin{enumerate}
\item $\bt$ is positive, or
\item $\bt=\sg_1^{-2}\sg_2^{k}$ for $k>0$, or
\item $\bt$ can be written as
\begin{eqn}\label{btkl_}
\bt(k_1,l_1)\dots\bt(k_n,l_n)\,,
\end{eqn}
where
\begin{eqn}\label{btkl}
\bt(k,l):=\sg_1^{-1}\sg_2^k\sg_1^{-1}\sg_2^{-1}\sg_1^l\sg_2^{-1}\,.
\end{eqn}
\end{enumerate}
\end{prop}

\proof The first two cases are direct to see, so assume
$\bt$ is not positive and is of Schreier vector length $>2$.
So the Schreier vector has isolated positive entries.
The form \eqref{btkl_} is then observed by directly drawing
the picture of the $A$-state. \qed

Since the form in \eqref{btkl} will be needed later,
let us call $\bt(k,l)$ \em{exceptional syllables} and
the words in \eqref{btkl_} \em{exceptional words}.
Note that by braid relations, an exceptional word
can be made to have only one of $k_1,l_1,\dots,k_n,l_n$,
say $l_n$, being non-zero, and can be rewritten
then as $\dl^{-2n}\sg_{1,2}^{l_n+2n}$. Here $\dl^2=
(\sg_1\sg_2)^3$ generates the center of $B_3$.

% \begin{rem}
% For links a bit more care is needed; see for example the connected sum
% of a positive and negative Hopf link.
% \end{rem}

We will use in a separate paper \cite{mwf}
semiadequacy to answer positively the
suggestive question: If a link is the closure of a positive braid,
and has braid index 3, is it always the closure of a positive 3-braid?
This question is motivated in particular by counterexamples of 4-braids.

Although we restricted ourselves to 3-braids, of course some
statements for general strand number are possible.

\begin{prop}
Let $\bt$ be a braid such that when writing
\[
\bt\,=\,\prod_{j=1}^k\sg_{p_j}^{q_j}
\]
with $p_j\ne p_{j+1}$ and $q_j\ne 0$, we have
$q_jq_k>0$ for $p_j=p_k$ (that is, $\bt$ is homogeneous)
and $q_j\ne -1$. Then $\hat\bt$ has non-trivial Jones polynomial.
\qed
\end{prop}

\begin{rem}
The knot $9_{47}$ has a homogeneous braid representation, but is not
weakly adequate, so that a general homogeneous braid representation
does not suffice to prove non-triviality of the Jones polynomial,
at least using the present approach.
\end{rem}

\subsection{3-braids with polynomials of $(2,\,.\,)$-torus links%
\label{S2q}}

In a separate paper \cite{mwf}, we explained how to find the closed
3-braids for a given Jones polynomial. This method will in general
not help to deal with infinitely many polynomials at a time.
A tool in that vein is probably too much to hope for, given
that the distribution of Jones polynomials on closed 3-braids
has surely some degree of difficulty (as became evident to
Birman \cite{Birman}). Dealing with a particular
interesting infinite family, here we will apply the formulas for
$V_{1,2}$, together with some representation and skein theoretic
arguments, to determine exactly which closed 3-braids have the
Jones polynomials $V_{(2,q)}$ of the $(2,q)$-torus links $T(2,q)$.

\begin{lemma}\label{lps}
Let $q>1$. Then there is no positive 3-braid different from
the $(2,q)$-torus link with the same Jones polynomial as $T(2,q)$.
\end{lemma}

\proof Since $2\md V=1-\chi$, we know that a positive 3-braid
$\bt$ with $V(\hat\bt)=V_{(2,q)}$ has $q+1$ crossings. Now we
know (see \cite{LickThis}) that
\begin{eqn}\label{sq}
\spn V(D)\,\le\,\frac 12\,\Bigl(c(D)+v(A(D))+v(B(D))-2\Bigr)\,.
\end{eqn}
Applying this to $D=\hat\bt$, we find that $\spn V=q=c(D)-1$ holds only
if the Schreier vector of $\bt$ has length $2$ or $4$. In former case
we are easily done, and in latter case we have additionally that
the inequality  \eqref{sq} is sharp, so that the extremal $B$-states
(in the language of \cite{Bmo}) have a non-cancelling contribution.
Using (for example) the formalism in \cite{Bmo}, one easily sees
that then either one $\sg_1$ and one $\sg_2$ syllable are trivial
(and then the braid reduces to a $(2,\,.\,)$-torus braid), or all
$4$ syllables are non-trivial. The we have an adequate diagram $D$,
and $|\bV_1|=1$ implies that the syllables have exponent $2$ except
exactly one. So we have $\sg_1^2\sg_2^2\sg_1^2\sg_2^k$ for $k\ge 3$,
and we can exclude these braids easily by induction on $k$
using the skein relation for $V$. \qed

Let $[\bt]$ be the exponent sum of $\bt$ and $\dl=\sg_1\sg_2\sg_1$.
We write $V(\bt)$ for $V(\hat\bt)$.
The \em{Birman dual} \cite{Birman} of $\bt\in B_3$ is $\bt^*=\bt^{-1}
\dl^{4p}$, where we assume $[\bt]=6p$ is divisible by 6. Birman
introduced this construction in order to obtain different 3-braid
links with the same Jones polynomial, proving $V(\bt^*)=V(\bt)$.

\begin{theorem}\label{tt}
Let $q\ge 0$. Then
\def\theenumi{\alph{enumi}}
\def\labelenumi{(\theenumi)}
\begin{enumerate}
\item If $q\ge 7$ and $q\equiv \pm 1\,(6)$, then except the
$(2,q)$-torus knot, there is exactly one other closed 3-braid
with the same Jones polynomial. It is obtained by taking the
Birman dual of $\sg_1^{q\mp 1}\sg_2^{\pm 1}$. 

\item For all other $q$, the $(2,q)$-torus link is the only
closed 3-braid with this Jones polynomial.
\end{enumerate}
\end{theorem}

\proof Modulo mirroring assume $\bt$ is an $A$-adequate 3-braid
with the polynomial of a $(2,q)$-torus link, where now $q\in \bZ$.
We can exclude $|q|\le 1$, since these are trivial polynomials
that were dealt with before.

First assume $q\ge -2$. Then $V_1=0$. By direct observation
(see proposition \reference{pexc}), $\bt$ is of the following
forms (after conjugation and braid relations):
\begin{itemize}
\item $\sg_1^{-2}\sg_2^k$, and $\sg_1^{-1}\sg_2^k$,
\item positive, or
\item $(\sg_1\sg_2^2\sg_1)^{-k}\sg_2^l$, with $k,l\ge 0$.
\end{itemize}

The first form is easy to handle, and the second one is handled
by lemma \reference{lps}, in both cases with the expected outcome.
So consider the third form. We can evaluate the Jones polynomial
using the trace of the (reduced) Burau representation $\psi$. We have
\[
\psi(\sg_1)\,=\,\left[\,\begin{array}{cc}-t & 1 \\ 0 & 1\end{array}
\right]\,,\qquad
\psi(\sg_2)\,=\,\left[\,\begin{array}{cc}1 & 0 \\ t & -t\end{array}
\right]\,.
\]
So for $k\in\bZ$
\[
\psi(\sg_1^k)\,=\,\left[\,\begin{array}{cc}(-t)^k &
  \ffrac{1-(-t)^k}{1+t} \\ 0 & 1\end{array}\right]
\,,\qquad
\psi(\sg_2^k)\,=\,\left[\,\begin{array}{cc}1 & 0 \\
  t \ffrac{1-(-t)^k}{1+t} & (-t)^k\end{array}\right]
\,.
\]
Then $\tr\,\psi(\sg_1^k)\,=1+(-t)^k$ and
\begin{eqn}\label{31.1}
\tr\,\psi(\sg_1^k\sg_2^l)\,=\,(-t)^k+(-t)^l+
t\,\frac{(1-(-t)^k)(1-(-t)^l)}{(1+t)^2}\,.
\end{eqn}
Also $\psi(\dl^2)=t^3\,\cdot\,Id$, where $\dl^2=(\sg_1\sg_2)^3$
is the center generator. Then (see \cite{Birman} or \cite{Jones})
\begin{eqn}\label{31.2}
V_{\hat\bt}(t)\,=\,\left(-\sqrt{t}\right)^{e-2}
\left[\,t\cdot \tr\psi\,+\,(1+t^2)\,\right]\,,
\end{eqn}
with $e=[\bt]$ the exponent sum. So with
$\bt=(\sg_1\sg_2^2\sg_1)^{-k}\sg_2^l$, we find 
\begin{eqn}\label{V}
V(\hat\bt)\,=\,V(\dl^{-2k}\sg_2^{l+2k})\,=\,
\left(-\sqrt{t}\right)^{l-4k-2}\,
\left[\,t\cdot t^{-3k}\,(1+(-t)^{l+2k})\,+\,(1+t^2)\,\right]\,.
\end{eqn}
We can assume $k,l>0$. Otherwise we have for $k=0$ a positive
braid, and can apply the lemma \ref{lps}, or for $l=0$ a negative braid.
In latter case, since we assumed $q\ge -2$, the only possibility
is the negative Hopf link polynomial, and since for a negative
braid $\Md V=-1/2$ determines $\chi(\hat\bt)$, we are easily done.

Now $k,l>0$. Then from \eqref{V} we have
\begin{eqn}\label{V2}
\md V\,=\,\frac{l-4k-2}{2}+1-3k\,=\,\frac{l}{2}-5k\,\ge\,-\frac 52\,,
\end{eqn}
the last inequality because of $q\ge -2$. Also, since $l+2k>2$, we find
\[
\Md V\,=\,\frac{l-4k-2}{2}+\left\{\,\begin{array}{c@{\quad}l}
1+l-k & 1+l-k>2 \\ 2 & 1+l-k<2  \\ 2 & 1+l-k=2 \mbox{ and $l$
even} \\ 0 & 1+l-k=2 \mbox{ and $l$ odd} \end{array} \right.
\]
In particular, $\spn V\ge l+2k-2$. But from \eqref{V} we see also
that $V$ has at most $4$ monomials, and by observing the shape of
$V_{(2,q)}$, we conclude $\spn V=q\le 4$, so that
\begin{eqn}\label{V3}
l+2k\le 6\,.
\end{eqn}
No integers $k,l>0$ satisfy simultaneously the inequalities in
\eqref{V2} and \eqref{V3}. This finishes the case $V_1=0$.

Now let $V_1=1$, that is, $q\le -3$. By a similar observation
we can determine the forms of $\bt$ (up to conjugacy and braid
relations):
\begin{itemize}
\item $\sg_1^k\sg_2^{-l}$, $l>2$, $k\ge -2$,
\item $\sg_1^k\sg_2^{-l}\sg_1^m\sg_2^{-n}$, $k,l,m,n>0$, or
\item $\sg_1^{-k}(\sg_1\sg_2^2\sg_1)^{-l}\sg_2^m$, with $k> 0$.
\end{itemize}
The first case is again easy to settle. The second case gives an
alternating (braid) diagram. Using the twist number as in
\cite{DL}, or the theorem in \cite{Thistle2}, one easily
verifies that among alternating links the $(2,q)$-torus links
are the only ones with their Jones polynomial. So
the third form is again the only interesting one.

We can assume again that $l,m>0$. Now we use proposition \reference{pv2}. Consider 
the $A$-adequate diagram $D=\hat\bt$ and its $A$-state. We have
$V_1=1$ and by direct observation,
\[
\begin{array}{l}
e_{++}=2l, \\ \dl=2l-1,
\end{array}\quad\mbox{and}\quad
\tg A(D)'=
\,\left\{\,\begin{array}{cl}
1 & \mbox{if\ } k=1 \\ 0 & \mbox{if\ } k>1\,\end{array}\right.
\,.
\]
So
\[
V_0V_2=\,\left\{\,\begin{array}{cl}
2 & \mbox{if\ } k>1\,\\
1 & \mbox{if\ } k=1\,\end{array}\right.\,.
\]
Since for $V_{(2,q)}$ we have $V_0V_2\in\{0,1\}$, we find $k=V_0V_2=1$
(and also $V_2\ne 0$, so that $q\ne -3$). Using \eqref{31.1} and
\eqref{31.2}, we have again
\begin{eqnarray}
\nonumber V_{\hat\bt}\,=\,V(\dl^{-2l}\sg_2^{m+2l}\sg_1^{-1}) & = &
\left(-\sqrt{t}\right)^{-4l+m-3}\,
\left[\,t\cdot t^{-3l}\,\left[
(-t)^{m+2l}-\frac{1}{t}+\frac{1-(-t)^{m+2l}}
{1+t}\,\right]\,+\,1+t^2\,\right]\,  \\
& = &
\label{*q} \left(-\sqrt{t}\right)^{-4l+m-3}\,
\Bigl[\,\underbrace{t^{1-3l}\,\sum_{i=-1}^{m+2l}(-t)^i}_{(*)}
\,+\,1+t^2\,\Bigr]\,.
\end{eqnarray}
If $m$ is even (and $\hat\bt$ is a 2-component link), then $V_0=1$,
which is not the case for any $V_{(2,q)}$ when $q\le -3$ (having
$V_0=-1$). So assume $m$ is odd (and $\hat\bt$ is a knot).
The term (*) in \eqref{*q} gives a coefficient list
$-1\,1\,\dots\,-1\,1\,-1$, with maximal polynomial degree, say,
$y$. In order to obtain a coefficient list
$-1\,1\,\dots\,-1\,1\,-1\,1\,0\,1$ by adding $t^x(1+t^2)$,
there are two options, $x=y+1$ and $x=y-2$. Then
$1-3l+m+2l\in\{-1,2\}$, so $m=l-2$ or $m=l+1$.

In these cases we have the polynomials $V_{(2,q)}$ for some
$q\le -2$ up to units. Then ``up to units'' can be removed,
either by direct check of the degree, or more easily by
the previous argument using $V(1)=1$ and $V'(1)=0$, since
we have knot polynomials. So we have the forms
\begin{itemize}
\item for $l=m+2$ odd, $\sg_1^{-1}(\sg_1\sg_2^2\sg_1)^{-l}
\sg_2^{l-2}$, and
\item for $l=m-1$ even, $\sg_1^{-1}(\sg_1\sg_2^2\sg_1)^{-l}
\sg_2^{l+1}$\,.
\end{itemize}
It is direct to check that the first family is the (mirrored)
Birman dual of $\sg_1^{6p-1}\sg_2$ with $l=2p-1$ and the
second family is the dual of $\sg_1^{6p+1}\sg_2^{-1}$ with $l=2p$.

Birman duality is trivial (i.e., the duals are equal) for $p=0$, or
in the first family if $p=1$. In the other cases one can use (as
in \cite{Birman}) the signature $\sg$ to distinguish the duals.
For the first family, first change a crossing in $\ap=\sg_1^{6p-1}
\sg_2$ to obtain an alternating braid $\ap_0=
\sg_1^{6p-1}\sg_2^{-1}$, and calculate $\sg(\ap_0^{-1}\dl^{4p})$
using Murasugi's formulas. Then observe that it differs for $p\ge 2$
at least by $4$ from $\sg(T(2,6p-1))=6p-2$, and that $\sg$ changes at
most by $2$ under the crossing change, that turns $\ap_0^{-1}\dl^{4p}$
into $\ap^{-1}\dl^{4p}=\ap^*$. \qed

The following consequence, due to El-Rifai, classifies 3-braids
whose Morton-Williams-Franks inequality for the braid index is inexact.
Here $MWF(\bt)=1+(\Md_lP-\md_lP)/2$ with $P=P(\hat\bt)$ is a notation
of \cite{Nakamura2} for this bound on (the closure of) a braid $\bt$.

\begin{corr}\label{CMWF}(El-Rifai)
A braid $\bt\in B_3$ has $MWF(\bt)\le 2$, if and only if $\bt$
is (up to mirroring and conjugacy) as specified in theorem
\reference{tt}.
\end{corr}

\proof It is an easy observation (using for example
\cite[proposition 21]{LickMil}; see remark 1.4 in \cite{Nakamura2})
that $MWF(\bt)\le 2$ implies that $P(\hat\bt)=P_{(2,q)}$
for some $q$, so in particular $V(\hat\bt)=V_{(2,q)}$. In the
reverse direction, note (as in \cite{Birman}) that Birman duality
preserves not only $V$ but also $P$. \qed

\begin{rem}
El-Rifai's result is obtained in his (unpublished) thesis \cite{El},
which I learned about from \cite{Nakamura2}. It was the interest in
this problem (and the inavailability of El-Rifai's solution) that
motivated theorem \reference{tt}.
\end{rem}

Theorem \reference{tt} and corollary \reference{CMWF} also
relate to some evidence we have that closed 3-braids with the
same Jones polynomial have in fact the same skein polynomial
$P$. From \cite{mwf} we know that for given $V$ at most three
different $P$ occur, and $P$ is unique if $\Md V\cdot \md V<0$.
This condition is not fulfilled, though, by torus link polynomials,
so it is justified to consider corollary \ref{CMWF} (still only)
as a consequence of theorem \reference{tt}.

\section{Upper bounds on the volume\label{vest}}

In this short section we explain how our formulas for $V_{1,2}$ relate
to the twist number of \S\ref{BP}, and volume. Another relationship
is shown for all $V_i$ for alternating links in \S\reference{SMM}.

\subsection{Adequate links}

\begin{defi}
Let $t_k(D)$ be the number of twist equivalence classes of
$k$ elements in a diagram $D$.
\end{defi}

\begin{prop}\label{p21}
For an adequate diagram $D$ of a knot or link $K$ we have
\begin{eqn}\label{six}
t(D)\,\ge\,|V_{1}|+|\bV_{1}|-2c(K)+2\spn V(K)\,\ge\,
-t_1(D)-t_2(D)\,,
\end{eqn}
and
\begin{eqn}\label{seven}
t(D)\,\ge\,2c(K)-2\spn V(K)\,.
\end{eqn}
\end{prop}

\proof We have as before $ |V_{1}|=e(G(A(D))')-v(G(A(D))')+1\,, $
and similarly for $\bV_1$ and $B(D)$.

Now we claim that
\begin{eqn}\label{ee}
c(D)-t_2(D)-t_1(D)\le e(G(A(D))')+e(G(B(D))')\le c(D)+t(D).
\end{eqn}
To see this note that we may assume (as in \cite{DL} or \cite{Lackenby})
that twist equivalent crossings form a twist tangle (i.e., that the
diagram is \em{twist reduced} in the terminology of \cite{Lackenby}).
Then in at least one of the $A$ and $B$-state of $D$, the traces of all
crossings in a twist equivalence class will connect the same two loops.

This explains the second inequality in
\eqref{ee}. For the first inequality, observe that in the other state, the 
splicings form a small (inner) loops between $p$ and $q$, and $p$ and $q$
connect the same loops only if there are no more twist equivalent
crossings to $p$ and $q$. (Only then the outer loop can be the same.)
Also no other crossings connect the same pairs of loops. The identification of 
crossings single in their twist equivalence class can not be controlled 
easily in $A(D)$ or $B(D)$. So
\[
t_2(D)\,+\,\sum_{k>2}k\cdot t_k(D)\,\le\,e(G(A(D))')+e(G(B(D))')\,,
\]
which is the left inequality in \eqref{ee}.

If $D$ is now ($A$ and $B$-)adequate, then by \cite{Thistle} we have
$c(D)=c(K)$ and by \cite{LickThis} that
\begin{eqn}\label{vv}
v(A(D))+v(B(D))=v(G(A(D))')+v(G(B(D))')= 2\spn V(D)-c(D)+2\,.
\end{eqn}
So
\[
|V_{1}|+|\bV_{1}|-2c(K)+2\spn V(K)\,=\,e(G(A(D))')+e(G(B(D))')-c(D)\,,
\]
and \eqref{six} follows. To see \eqref{seven}, we use \eqref{vv}
and must show that at least $c(D)-t(D)+2$ loops exist in $A(D)$
and $B(D)$ together. Now each equivalence class of $k$ crossings
gives rise to $k-1$ (inner) loops in one of the states, and each
state has at least one more loop to connect the outer splicings.
\qed

\begin{corr}\label{cor31}
If $K$ is adequate, for any adequate diagram $D$ of $K$ we have
(with $v_0=\vol(4_1)/2$ the ideal tetrahedral volume)
\[
\vol(K)\le 10v_0\,\bigl(\,2\spn V(K)+t_1(D)-c(K)-1\,\bigr)\,.
\]
\end{corr}

\proof Use $c(D)-t(D)+t_1(D)\ge t(D)$. The claim follows using 
the inequality $\vol(K)\le 10v_0(t(D)-1)$ 
of Lackenby-Agol-Thurston \cite{Lackenby}. \qed

\begin{rem}
Alternating diagrams are adequate with the additional property that
the $A$- and $B$-state loops are connected by crossing traces
only from outside (modulo change of the infinite region). Then it
is easy to see that crossings connecting the same vertices in
$A(D)$ or $B(D)$ are twist equivalent, and so the right inequality
of \eqref{ee} is exact. Then we have a different proof of the
Dasbach-Lin lemma \reference{lDL}.
\end{rem}

\begin{rem}
I do not know whether $V(K)$ determines $c(K)$ for adequate $K$.
One should expect that it does not, but among adequate knots
of $\le 16$ crossings there is no pair with different crossing
number but same Jones polynomial. Thistlethwaite \cite{Thistle}
proved that the Kauffman polynomial determines the crossing number
for adequate links. Also, the obvious fact that the left hand-side
of \eqref{vv} is at least 2 implies that $\spn V(K)\ge c(K)/2$, so
that $V(K)$ at least bounds above $c(K)$.
\end{rem}

\begin{exam}
It is in general not true that $t(D)$ is an invariant of adequate
diagrams $D$ of $K$. The simplest counterexample is $11_{440}$,
which has adequate diagrams of 6 and 7 twist equivalence classes.
(This corrects my initial blunder to believe the right inequality
of \eqref{ee} is exact also for general $A$-adequate diagrams,
and then that so would be the inequality in the corollary.) Also
$T_1(K)=|V_1|+|\bV_1|\le t(D)$ is not always true, as show the
examples $15_{253246}$ and $15_{253273}$.
\end{exam}
% 
% This allows to extend the definition of $t_a$ to adequate knots $K$.

The terms $t_{1,2}(D)$ seem a bit artificial, but they are necessary
at least for estimating the twist numbers, since in general 
adequate diagrams a big discrepancy between $T_1(D)$
and $t(D)$ may occur. To illustrate this, we give a few constructions
of such examples. Call the replacement of a crossing in a link
diagram by a clasp of the same checkerboard sign a \em{clasping}.
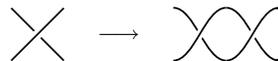
\[
\diag{7mm}{1}{1}{
  \picmultiline{0.18 1 -1.0 0}{0 1}{1 0}
  \picmultiline{0.18 1 -1.0 0}{0 0}{1 1}
}\quad\lra\quad
\diag{7mm}{2}{2}{
  \picPSgraphics{0 setlinecap}
  \pictranslate{1 1}{
    \picrotate{-90}{
      \rbraid{0 -0.5}{1 1}
      \rbraid{0 0.5}{1 1}
      \pictranslate{-0.5 0}{
      \picline{0.02 .95}{0 1}
      \picline{0.98 .95}{1 1}
    }
    }
  }
}
\]
(This may preserve or alter the number of components dependingly
on the crossing and orientation of the clasp.)

\begin{exam}\label{x0}
The diagrams $D_k$ given by the closures of the $n$-braids 
$(\sg_1\sg_2\dots\sg_{n-1}\sg_{n-1}\dots\sg_1)^k$
for $n>3$ fixed and $k\to\infty$ have $t_{1,2}\to\infty$,
but the contribution of the crossings with no twist equivalent one
to $e(G(A(D_k))')+e(G(B(D_k))')$ is bounded. Also, the crossings
in equivalent pairs are identified in both $G(A(D_k))'$ and
$G(B(D_k))'$. Thus the left inequality in \eqref{ee} is exact up to
a bounded difference. If one likes knot diagrams, apply
claspings at the last $n-1$ braid crossings, making $\sg_{n-1}
\dots\sg_1$ into $\sg_{n-1}^2\dots\sg_1^2$. So the $t_{1,2}$ terms
on the left of \eqref{ee} are inavoidable. Whether the $t_{1}(D)$
term is indeed essential in corollary \reference{cor31} is not
evident, though, since our examples have bounded volume.
\end{exam}

\begin{exam}\label{xdf}
Let $K_1$ be an adequate knot, and $K_{n+1}$ be constructed
from $K_n$ as follows. Take a disconnected twisted with blackboard
framing 2-cable of an adequate diagram $D_n$ of $K_{n}$, and
make a clasping at one of the mixed crossings of the resulting
2-component link diagram, so that the result is a knot and not a
cable. We have the diagram $D_{n+1}$ of $K_{n+1}$, which is easily
seen to be adequate (similarly to \cite{LickThis}). One can
check using lemma \reference{lcab} that cabling preserves
(adequacy and) $T_1=|V_1|+|\bV_{1}|$, while clasping augments it
by one. Thus $T_1(K_{n+1})=T_1(K_{n})+1$, while $c(K_{n+1})=4c
(K_{n})+1$, and the twist numbers of $D_n$  grow similarly. 
Thus $T_1(D_n)=O(\log_4 t(D_n))$.
% So $c(D)$ and $t(D)$ can grow exponentially compared to $T_1(D)$
% for adequate \em{knot} diagrams $D$. 
\end{exam}

% We can try to obtain inequalities for $\vol(K)$ that do not depend on
% a particular adequate diagram $D$ of $K$, or on $c(K)$, but just on
% $V(K)$, and cancel the $t_{1}(D)$ term with $c(D)$. However, then
% the estimate in proposition \reference{p21} is trivial, since $\spn
% V(K)\ge c(K)/2$. Then one may seek to estimate $e(A(D)')+e(B(D)')$
% from below increasingly using $c(D)$ or $t(D)$ only. However, we can 
% see from example \reference{x0} that in general no such estimate is 
% possible (including for knot diagrams). %; apply a few claspings). 

Still one can ask (with little hope, though)

\begin{question}\label{q1}
Are there only finitely many adequate or hyperbolic knots of given
$T_1$?
\end{question}

\begin{exam}
The existence of $T_1=0$ knots immediately implies (by connected sums)
that for knots which are non-hyperbolic and not adequate the answer
is negative. If the number of link components may grow unboundedly, the
following sequence of links gives a negative
answer: take a trivial 3-braid and place circles around strands $1,2$
and $2,3$ interchangingly (as in fig.~13 of \cite{CR}). It is now easy
to check that these links are adequate, with $T_1=0$, and their volume
grows unboundedly.
\end{exam}

We can thus also answer negatively the question of
Dasbach-Lin whether $T_1$ bounds (above) the volume at least for
links. For knots, it would be interesting to know the behaviour
of $\vol(K_{n})$ in example \reference{xdf}. It might show that
for an adequate link an upper bound on the volume, if it at all
exists, no better than exponential in terms of $T_1$ is possible
(in contrast to Dasbach-Lin's linear bound for alternating links).

% , but the previous construction
% similarly gives good countercandidates at least for links:
% take the disconnected (blackboard) framed $n$-cable of a fixed 
% adequate diagram, and apply one clasping. Then let $n\to\infty$.
%
% The $n$-cabling argument %in particular also shows that the answer
% to question \reference{q1} is negative if we allow links of an
% unbounded number of components. (In fact, the resulting
% diagrams have also unbounded twist number, so that it cannot
% be controlled by $T_1$ alone.) It
% also shows (as mentioned above)
% that $e(A(D)')+e(B(D)')$ can be as small as $C\cdot\sqrt{c(D)}$
% (even for knot diagrams; apply a few claspings).

% For an adequate knot $K$ we have
% \[
% \vol(K)\le 10v_0(|a_{1}|+|a_{d-1}|+2\spn V(K)-\myfrac{3}{2}c(K)-1)\,.
% \]
% and (unless $K$ is the unknot or trefoil)
% \[
% \vol(K)\le 10v_0(|a_{1}|+|\bV_{1}|+\myfrac{1}{2}\spn V(K)-3)\,.
% \]
% \end{corr}
% 
% \proof This is a consequence of theorem \reference{theo1}, since
% $t(D)=t(K)$, and so we have $t_2(D)\le c(K)/2$. Then the first
% inequality follows. The second inequality follows from the first
% since $\spn V(K)\le c(K)-1$, unless $K$ is alternating (see
% \cite{Kauffman}). For $K$
% alternating the Dasbach-Lin result implies the inequality unless
% $\spn V(K)=c(K)<4$. \qed

\subsection{3-braid and Montesinos links}

Concerning 3-braids, the finiteness theorem in \cite{mwf} implies
the existence of some upper bound on the volume in terms of the
Jones polynomial. We can make now a more concrete estimate:

\begin{prop}\label{VOL3}
Let $L$ be a 3-braid link. Then there is a constant
$C$ independent on $L$ such that
$\vol(L)\le C\,(T_1'(L)+1)$, where $T_1'(L)=|V_1'|+|\bV_{1}'|$,
and 
\[
V_1'\,=\,\left\{\,\begin{array}{c@{\ \ \mbox{if}\ }c}
V_1 & V_1\ne 0 \\ \md V & V_1=0 \end{array}\,\right.\,,
\quad
\bV_{1}'\,=\,\left\{\,\begin{array}{c@{\ \ \mbox{if}\ }c}
\bV_{1} & \bV_{1}\ne 0 \\ \Md V & \bV_{1}=0 \end{array}\,\right.\,.
\]
\end{prop}

\proof Assume $L$ has an $A$-adequate 3-braid representation $\bt$.
The links with Schreier vector length $2$ are non-hyperbolic.
If $A$-adequacy of $\bt$ comes from a positive 3-braid representation,
the number of crossings of $\bt$ is bounded linearly by the degree of
$V$. So assume $A$-adequacy comes from isolated positive Schreier
vector entries. From a straightforward generalization of the
argument proving proposition \reference{pexc} it follows that $V_1$
bounds linearly the length of the Schreier vector entries in $\bt$,
that do not belong to exceptional subwords. Now, we observed after
proposition \reference{pexc} that exceptional subwords can be
made to a single syllable by factoring out a power of $\dl^2$.
Such powers can be collected at the end of the braid, because
$\dl^2$ is central. So $\bt=\bt'\dl^{2n}$, and the Schreier
vector length of $\bt'$ is bounded linearly by $V_1$. The rest
follows by Thurston's hyperbolic surgery theorem.

% If an $A$-adequacy comes from isolated positive Schreier vector
% entries, then by direct verification $V_1\ne 0$ (see the proof
% of corollary \reference{c5.4}). So
% when $V_1=0$, then the 3-braid representation is positive.
% Then the number of crossings is bounded linearly by the degree of $V$.
% % orientatable. 
% % By checking the $A$-graph, we see that it is actually positive. 
% Otherwise $V_1$ bounds the length of the Schreier vector, and
% hence, the volume, linearly.
% 
A similar argument applies for
$B$-adequate 3-braid representations. \qed
% 
% A simple check shows that $|V_1|$ counts the number of negative
% Schreier vector entries (and the positive ones are not more than
% the negative ones), and so the length of the Schreier vectors of
% $\bt_i$ are bounded. Then the effect of adding positive crossings
% to a single $\sg_i^k$
% is to augment $\md V$ by $1/2$, while the effect of adding negative

The case of Montesinos links is more difficult, and the additional
work on $V_2$ is needed. We resume the notation of \S\reference{Ml}.

% here a complete
% result is not possible. We have, though, fairly strong restrictions on
% degeneracies that can occur.
% 
% \begin{defi}
% We call a sequence $L_i$ of Montesinos links \em{pretzel-like} if 
% after possible mirroring,
% \begin{enumerate}
% % \item the length $n_i\to\infty$,
% \item $L_i$ have a representation 
% \[
% L_i=M(p_{i,1}/q_{1,i}\,\dots,p_{i,n_i-1}/q_{i,n_i-1},
% -p_{i,n_i}/q_{i,n_i})
% \]
% where \em{after possible permuation of $p_{i,k}/q_{i,k}$}, we have
% $0<p_{i,k}/q_{i,k}< 1/2$ for $1\le k\le n_i-1$, and
% $0<p_{i,n_i}/q_{i,n_i}\le 1/2$, and
% \item when $i\to\infty$, only a bounded number of $p_{i,k}\ne 1$.
% \end{enumerate}
% \end{defi}

\begin{prop}\label{VOLM}
If $L_i$ are Montesinos links and have the same Jones polynomial,
then $\vol(L_i)$ are bounded. More exactly, for every Montesinos link
$L$, we have $\vol(L)\le C'\cdot(T_1(L)+T_2(L))$, with a constant
$C'$ independent on $L$.
\end{prop}

\proof Choose the general form of representation 
\[
L_i=M(p_{i,1}/q_{i,1}\,\dots,p_{i,n_i-1}/q_{i,n_i-1},
p_{i,n_i}/q_{i,n_i},e_i)\,,
\]
with $p_{i,j},q_{i,j}\ge 1$. W.l.o.g. assume $L_i$ are $B$-adequate,
so that $e_i\ne -1$.

First observe that, from part (1) of proposition \reference{pmnt},
if $L_i$ are adequate, i.e. $e_i\ne 1-n_i$, then the twist number is
bounded from a multiple of $T_1$. So we can assume that $e_i= 1-n_i$.

% If $e_i\ne 1-n_i$ or
% $
% \#\,\{\,k\,:p_{i,k}/q_{i,k}\le 1/2\,\}\ne 1
% $
% then observe that $L_i$ are adequate, or apply part (2) of
% proposition \reference{pmnt} on $!L_i$, to deduce that
% $c(L_i)$ is bounded from $\spn V(L_i)$. 
% 
% So $e_i=1-n_i$ and $\{\,k\,:p_{i,k}/q_{i,k}\le 1/2\,\}=\{k_i\}$.

By looking at $\bV_{1}(L_i)$ and using part (1) of proposition
\reference{pmnt} and proposition \reference{prat}, we see that
for fixed $\bV_{1}$ the twist number may go to infinity only if
the lengths $n_i$ go to infinity and 
almost all $q_{i,k}=p_{i,k}+1$ (i.e. for given $i$ the number of
$k$ not satisfying this condition remains bounded as a sequence
over $i$ when $i\to\infty$). Otherwise $\bV_1$ bounds the twist
number.

Now consider $\bV_{2}(L_i)$. We need to study the behaviour of
the terms independent from $\bar V_1$, which are $\tg B(D)'$,
$e_{++}(B(D))$ and $\dl B(D)'$. The first quantity remains
bounded. Every $q_{i,k}=p_{i,k}+1$ contributes one to $e_{++}
(B(D))$. This contribution can be equilibrated by $\dl B(D)'$
only if $e_i=-2$. (Otherwise there is no intertwined pair of
connections in $B(D)$.) Then, however, for length $n_i>3$, the
$L_i$ are adequate. Thus a non-zero multiple of all twists will
be detected by $|\bV_{2}|+|\bV_{1}|$, and we are done.
% 
% Then incorporate the $n_i-1$ negative $e_i$-twists into the tangles
% $p_{i,k}/q_{i,k}>1/2$ for $k\ne k_i$, making $p_{i,k}/(p_{i,k}+1)$
% into $-1/(p_{i,k}+1)$, put after permuting the $n_i-1$ now negative
% tangles first, and go over to the mirror image.
% 
% We obtain that $L_i$ are pretzel-like. 
\qed

\begin{rem}\label{remV3}
Unfortunately, I still don't know how to prove that
there are only finitely many Montesinos links with the same
Jones polynomial, since controlling the crossing number from
$\spn V$ is non-trivial. There are situations, for example the
$(-3,k,l)$-pretzel links with $k,l\ge 4$, where in the almost
alternating representation with $e=-1$, the third terms
in $A\ag{D_A}$ and $A^{-1}\ag{D_B}$ of the proof of part
(2) of proposition \reference{pmnt} cancel also.
In the paper \cite{HTY} it was shown that the difference
between $\spn V(L)$ and $c(L)$ for certain pretzel links
can actually be made arbitrarily large, so that even a study
of $V_3$ etc.\ may not be helpful.
% It is possible that this is the worst case, and $c(L)\le
% \spn V(L)+4$, but a study of $V_3$ appears extremely
% complicated, even only for the cases we need.
\end{rem}

\section{Mahler measure and twist numbers\label{SMM}}

\subsection{Estimating the norm of the Jones polynomial}

We use the notation $t(K)$ and $t_a(K)$ from \S\reference{BP}.

\begin{theorem}\label{theo1}
For all numbers $n\in\bN$ there are numbers $c_n,d_n$ such that
if $t(K)\le n$, then $(t+1)^nV_K(t)$ has at most $c_n$ non-zero
coefficients, each of which has absolute value at most $d_n$.
\end{theorem}

In other words, the $1$-norm $||X||_1=\sum_{i}|X_i|$ (with $X_i$
the coefficients of $X$) of $X=(t+1)^nV_K(t)$ is bounded when
$t(K)\le n$.

\begin{rem}\label{r1}
One can write down explicit values for $c_n$ and $d_n$, which would
be exponential in $n$.
\end{rem}

In simultaneous joint work with Dan Silver and Susan Williams
\cite{STW} we give similar results to theorem \ref{theo1}
for the Alexander and skein polynomial. Latter result implies
a different, but more technical, proof of theorem \ref{theo1}.
%, in which the bounds $d_n$ on the coefficients are not stated.)

Dan Silver also informed me of recent related work by I.~Kofman
and A.~Champanerkar \cite{CK}. While they work in a more general
situation concerning (multi-strand) twisting and (colored) Jones
polynomial, our proofs are more conceptual and less computational,
and owe to the fact that we have done some of the relevant work
previously, in \cite{gwg,gbi2,coeff}. They should also emphasize
the relation to the questions in \cite{DL}.

\begin{corr}
$||(t+1)^nV_K||_{L^2(S^1)}$ is bounded for $t(K)\le n$.
\end{corr}

Here $S^1$ is understood as the set of unit norm complex numbers.

\proof By \cite[\S 6]{coeff}, we have $||\,X\,||^2_{L^2}\,=\,\sum
\,|X_i|^2$, where $X_i$ are the coefficients of $X$. \qed

The Mahler measure $M(P)$ of a polynomial $P$ is defined as the product
of norms of all (complex) roots outside the unit circle $S^1$ and the
norm of the leading coefficient. For the relevance of this concept
and further discussion see for example \cite{GH,SW}.

\begin{corr}
The Mahler measure $M(V_K)$ is bounded for $t(K)\le n$.
\end{corr}

\proof Apply \cite[lemma 6.1]{coeff}. \qed

\begin{corr}
The Mahler measure $M(V_K)$ bounds below increasingly $\vol(K)$ when
$K$ is alternating. (That is, if $M(V_K)$ grows to infinity, then
so does $\vol(K)$.)
\end{corr}

\proof Use Lackenby's inequality \cite{Lackenby},
$\vol(K)\ge v_0(t_a(K)-2)/2$ for $K$ alternating. \qed

\begin{corr}
Let $V_K$ be the Jones polynomial represented as in \eqref{Vc}.
Then for all $l$ and $n$ there are numbers $c_{n,l}$ such that if
$t(K)\le n$, then $|V_{l}|\le c_{n,l}$.
\end{corr}

\proof The possible values for $a_0t^k+a_1t^{k+1}+\dots+a_lt^{k+l}$
are determined by the first $l+1$ (with respect to the minimal
degree) coefficients of $(t+1)^nV_D$, and there are finitely many
such choices. \qed

\begin{corr}\label{ZZ}
Any coefficient $|V_{l}|$ in the representation \eqref{Vc}
(that is, fixed with respect to the minimal degree of $V$)
bounds below increasingly $\vol(K)$ when $K$ is alternating. \qed
\end{corr}

Similarly, so do $T_l=|V_{l}|+|\bV_{l}|$, which explains rigorously
some of the experimental observations of Dasbach-Lin \cite{DL}.
By remark \reference{r1}, explicit estimates obtained this way
are not optimal, so we clearly lose some efficiency for the sake
of generality. Possibly one can prove better estimates with more
effort from the graph theoretical setting of \cite{DL}, but at least
for non-alternating knots, where the graph approach fails, we have
new information from our main theorem. A few more applications follow
after its proof.

\proof[of theorem \reference{theo1}]
Let $T$ be a diagram \em{template}. It is obtained from a link diagram
by replacing some crossings by slots
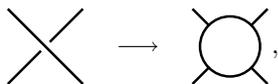
\[
\diag{1cm}{1}{1}{
  \picmultiline{-5 1 -1 0}{0 0}{1 1}
  \picmultiline{-5 1 -1 0}{0 1}{1 0}
}\quad\lra\quad
\diag{1cm}{1}{1}{
  \picline{0 0}{1 1}
  \picline{1 0}{0 1}
  \picfilledcircle{0.5 0.5}{0.4}{}
}\,,
\]
as was done for tangles in
\cite{SunThis}. We represent any twist sequence of link diagrams
by a template $T$, into whose slots twist tangles are inserted.
We say that we \em{associate} a diagram to the template in this way.
We meet the convention that twist tangles inserted into template slots
have an even number of half-twists (possibly leave a crossing outside
the slot).

We assign a sign to each slot of $T$, which is positive if the
twists in the slot are parallel and negative if the twists are
reverse. (This means that not necessarily any sign choice can be
realized by proper link component orientation.)

Let $s(T)$ be the number of slots of $T$, $c(T)$ the number of
crossings (which do not include slots) and $s_p(T)$ the number of
positive (parallel twist) slots of $T$. 

\begin{lemma}
For each template $T$ there are numbers $c_T$ and $d_T$ such that if
$D$ is a diagram associated to $T$, then $(t+1)^{s(T)}V_D(t)$ has
at most $c_T$ non-zero coefficients, each of which has absolute
value at most $d_T$.
\end{lemma}

Since for diagrams $D$ with $t(D)\le n$ we need to consider
templates $T$ with $s_p(T)\le s(T)\le n$ and $c(T)\le n$ (the
crossings left outside the slots to adjust half-twist parity),
which are finitely many, the theorem follows.

We prove the lemma inductively on $s_p(T)$. 

First the lemma is easy to verify when $s(T)=1$ and $c(T)=0$
(the $(2,n)$-torus links). For $s_p(T)=0$ it follows from 
formula (8) of \cite{gwg} (as explained in \cite{gbi2},
including the correction of the misprint in the original version).

Now let $s_p(T)>0$ and consider a parallel twist slot
\begin{eqn}\label{3}
\diag{1cm}{1}{1}{
  \picline{0 0}{1 1}
  \picline{1 0}{0 1}
  \picfilledcircle{0.5 0.5}{0.4}{}
}\,,
\end{eqn}
of $T$, in which we assume the twists to be vertical.
Replace $V(D)$ by the Kauffman bracket $[D]$ of $D$ (see
\cite{Kauffman}). We have
\[
V_D(t)\,=\,\left(-t^{-3/4}\right)^{-w(D)}\,[D]
\raisebox{-0.6em}{$\Big |_{A=t^{-1/4}}$}\,,
\]
so the claims for $(t+1)^nV(D)$ and $(A^4+1)^n[D]$ are equivalent.
The bracket of \eqref{3} is a linear combination of those of
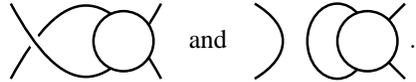
\[
\diag{1cm}{2}{1}{
  \picmulticurve{-5 1 -1 0}{0 0}{0.7 1.3}{1.3 1.3}{2 0}
  \picmulticurve{-5 1 -1 0}{0 1}{0.7 -0.3}{1.3 -0.3}{2 1}
  \picfilledcircle{1.5 0.5}{0.4}{}
}\quad\mbox{and}\quad
\diag{1cm}{2}{1}{
  \picline{2 1}{1.5 0.5}
  \picline{2 0}{1.5 0.5}
  \picellipsearc{1.1 0.5}{0.4 0.5}{45 -45}
  \picfilledcircle{1.5 0.5}{0.4}{}
  \piccurve{0 1}{0.5 0.7}{0.5 0.3}{0 0}
}\,.
\]
The first diagram has smaller $s_p$ (since one of the components
in the twist slot has now reverse orientation), while the other one
is a connected sum of a $(2,n)$-torus link and a diagram of smaller
$s$ and $s_p$, and the claim follows easily. \qed

\subsection{Link families with growing twist numbers}

% Here is another small application 
We have now a simpler proof of

\begin{prop}
Let $V,X\in\bZ[t,t^{-1}]$ with $X\ne 0$ and $V\ne 1$, and $K_n$
be \em{knots} (not links) with $V(K_n)=XV^n$. Then $t(K_n)\to\infty$.
\end{prop}

\proof If $V$ is zero or a unit in $\bZ[t,t^{-1}]$, then for all
but finitely many $n$ there is no knot with Jones polynomial
$XV^n$. (If the Jones polynomial of a \em{knot} is determined up
to units, it is uniquely determined, because $V(1)=1$ and $V'(1)=0$.)
Thus $||V||_{L^2}>1$, and so (see \cite[\S 6]{coeff})
$\max_{|t|=1}|V(t)|>1$. Then $|V|>1+\eps$ and $|\hat X|>\eps$ on an arc
of $S^1$ for any non-zero polynomial $\hat X$, so that $||\hat XV^n||_
{L^2}\to\infty$ when $n\to\infty$. Now taking $\hat X=(t+1)^pX$
we would have a contradiction to $t(K_{n_i})\le p$ for any
subsequence of $K_n$. \qed

In particular, taking $V$ to be the polynomial of a $T_1=0$ knot,
we have

\begin{corr}
There exist $T_1=0$ knots of arbitrarily large twist number. \qed
\end{corr}

Taking $V$ to be the polynomial of the figure-8-knot, we have

\begin{corr}
There exist (even alternating) knots with $M(V_K)=1$ and arbitrarily
large twist number (not only alternating twist number).  \qed
\end{corr}

More generally all the work of \cite{gwg} done on the relation
between Jones polynomial and canonical (weak) genus, can be extended
to twist numbers. Among others, we can now state:

\begin{corr}
Let $K$ be knot that has a (not necessarily untwisted) Whitehead double
$W_K$, whose Jones polynomial $V(W_K)$ has an absolute coefficient
$|V_i|\ge 4$ for some $i$. Then any sequence of Whitehead doubles
of $\#^nK$ has unbounded twist number when $n\to\infty$. \qed
\end{corr}

Note that all the preceding examples are non-hyperbolic, so that
the volume cannot be used to estimate the twist number.

{}From \cite{gwg} we obtain also

\begin{corr}
Let $|t|=1$ and $t(K)\le n$. Then for $t\ne -1$ we have $|V_K(t)|\le
C_{t,n}$ for some constant $C_{t,n}$ depending only on $t$ and $n$.
Also $|V_K(-1)|\le C_n(\spn V_K)^n$ for some constant $C_n$ depending
only on $n$. \qed
\end{corr}

Latter property can be used for example to easily deduce

\begin{prop}
For $n\ge 2$ there exist $n$-component links with
trivial Jones polynomial and arbitrarily large twist number.
\end{prop}

\proof By \cite{EKT}, there are such links that have as a component
a 2-bridge knot whose alternating Conway notation contains an
iterative subsequence. Clearly the twist number of a link is not
smaller than the one of its components. But it is easy to see that
the determinant $|V(-1)|$ of such 2-bridge knots grows exponentially
(with the number of occurrences of the subsequence), while their
crossing number (or span of the Jones polynomial) grows linearly. \qed

\begin{question}\label{q71}
Are there infinitely many \em{knots} with the same Jones polynomial
that have arbitrarily large twist number?
\end{question}

\begin{rem}
Kanenobu \cite{Kanenobu} found infinitely many knots with the
same Jones polynomial, but these knots have bounded twist number.
\end{rem}

The following questions arise (and were posed also by others) as a
natural follow-up to theorem \reference{theo1} and its applications.

\begin{question}\label{qr}
Are there knots/links of unbounded volume, which
\def\theenumi{\arabic{enumi}}
\def\labelenumi{(\theenumi)}
\begin{enumerate}
\item \label{(S)} (Dan Silver) have Jones polynomial of bounded
	Mahler measure? Or, as possible improvements,
\item \label{(a)} are there alternating such knots or 
\item \label{(K)} (Efstratia Kalfagianni \cite{Kalfagianni}) are
      there knots with Jones polynomial of bounded $1$-norm $\sum
      |V_i|$, or
\item \label{(b)} even with the same polynomial?
\end{enumerate}
\end{question}

\begin{rem}
Since the $1$-norm bounds the Mahler measure (see \cite{coeff}), a
positive answer to part \ref{(K)} implies such for part \ref{(S)}.
Part \reference{(a)} is related to the question how to find \em{prime}
alternating knots whose Jones polynomial is a power of $V(4_1)$ (or
more generally cyclotomic; for this problem see \cite{CK2}). Since
for alternating knots the $1$-norm of the Jones polynomial is equal
to the determinant, a series of knots answering positively part
\reference{(K)} (or even \reference{(b)}) are \em{a forteriori}
non-alternating. Since bounded twist number implies bounded volume, 
a positive answer to part \reference{(b)} for knots will imply such to
question \reference{q71}. In particular, again Kanenobu's series
can not be used. Similarly fails now also the construction
in \cite{EKT}, which yields non-hyperbolic links. 
It is worth remarking that for the Alexander polynomial most
of these questions can be effectively decided \cite{Kalfagianni,%
STW,aparb}.
\end{rem}

\section{Odd crossing number achiral knots\label{mpe}}

\subsection{Theorem and initial remarks}

We conclude with the most substantial application of the work
in this paper. It relates, as explained, to one of the oldest
problems in knot theory, concerning the crossing number of
achiral knots, dating back to more than a century ago, when
it was observed by Tait in his pioneering work on tabulating
the simplest knots. For alternating knots, this problem was
solved in \cite{Kauffman,Murasugi,Thistle2} as a consequence
of the proof of another of Tait's conjectures, the minimality
of reduced alternating diagrams. This was one of the celebrated
achievements of the Jones polynomial. Contrarily, we can state now

\begin{theorem}\label{th15}
For all numbers $n=11+4k$, $k>0$, there exists a prime achiral knot
$K_k$ of crossing number $n$. Specifically, $K_k$ is the closure of
the $5$-braid $-12^23^{2k}4-32-1(-2)^{2k}(-3)^24-23$. 
\end{theorem}

(In the braid word, $i>0$ stands for the Artin generator $\sg_i$,
and $-i$ for $\sg_i^{-1}$.) One easily sees that the given knots
are negative amphicheiral. The case $k=1$ is $K_1=15_{224980}$,
Thistlethwaite's example of an amphicheiral 15 crossing knot
(see \cite{HTW}).
% (which also likely triggered the question in Kirby's book).
It apparently came up in routine knot tabulation, and 
the insight that this knot has crossing number 15
is a result of generating all diagrams of fewer crossings,
and being able to distinguish $K_1$ from these diagrams. Even
though candidates to generalize such an example are straightforward,
the major problem one faces when the crossing number increases is
that exhaustive diagram verification quickly becomes impossible.

So one is led to seek general tools to rule out all potential
fewer crossing diagrams. Our formulas for the coefficients
of $V$ seem the first tool that is powerful enough to accomplish
the main part of the work, although a lot remains to be done,
and again computer aid at some point seems very helpful. (A
careful reading will reveal that a minor part of our arguments
is in fact redundant. Even aware of this, we preferred not to
omit them, since the complexity of proof makes the danger
serious to err on the opposite side and introduce a gap.)

Our result settles the first half of theorem \reference{thm}
and bases on Thistlethwaite's example. We will later deal with
the other half of odd crossing numbers by the same method,
though that case poses added difficulty. (The proof of theorem
\ref{th15} % the exposition would certainly become too long. 
will explain why; see remark \reference{r9,1} after it.) 
% To some extent it is more important to exhibit the approach
% that now allows, in principle, to attack such a question.
% A similar remark applies
% on the word ``almost''. The finite number of diagrams we will
% leave out can be explicitly tracked down during the proof. But
% I have so far no way to treat them except on a case by case basis.
% This turned out to be too tedious and technically difficult
% (both by hand and by computer), and would blow up the proof
% quantitatively, without that we gain much qualitatively.

For the proof of theorem \reference{th15} we need most of the
preparations, in particular the invariants in \S\reference{Cab}.
It is useful to assume below that $k>1$. That is, we waive on
unnecessarily reproving Thistlethwaite's example, for which our
argument would work, but would need slight modification.

A few simple observations are useful to collect in advance.

\begin{lemma}
$K_k$ are non-alternating.
\end{lemma}

\proof We have $V_0V_1=\bV_0\bV_1=-1$, so if alternating,
these knots must have twist number two. Then they are either
rational of genus $1$, or positive. Former is excluded looking
at the Alexander polynomial, and latter, for example, using
the signature. \qed

\begin{lemma}
$K_k$ are prime.
\end{lemma}

\proof $K_k$ have braid index $b\le 5$. We use the Birman-Menasco 
result \cite{BirMen2} on the (1-sub)additivity of $b$ under connected
sum. If composite, $K_k$ must be either connected sum of \\
\hbox{\vbox{%
\begin{tabbing}
\kern1cm\= \kill \\[-9mm]
(i)   \> four braid index $2$ factors, \\
(ii)  \> an amphicheiral $b=3$ knot and two braid index $2$ factors, \\
(iii) \> two amphicheiral $b=3$ knots, or \\
(iv)  \> a $b=3$ knot with its mirror image.\\[-8mm]
\end{tabbing}
}}

By Birman-Menasco \cite{BirMen}, or the work in \S\ref{S3B} here,
amphicheiral $b=3$ knots are alternating, and 
since braid index $2$ knots are too, but $K_k$ are not, we can have 
only option (iv). In that case look at the skein polynomial $P$.
Its second highest $m$-coefficient $X=[P]_{m^{\Md_mP-2}}$ must 
be of the form $f(l)+f(l^{-1})$ for some $f\in\bZ[l]$, in particular
$[X]_{l^0}$ must be even. But by direct calculation it is an
odd integer. \qed

\subsection{Start of proof of theorem \reference{th15}: basic setup}

The braid representation in the theorem gives a diagram of $K_k$ of
$12+4k$ crossings, that reduces (by a generalized version of the
move for Thistlethwaite's example) 
to a semiadequate diagram $\tl D_k$ of writhe
$\pm 1$ and $11+4k$ crossings. By \cite{Thistle}, we have $c(K_k)\ge
10+4k$, and if $c(K_k)=10+4k$, then $K_k$ is adequate. We assume, fixing
$k$, that $D_k=D$ is an adequate diagram of $K_k$. The proof will
consist in successively ruling out all possibilities for $D$.
% (up to finitely many).

First we remark that $D$ is prime, because $K_k$ is prime and
an adequate (connected sum factor) diagram represents a
non-trivial knot.

We use the invariants of \S\reference{Cab}. Either by looking at the
semiadequate diagrams $\tl D_k$, or (more tediously) by calculation of
the Jones polynomial and its 2-cable, we find that $D_k$ must satisfy
\begin{eqn}\label{tstcnd}
\tg=0,\quad \chi(G)=0,\quad \chi(IG)=0\ \mbox{ for
\em{both} $A$- and $B$-state,\quad and} \quad \gm=4\,.
% \parbox{0.8\textwidth}{\vbox{
% \begin{itemize}
% \item $\tg=0$
% \item $\chi(G)=0$
% \item $\chi(IG)=0$
% \item $\gm=4$
% \end{itemize}
% }}
\end{eqn}

These properties will be used throughout the whole proof.
The quantity $\gm$ refers to \eqref{gm}. 

We will use a similar approach as in \S\reference{G3} and again
construct the $A$-state of $D$ by successively attaching loops as in
\eqref{*}. A difference to \S\reference{G3} is that now $\chi(G)=0$,
so $G$ contains (exactly) one cycle, and we must start building
$A(D)$ with this cycle rather than a single loop.
We call the cycle of loops in $A(D)$ the \em{outer cycle}\footnote
{Beware in particular of confusing the words `loop' and `cycle';
the meanings they are used in here are very different! A loop is
a piece of a diagram obtained after splicing all crossings, and a
cycle in a set of loops connected by splicing traces in an
appropriate manner.}.
Since $D$ is not positive, the length $n$ of the outer cycle is odd,
and since $\tg=0$, we have $n>3$.
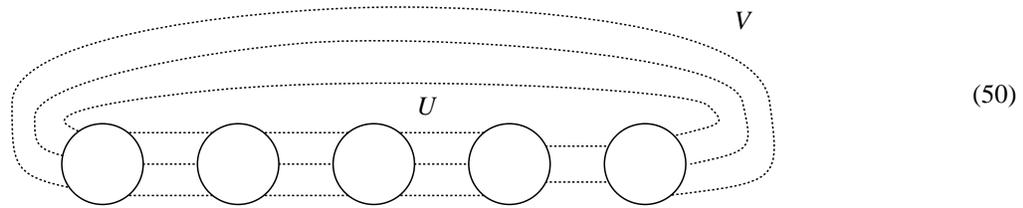
\begin{eqn}\label{oaz}
\diag{6mm}{17}{5}{
  {\piclinedash{0.05 0.05}{0.01}
    \picstroke{
      \opencurvepath{14.5 1.6}{15.8 1.9}{15 2.7}{8 3}{1 2.3}{1 2}
	{1.5 1.7}{2 1.7}{}
      \opencurvepath{15 1}{16.4 1.3}{15.9 3.5}{9 3.9}
        {4 3.5}{0.5 2.8}{0.5 1.3}{2 1}{}
      \opencurvepath{14.5 0.3}{17.0 .6}{16.5 4.2}{11 4.5}{6 4.5}
	{0 4}{0 0.7}{2.2 0.3}{}
    }
   \picmultigraphics{3}{3 0}{
     \picmultigraphics{3}{0 0.7}{
       \picline{2.5 0.3}{4.5 0.3}
     }
   }
   \picline{11.5 0.6}{13.5 0.6}
   \picline{11.5 1.4}{13.5 1.4}
  }
  \picmultigraphics{5}{3 0}{
    \picfilledcircle{2 1}{0.9}{}
  }
  \picputtext{9.2 2.3}{$U$}
  \picputtext{16.2 4.2}{$V$}
} %Fig 0
\end{eqn}
Let us for convenience assume $D$ drawn so that the loops
in the outer cycle are always connected from the outside,
i.e. the unbounded region in their complement. We define
\em{interor} and \em{exterior} of these loops according to
this convention. 

Every loop $l$ which is not in the outer cycle is attached to a
single previous loop $m$. This means that all traces of crossings,
that connect $l$ at one end connect to $m$ on the other. The
parallel equivalence classes of such traces are called below
\em{legs}. So legs are edges in $A(D)$ connecting an attached
loop. The number of legs of $l$ is called the \em{valence} of
$l$. As in \eqref{parallel}, and unless stated clearly otherwise,
we group parallel traces into a single one (with multiplicity
indicated, or explained from the context) when drawing loop
diagrams. So a dashed line in a diagram starting from a loop
attached by \eqref{*} (usually) stands for a leg.

Note that $\gm$ does not depend on the size (multiplicity) of
parallel equivalence classes of traces in $A(D)$ or $B(D)$.
Again the addition of legs, and consequently the move \eqref{*},
never reduces $\gm$, so if we reach a stage that $\gm>4$, or we can
delete legs and see that $\gm>4$, we can discard any (diagram 
obtained by) continuation of loop attachments.

It is clear that the legs of each attached loop form a multiple
connection of traces. So the existence of attachments forces
the intertwining graph $IG(A)$ to be non-empty (i.e. have
at least one vertex). Then the condition $\chi(IG)=0$ means
that we must have a cycle in $IG(A)$, and similarly in $IG(B)$.
The existence of this cycle will be helpful at several places below.

\begin{conv}
We also meet the convention that we consider diagrams \em{up
to mutations}. All our tests involve (quantities and properties
determined by) invariants that do not change under mutations
(see in particular corollary \reference{cormut}), so
this restriction is legitimate.
\end{conv}

Note that we have the same conditions for $B(D)$. We will
try to use these conditions to restrict the possible ways
in which we can build $A(D)$ successively by attaching loops.

Now we define the \em{depth} of loops in $A(D)$ inductively over the
order of attachments as follows. The depth of the loops in the outer
cycle, and of loops attached to them from the exterior, is $0$. Assume
next a loop $M''$ is attached to a loop $M$. If $M$ has depth
$k$ and is connected to a loop $M'$ of depth $k'\le k$, then any
loop $M''$ connected to $M$ from the opposite side to $M'$ has
depth $k+1$. If $M''$ is attached to $M$ from the same side as
$M'$, then we set the depth of $M''$ to be $k$. For example
in part (a) of figure \reference{Fig19} (page \pageref{Fig19}),
we have an outer cycle loop $L$ (of depth 0), other depth-0 loops
$E,F,G$ (not on the outer cycle), and depth-1 loops $A,B,C,D$.

If $k>0$, we call the side of a loop that contains 
loops of not larger depth the \em{exterior} of $M$. The other
region is the \em{interor} of $M$. If the interor contains
no further loops, we call $M$ \em{empty}.

We call $a$ \em{attached} to $b$ if there is a leg of $a$
connected to $b$. If $a$ is attached inside $b$ then $a$
has higher depth, and if $a$ is attached outside, it has the
same depth, but a higher distance to a loop of smaller depth.
 
By mutations we assume that the total depth of loops is the
smallest possible. In particular, the following flypes/mutations
are applied to reduce the total depth:
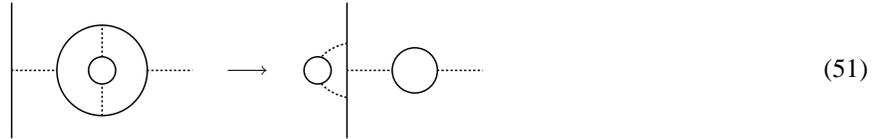
\begin{eqn}\label{Fig1}
\diag{6mm}{4}{3}{
  \picline{0 0}{0 3}
  \picfilledcircle{2 1.5}{1}{}
  {\piclinedash{0.05 0.05}{0.01}
   \picline{0 1.5}{1 1.5}
   \picline{3 1.5}{4 1.5}
   \picline{2 0.5}{2 2.5}
  }
  \picfilledcircle{2 1.5}{0.3}{}
}
\quad\lra\quad
\diag{6mm}{4}{3}{
  \picline{1 0}{1 3}
  {\piclinedash{0.05 0.05}{0.01}
  \piccurve{1 2.1}{0.1 1.9}{0.1 1.1}{1 0.9}
  \picline{1 1.5}{4 1.5}
  }
  \picfilledcircle{0.35 1.5}{0.3}{}
  \picfilledcircle{2.5 1.5}{0.5}{}
}
% Fig 1
\end{eqn}

We call a diagram \em{flat} if it has minimal total depth
(sum of depth of its loops) in its mutation equivalence class.
\em{So we assume from now on that $D$ is flat.}

The next assumption we can make is that, among diagrams
of the same total depth, we choose the one with the
fewest edges, i.e. \em{we assume $D$ is edge-reduced}.

\begin{lemma}\label{lmdf}
The two legs of a valence $2$-loop $L$ are not attached to the
same region, i.e. $A\ne B$ in
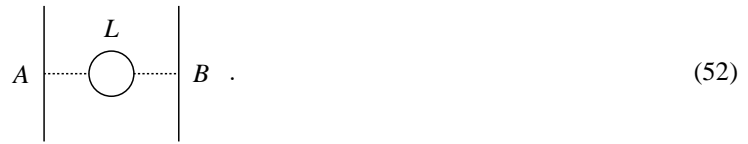
\begin{eqn}\label{Fig2}
\diag{6mm}{5}{3}{
  \picline{4 0}{4 3}
  \picline{1 0}{1 3}
  {\piclinedash{0.05 0.05}{0.01}
   \picline{1 1.5}{4 1.5}
  }
  \picfilledcircle{2.5 1.5}{0.5}{}
  \picputtext{0.5 1.5}{$A$}
  \picputtext{4.5 1.5}{$B$}
  \picputtext{2.5 2.5}{$L$}
}\,.
% Fig 2
\end{eqn}
\end{lemma}

\proof Straightforward, because $D$ is prime.  \qed

In the following we say in \eqref{Fig2} that $L$ \em{connects}
or \em{identifies} regions $A$ and $B$. The meaning comes from
looking at $B(D)$. There $A$ and $B$ correspond to loops, and the
attachment of $L$ has the effect of joining these loops.

\subsection{Controlling the number of attachments}

We need to extend lemma \reference{lmdf} to higher
valence attachments. For this we decisively use $\gm=4$ in
\eqref{tstcnd}. We need also the following piece of information. 

\begin{defi}
Assume $D$ is $A$-adequate. A loop in $A(D)$ is separating
if it is connected by crossing traces from both sides (see
\S\reference{S3c}). Define a
\em{block component} of $A(D)$ to be a region of the complement
of the set of loops in $A(D)$ together with the loops in its
boundary and all crossing traces within that region.

A connected (sum) component of a block component is defined
as for link diagrams: a block component decomposes as connected
sum along a curve intersecting it in two points on the loops,
and in no point on crossing traces. An $A$-state \em{atom} of
$D$ is defined to be a connected component of a block component
of $D$. We write $a_A(D)$ for the ($A$-state) \em{atom number}
of $D$, the number of its $A$-state atoms. Obviously, the same
definitions can be set up for the $B$-state. We will sometimes
write $a(D)$ if it is clear from the context if the $A$-state or
$B$-state is meant.
\end{defi}

A paraphrasing of $a_A(D)$ in terms of the $A$-state graph
$G=G(A(D))$ is that $a_A(D)=a(G)$, where $a(G)$ is the atom
number from definition \reference{SDE}.

A first illustration is the case when $D$ is positive.
Then the (separating) $A$-state loops are the (separating)
Seifert circles, the notion of block( component)s coincides
with the one of Cromwell \cite{Cromwell}, and ($A$-state)
atoms are connected components of blocks.

We visualize the notion also with the special case we will use
in our proof: the outer cycle in \eqref{oaz} is an ($A$-state)
atom of $D$, and each attachment \eqref{*} adds one new atom.
So in our case the atom number counts the number of attachments
(plus one). A decisive merit of this quantity is that it is
controllable from a link invariant. The following is proved in
\cite{adeq} using Thistlethwaite's work on the Kauffman polynomial%
\footnote{Note, though, that $a_A(D)$ is called $m(D)$ in
\cite{adeq}, and $a(D)$ is used there with a different meaning.}.

\begin{theorem}(see proposition 5.3 in \cite{adeq})\label{ozp}
For $A$-adequate diagrams $D$ of a link $L$, the number
$c(D)-a_A(D)$ is an invariant of $L$.
\end{theorem}

Then calculation shows

\begin{corr}\label{cr3l}
To construct $D_k$, we must perform exactly 3 attachments
\eqref{*} on the outer cycle loop. \qed
\end{corr}

This information will be used for the crossing numbers $15+4k$
only secondarily. A banal reason was that the author became
aware of it only after most of this proof was completed, and
corollary \ref{cr3l} is used only for minor simplifications
and fixes of the argument. However, the corollary became the
main motivation for attempting the case $17+4k$. Indeed, this
corollary becomes a major tool in that part, for without it
the proof would have been irrecordably complicated.

In particular, we have

\begin{corr}\label{cr3l'}
There are at least two attachments (of depth 1) inside an outer cycle
loop. In particular, there is at most one attachment outside an
outer cycle loop.
\end{corr}

\proof Were there no two attachments inside an outer cycle loop, 
we would have no cycle in $IG(A(D))$. \qed

The next simplifications need a definition.

\begin{defi}\label{cfiq}
Call a crossing of a diagram $D$ to be $A$-inadequate if its
trace in $A(D)$ connects the same loop. (So $D$ is $A$-adequate
if it has no $A$-inadequate crossing.) Similarly one defines
$B$-inadequate crossings.

Let $D'$ be obtained from $D$ by making all edges in $A(D)$
simple (i.e. of multiplicity 1).
We call an edge in $A(D)$ to be \em{inadequate} if its
trace in $D'$ comes from a $B$-inadequate crossing (of $D'$).
A loop of valence $2$, attached by \eqref{*} in $A(D)$, is
inadequate if its both legs are inadequate.
\end{defi}

An example, and simultaneously the important special case
we will use below, is as follows. Attach inside a loop $L$
another loop $M$, and let $R$ be a region outside $L$ touched
by a leg basepoint of $M$. If no other loop inside $L$
has a leg basepoint touching $R$, and no loop outside
$L$ has non-empty interior, then $M$ is inadequate.

\begin{lemma}\label{lmio'}
If $e$ has an inadequate edge of $A(D)$, then $e$ has multiplicity 2.
\end{lemma}

\proof Clearly $e\ge 2$ by $A$-adequacy of $D$. If $e$ has
multiplicity $\ge 3$, we have a cycle of $B(D)$ that consists
of simple connections between the loops only (dual to the
traces of the edge $e$). This must be then the only cycle.
However, then with corollary \reference{cr3l}, we see that by
3 attachments we cannot create then a cycle in $IG(B(D))$. \qed

\begin{lemma}\label{lmio}
If $D$ has an inadequate loop $P$, then $D$ cannot be one of our
hypothetic diagrams $D_k$.
\end{lemma}

\proof Let $x,y$ be the multiplicities of the legs of $P$ in $D$.
Clearly $x,y\ge 2$, and by connectivity, one is at least 3. But
this contradicts lemma \reference{lmio'}. \qed
% means that $B(D)$ has a cycle made up entirely of simple traces/edges.
% If $D=D_k$, then this would be the sole cycle. But one cannot
% attach 3 loops to such an outer cycle creating a cycle in $IG(B(D))$.

\begin{lemma}\label{lmZ}
There is no region touched in two different segments by legs of
attachments \eqref{attach} on the opposite side.
\end{lemma}

\proof The situation in which two such legs could occur is that
they connect to the same region on both sides of an outside
attachment, like 
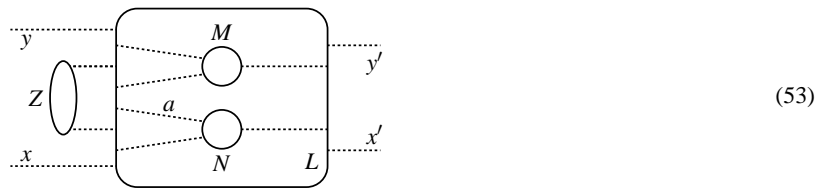
\begin{eqn}\label{pecu}
\small
\diag{7mm}{7}{3.8}{
 \pictranslate{0 0.4}{
  \picovalbox{4 1.5}{4 3.4}{0.4}{}
  {\piclinedash{0.05 0.05}{0.01}
   \picline{0 0.2}{2 0.2}
   \picline{0 2.8}{2 2.8}
   \picline{1 2.1}{2 2.1}
   {
   \picline{1 0.9}{2 0.9}
   \picline{2 0.5}{4 0.8}
   \picline{2 1.3}{4 1.0}
   \picline{6 0.9}{4 0.9}
   }
   \pictranslate{0 3}{
     \picscale{1 -1}{
       \picline{1 0.9}{2 0.9}
       \picline{2 0.5}{4 0.8}
       \picline{2 1.3}{4 1.0}
       \picline{6 0.9}{4 0.9}
     }
   }
   \picline{6 0.5}{7 0.5}
   \picline{6 2.5}{7 2.5}
  }
  \picfilledellipse{1 1.5}{0.25 0.7}{}
  \picfilledcircle{4 0.9}{0.37}{}
  \picfilledcircle{4 2.1}{0.37}{}
  \picputtext{0.3 0.4}{$x$}
  \picputtext{6.9 2.2}{$y'$}
  \picputtext{6.9 0.8}{$x'$}
  \picputtext{0.3 2.6}{$y$}
  \picputtext{0.49 1.5}{$Z$}
  \picputtext{3.0 1.35}{$a$}
  \picputtext{4 0.25}{$N$}
  \picputtext{4 2.75}{$M$}
  % \picputtext{2.2 0.2}{$L$}
  \picputtext{5.7 0.3}{$L$}
 }
}
\end{eqn}
If both segments of $L$ left and right of $Z$ are touched
by a leg of the same loop $M$ only, then we could reduce the
number of legs by a mutation, or an attachment inside $M$ that
prevents us from that would make $\gm>4$. Thus we can assume
they are touched by at least two, and then exactly two by
corollary \reference{cr3l}, loops $M,N$ inside $L$. This
also implies, by corollary \reference{cr3l'}, that the
analogous case when $Z$ is inside $L$, and $M,N$ would be
outside, does not occur. So we have \eqref{pecu}, with
the option that some of $a,x',y'$ may be zero.

Since we have an outside attachment $Z$, the other
two attachments are inside by corollary \reference{cr3l'}.
So $Z$ is the sole outside attachment, and there are no
attachments of depth $>1$. Note that $Z$ cannot have valence
more than 2. Otherwise we could not keep $\gm\le 4$ and
simultaneously make the diagram prime with the inside
attachments. But from \eqref{pecu} we see that $Z$ would
be inadequate, contradicting lemma \ref{lmio}.
\qed

In particular with $\gm=4$, we have

\begin{lemma}
There are no attached loops of valence different from 2 or 3.
\end{lemma}
 
\proof Valence 1 means a composite diagram. Also, we can assume
that traces of attached loops that connect to the same region
are parallel (because the diagram is prime). So valence
$\ge 4$ for a loop attached inside some other loop will
account in $\gm\ge 6$. If a loop $M$ of valence at least 4 is
attached outside to an outer cycle loop $L$, there are at least
4 intervals of $L$ to which traces of $M$ are attached from inside
(3 between the legs of $M$, and one more to make the diagram
prime). Then we have again $\gm\ge 6$.
\qed

\subsection{Properties of edge multiplicities}

As in \S\reference{G3}, we can assume up to flypes that 
\em{$D$ has at most one $\sim$-equivalence class of more than
2 crossings}. Since for the time being we trace only the number
of components, adequacy, and the invariants in \eqref{tstcnd},
we change the multiplicity of this class to $4$ or $5$ crossings,
preserving parity. The only possible change in the invariants in
\eqref{tstcnd} can occur if the multiplicity $3$ is increased to
$5$, and a triangle in $B(D)$ is destroyed. But we know that there
are no triangles anyway. So, for the sake of working with the
invariants in \eqref{tstcnd}, we can assume that
\begin{eqn}\label{lmt}
\mbox{\em{parallel equivalence classes have $1$, $2$, $4$ or
$5$ elements, and at most one class has $4$ or $5$ elements}.}
\end{eqn}
Only at the very final stage of the proof, in \S\reference{S9.6},
we will change $4$ and $5$ back to $2n$ and $2n+1$.

\begin{conv}\label{cthi}
Below we will draw in diagrams the $B$-state in thicker
(solid for loops and dashed for traces) lines than the
$A$-state. Often drawing $A$- and $B$-state in the same
diagram makes it hard to parse, though, and drawing separate
diagrams takes extra space we like to save. So we will several
times waive on drawing the $B$-state even when we need to
argue with it. The indicated examples should advise how to
obtain it in the other cases. In other situations, to enhance
legibility, we will draw only the loops of the $B$-state, but
not the traces, which are easy to reconstruct.
\end{conv}

\begin{lemma}\label{lemmaW}
Assume $D$ has an empty loop of valence $2$. 
Then the leg valences are $(1,1)$, $(1,4)$ or $(2,5)$.
\end{lemma}

\proof Multiplicities $(2,2)$ and $(2,4)$ are ruled out because of
connectivity. (The tangle has a closed component, so the diagram is
not of a knot.) $(1,5)$ is ruled out because in
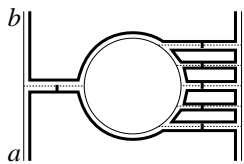
\[
\diag{4mm}{8}{5}{
    \picline{0.4 0}{0.4 5}
    \picline{7.6 0}{7.6 5}
    {\piclinedash{0.05 0.05}{0.01}
     \picmultigraphics{5}{0 0.68}{
       \picline{7.6 1.15}{4 1.15}
       {\piclinewidth{10}\piclinedash{0.25 0.25}{0.01}
	\picline{6.3 1.05}{6.3 1.25}
       }
     }
     \picline{0.4 2.5}{2.4 2.5}
    }
    {\piclinewidth{10}\piclinedash{0.25 0.25}{0.01}
     \picline{1.5 2.3}{1.5 2.7}
    }
    \picfilledcircle{4 2.5}{1.6}{}
    \picputtext{0.1 4.8}{$b$}
    \picputtext{0.1 0.2}{$a$}
    \piclinewidth{10}
    \picstroke{
      \picline{7.4 5}{7.4 3.99}
      \piclineto{5.0 3.99}
      \piccurveto{4.3 4.5}{2.6 4.5}{2.2 2.7}
      \piclineto{0.6 2.7}
      \piclineto{0.6 5}
      \picline{7.4 3.7}{7.4 3.3}
      \piclineto{5.6 3.3}
      \piclineto{5.3 3.7}
      \picPSgraphics{cp}
      \picline{7.4 3.1}{7.4 2.65}
      \piclineto{5.8 2.65}
      \piclineto{5.7 3.1}
      \picPSgraphics{cp}
      \picline{7.4 3.7 5 x -}{7.4 3.3 5 x -}
      \piclineto{5.6 3.3 5 x -}
      \piclineto{5.3 3.7 5 x -}
      \picPSgraphics{cp}
      \picline{7.4 3.1 5 x -}{7.4 2.65 5 x -}
      \piclineto{5.8 2.65 5 x -}
      \piclineto{5.7 3.1 5 x -}
      \picPSgraphics{cp}
      \picline{7.4 5 5 x -}{7.4 3.99 5 x -}
      \piclineto{5.0 3.99 5 x -}
      \piccurveto{4.3 4.5 5 x -}{2.6 4.5 5 x -}{2.2 2.7 5 x -}
      \piclineto{0.6 2.7 5 x -}
      \piclineto{0.6 5 5 x -}
    }
}
% Fig 9
\]
if $a=b$, the diagram is not adequate, else $B(D)$ has an even length
cycle.  \qed 

\begin{conv}
We denote an edge (parallel equivalence class of traces)
and its multiplicity by the same letter. For example `$x=2$'
means that the edge $x$ has multiplicity 2.
\end{conv}

\subsection{Restricting parallel loops}

We call loops \em{parallel} if they have both valence $2$ and
their legs connect the same regions.

\begin{lemma}\label{lemmaX}
Assume $D$ has a pair of parallel loops. Assume both loops are empty.
Then their valencies are $(1,1)$ and $(1,4)$ resp.
\end{lemma}

\proof By mutations we can assume that the box $B$ in 
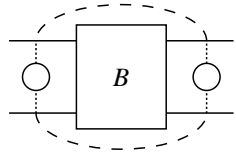
\[
\diag{6mm}{5}{3}{
  \picline{0 0.7}{5 0.7}
  \picline{0 2.3}{5 2.3}
  {\piclinedash{0.05 0.05}{0.01}
  \picline{0.6 0.7}{0.6 2.3}
  \picline{4.4 0.7}{4.4 2.3}
  }
  \picfilledbox{2.5 1.5}{2 2.3}{$B$}
  \picfilledcircle{0.6 1.5}{0.3}{}
  \picfilledcircle{4.4 1.5}{0.3}{}
  {\piclinedash{0.25 0.25}{0.01}
   \picellipsearc{2.5 2.3}{1.9 0.8}{0 180}
   \picellipsearc{2.5 0.7}{1.9 0.8}{180 0}
  }
}
% Fig 3
\]
is empty. (The dashed lines on top and bottom should indicate
that one can connect the basepoints of the legs of the two
loops without crossing any other traces or loops.)
Let $(a,b)$, $(c,d)$ be the multiplicities of the two leg pairs.
(We assume $a,b,c,d\in\{1,2,4,5\}$.) For any of both pairs, $(1,2)$
is ruled out:
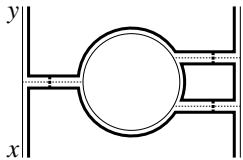
\[
\diag{4mm}{8}{5}{
    \picline{0.4 0}{0.4 5}
    \picline{7.6 0}{7.6 5}
    {\piclinedash{0.05 0.05}{0.01}
     \picline{0.4 2.5}{5 2.5}
     \picline{7.6 1.7}{5 1.7}
     \picline{7.6 3.3}{5 3.3}
    }
    \picfilledcircle{4 2.5}{1.6}{}
    \picputtext{0.1 4.8}{$y$}
    \picputtext{0.1 0.2}{$x$}
    \piclinewidth{10}
    \picstroke{
      \picline{7.4 5}{7.4 3.5}
      \piclineto{5.45 3.5}
      \picarcto{2.25 2.7}{1.75}
      \piclineto{0.6 2.7}
      \piclineto{0.6 5}
      \picline{7.4 3.1}{7.4 1.9}
      \piclineto{5.65 1.9}
      \picarcto{5.65 3.1}{1.7}
      \piclineto{7.4 3.1}
      \picPSgraphics{cp}
      \picline{7.4 5 5 x -}{7.4 3.5 5 x -}
      \piclineto{5.45 3.5 5 x -}
      \picarcto{2.25 2.7 5 x -}{-1.75}
      \piclineto{0.6 2.7 5 x -}
      \piclineto{0.6 5 5 x -}
    }
    {\piclinedash{0.1 0.05}{0.01}
     \picline{1.3 2.3}{1.3 2.7}
     \picline{6.7 3.1}{6.7 3.5}
     \picline{6.7 1.5}{6.7 1.9}
    }
}
% Fig 8
\]
if $x=y$, the diagram is not adequate, else $\tg>0$. Then
if $(a,b),(c,d)\ne (1,1)$, then any of those creates a cycle
in $B(D)$, so $\chi(B(D))<0$. So we can assume that w.l.o.g.
$(a,b)=(1,1)$, and we need to sort out only the possibilities
for $(c,d)$. The cases $(c,d)=(2,2),(2,4),(1,5)$, are ruled out
as before. The case $(2,5)$ is ruled out because in that case
$(a,b)=(1,1)$ accounts for the loops $x,y$ in
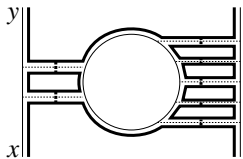
\[
% Fig 10
\diag{4mm}{8}{5}{
    \picline{0.4 0}{0.4 5}
    \picline{7.6 0}{7.6 5}
    {\piclinedash{0.05 0.05}{0.01}
     \picmultigraphics{5}{0 0.68}{
       \picline{7.6 1.15}{4 1.15}
       {\piclinewidth{10}\piclinedash{0.1 0.05}{0.01}
	\picline{6.3 1.05}{6.3 1.25}
       }
     }
     \picline{0.4 3.0}{2.4 3.0}
     \picline{0.4 2.0}{2.4 2.0}
    }
    {\piclinewidth{10}\piclinedash{0.1 0.05}{0.01}
     \picline{1.5 1.8}{1.5 2.2}
     \picline{1.5 2.8}{1.5 3.2}
    }
    \picfilledcircle{4 2.5}{1.6}{}
    \picputtext{0.1 4.8}{$y$}
    \picputtext{0.1 0.2}{$x$}
    \piclinewidth{10}
    \picstroke{
      \picline{7.4 5}{7.4 3.99}
      \piclineto{5.0 3.99}
      \piccurveto{4.3 4.5}{2.8 4.4}{2.4 3.2}
      \piclineto{0.6 3.2}
      \piclineto{0.6 5}
      \picline{7.4 3.7}{7.4 3.3}
      \piclineto{5.6 3.3}
      \piclineto{5.3 3.7}
      \picPSgraphics{cp}
      \picline{7.4 3.1}{7.4 2.65}
      \piclineto{5.8 2.65}
      \piclineto{5.7 3.1}
      \picPSgraphics{cp}
      \picline{7.4 3.7 5 x -}{7.4 3.3 5 x -}
      \piclineto{5.6 3.3 5 x -}
      \piclineto{5.3 3.7 5 x -}
      \picPSgraphics{cp}
      \picline{7.4 3.1 5 x -}{7.4 2.65 5 x -}
      \piclineto{5.8 2.65 5 x -}
      \piclineto{5.7 3.1 5 x -}
      \picPSgraphics{cp}
      \picline{7.4 5 5 x -}{7.4 3.99 5 x -}
      \piclineto{5.0 3.99 5 x -}
      \piccurveto{4.3 4.5 5 x -}{2.8 4.4 5 x -}{2.4 3.2 5 x -}
      \piclineto{0.6 3.2 5 x -}
      \piclineto{0.6 5 5 x -}
      \picline{0.6 2.2}{2.3 2.2}
      \piccurveto{2.25 2.4}{2.25 2.6}{2.3 2.8}
      \piclineto{0.6 2.8}
      \picPSgraphics{cp}
    }
}
\]
being different, and so again $B(D)$ has two different cycles. \qed

\begin{lemma}\label{l9.10}
Assume $D$ has a pair of parallel loops. Then both loops are empty.
\end{lemma}

\proof We can assume again that we have
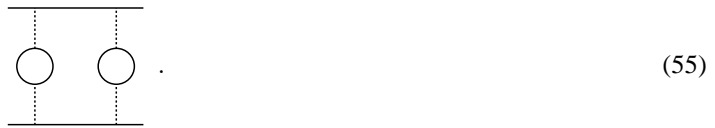
\begin{eqn}\label{Fig4}
\diag{6mm}{3}{4}{
  \picline{0 0.7}{3 0.7}
  \picline{0 3.3}{3 3.3}
  {\piclinedash{0.05 0.05}{0.01}
  \picline{0.6 0.7}{0.6 3.3}
  \picline{2.4 0.7}{2.4 3.3}
  }
  \picfilledcircle{0.6 2.0}{0.4}{}
  \picfilledcircle{2.4 2.0}{0.4}{}
}\ \ .
% Fig 4
\end{eqn}
We can assume now using corollary \ref{cr3l} that at least
one of the loops is non-empty, and that \eqref{Fig4} is a
pair of such loops of maximal depth.
% We consider a pair of such loops of maximal depth.

Now if a loop is non-empty, we can flype its interor out by
\eqref{Fig1} (and have a contradiction to the flatness condition),
unless we have a super-configuration of
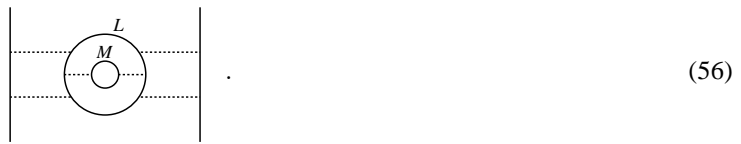
\begin{eqn}\label{Fig5}
{\scriptsize
\diag{6mm}{5}{3}{
  \picline{0.4 0}{0.4 3}
  \picline{4.6 0}{4.6 3}
  {\piclinedash{0.05 0.05}{0.01}
   \picline{1 0.4 x}{1 4.6 x}
   \picline{2 0.4 x}{2 4.6 x}
  }
  \picfilledcircle{2.5 1.5}{0.9}{}
  {\piclinedash{0.05 0.05}{0.01}
   \picline{1.6 1.5}{3.4 1.5}
  }
  \picfilledcircle{2.5 1.5}{0.3}{}
  \picputtext{2.5 2.0}{$M$}
  \picputtext{2.8 2.6}{$L$}
}}\ .
% Fig 5
\end{eqn}
($X$ is a super-configuration of $Y$, if $Y$ arises from
$X$ be deleting some number of, possibly no, legs.)
If such a situation occurs in both loops in \eqref{Fig4},
then we delete one of the $4$ legs in each copy of \eqref{Fig5},
flype the interor out, and see $\gm>4$, since the legs of the
(low empty loops) $L$, $L'$ identify $4$ different regions
$A$, $B$, $C$, $D$.
{\small
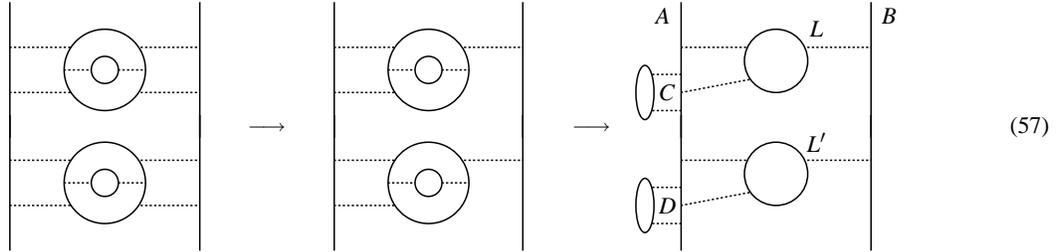
\begin{eqn}\label{Fig6}
\diag{6mm}{5}{5.5}{
  \picmultigraphics{2}{0 2.5}{
    \picline{0.4 0}{0.4 3}
    \picline{4.6 0}{4.6 3}
    {\piclinedash{0.05 0.05}{0.01}
     \picline{1 0.4 x}{1 4.6 x}
     \picline{2 0.4 x}{2 4.6 x}
    }
    \picfilledcircle{2.5 1.5}{0.9}{}
    {\piclinedash{0.05 0.05}{0.01}
     \picline{1.6 1.5}{3.4 1.5}
    }
    \picfilledcircle{2.5 1.5}{0.3}{}
  }
}
\quad\lra\quad
\diag{6mm}{5}{5.5}{
  \picmultigraphics{2}{0 2.5}{
    \picline{0.4 0}{0.4 3}
    \picline{4.6 0}{4.6 3}
    {\piclinedash{0.05 0.05}{0.01}
     \picline{1 0.4 x}{1 2.5 x}
     \picline{2 0.4 x}{2 4.6 x}
    }
    \picfilledcircle{2.5 1.5}{0.9}{}
    {\piclinedash{0.05 0.05}{0.01}
     \picline{1.6 1.5}{3.4 1.5}
    }
    \picfilledcircle{2.5 1.5}{0.3}{}
  }
}
\quad\lra\quad
\diag{6mm}{5}{5.5}{
  \pictranslate{0.5 0}{
  \picmultigraphics{2}{0 2.5}{
    \picline{0.4 0}{0.4 3}
    \picline{4.6 0}{4.6 3}
    {\piclinedash{0.05 0.05}{0.01}
     \picline{1 0.4 x}{1.4 2.5 x}
     \picline{2 0.4 x}{2 4.6 x}
    }
    \picfilledcircle{2.5 1.7}{0.7}{}
    \pictranslate{-0.6 -0.5}{
      {\piclinedash{0.05 0.05}{0.01}
       % \piccurve{1 2.1}{0.1 1.9}{0.1 1.1}{1 0.9}
       \picline{1 1.9}{0.2 1.9}
       \picline{1 1.1}{0.2 1.1}
      }
      \picfilledellipse{0.2 1.5}{0.2 0.6}{}
    }
  }
  \picputtext{3.4 2.4}{$L'$}
  \picputtext{3.4 4.9}{$L$}
  \picputtext{0.1 1.0}{$D$}
  \picputtext{0.1 3.5}{$C$}
  \picputtext{0.0 5.2}{$A$}
  \picputtext{5.0 5.2}{$B$}
  }
}
% Fig 6
\end{eqn}
}
Thus we can assume that exactly one of the loops in \eqref{Fig4}
is a super-configuration of \eqref{Fig5}. Deleting a leg, and
flyping the interor out, as in \eqref{Fig6}, we see that $\gm=4$.
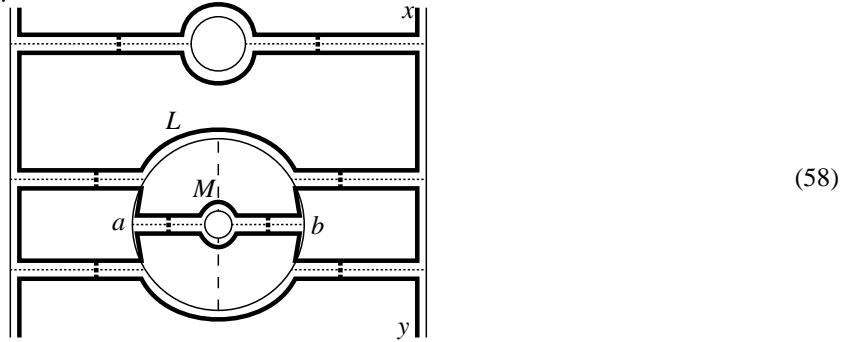
\begin{eqn}\label{Fig7}
\diag{6mm}{10}{7}{
    \picline{0.4 0}{0.4 7.3}
    \picline{9.6 0}{9.6 7.3}
    {\piclinedash{0.05 0.05}{0.01}
     \picline{1.5 0.4 x}{1.5 9.6 x}
     \picline{3.5 0.4 x}{3.5 9.6 x}
     \picline{6.5 0.4 x}{6.5 9.6 x}
    }
    \picfilledcircle{5 2.5}{1.9}{}
    {\piclinedash{0.05 0.05}{0.01}
     \picline{3.1 2.5}{6.9 2.5}
    }
    {\piclinedash{0.25 0.25}{0.01}
     \picline{5 0.6}{5 4.4}
    }
    \picfilledcircle{5 2.5}{0.3}{}
    \picfilledcircle{5 6.5}{0.6}{}
    \picputtext{4.0 4.8}{$L$}
    \picputtext{4.7 3.3}{$M$}
    \picputtext{2.8 2.5}{$a$}
    \picputtext{7.2 2.5}{$b$}
    \picputtext{9.2 7.2}{$x$}
    \picputtext{9.1 0.2}{$y$}
    \piclinewidth{10}
    \picstroke{
      \picline{0.6 7.3}{0.6 6.7}
      \piclineto{4.2 6.7}
      \piccurveto{4.3 7.6}{5.7 7.6}{5.8 6.7}
      \piclineto{9.4 6.7}
      \piclineto{9.4 7.3}
      \picline{0.6 3.7}{0.6 6.3}
      \piclineto{4.2 6.3}
      \piccurveto{4.3 5.4}{5.7 5.4}{5.8 6.3}
      \piclineto{9.4 6.3}
      \piclineto{9.4 3.7}
      \piclineto{6.7 3.7}
      \piccurveto{6.1 4.9}{3.9 4.9}{3.3 3.7}
      \piclineto{0.6 3.7}
      \picPSgraphics{cp}
      \picline{0.6 0}{0.6 1.3}
      \piclineto{3.3 1.3}
      \piccurveto{3.9 0.1}{6.1 0.1}{6.7 1.3}
      \piclineto{9.4 1.3}
      \piclineto{9.4 0}
      \picline{0.6 1.7}{0.6 3.3}
      \piclineto{3.3 3.3}
      \piclineto{3.2 2.7}
      \piclineto{4.6 2.7}
      \piccurveto{4.8 3.1}{5.2 3.1}{5.4 2.7}
      \piclineto{6.8 2.7}
      \piclineto{6.7 3.3}
      \piclineto{9.4 3.3}
      \piclineto{9.4 1.7}
      \piclineto{6.7 1.7}
      \piclineto{6.8 2.3}
      \piclineto{5.4 2.3}
      \piccurveto{5.2 1.9}{4.8 1.9}{4.6 2.3}
      \piclineto{3.2 2.3}
      \piclineto{3.3 1.7}
      \piclineto{0.6 1.7}
      \picPSgraphics{cp}
    }
    {\piclinedash{0.1 0.05}{0.01}
     \picline{2.8 6.3}{2.8 6.7}
     \picline{7.2 6.3}{7.2 6.7}
     \picline{2.3 3.3}{2.3 3.7}
     \picline{7.7 3.3}{7.7 3.7}
     \picline{2.3 1.3}{2.3 1.7}
     \picline{7.7 1.3}{7.7 1.7}
     \picline{3.9 2.3}{3.9 2.7}
     \picline{6.1 2.3}{6.1 2.7}
    }
}
% \qquad
% \diag{6mm}{10}{7}{
%     \picline{0.4 0}{0.4 5.3}
%     \picline{9.6 0}{9.6 5.3}
%     {\piclinedash{0.05 0.05}{0.01}
%      \picline{1.5 0.4 x}{1.5 9.6 x}
%      \picline{3.5 0.4 x}{3.5 9.6 x}
%      \picline{3.2 2.9}{2.5 2.9}
%      \picline{3.2 2.1}{2.5 2.1}
%     }
%     \picfilledcircle{5 2.5}{1.9}{}
%     {\piclinedash{0.05 0.05}{0.01}
%      \picline{3.1 2.5}{6.9 2.5}
%     }
%     \picfilledcircle{5 2.5}{0.3}{}
%     \picfilledellipse{2.5 2.5}{0.2 0.6}{}
% }
% Fig 7
\end{eqn}
Adding a leg at some of the vertical dashed positions,
joins two loops in $B(D)$, so $\gm$ gets $>4$.
% (Note that the
% third to-attach loop may be a valence-2 loop parallel to $M$.
% In that case, its addition does not reduce $\gm$, so the argument
% that $\gm>4$ remains correct.)

By corollary \reference{cr3l} we may
assume that the interor of $M$ is empty and $M$ has no
parallel loop, and no loop is attached outside to $L$.
Let $a$, $b$ be the multiplicities of the legs.
If one of $a$ or $b$ is single, then $D$ is not adequate.
If $a=b=2$, then $D$ is not a knot diagram. So $\max(a,b)\ge 3$.
But then \eqref{Fig7} contains one cycle in $B(D)$, and there
must be another one, that links $x$ and $y$ via the rest of
the diagram ($x\ne y$ because $\tg=0$), because $D$ is prime.
So $B(D)$ has at least two cycles, a contradiction.
% 
% The right fragment in \eqref{Fig7} is excluded with the same argument.
% 
% Assume second $L$ has one parallel loop $L'$. By maximal depth
% assumption in \eqref{Fig4}, both $L$ and $L'$ have empty interor. Then
% by lemma \reference{lemmaX}, again one of $L$ and $L'$ has valency
% $(1,4)$, so again $B(D)$ has a cycle inside \eqref{Fig7}, with the
% same contradiction. 
% 
% Assume finally $L$ has non-empty interor. Then one has a nested
% situation like \eqref{Fig5}, which is easily observed to have $\gm>4$.
\qed

\begin{lemma}
There is no triple of parallel (empty) loops, and no two different
pairs.
\end{lemma}

\proof Use $\chi(G(B(D)))=0$ and corollary \reference{cr3l}. \qed

\begin{lemma}\label{lempar}
There exist no parallel (empty) loops in $D$.
\end{lemma}

\proof In order to create a cycle in $IG(A)$ with three attachments
(corollary \reference{cr3l}), a parallel pair $X,Y$ of loops in $D$
must be attached (with depth 1) inside an outer cycle loop $L$, and
intertwine at least two multiple edges. (We assume $X,Y$ are
empty by lemma \reference{l9.10}.) In particular, $X,Y$ do not
identify any of the regions $U$ or $V$ in \eqref{oaz}.
Now the third attachment is supposed to establish $\gm=4$ and
not be an inadequate loop (by lemma \reference{lmio}). 

Since $X,Y$ do not identify $U$ or $V$, this means that a loop
of valence $2$ attached within an outer cycle loop $L'$ different
from $L$ would be inadequate. The condition $\gm=4$ also
rules out the option of attaching the third loop outside to $L$,
or attaching a loop of valence 3 inside some $L'\ne L$.
Then the only option is
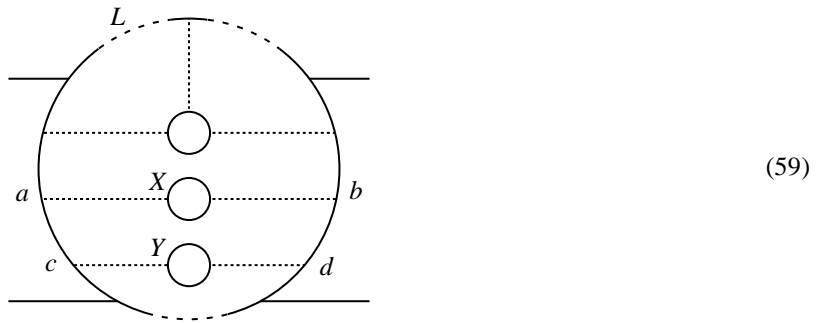
\begin{eqn}\label{lDD}
\diag{8mm}{6}{5}{
  \pictranslate{3 2.5}{
    \picline{-3 1.5}{3 1.5}
    \picline{-3 -2.2}{3 -2.2}
    \picfillnostroke{\piccircle{0 0}{2.5}{}}
    \picclip{\piccircle{0 0}{2.5}{}}{
      {\piclinedash{0.05 0.05}{0.01}
       \picline{0 1}{0 2.5}
      }
      \picmultigraphics{3}{0 1.1}{
      {\piclinedash{0.05 0.05}{0.01}
       \picline{-2.5 -1.6}{2.5 -1.6}
      }
      \picfilledcircle{0 -1.6}{0.35}{}
      }
   }
   \picputtext{2.8 -7 polar}{$b$}
   \picputtext{2.8 -35 polar}{$d$}
   \picputtext{2.8 189 polar}{$a$}
   \picputtext{2.8 215 polar}{$c$}
   \picputtext{2.8 115 polar}{$L$}
   \picputtext{-.5 -1.3}{$Y$}
   \picputtext{-.5 -0.2}{$X$}
   \piccirclearc{0 0}{2.5}{-75 53}
   \piccirclearc{0 0}{2.5}{127 255}
   \piccirclearc{0 0}{2.5}{83 97}
   {\piclinedash{0.1 0.15}{0.05}
    \piccirclearc{0 0}{2.5}{53 83}
    \piccirclearc{0 0}{2.5}{97 127}
    \piccirclearc{0 0}{2.5}{255 285}
   }
 }
}
\end{eqn}
(The dashed pairs of the loop indicate that something may
be attached to it there from the outside.)
Now, by connectivity reasons, not all of $a,b,c,d$ are $1$,
and then we have two cycles in $G(B(D))$. \qed

\subsection{Ruling out nesting of loops}

\begin{lemma}
All valence 2-loops of depth $1$, or such of depth $0$
outside the outer cycle, are empty.
\end{lemma}

\proof Assume $L$ is a valence 2-loop with non-empty interior.
As before, we have a similar picture to \eqref{Fig5}: 
%We draw it again:
\begin{eqn}\label{Fig18}
\diag{6mm}{10}{5}{
    \picline{0.4 0}{0.4 5}
    \picline{9.6 0}{9.6 5}
    {\piclinedash{0.05 0.05}{0.01}
     \picline{1.5 0.4 x}{1.5 9.6 x}
     \picline{3.5 0.4 x}{3.5 9.6 x}
    }
    \picfilledcircle{5 2.5}{1.9}{}
    {\piclinedash{0.05 0.05}{0.01}
     \picline{3.1 2.5}{6.9 2.5}
    }
    \picfilledcircle{5 2.5}{0.3}{}
    \picputtext{4.8 4.9}{$L$}
    \picputtext{4.8 3.3}{$M$}
    \picputtext{2.9 2.5}{$\ap$}
    \picputtext{7.1 2.5}{$\bt$}
    \picputtext{0.1 4.8}{$a$}
    \picputtext{0.1 0.2}{$c$}
    \picputtext{9.1 2.7}{$b$}
    \piclinewidth{10}
    \picstroke{
      \picline{9.4 5}{9.4 3.7}
      \piclineto{6.7 3.7}
      \piccurveto{6.1 4.9}{3.9 4.9}{3.3 3.7}
      \piclineto{0.6 3.7}
      \piclineto{0.6 5}
      \picline{0.6 0}{0.6 1.3}
      \piclineto{3.3 1.3}
      \piccurveto{3.9 0.1}{6.1 0.1}{6.7 1.3}
      \piclineto{9.4 1.3}
      \piclineto{9.4 0}
      \picline{0.6 1.7}{0.6 3.3}
      \piclineto{3.3 3.3}
      \piclineto{3.2 2.7}
      \piclineto{4.6 2.7}
      \piccurveto{4.8 3.1}{5.2 3.1}{5.4 2.7}
      \piclineto{6.8 2.7}
      \piclineto{6.7 3.3}
      \piclineto{9.4 3.3}
      \piclineto{9.4 1.7}
      \piclineto{6.7 1.7}
      \piclineto{6.8 2.3}
      \piclineto{5.4 2.3}
      \piccurveto{5.2 1.9}{4.8 1.9}{4.6 2.3}
      \piclineto{3.2 2.3}
      \piclineto{3.3 1.7}
      \piclineto{0.6 1.7}
      \picPSgraphics{cp}
    }
    {\piclinedash{0.1 0.05}{0.01}
     \picline{2.3 3.3}{2.3 3.7}
     \picline{7.7 3.3}{7.7 3.7}
     \picline{2.3 1.3}{2.3 1.7}
     \picline{7.7 1.3}{7.7 1.7}
     \picline{3.9 2.3}{3.9 2.7}
     \picline{6.1 2.3}{6.1 2.7}
    }
}\,.
% Fig 18
\end{eqn}
(Here $\ap,\bt$ denote the leg multiplicities of $M$.)
If $a\ne c$, then we argued to have two cycles in $B(D)$
(see the proof of lemma \reference{lemmaX}). So $a=c$.
Then again $\gm=4$.

Assume first $M$ has no parallel loop inside $L$.
Now, regardless of the way the $B(D)$ state loops are joined
by possible $A(D)$-state loops attached to $M$,
for an adequate diagram we must have $\ap,\bt\ge 2$.
For connectivity reasons then again $\max(\ap,\bt)\ge 3$,
so we have one cycle in $B(D)$. Now (again for both ways
the loops are connected inside $M$), for connectivity reasons,
one of the outer 4 legs of $L$ must be multiple (else
the tangle in \eqref{Fig18} contains a closed component).
But then, since $b\ne a=c$, we have a second cycle in $B(D)$.

If $M$ has a parallel loop inside $L$, we have by lemma
\reference{lemmaX} an edge of multiplicity $4$ inside
$L$, which accounts for one cycle in $B(D)$, and the
argument with the outer legs of $L$ applies as before
to yield a second cycle in $B(D)$. \qed

\begin{lemma}
Assume $M$ is a valence 3-loop of depth $\le 1$
(but not an outer cycle loop if of depth 0). Then $M$ is empty.
\end{lemma}

\proof Assume $M$ has a non-empty interior.
By flatness we can assume we have a loop $N$ inside $M$ like
\begin{eqn}\label{Fig16}
\raisebox{4mm}{
\diag{8mm}{4}{3}{
  \pictranslate{2 1}{
    {\piclinedash{0.05 0.05}{0.01}
     \picmultigraphics[rt]{3}{120}{
       \picline{1 -30 polar}{2 -30 polar}
       \picline{1 30 polar}{0.25 30 polar}
     }
    }
    \piccircle{0 0}{1.0}{}
    \piccircle{0 0}{0.25}{}
    \picputtext{-0.3 -0.4}{$N$\,}
    \picputtext{1.2 0.3}{\,$M$}
  }
}}\ .
% Fig 16
\end{eqn}

Observe that
\begin{eqn}\label{Fig17}
\gm\ \left(\ \ 
\diag{5mm}{4}{3.5}{
  \pictranslate{2 1.5}{
    {\piclinedash{0.05 0.05}{0.01}
     \picmultigraphics[rt]{3}{120}{
       \picline{1 -30 polar}{2 -30 polar}
       \picline{1 30 polar}{0.25 30 polar}
     }
    }
    \picmultigraphics[rt]{3}{120}{
      \picline{2.1 -50 polar}{2.1 -10 polar}
    }
    \piccircle{0 0}{1.0}{}
    \piccircle{0 0}{0.25}{}
  }
}
\ \ \right)\quad=\quad
\gm\ \left(\ \ 
\diag{5mm}{4}{3.5}{
  \pictranslate{2 1.5}{
    {
     \picmultigraphics[rt]{3}{120}{
       \picline{2.1 -50 polar}{2.1 -10 polar}
     }
    }
  }
}
\ \ \right)\ +\ 4\,,
% Fig 17
\end{eqn}
and remember that $\gm(D)\le 4$.

If $M$ has depth 0, then removing it with its interor results
in a diagram that still has a depth 1 loop (by primeness), and
\eqref{Fig17} shows $\gm\ge 6$. So assume that $M$ is a depth-one loop.

First we argue that we have no loop attached to $M$ from
the exterior. Namely, such a loop must be intertwined with
some loop in $M$'s interior (because $D$ is prime), but then
$\gm\ge 6$. 

We claim next that all loops attached inside $M$ are empty.
If a valence $2$-loop $N'$ attached inside $M$ is non-empty.
then (by deleting some legs if necessary), we flype its interior
out, and have a loop attached outside to $M$, and $\gm\ge
6$. If we have a non-empty valence $3$-loop attached inside
$M$, we find by flatness a fragment \eqref{Fig16} inside $M$.
Using \eqref{Fig17} and that $M$ remains of depth $1$, we have
$\gm\ge 6$.

This means that $IG(A)$ has no cycle resulting from edges inside
$M$.

Also, from \eqref{Fig17} we see that removing $M$ and its
interior we have a diagram with $\gm=0$, i.e. no loops of
non-zero depth.  Since $D$ is non-alternating and has loops of
positive depth, in particular, $M$ must be the only loop attached
inside an outer cycle loop. But then one easily sees that
$IG(A)$ is a forest, and so $\chi(IG(A))>0$. \qed

We should note the following important implications of the
last two lemmas.

\begin{corr}
There are no loops attached (neither from the interior, nor from the
exterior) to depth 1 loops, or to depth 0 loops
which are not outer cycle loops.
\end{corr}

\proof For the interior this is proved in the previous two
lemmas. For the exterior this follows then because $D$ is prime. \qed

\begin{corr}
There are no loops of depth $>1$ in $D$. \qed
\end{corr}

% Every 
\begin{lemma}\label{l1ne}
Only one loop of the outer cycle is non-empty.
\end{lemma}

\proof % We proved that we can achieve this if we have parallel
% loops. So assume we have no such.
Assume two loops $X$ and $Y$ of the outer cycle are non-empty.
Let $L$ and $M$ be loops attached to $X$ and $Y$ inside.
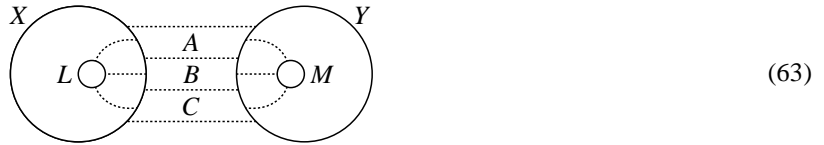
\begin{eqn}\label{eqw}
\diag{6mm}{8}{3}{
    {\piclinedash{0.05 0.05}{0.01}
      \picmultigraphics{4}{0 0.7}{
	\picline{1.5 0.45}{6.5 0.45}
      }
    }
    \picfilledcircle{1.5 1.5}{1.5}{}
    \pictranslate{4 1.5}{
      \picmultigraphics[S]{2}{-1 1}{
	\picfilledcircle{2.5 0}{1.5}{}
	\pictranslate{2.5 0}{
	  \piclinedash{0.05 0.05}{0.01}
	  \piccurve{1.5 150 polar}{-0.0 0.9}{-0.0 -0.9}{1.5 210 polar}
	  \picline{1.5 180 polar}{-0.1 0}
	}
	\picfilledcircle{2.2 0}{0.3}{}
      }
      \picputtext{0 0}{$B$}
      \picputtext{0 0.7}{$A$}
      \picputtext{2.9 0}{$M$}
      \picputtext{-2.8 0}{$L$}
      \picputtext{0 -0.7}{$C$}
      \picputtext{3.8 1.3}{$Y$}
      \picputtext{-3.8 1.3}{$X$}
    }
}
% Fig 13
\end{eqn}
If, say, $L$ has valence $\ge 3$, then it is $3$, and we have
$\gm=4$ already after attaching $L$. Let $A,B,C$ be the regions
connected by the legs of $L$. Then $M$ must connect the
same triple (if valence 3) or 2 of the three (if valence 2)
of $A,B,C$. If $M$ has valence 2, one can mutate it into
$X$. If $M$ has valence 3, then $L$ and $M$ are both empty
inside. Else one could as before, after possibly deleting some
legs, flype out their interor and see that they do not connect
the same triple of regions, so $\gm>4$. Also, since the length
of the outer cycle is $>2$, the coincidence of the 3 regions
connected forces $L$ and $M$ to be intertwined with at most one
multiple edge in the outer cycle. Then $IG(A)$ is a forest,
and $\chi(IG(A))>0$, a contradiction. The same argument applies
if inside $X$ or $Y$ several valence-2 loops are attached that
identify three different regions $A,B,C$.

So we can assume that loops attached inside $X$ and $Y$
identify only two pairs of regions. In particular there is
only one such loop for each $X$ and $Y$, because two such
are parallel, which we excluded with lemma \reference{lempar}.
So $L$ and $M$ are the single 
loops attached inside $X$ and $Y$, and have both valence 2.
Let $(A,B)$ and $(C,D)$ be the region pairs their legs connect.
(Pairs are unordered.) If $(A,B)=(C,D)$ we can move by mutation
$L$ and $M$ into $X$, and see that $L$ and $M$ are parallel,
again in contradiction to lemma \reference{lempar}. 

So assume $|\{A,B,C,D\}|\ge 3$.
Then $\gm=4$ after installing $L$ and $M$. If $L$ or $M$ are
non-empty, flatness forces again fragments like \eqref{Fig5}.
Then delete a leg and flype to see that $\gm>4$. So $L$ and $M$
are empty. Then we have some subconfiguration of (i.e. obtained by
possibly deleting some depth-1-loops from) one of the following:
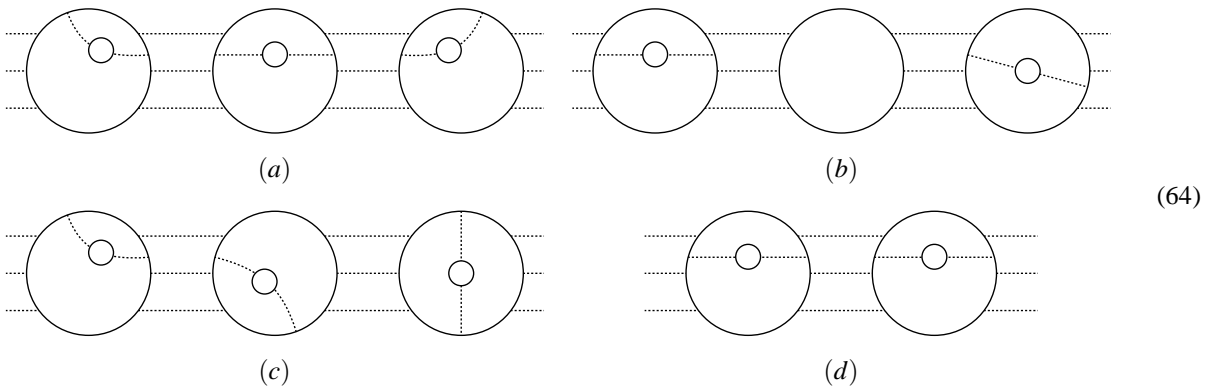
\begin{eqn}\label{Fig14}
\begin{array}{cc}
\diag{5.5mm}{13}{3}{
  \pictranslate{0.5 0}{
    {\piclinedash{0.05 0.05}{0.01}
      \picmultigraphics{3}{0 0.9}{
	\picline{-0.5 0.6}{12.5 0.6}
      }
    }
    \picmultigraphics{3}{4.5 0}{
      \picfilledcircle{1.5 1.5}{1.5}{}
    }
    \pictranslate{1.5 1.5}{
      \gdef\A{
      {\piclinedash{0.05 0.05}{0.01}
       \piccurve{1.5 110 polar}{-0.2 0.6}{0.3 0.3}{1.5 15 polar}
      }
      \picfilledcircle{0.3 0.5}{0.3}{}
      }\A
    }
    \pictranslate{6 1.5}{
      \gdef\B{
      {\piclinedash{0.05 0.05}{0.01}
       \picline{1.5 165 polar}{1.5 15 polar}
      }
      \picfilledcircle{0 0.4}{0.3}{}
      }\B
    }
    \pictranslate{10.5 1.5}{
      \gdef\C{\picscale{-1 1}{
	\A
      }}\C
    } 
  }
} &      
\diag{5.5mm}{13}{3}{
  \pictranslate{0.5 0}{
    {\piclinedash{0.05 0.05}{0.01}
      \picmultigraphics{3}{0 0.9}{
	\picline{-0.5 0.6}{12.5 0.6}
      }
    }
    \picmultigraphics{3}{4.5 0}{
      \picfilledcircle{1.5 1.5}{1.5}{}
    }
    \pictranslate{1.5 1.5}{
      \B
    }
    \pictranslate{10.5 1.5}{
     {\piclinedash{0.05 0.05}{0.01}
      \picline{1.5 165 polar}{1.5 -15 polar}
     }
     \picfilledcircle{0 0}{0.3}{}
    } 
  }
} \\[9mm]
(a) & (b) \\[3mm]
\diag{5.5mm}{13}{3}{
  \pictranslate{0.5 0}{
    {\piclinedash{0.05 0.05}{0.01}
      \picmultigraphics{3}{0 0.9}{
	\picline{-0.5 0.6}{12.5 0.6}
      }
    }
    \picmultigraphics{3}{4.5 0}{
      \picfilledcircle{1.5 1.5}{1.5}{}
    }
    \pictranslate{1.5 1.5}{
      \A
    }
     \pictranslate{6 1.5}{
       \picscale{-1 -1}{
         {\piclinedash{0.05 0.05}{0.01}
          \piccurve{1.5 110 polar}{-0.2 0.4}{0.3 -0.1}{1.5 -15 polar}
         }
         \picfilledcircle{0.25 0.2}{0.3}{}
       }
     }
    \pictranslate{10.5 1.5}{
      {\piclinedash{0.05 0.05}{0.01}
       \picline{1.5 90 polar}{1.5 -90 polar}
      }
      \picfilledcircle{0 0}{0.3}{}
    } 
  }
} &      
\diag{5.5mm}{9.5}{3}{
  \pictranslate{0.5 0}{
    {\piclinedash{0.05 0.05}{0.01}
      \picmultigraphics{3}{0 0.9}{
	\picline{-0.5 0.6}{9 0.6}
      }
    }
    \picmultigraphics{2}{4.5 0}{
      \picfilledcircle{2 1.5}{1.5}{}
    }
    \pictranslate{2 1.5}{
      \B
    }
    \pictranslate{6.5 1.5}{
      \B
    } 
  }
} \\[9mm]
(c) & (d) \\
\end{array}
% Fig 14
\end{eqn}
where we may add a depth-$0$ loop for some of the depth-$1$ loops like
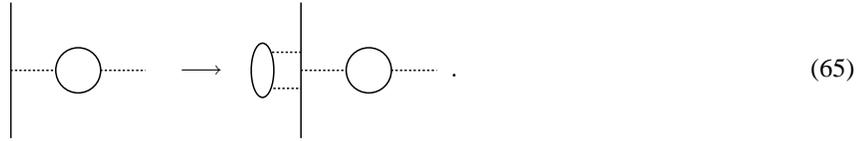
\begin{eqn}\label{Fig15}
\diag{6mm}{3}{3}{
  \picline{0 0}{0 3}
  {\piclinedash{0.05 0.05}{0.01}
   \picline{0 1.5}{3 1.5}
  }
  \picfilledcircle{1.5 1.5}{0.5}{}
}
\quad\lra\quad
\diag{6mm}{4}{3}{
  \picline{1 0}{1 3}
  {\piclinedash{0.05 0.05}{0.01}
  \picline{1 1.9}{0.1 1.9}
  \picline{0.1 1.1}{1 1.1}
  \picline{1 1.5}{4 1.5}
  }
  \picfilledellipse{0.15 1.5}{0.25 0.6}{}
  \picfilledcircle{2.5 1.5}{0.5}{}
}\ \ .
% Fig 15
\end{eqn}
However, we see that in all cases of \eqref{Fig14}, we have
$\chi(IG(A))>0$, and \eqref{Fig15} does not change
$\chi(IG(A))$. \qed
% loop into 

\begin{lemma}\label{l9,15}
The non-empty outer cycle loop has at least two loops attached inside.
\end{lemma}

\proof The graph $IG(A)$, which most be non-empty, must have a cycle.
There is no way of creating a cycle with edges of depth $>1$ (because
there are no such), and if only one loop is attached inside the
non-empty outer cycle loop, then $IG(A)$ cannot have a cycle at all.
\qed

\subsection{Diagram patterns\label{SPAT}}

The preparations so far yield a very restricted type of
combinatorial shapes of $A(D)$, which we call below \em{patterns}.
Formally, a pattern is a set of $A$-state diagrams, which are
obtainable from a given one by replacing (traces of) simple edges
$e$ by a number of parallel traces (i.e. a multiple edge),
including 0 (i.e. deletion of $e$).

\begin{figure}
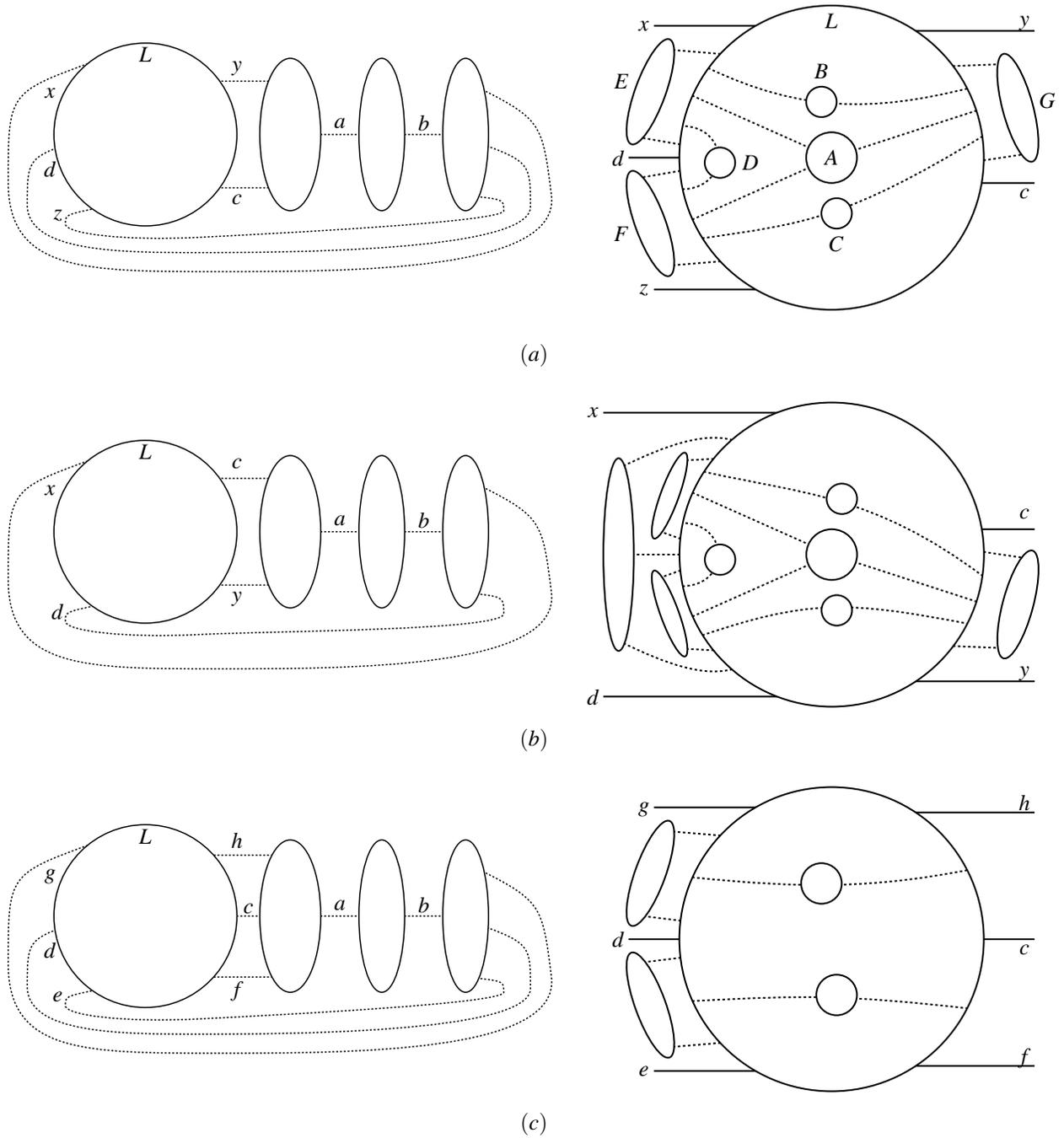
%[h]
\begin{center}
\[
\begin{array}{cc}
\diag{6mm}{15}{6}{
  {\piclinedash{0.05 0.05}{0.01}
    \picstroke{
      \opencurvepath{12 2}{13 1.9}{13 1.3}{7 1}{3 0.8}{1.5 1}{1.5 1.5}
        {2.5 1.7}{}
      \opencurvepath{12 3.5}{13.7 2.9}{13.7 0.5}{7 0.5}
        {3 0.5}{0.5 0.5}{0.5 3.0}{2 3.5}{}
      \opencurvepath{12 5}{14.0 4}{14.5 0}{9 0}{2 0}{0 0}{0 4.7}
        {2.2 5.5}{}
      \picline{4 5}{7.6 5}
      \picline{4 2.2}{7.6 2.2}
      \picline{8 3.6}{9.6 3.6}
      \picline{10 3.6}{11.6 3.6}
    }
  }
  \pictranslate{3.6 3.6}{
    \picfilledcircle{0 0}{2.4}{}
    \picfilledellipse{3.8 0}{0.8 2.0}{}
    \picfilledellipse{6.2 0}{0.6 2.0}{}
    \picfilledellipse{8.4 0}{0.6 2.0}{}
  }
  \picputtext{1.1 4.7}{$x$}
  \picputtext{1.1 2.7}{$d$}
  \picputtext{1.3 1.5}{$z$}
  \picputtext{6 5.4}{$y$}
  \picputtext{6 1.9}{$c$}
  \picputtext{8.7 3.9}{$a$}
  \picputtext{10.9 3.9}{$b$}
  \picputtext{3.6 5.7}{$L$}
} & 
\diag{8mm}{8}{6}{
  \picline{0.5 5.6}{3 5.6}
  \picline{0.5 0.4}{3 0.4}
  \picline{0 3.0}{3 3.0}
  \picline{5.5 5.5}{8 5.5}
  \picline{6 2.5}{8 2.5}
  \picputtext{0.3 5.6}{$x$}
  \picputtext{-0.2 3.0}{$d$}
  \picputtext{0.3 0.4}{$z$}
  \picputtext{7.8 2.3}{$c$}
  \picputtext{7.8 5.7}{$y$}
  \pictranslate{4 3}{
    \picfilledcircle{0 0}{3}{}
    {\piclinedash{0.05 0.05}{0.01}
      \picline{3 137 polar}{3.8 146 polar}
      \picline{3 175 polar}{3.8 174 polar}
      \picline{3 185 polar}{3.8 186 polar}
      \picline{3 223 polar}{3.8 214 polar}
      \picline{3 -1 polar}{3.8 1 polar}
      \picline{3 37 polar}{3.8 29 polar}
      \picline{3 156 polar}{0 0}
      \picline{3 204 polar}{0 0}
      \picline{3 18 polar}{0 0}
      \piccurve{3 144 polar}{-0.6 0.9}{0.8 0.9}{3 27 polar}
      \piccurve{3 192 polar}{-2.0 -0.6}{-2.0 0.6}{3 168 polar}
      \piccurve{3 212 polar}{-0.5 -1.3}{0.3 -1.2}{3 8 polar}
    }
    \pictranslate{15 3.8 x polar}{
      \picrotate{15}{
	\picfilledellipse{0 0}{0.3 1.1}{}
      }
    }
    \picputtext{15 4.4 x polar}{$G$}
    \pictranslate{160 3.8 x polar}{
      \picrotate{160}{
	\picfilledellipse{0 0}{0.3 1.1}{}
      }
    }
    \picputtext{160 4.4 x polar}{$E$}
    \pictranslate{200 3.8 x polar}{
      \picrotate{200}{
	\picfilledellipse{0 0}{0.3 1.1}{}
      }
    }
    \picputtext{200 4.4 x polar}{$F$}
    \picfilledcircle{0 0}{0.5}{$A$}
    \picfilledcircle{-0.2 1.1}{0.3}{}
    \picputtext{-0.2 1.7}{$B$}
    \picfilledcircle{0.1 -1.1}{0.3}{}
    \picputtext{0.1 -1.7}{$C$}
    \picfilledcircle{-2.2 -0.1}{0.3}{}
    \picputtext{-1.6 -0.1}{$D$}
    \picputtext{0 2.7}{$L$} 
  }
}  \\[27mm] \multicolumn{2}{c}{(a)} \\[5mm]
\diag{6mm}{15}{6}{
  {\piclinedash{0.05 0.05}{0.01}
    \picstroke{
      \opencurvepath{12 2}{13 1.9}{13 1}{7 1}{3 0.8}{1.5 1}{1.5 1.5}
        {2.5 1.7}{}
      \opencurvepath{12 5}{14.0 4}{14.5 0}{9 0}{2 0}{0 0}{0 4.7}
        {2.2 5.5}{}
      \picline{4 5}{7.6 5}
      \picline{4 2.2}{7.6 2.2}
      \picline{8 3.6}{9.6 3.6}
      \picline{10 3.6}{11.6 3.6}
    }
  }
  \pictranslate{3.6 3.6}{
    \picfilledcircle{0 0}{2.4}{}
    \picfilledellipse{3.8 0}{0.8 2.0}{}
    \picfilledellipse{6.2 0}{0.6 2.0}{}
    \picfilledellipse{8.4 0}{0.6 2.0}{}
  }
  \picputtext{1.1 4.7}{$x$}
  \picputtext{1.3 1.5}{$d$}
  \picputtext{6 5.4}{$c$}
  \picputtext{6 1.9}{$y$}
  \picputtext{8.7 3.9}{$a$}
  \picputtext{10.9 3.9}{$b$}
  \picputtext{3.6 5.7}{$L$}
} &
\diag{8mm}{9}{6}{
  \pictranslate{0.5 0}{
  \picline{-0.5 5.8}{3 5.8}
  \picline{-0.5 0.2}{3 0.2}
  \picline{5.5 0.5}{8 0.5}
  \picline{6 3.5}{8 3.5}
  \picputtext{-0.7 5.8}{$x$}
  \picputtext{-0.7 0.2}{$d$}
  \picputtext{7.8 3.8}{$c$}
  \picputtext{7.8 0.7}{$y$}
  }
  \pictranslate{4.5 3}{
    \picfilledcircle{0 0}{3}{}
    {\piclinedash{0.05 0.05}{0.01}
      % \picline{-3 0}{0 0}
      \picmultigraphics[S]{2}{1 -1}{
        % \piccurve{3 125 polar}{-2.2 2.5}{-3.5 2.5}{-5 1.9}
        \piccurve{3 131 polar}{-2.5 2.4}{-3.2 2.2}{-4.2 1.7}
      }
      \picline{3 180 polar}{4.1 180 polar}
      \picline{3 141 polar}{3.5 148 polar}
      \picline{3 174 polar}{3.5 172 polar}
      \picline{3 186 polar}{3.5 188 polar}
      \picline{3 219 polar}{3.5 212 polar}
      \picline{3 1 polar}{3.8 -1 polar}
      \picline{3 -37 polar}{3.8 -29 polar}
      \picline{3 156 polar}{0 0}
      \picline{3 204 polar}{0 0}
      \picline{3 -18 polar}{0 0}
      \piccurve{3 147 polar}{-0.6 1.2}{0.8 1.2}{3 -8 polar}
      \piccurve{3 192 polar}{-2.0 -0.6}{-2.0 0.6}{3 168 polar}
      \piccurve{3 212 polar}{-0.5 -0.9}{0.3 -0.9}{3 -27 polar}
    }
    \picfilledellipse{-4.2 0}{0.3 1.9}{}
    % \picfilledellipse{-5 0}{0.3 2.3}{}
    \pictranslate{-15 3.8 x polar}{
      \picrotate{-15}{
	\picfilledellipse{0 0}{0.3 1.1}{}
      }
    }
    \pictranslate{160 3.4 x polar}{
      \picrotate{160}{
	\picfilledellipse{0 0}{0.18 0.9}{}
      }
    }
    \pictranslate{200 3.4 x polar}{
      \picrotate{200}{
	\picfilledellipse{0 0}{0.18 0.9}{}
      }
    }
    \picfilledcircle{0 0}{0.5}{}
    \picfilledcircle{0.2 1.1}{0.3}{}
    \picfilledcircle{0.1 -1.1}{0.3}{}
    \picfilledcircle{-2.2 -0.1}{0.3}{}
  }
}  \\[25mm] \multicolumn{2}{c}{(b)} \\[5mm]
\diag{6mm}{15}{6}{
  {\piclinedash{0.05 0.05}{0.01}
    \picstroke{
      \opencurvepath{12 2}{13 1.9}{13 1.3}{7 1}{3 0.8}{1.5 1}{1.5 1.5}
        {2.5 1.7}{}
      \opencurvepath{12 3.5}{13.7 2.9}{13.7 0.5}{7 0.5}
        {3 0.5}{0.5 0.5}{0.5 3.0}{2 3.5}{}
      \opencurvepath{12 5}{14.0 4}{14.5 0}{9 0}{2 0}{0 0}{0 4.7}
        {2.2 5.5}{}
      \picline{4 5.2}{7.6 5.2}
      \picline{4 3.6}{7.6 3.6}
      \picline{4 2.0}{7.6 2.0}
      \picline{8 3.6}{9.6 3.6}
      \picline{10 3.6}{11.6 3.6}
    }
  }
  \pictranslate{3.6 3.6}{
    \picfilledcircle{0 0}{2.4}{}
    \picfilledellipse{3.8 0}{0.8 2.0}{}
    \picfilledellipse{6.2 0}{0.6 2.0}{}
    \picfilledellipse{8.4 0}{0.6 2.0}{}
  }
  \picputtext{1.1 4.7}{$g$}
  \picputtext{1.1 2.7}{$d$}
  \picputtext{1.3 1.5}{$e$}
  \picputtext{6 5.6}{$h$}
  \picputtext{6.3 3.8}{$c$}
  \picputtext{6 1.7}{$f$}
  \picputtext{8.7 3.9}{$a$}
  \picputtext{10.9 3.9}{$b$}
  \picputtext{3.6 5.7}{$L$}
} &
\diag{8mm}{8}{6}{
  \picputtext{0.3 5.6}{$g$}
  \picputtext{-0.2 3.0}{$d$}
  \picputtext{0.3 0.4}{$e$}
  \picputtext{7.8 2.8}{$c$}
  \picputtext{7.8 5.7}{$h$}
  \picputtext{7.8 0.7}{$f$}
  \picline{0.5 5.6}{3 5.6}
  \picline{0.5 0.4}{3 0.4}
  \picline{0 3.0}{3 3.0}
  \picline{5.5 5.5}{8 5.5}
  \picline{5.5 0.5}{8 0.5}
  \picline{6 3.0}{8 3.0}
  \pictranslate{4 3}{
    \picfilledcircle{0 0}{3}{}
    {\piclinedash{0.05 0.05}{0.01}
      \picline{3 137 polar}{3.8 146 polar}
      \picline{3 173 polar}{3.8 173 polar}
      \picline{3 187 polar}{3.8 187 polar}
      \picline{3 223 polar}{3.8 214 polar}
      \piccurve{3 156 polar}{-0.6 1.0}{0.8 1.0}{3 27 polar}
      \piccurve{3 204 polar}{-0.5 -1.1}{0.3 -1.1}{3 -27 polar}
    }
    \pictranslate{200 3.8 x polar}{
      \picrotate{200}{
	\picfilledellipse{0 0}{0.3 1.1}{}
      }
    }
    \pictranslate{160 3.8 x polar}{
      \picrotate{160}{
	\picfilledellipse{0 0}{0.3 1.1}{}
      }
    }
    \picfilledcircle{-0.2 1.1}{0.4}{}
    \picfilledcircle{0.1 -1.1}{0.4}{}
  }
}  \\[25mm]  \multicolumn{2}{c}{(c)} 
\end{array}
\]
\end{center}
\caption{\label{Fig19}Diagram patterns. The left pictures show the
outer cycle loops and the edges connecting them. The right parts show
the non-empty outer cycle loop $L$ and the ways other loops can be
attached to it, up to omissions and parallel loops.}
\end{figure}

Depending on the choice of regions
(a triple, a pair, or two pairs) identified by the loops attached 
inside the outer cycle loop $L$ we have three types of patterns we
refer to below as (a), (b) and (c), see figure \reference{Fig19}.
% (There are a few modifications of patterns (a) and (b), whose
% treatise is postponed to \S\reference{Sfinal}.)
Each edge in figure \ref{Fig19} stands for a group of parallel
traces. In the following we identify an edge with its multiplicity.
(We allow in general this multiplicity to be 0, i.e. to delete
the edge. Then we delete also loops becoming isolated, i.e. when
all their legs are 0.) We assign to an edge a letter as indicated
in figure \reference{Fig19}. 

Now, by lemma \reference{lempar},
a pair of parallel loops is excluded in all patterns. % (b)~- and
% only there, because by lemma \reference{lempar}, loops attached inside
% the outer cycle loop identify only two regions.
Also in patterns (c)
the only choice of two pairs of regions so that $IG(A)$ contains a
cycle is as indicated. We allow in patterns edges that connect $L$ to a
neighboring loop in the outer cycle to be of multiplicity 0 (i.e. to
be deleted) when generating diagrams $D'$. Still many of these cases can
be discarded because of repetitions, or because $D'$ (and so $D$)
becomes a composite diagram. In particular in figure \reference{Fig19}
we chose the naming of edges in all 3 cases so that initial alphabet
letters, $a$ to $h$, stand for multiplicities that can be assumed to
be at least $1$ (by symmetry, because $D$ is prime, or because, in case
$(c)$, a cycle must occur in $IG(A)$). The names $a$ and $b$ always
stand for edges that connect two empty outer cycle loops. The
omitable (i.e. possibly of multiplicity 0) edges are named by
end-alphabet letters $x,y,z$. In case (c), if $c$ or $d$ is
missing, we have a special case of (b).

It is possible that all regions identified by loops attached
inside $L$ are on the same side, that is, between traces/edges that
connect $L$ to the same neighboring loop $L'$ in the outer cycle.
One easily sees that in the pattern (c) both loops of depth 1
(which identify two disjoint pairs of regions) must intertwine two
connections of the outer cycle. So we need to look at patterns
(a) and (b), and we call the modified patterns of type (a$'$) and
(b$'$). Here are examples for either case (keep in mind corollary
\reference{cr3l}):
\begin{eqn}\label{exss}
\begin{array}{c@{\kern2cm}c}
\diag{7mm}{6.5}{6}{
    {\piclinedash{0.05 0.05}{0.01}
      \picline{0 0.9}{3 0.9}
      \picline{0 2.3}{3 2.3}
      \picline{0 3.7}{3 3.7}
      \picline{0 5.1}{3 5.1}
      \picline{3 3}{6.5 3}
      \picline{0.5 4.0}{3.5 4.0}
      \picline{0.5 4.8}{3.5 4.8}
    }
  \pictranslate{3.4 3}{
    \picfilledcircle{0 0}{2.5}{}
    \picfilledellipse{154 3.2 x polar}{0.25 0.55}{}
    {\piclinedash{0.05 0.05}{0.01}
      \piccurve{220 2.5 x polar}{-0.6 -0.1}{-0.5 .6}{139 2.5 x polar}
      \piccurve{210 2.5 x polar}{-1.4 -0.1}{-1.3 .4}{149 2.5 x polar}
      \picline{185 1.6 x polar}{180 2.5 x polar}
    }
    \picfilledcircle{170 0.9 x polar}{0.27}{}
    \picfilledcircle{185 1.6 x polar}{0.27}{}
  }
  \picputtext{3.4 5.2}{$L$}
  \picputtext{0.4 5.4}{$e$}
  \picputtext{0.4 1.2}{$e'$}
}
& 
\diag{7mm}{6.5}{6}{
    {\piclinedash{0.05 0.05}{0.01}
      \picline{0 0.9}{3 0.9}
      \picline{0 3}{3 3}
      \picline{0 5.1}{3 5.1}
      \picline{3 3}{6.5 3}
      \picline{0.5 4.0}{3.5 4.0}
      \picline{0.5 4.8}{3.5 4.8}
    }
  \pictranslate{3.4 3}{
    \picfilledcircle{0 0}{2.5}{}
    \picfilledellipse{154 3.2 x polar}{0.25 0.55}{}
    {\piclinedash{0.05 0.05}{0.01}
      \piccurve{220 2.5 x polar}{-0.6 -0.1}{-0.5 .6}{139 2.5 x polar}
      \piccurve{210 2.5 x polar}{-1.4 -0.1}{-1.3 .4}{149 2.5 x polar}
      \picline{185 1.6 x polar}{170 2.5 x polar}
    }
    \picfilledcircle{170 0.9 x polar}{0.27}{}
    \picfilledcircle{185 1.6 x polar}{0.27}{}
  }
  \picputtext{3.4 5.2}{$L$}
  \picputtext{0.4 5.4}{$e$}
  \picputtext{0.4 1.2}{$e'$}
} \\
\ry{1.1em}(a') & (b')
\end{array}
\end{eqn}
% 
% In patterns (a) and (b) the choices of
% regions could be all of the same side (that is, between edges that
% connect $L$ to the same neighboring loop in the outer cycle).
Now, by a mutation, we can move down the uppermost (possibly
multiple) edge $e$ between $L$ and $L'$ (the loop left of $L$,
which is not drawn), making $e$ parallel to $e'$. This way we
convert patterns as in \eqref{exss} into such included in figure
\ref{Fig19} (with some outer cycle connection of multiplicity 0).

There is one more peculiar instance of patterns of types (b) and
(b$'$), in which one may attach inside $L$ two loops of valence three,
as shown in \eqref{pecu}.
(This diagram is for patterns (b); in the analogous picture for
patterns (b$'$), one pair of edges $(x,x')$ or $(y,y')$ may have the
same two loops connected, or even unify to the same edge.) However,
such patterns were ruled out by lemma \reference{lmZ} (because the
loop $Z$ is inadequate).

Next observe that patterns (c) are excluded by lemma \ref{lmio}. We
need to test only types (a) and (b). While one can rule out by hand
many of these cases, still too many different possibilities remain
to make a ``manual'' continuation worthwhile. It appears more helpful
to check these patters with a computer. We found no (technically)
easy way to calculate the $A$- and, particularly, the $B$-state
invariants from the patters directly, so we relied on a calculation
of their Jones polynomials. To perform this calculation, a few more
preparations are necessary. The first is a technical help. 

\begin{lemma}\label{lmoc}
There exists at most one multiple edge between empty loops of the
outer cycle. Consecutively, there exist at least two simple edges.
\end{lemma}

\proof That the second claim follows from the first claim is
a consequence of the fact that the outer cycle has length at least
$5$, and only one non-empty loop. To see the first claim, assume
there were at least two multiple edges. They form then isolated
vertices in $IG(A)$. However, by three attachments it is impossible
to create three (homologically independent) cycles in $IG(A(D))$,
and so $\chi(IG(A(D)))$ would remain positive. \qed
% By connectivity reasons
% there must be at least two of multiplicity at least 3. However,
% since the loops connected are empty, such edges are parallel.
% We argued, though, that we have at most one parallel equivalence
% class of more than $2$ edges. \qed

We can now \em{desplice} a single edge between two empty
outer cycle loops and join them.
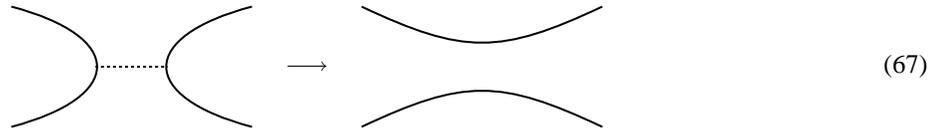
\begin{eqn}\label{desp}
\diag{8mm}{4}{2}{
   \pictranslate{2 1}{
       \picmultigraphics[S]{2}{-1 1}{
           \piccurve{2 -1}{0.1 -0.5}{0.1 0.5}{2 1}
       }
       {\piclinedash{0.05 0.05}{0.01}
        \picline{0 -0.6 x}{0 0.6 x}
       }
   }
}
\quad\lra\quad
\diag{8mm}{4}{2}{
   \pictranslate{2 1}{
       \picmultigraphics[S]{2}{1 -1}{
           \piccurve{-2 1}{-0.3 0.2}{0.3 0.2}{2 1}
       }
       {\piclinedash{0.05 0.05}{0.01}
        %\picline{0 -0.4}{0 0.4}
       }
   }
}
\end{eqn}
Then we continue until we have
a diagram $D'$ with an outer cycle of length $4$. Thus (and this
is to be kept in mind) the number of desplicings that turn $D$ into
$D'$ is odd. We can recover $D$ by choosing an outer cycle loop and
\em{splicing} it again (the reverse operation to \eqref{desp}).

The reason for the change to $D'$ is computational. For the
calculation of the Jones polynomials, we tried to apply the program
of Millett-Ewing in \cite{KnotScape} that uses a skein algorithm.
It requires a choice of orientation on all components, but this
choice is not natural for an odd length outer cycle. (The way
components are connected varies drastically for the different 
combinations of edge multiplicities.) In opposition, for even
length we have a canonical positive orientation on $D'$.

In our situation the (de)splicings between $D'$ and $D$ do not change
adequacy and the invariants of $A(D)$. Let $a,b$ be the (multiplicities
of) the outer loop edges that do not connect the loop with non-empty
interior. By mutations we can assume that $a=1$.

The effect of \eqref{desp} on $B(D)$ is that (up to mutations that
preserve all $B$-state invariants) an edge is reduced to multiplicity
$1$ (if $b>1$) or $2$ (if $b=1$). This does not change $\chi(G(B(D)))$
and $\tg(B(D))$, and the only effect it can have on $IG(B(D))$ is that
an isolated vertex is removed if $b>1$. So we can conclude:
\begin{itemize}
\item The $A$-state invariants (in particular also
  both $V_{1,2}$) coincide on $D,D'$.
\item We have $\chi(B(D'))=\chi(B(D))=0$, and $\tg(B(D'))=\tg(B(D))=0$.
%  and $\chi(IG(B(D')))$ is $0$ or $-1$.
% \item In the case $\chi(B(D'))=\chi(B(D))+1=1$, the despliced edge is
%   the only single edge in the outer loop, so in $D'$ we have $a,b>1$.
\item $\chi(IG(B(D')))$ is $-1$ or $0$ (so $\bV_0\bV_2=0,1$), and if
  $b=1$, then $\chi(IG(B(D')))=0$ (and $\bV_0\bV_2=1$).
\end{itemize}

% IN CALC WHEN IN DGS OF page 10 
% if some of a,b>1 check if \bar V_1=-1
% if both a,b>1 check if \bar V_1=-1,0
% if a=b=1 check if \bar V_1=-1, \bar V_2=-1 as usual

The generation process of diagrams in the above patterns
follows the general rules:
\begin{itemize}
\item a loop except the outer cycle loops may or may not be there,
\item exactly 3 (non-outer cycle) 
  loops are attached in total (corollary \ref{cr3l}),
\item no pair of parallel loops occurs (lemma \reference{lempar}),
\item there are at least two depth-one loops (lemma \reference{l9,15}),
\item leg multiplicities for valence two loops
  are given in lemma \reference{lemmaW}, and for
% \item leg multiplicities
  valence three loops are not all even (for connectivity reasons),
% \item in case (b) exactly one valence 2-loop can be a parallel
  % pair, with leg multiplicities given in lemma \reference{lemmaX},
\item leg multiplicities follow \eqref{lmt},
\item every depth-0 loop which is not an outer cycle loop
  is intertwined with some depth-one loop (else the diagram is
  not prime), and
\item $c(D')$ is odd (because $c(D)$ is even and the number of
  desplicings is odd).
\end{itemize}

There is a large but finite number of cases for $D'$. We wrote
a computer program to generate them. In case (b) we simultaneously
calculated the invariant $\chi(IG(A))$ from the data (since now this
is possible directly without calculation of the whole polynomial!),
to minimize the number of cases in advance. (The values of $\gm$,
$\chi(G(A))$ and $\tg(A)$ are correct by construction.) This way
the list reduces to about 120,000 diagrams $D'$, the most
complicated ones having 27 crossings.

We calculate the Jones polynomial of these diagrams (it takes a few
minutes) and check that we have a diagram $D'$ with 
\begin{enumerate}
\item $\spn V=c(D')-2$ (still useful to test, since despite
  $\gm=4$, we do not know if $D'$ is $B$-adequate), 
\item $V_0V_1=\bV_0\bV_1=-1$, and $V_0V_2=1$.
% \item If in the outer loop of $D'$ one of $a$ or $b$ is $1$,
%   then we can test $\chi(G(B))=0$.
\item If in the outer loop of $D'$ we have $(a=)b=1$, then we can test
  also $\bV_0\bV_2=1$, otherwise this product is $0$ or $1$.
\item $D'$ has one or two components. (This was for technical
  reasons easier to check using the value $V(1)$ rather than at
  the time of generating the diagrams.)
\end{enumerate}

Thus it is legitimate to \em{discard} during the generation process
outer
loop edge multiplicities that satisfy some of the following properties:
\begin{enumerate}
\item Type (a): $a>1$, which can be avoided by mutations, or $x=y=0$,
  because then regions connected by depth-1 loops identify, and
  we have a special case of type (b).
  % $0=d\%2+x\%2+z\%2+b\%2$, $0=y\%2+c\%2+a\%2$,
  % $0=d\%2+x\%2+z\%2+y\%2+a\%2$, or $0=a\%2+b\%2$ (then disconnected)
\item Type (b): 
  % \begin{itemize}
    % \item
      $a>1$, or $(a=)b=1$ and $c+y>d+x$ (then can flip by 180 degree
      or do mutations)
    % \item no pair of parallel edges occurs, %then discard diagram 
      % \em{unless} $a=b=c=1$, $y=0$, $d=2$, $x=1$ (because of lemma
      % \reference{lemmaX})
    % \item $0=a\%2+y\%2+c\%2$, $0=x\%2+d\%2+b\%2$,
      % $0=x\%2+d\%2+c\%2+y\%2$, or $0=a\%2+b\%2$ (then disconnected)
  % \end{itemize}
\item Types (a) and (b): for some outer cycle loop all edges
  connecting it to neighboring outer cycle loops are of even
  multiplicity (then $D'$ and $D$ are disconnected)
\end{enumerate}

% IN CALC WHEN IN DGS OF page 10 
% if some of a,b>1 check if \bar V_1=-1
% if both a,b>1 check if \bar V_1=-1,0
% if a=b=1 check if \bar V_1=-1, \bar V_2=-1 as usual

% !!!!!!!!!!!!!!!!!!!!!!!!!
% this is a mistake of v_ktest3a.C that is being corrected
% in v_ktest3c.C

% Type (c) patterns are the simplest, in that the fewest combinations
% are possible.
Most combinations are ruled out by the Jones polynomial conditions:
% For patterns of type (a) and (b), we have more work:
only 40 patterns of diagrams of type (a) and 2 patterns of diagrams of
type (b) remain. %and their discussion becomes too tedious.
The next part of the proof handles these cases.

\subsection{Concluding part: Generic and sporadic cases\label{S9.6}}

Let us call in mind how one recovers all possible diagrams $D$ from
$D'$. First one can, reversing \eqref{desp}, by (iterated) splicing
create single edges in the
outer cycle of $A(D)$. Then remember our convention from the
beginning of the proof that we worked with parallel equivalence
classes of $4$ or $5$ edges in $A(D)$. This was done so as to have only
a finite number of diagrams to test with the Jones polynomial,
and because the test is equivalent for all the others. Now we must
reconvert these multiplicities to $2n$ or $2n-1$ (for $n\ge 2$).
So if $D'$ has an edge of $4$ or $5$ traces,
we have two sites to twist at to generate diagrams $D$, while
otherwise we have only one (the outer cycle).

% Here comes the situation where extra cases occur.

Let us call a diagram $D'$ \em{sporadic} if among the
diagrams $D$ that result from $D'$ there is only a finite
number with $w(D)=0$ and $c(D)-2g(D)\le n$ for any fixed
number $n$. (Here $g(D)$ refers to the canonical genus
of the diagram.) Otherwise call $D'$ \em{generic}.

The (hypothetic) diagrams $D=D_k$ of $K_k$ satisfy these
conditions for $n=6$ (by direct calculation of the
Alexander polynomial). So generic diagrams $D'$ give infinitely
many candidates for diagrams $D$, while sporadic $D'$ give
only finitely many. %allowing for a finite number
% of exceptional diagrams $D$, it is enough to deal with generic
% diagrams $D'$ only.
The assumption $D'$ to be generic leads to
further significant computational simplifications.

\begin{lemma}
If $D'$ is generic then it has one (and then only one)
edge of $4$ or $5$ traces.
\end{lemma}

\proof Twisting at only one site
(the outer loop in $A(D)$ when recovering $D$ from $D'$)
can give at most one diagram $D$ of writhe 0. \qed

\begin{lemma}\label{l9.18}
If $D'$ is generic, we can assume that $D$ (or $D'$)
has no multiple edges between empty outer cycle loops.
\end{lemma}

\proof We can assume henceforth, because $w(D)=0$, that $D'$
has an edge of $4$ or $5$ traces, and the
recovery of $D$ requires twisting at two different
sites, with the number of twists having a bounded difference.
By lemma \reference{lmoc} we know already that
all but finitely many edges between empty outer cycle loops
in $D$ are simple. Since the twists added by these edges
must be parallel (else $g(D)$ will grow at most like $c(D)/4$),
all multiplicities are odd. Since at most one edge
has multiplicity $\ge 3$, at most one
edge is multiple. But if it is multiple, then it is
the only edge in $A(D)$ of multiplicity $\ge 3$.
So in recreating $D$ we would have to twist at that
site. Now one easily observes that if the twists 
that correspond to edges in the outer cycle are
parallel, those that correspond to the multiple edge
are reverse. So again $g(D)$ will grow at most like $c(D)/4$,
and not $c(D)/2$. \qed

Then one easily observes:

\begin{corr}
If $D'$ is generic, the change \eqref{desp} from $D$ to $D'$ 
(or vice versa) preserves all $A$- and $B$-state invariants.
That is, it is legitimate to test the conditions \eqref{tstcnd}
on $D'$. \qed
\end{corr}

With the same argument as in lemma \ref{l9.18} we have

\begin{lemma}\label{lmod}
If $D'$ is generic, then for each of both neighbors $L'$ of the
non-empty loop $L$ in the outer cycle, there is an edge of
odd multiplicity between $L$ and $L'$. \qed
\end{lemma}

Recall also that $D'$ must have at most two components.

The conditions on generic $D'$ we just put together rule out the
2 remaining possibilities of patters (b), and leave (from the
40) only 24 diagrams $D'$ of patters (a). All these diagrams
have two components. The vectors describing them are below:

\begin{figure}
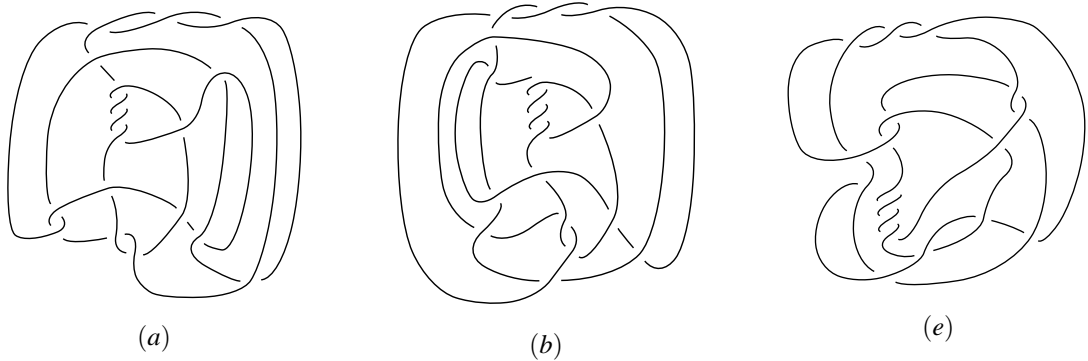
%[htb]
\begin{center}
\begin{eqnarray*}%{*4c}
&\rx{-7mm} \vis{(a)}{t1-oddach_case_a} \qquad \vis{(b)}{t1-oddach_case_b} \qquad
  \vis{(e)}{t1-oddach_case_e} & \\ [0mm]
% &\rx{-7mm} \vis{(d)}{t1-oddach_case_d} \qquad \vis{(e)}{t1-oddach_case_e} \qquad
% \vis{(f)}{t1-oddach_case_f} & \\ [3mm]
% &\rx{-7mm} \vis{(g)}{t1-oddach_case_g} \qquad \vis{(h)}{t1-oddach_new1} \qquad
% \vis{(k)}{t1-oddach_new2} & \\ [3mm]
% &\rx{-7mm} \vis{(l)}{t1-oddach_new3} \qquad \vis{(m)}{t1-oddach_new4}
% \qquad  \vis{(n)}{t1-oddach_case_h}  & \\ [3mm]
\end{eqnarray*}
\end{center}
\caption{\label{Fig20}The diagrams $\hD_2$ in the remaining
3 flype-inequivalent generic sequences.}
\end{figure}

{\small
\begin{eqn}\label{css}
\begin{array}{*8{c|}c}
  A   & B   & C   & D   & E   & F   & G   & \mbox{outer cycle} \\[1mm]
  \hline
\ry{1.4em}%
1\,1\,2 & 1\,1 & 0\,0 & 0\,0 & 0\,0 & 0\,0 & 1\,1 & 1\,1\,5\,2\,0\,1\,1 & * \\
1\,1\,2 & 1\,1 & 0\,0 & 0\,0 & 0\,0 & 0\,0 & 1\,4 & 1\,1\,1\,2\,1\,0\,0 & * \\
1\,1\,2 & 1\,1 & 0\,0 & 0\,0 & 1\,1 & 0\,0 & 0\,0 & 1\,1\,1\,1\,4\,0\,1 & * \\
1\,1\,2 & 1\,1 & 0\,0 & 0\,0 & 1\,4 & 0\,0 & 0\,0 & 1\,1\,1\,1\,1\,0\,1 & * \\
1\,1\,2 & 1\,1 & 0\,0 & 0\,0 & 1\,4 & 0\,0 & 0\,0 & 1\,1\,2\,1\,0\,1\,0 & * \\
1\,1\,2 & 1\,4 & 0\,0 & 0\,0 & 0\,0 & 0\,0 & 1\,1 & 1\,1\,1\,2\,1\,0\,0 & (a) \\
1\,1\,2 & 1\,4 & 0\,0 & 0\,0 & 1\,1 & 0\,0 & 0\,0 & 1\,1\,1\,1\,1\,0\,1 & (b) \\
1\,1\,2 & 1\,4 & 0\,0 & 0\,0 & 1\,1 & 0\,0 & 0\,0 & 1\,1\,2\,1\,0\,1\,0 & (c) \\
1\,2\,1 & 0\,0 & 1\,1 & 0\,0 & 0\,0 & 0\,0 & 1\,1 & 1\,1\,1\,2\,1\,5\,0 & \# \\
1\,2\,1 & 0\,0 & 1\,1 & 0\,0 & 0\,0 & 1\,4 & 0\,0 & 1\,1\,1\,1\,0\,2\,0 & \# \\
1\,2\,1 & 0\,0 & 1\,4 & 0\,0 & 0\,0 & 1\,1 & 0\,0 & 1\,1\,1\,1\,0\,2\,0 & (d) \\
1\,2\,4 & 0\,0 & 1\,1 & 0\,0 & 0\,0 & 1\,1 & 0\,0 & 1\,1\,1\,1\,0\,2\,0 & (e) \\
\end{array}
\quad
\begin{array}{*8{c|}c}
  A   & B   & C   & D   & E   & F   & G   & \mbox{outer cycle} \\[1mm]
  \hline
\ry{1.4em}%
1\,4\,2 & 1\,1 & 0\,0 & 0\,0 & 1\,1 & 0\,0 & 0\,0 & 1\,1\,1\,1\,1\,0\,1 & (f) \\
1\,4\,2 & 1\,1 & 0\,0 & 0\,0 & 1\,1 & 0\,0 & 0\,0 & 1\,1\,2\,1\,0\,1\,0 & (g) \\
2\,1\,1 & 0\,0 & 0\,0 & 1\,1 & 0\,0 & 1\,1 & 0\,0 & 1\,1\,1\,1\,0\,2\,5 & \$ \\
2\,1\,1 & 0\,0 & 0\,0 & 1\,1 & 0\,0 & 1\,4 & 0\,0 & 1\,1\,1\,1\,2\,0\,0 & \$ \\
2\,1\,1 & 0\,0 & 0\,0 & 1\,1 & 1\,1 & 0\,0 & 0\,0 & 1\,1\,1\,4\,1\,0\,1 & \$ \\
2\,1\,1 & 0\,0 & 0\,0 & 1\,1 & 1\,1 & 0\,0 & 0\,0 & 1\,1\,2\,1\,5\,1\,0 & \$ \\
2\,1\,1 & 0\,0 & 0\,0 & 1\,1 & 1\,4 & 0\,0 & 0\,0 & 1\,1\,1\,1\,1\,0\,1 & \$ \\
2\,1\,1 & 0\,0 & 0\,0 & 1\,4 & 0\,0 & 1\,1 & 0\,0 & 1\,1\,1\,1\,2\,0\,0 & (h) \\
2\,1\,1 & 0\,0 & 0\,0 & 1\,4 & 1\,1 & 0\,0 & 0\,0 & 1\,1\,1\,1\,1\,0\,1 & (k) \\
2\,1\,4 & 0\,0 & 0\,0 & 1\,1 & 0\,0 & 1\,1 & 0\,0 & 1\,1\,1\,1\,2\,0\,0 & (l) \\
2\,4\,1 & 0\,0 & 0\,0 & 1\,1 & 1\,1 & 0\,0 & 0\,0 & 1\,1\,1\,1\,1\,0\,1 & (m) \\
4\,1\,2 & 1\,1 & 0\,0 & 0\,0 & 0\,0 & 0\,0 & 1\,1 & 1\,1\,1\,2\,1\,0\,0 & (n) \\
\end{array}
\end{eqn}
}
These vectors give the multiplicities of edges in $D'$.
The legs of the loop A are listed in counterclockwise order
with the first leg connecting to the right side in figure \ref{Fig19}
(between the legs of $G$). The edges of the outer cycle are listed
in the order $a,b,c,d,x,y,z$. Up to mutations (actually flypes)
this information specifies the diagrams. Twelve of the 24 diagrams $D'$
do not give a knot diagram $D$ by undoing the desplicings \eqref{desp}
(i.e. crossings that connect the empty outer cycle loops of $A(D')$
do not involve both components). These patterns have leg multiplicities
$A=(1,1,2)$ and $B=(1,1)$ (the vectors of these cases are marked with
an asterisk), $A=(1,2,1)$ and $C=(1,1)$ (marked with a double
cross), or $A=(2,1,1)$ and $D=(1,1)$ (marked with `\$'). Consider the
other twelve diagrams, to which we assign labels (a) to (n) as above.

For all twelve $D'$ we obtain (indeed) diagrams $D$ of achiral
adequate knots that satisfy all criteria we could impose so far
coming from the Jones polynomial (and as well the Kauffman polynomial
in corollary \ref{cr3l}). Let  us call these diagrams $\hD_k$,
reintroducing the twist parameter $k$ from the beginning of the
proof.

The 12 initial diagrams $\hD_2$ for the series (a) to (n) have 18
crossings, and they have exactly two (parallel) twists of 3 or 4
crossings. Each next diagram $\hD_{k+1}$ arises from the previous
one $\hD_k$ by a adding two crossings each at its two twists. Now
it turns out that up to mirroring and flypes from these 12 series
only 3 are inequivalent. Since flypes persist under twists, it is
enough to deal with these 3. They are shown in figure \ref{Fig20}.
In fact, we have only two sequences of knots, since the series
(b) and (e) are not flype equivalent but still represent the same
knot. (The division of the 12 series between the two knots is
evident from \eqref{65}.)

% The 12 initial 18 crossing diagrams $\hD_2$ of each series
% are shown in figure \ref{Fig20}. (They are obtained from $D'$ by
% one splicing reverse to \eqref{desp}.) Each next diagram $\hD_{k+1}$
% arises from the previous one $\hD_k$ by a pair of (parallel) twists,
% one each at the two (obvious) sites where a twist of at least 3
% crossings occurs.

Now we must use some further information to distinguish
these series of $\hD_k$ from the hypothetic diagrams $D_k$
of our odd crossing number knot candidates. We use the
degree-2 Vassiliev invariant $v_2$ (from \eqref{v23}).

It is well-known (see \cite{bseq}) that in a sequence of diagrams
$\hD_k$, where $\hD_{k+1}$ is obtained from $\hD_{k}$ by a
pair of parallel twists at the same sites, $v_2(\hD_{k})$
is a(n at most) quadratic polynomial in $k$. So to determine
$v_2(\hD_{k})$ for all $k\ge 0$, it is enough to evaluate
$v_2$ on the first three diagrams. Actually one could also
switch crossings in the twists, as we did, to use simpler
diagrams instead. (This gives a meaning to $\hD_k$ and $D_k$
also when $k<1$.) The same behavior is exhibited by $v_2$
on our candidate knots. The calculation shows 
\begin{eqn}\label{65}
\begin{array}{l*4{|c}}
\multicolumn{1}{c|}{k} & -1 & 0 & 1 & 2 \\[1mm]
\hline
\ry{1.4em}
v_2(\hD_k) \mbox{ for series (a),(c),(d),(f),(k) and (m)}
  & -5 & -1 & 5 & 13 \\
v_2(\hD_k) \mbox{ for series (b),(e),(g),(h),(l) and (n)}
  & -4 &  0 & 6 & 14 \\
v_2(D_k)  & -3 &  1 & 7 & 15 \\
\end{array}
\end{eqn}
So the polynomials for $v_2$ for all eight sequences of
$\hD_{k}$ do not coincide with those for our hypothetic
diagrams $D_k$ (of the same crossing number). As explained,
the coincidences of $v_2$ on the two groups of 6 series are
because there are  only two series of knots.
% not accidental: these are in fact only two
% series of knots, and also some of the
% diagrams within the groups of 6 are equivalent up to flypes.

Now turn to the sporadic cases. We obtain by calculation 16 diagrams
$D'$ of pattern (a) in figure \reference{Fig19}
and 2 of pattern (b). (This is the rest of the
42 cases from the end of \S\reference{SPAT} after removing the
24 generic ones we just treated.) 

There are 4 sporadic diagrams $D'$, where $L$ is connected to an
outer cycle loop only by even valence edges (i.e. the conclusion
of lemma \reference{lmod} is violated). Since all these diagrams
have two components, and $a=b=1$, it is easy to see that 
(even an odd number of) splicings reverse to \eqref{desp} 
yield no knot diagrams.

Consider the other 14 sporadic diagrams $D'$. We can handle them
similarly to the described procedure for the generic ones, but
there is a (technically) more convenient way. As it turns out,
we have that for all patterns $a=b=1$, and $c(D')\in \{11,15\}$.
This means that one desplicing should give an adequate knot diagram
$D$, with $c(D)\in \{12,16\}$, $\gm(D)=4$ and all $A$- and $B$-state
invariants as for $K_k$. We can thus check in the tables of
\cite{KnotScape} for adequate knots $K$ with $\spn V(K)=c(K)-2$,
$V_0V_1=\bV_0\bV_1=-1$ and $V_0V_2=\bV_0\bV_2=1$. There are five 14
crossing knots, which are not relevant, and three 16 crossing knots,
%($14_{31449}$, $14_{33181}$, $14_{41330}$, $14_{41609}$ $14_{45653}$),
$16_{1166743}$, $16_{1223228}$, $16_{1303786}$. However, beside
the first and last three coefficients, we have one more piece of
information (which we ignored so far for the sake of testing
$V$ directly), say, the triangle number. By simply drawing a 16
crossing diagram of each of these three knots, one observes that
one of the $A$ or $B$ state has a triangle. 

This argument shows that
the diagrams $D$ coming from our 18 sporadic diagrams $D'$
either have multiple components, or $B(D)$ has a triangle.
We can thus rule out all sporadic diagrams $D'$ too.
The exclusion of these final remaining cases completes the
proof of theorem \reference{th15}. %\qed

\begin{rem}\label{r9,1}
The reason for the insufficiency of this proof in crossing numbers
$17+4k$ is that I found so far no generalization of Thistlethwaite's
example with $\chi(G(A))=0$, while for $\chi<0$ the combinatorics of
the argument seems to become considerably more complicated. The last
section completes the work by dealing with this difficulty using
all tools we have set up so far.
\end{rem}

% The reason for this 
% a cycle of length $4$.
% {}From now one we distinguish two main cases.
% 
% CASE I. We have parallel loops with the condition on $A(D)$ as
% in the lemma.
% 
% CASE II. We have no parallel loops.
% 
% We deal first with case II. Note that addition of a parallel loop
% does not change any of the invariants of $A(D)$ and it may only
% increase $\gm$. So we need to add a parallel loop only in the
% situation in CASE II which we cannot reult out by looking at $A(D)$.
% 
% CASE II. 

\section{Proof of main result\label{S10}}

\subsection{Examples}

To complete the proof of theorem \ref{thm}, we need to settle
crossing numbers $17+4k$. The examples we find are as follows.
Let $\dl_i$ be the tangle (with $n$ strings) obtained by 
(ignoring orientations and) replacing the crossing $\Pos{1.7em}{
\piclinewidth{50}}$ in the braid generator tangle $\sg_i\in B_n$
by 
\vcbox{
  \rbox{
    \parbox{1.3em}{
       \diag{1em}{1.3}{2}{
          \piccirclearc{1.3 1}{1}{90 180}
          \piccirclearc{0 1}{1}{-90 0}
          \picmulticirclearc{-9 1 -1 0}{1.3 1}{1}{180 -90}
          \picmulticirclearc{-9 1 -1 0}{0 1}{1}{0 90}
       }
    }
  }
}, and write $\bar\dl_i$ for the mirror image of $\dl_i$,
containing a
\vcbox{
  \rbox{
    \parbox{1.3em}{
       \diag{1em}{1.3}{2}{
          \piccirclearc{1.3 1}{1}{180 -90}
          \piccirclearc{0 1}{1}{0 90}
          \picmulticirclearc{-9 1 -1 0}{1.3 1}{1}{90 180}
          \picmulticirclearc{-9 1 -1 0}{0 1}{1}{-90 0}
       }
    }
  }
}\,. (This is no longer an inverse in a group theoretic sense.)

\begin{theo}\label{th10.1}
The knots $K_k$, given by the (braid-like, but 
now unoriented) closure of the tangles $T_k=
-12^23^{2k+2}4-32-1(-2)^{2k+2}3^{-2}4\bar\dl_2\dl_3$
are prime, amphicheiral, and have crossing number $17+4k$.
\end{theo}

Below on the left is shown the diagram of the theorem for $k=2$:
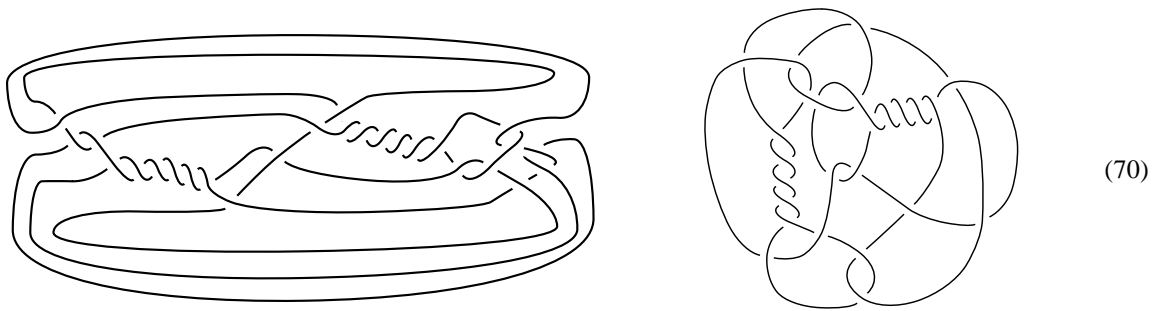
\begin{eqn}\label{Ti}
\kern5mm
\vcbox{
\hbox to 0.5\textwidth{\hss
\rbox{
  \epsfsv{6cm}{t1-17+4k_ach_tang}
}
\hss}
}
\kern10mm
\raisebox{0.15cm}{
  \epsfsv{4.2cm}{t1-17+4k_ach_tang_25}
}
\kern5mm
\end{eqn}
Amphicheirality of $K_k$ is again evident, and it is not too
hard to see that the $18+4k$ crossing diagrams as closure of the
specified tangle reduce by one crossing, as in the right diagram
of \eqref{Ti}. The minimality of
this diagram is again the object of difficulty.

The two replacements of braid crossings by the $\dl$ and $\bar\dl$
preserve the connectivity, but alter some geometric properties.
For example, the knots from the previous section all have the
same braid index (5), in opposition to our $K_k$, which in turn
have the same genus. The genus (equal to 6) can be inferred from
\cite{achpol}, using that the knots are semi-homogeneous in the
sense introduced there, as their diagrams from the left of \eqref{Ti}
are Murasugi sums of (special) alternating diagrams.

The $17+4k$ crossing diagrams of $K_k$ (of writhe $\pm 1$) are 
again semiadequate. So again a $16+4k$ crossing diagram $D_k$, 
which is to be ruled out, must be adequate. The Jones invariants
(for both $A$- and $B$-state of $D_k$) become as follows:
\begin{eqn}\label{71} 
\chi(G)=-1\,,\quad\chi(IG)=-1\,,\quad\tg=1\,,
\end{eqn} 
except again for $k=0$, where $\tg=2$. These changes to the previous 
case will lead to almost entirely different proof, even though its 
spirit is similar. First note that the triangle number becomes 
important, in that it controls the length of one of the two cycles 
in the state graph $G$. In contrast, the other cycle 
must be of length $\ge 4$ (for $k>0$). There is no parity condition
on this length anymore, since the triangle spoils the bipartacy of 
$G$. We call again the longer cycle the \em{outer cycle}, while the
shorter will still be the \em{triangle}. 

Furthermore, again we have 
\begin{eqn}\label{gm4}
\gm\,:=\,c(D_k)+2-s(A(D_k))-s(B(D_k))\,=\,4\,.
\end{eqn}
Let us point out also the following, which was not (needed to be)
used before. Since we must have writhe $w(D_k)=0$, and the Jones
polynomial is reciprocal, the $A$-and $B$-state of $D_k$ must 
have equal number of loops, and from \eqref{gm4} we have
\begin{eqn}\label{ABl}
s(A(D_k))=s(B(D_k))=c(D_k)/2-1=7+2k\,.
\end{eqn}
This property is easy to test from the $A$-state directly,
and we use it in both our computations (see \S\reference{Scmp})
and arguments (see lemma \reference{UV}). Contrarily, the
(vanishing of) the writhe can be tested only from the knot diagram.
We thus postpone its application to the very end of the proof
in \S\reference{SSI}, where all semiadequacy invariant tests
are exhausted, and we recover the knot diagram.

Even if just a triangle, the second cycle in $G(A)$ nonetheless 
complicates things drastically. The length of the proof results
mainly from the extraordinary care that is needed not to leave
out some cases how the attachments are performed (similar to
\S\reference{mpe}, but now with the new option \eqref{tatt}, as
we discuss below). Once the cases are identified, the invariants we
have exercise a control strong enough to easily rule out most of them;
only a few must be entrusted to the computer.

To restrict the attachments, we will rely now more heavily
on the atom number invariant of theorem \reference{ozp}.
By verification, we see that corollary \ref{cr3l} holds 
further for our present family. Another circumstance that comes 
extremely to our help is that $IG$ has now also two cycles.

The atom number again
easily implies that $K_k$ are non-alternating. In contrast, due
to the different geometric picture, the primeness proof requires
a different, and somewhat longer, discussion. We start with
this. (The primeness proof now does not use, but shows
non-alternation in an alternative way.)

\subsection{Primeness}

\begin{lemma}\label{lpr}
The knots $K_k$ are prime. Hence so are their hypothetic diagrams
$D_k$.
\end{lemma}

Let us make the basic assumption that tangles are considered here
to have an arbitrary number $n\ge 2$ of strings (intervals properly 
embedded in a ball), but no closed components. Moreover, strands are 
unoriented. Tangles are considered equivalent up to homeomorphisms 
of the ball that keep its boundary sphere fixed.

We lean closely on the work of Kirby and Lickorish \cite{KL}.
The notion of a rational tangle is explained there.  

\begin{defi}\label{DF1}
An $n$-string tangle $T$ is \em{prime} if
\def\theenumi{\alph{enumi})}
\def\labelenumi{\theenumi}
\begin{enumerate}
\item \label{CDA}
there is no ball inside $T$ intersecting $T$ in a \em{single}
and knotted arc, and
\item $T$ is not a generalized rational tangle\,.
\end{enumerate}
\end{defi}

\begin{defi}\label{DF2}
An $n$-string tangle $T$ is a \em{generalized rational tangle}
if there is a division $C\cup D$ of the strings of $T$ into two
subsets (with $C\ne\vn\ne D$, and $C\cap D=\vn$), such that each
2-string subtangle of $T$ made up of a string from $C$ and one 
from $D$ is a rational tangle, with possible connected sum of some 
string with non-trivial knots $K'$.
\end{defi}

% 2 %
Below on the left is an example of 4 strings (with $|C|=1$, $|D|=3$),
that visualizes also the meaning of $K'$ and the occurrence of
different rational subtangles.
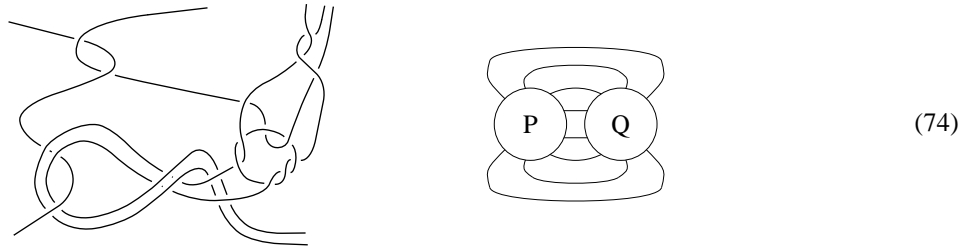
\begin{eqn}\label{62.5}
\epsfsv{3.6cm}{t1-17+4k_ach_tang5}
\kern2cm
\diag{6mm}{4}{3}{
  \pictranslate{0 0.5}{
    \piclinewidth{50}
    \opencurvepath{3.5 1.5}{4 2}{3.8 2.8}{0.2 2.8}{0 2}{0.5 1.5}{}
    \opencurvepath{3.5 0.5}{4 0}{3.8 -0.8}{0.2 -0.8}{0 0}{0.5 0.5}{}
    \opencurvepath{3 1.5}{3.2 2}{2.7 2.35}{1.3 2.35}{0.8 2}{1 1.5}{}
    \opencurvepath{3 0.5}{3.2 0}{2.7 -0.35}{1.3 -0.35}{0.8 0}{1 0.5}{}
    \piccirclearc{2 0.5}{1.3}{45 135}
    \piccirclearc{2 1.5}{1.3}{-135 -45}
    \picline{1 1.3}{3 1.3}
    \picline{1 0.7}{3 0.7}
    \pictranslate{1 1}{
      \picrotate{0}{
        \picscale{1 1}{
	  \picfilledcircle{0 0}{0.8}{P}
        }
      }
    }
    \pictranslate{3 1}{
      \picrotate{0}{
        \picscale{1 1}{
	  \picfilledcircle{0 0}{0.8}{Q}
        }
      }
    }
  }
} 
\end{eqn}
Again, we have, as for 2 strings, the sum of two tangles
(shown above on the right for $n=4$).

\begin{rem}
The local connected sum factors $K'$ at a string $s$ in definition
\ref{DF2} can be often excluded, for example if $s$ is unknotted.
Also, if $|C|=1$, then condition \ref{CDA} in definition \ref{DF1} 
prohibits $K'$ for the string $s\in C$. In particular, for 2-string
tangles $T$, where $|C|=|D|=1$, the $K'$ can be ignored when 
defining a prime tangle, and our definition coincides with 
the one of Kirby and Lickorish. (They allow for closed tangle
components, but for tangles sums which are supposed to be knots,
such components do not occur, so our initial constraint makes sense.)
\end{rem}

The following is rather easy to prove:

\begin{corr}\label{io}
Let $T$ be an $n$-string tangle. Assume for each two strings $s,s'$
of $T$ there is a sequence of strings $s=s_1,s_2,\dots,s_k=s'$ such
that for all $1\le i\le k-1$ the $2$-string subtangle of $T$ made of
$s_i$ and $s_{i+1}$ is a prime tangle (in the sense of \cite{KL}).
Then $T$ is a prime tangle (in the sense of definition \ref{DF2}).
\qed
\end{corr}

\begin{theo}
A knot $K$, which is the sum of two prime tangles $P,Q$ as
on the right of \eqref{62.5}, is prime.
\end{theo}

\proof This is an easy modification of Kirby-Lickorish's proof. 
One needs to look again at how a sphere that determines a factor
decomposition of $K$ would intersect the sphere of its tangle sum
decomposition. The primeness of tangles was defined so as to 
exclude all possibilities for this intersection. \qed

\proof[of lemma \ref{lpr}]
When the tangles $T_k$, whose closure is shown on the left of
\eqref{Ti}, are separated as by the dashed line,
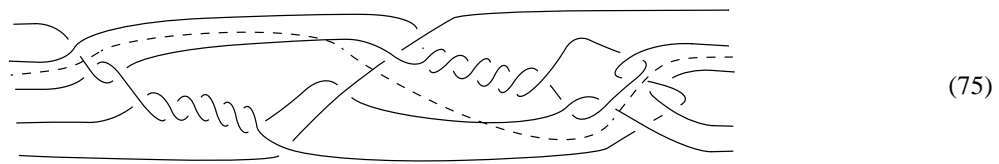
\begin{eqn}%\label{Ti}
\diag{8mm}{12}{4}{
  \picputtext{6 2}{\rbox{\epsfs{4.1cm}{t1-17+4k_ach_tang3}}}
  \piclinewidth{50}
  \piclinedash{0.15 0.15}{0.11}
  \opencurvepath{0 2.2}{1 2.3}{2 2.9}{5 3.1}{6 2.2}{9.2 0.4}
    {10.8 2.5}{12 2.5}{}
}
\end{eqn}
we obtain a sum decomposition of the closure $K_k$ into a
$3$-string tangle $T'_k$ and its mirror image. 
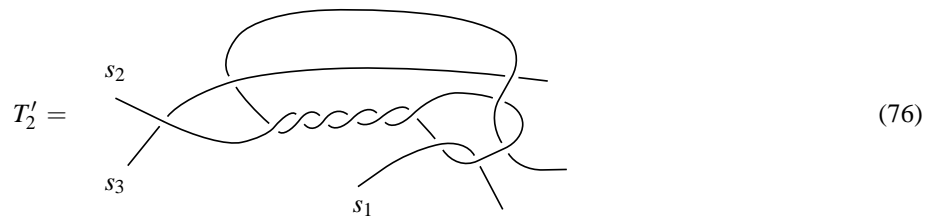
\begin{eqn}\label{Tpk}
T'_2\,=\quad
\diag{1cm}{6.5}{3}{
  \pictranslate{0 -0.4}{
    \picputtext{3.2 2}{\epsfs{4cm}{t1-17+4k_ach_tang4}}
    \picputtext{1.0 0.2 x}{$s_3$}
    \picputtext{2.5 0.2 x}{$s_2$}
    \picputtext{0.7 3.5 x}{$s_1$}
  }
}
\end{eqn}
So it suffices to prove that
$T'_k$ are prime. For this we use corollary \ref{io}. With the
ordering of the strings of $T'_k$ as in \eqref{Tpk}, it is easy
to verify, using \cite{KL}, that $(s_1,s_2)$ and $(s_2,s_3)$ are 
prime 2-string tangles. This completes the primeness proof for
$K_k$. That for $D_k$ follows again at once from the non-triviality
of (semi)adequate links. \qed

\subsection{Basic cases}

For the crossing number proof, we need again several preparations.

\subsubsection{The 17 crossing knot}

Let us first deal with the case $k=0$. We must rule out that
$K_0$ is an adequate non-alternating prime 16 crossing knot.
A check of the Jones polynomial (here $V_2=\bV_2=0$ because 
we have two triangles) of these knots, obtainable from the tables
of \cite{KnotScape}, shows only one knot with matching polynomial, 
$16_{934760}$. But the skein polynomial distinguishes it 
from our (therefore) 17 crossing example $K_0$. We thus assume
for the rest of the discussion that we have diagrams $D_k$ of
$\ge 20$ crossings, for $k>0$, and so $\tg=1$.

Again our proof will consist in trying to construct $A(D_k)$ 
successively by loop attachments, and seeing that we can never do 
it right. The two cycles change, though, the configuration we start 
with. Here we should introduce a case distinction, which will 
accompany us though most of the proof. It bases on the mutual 
position of the two cycles in $A(D_k)$.

Let us again draw the outer cycle in the plane so that all its 
edges lie outside the loops, i.e. in the region of their complement
that contains infinity. This position then defines a notion of
\em{interior} and \em{exterior} (resp. \em{in/outside}) for all 
loops, incl. those to be subsequently attached. For example, if a 
loop $M$ is attached inside (i.e. in the interior) of $L$, then
$L$ lies in the exterior of $M$. A loop is \em{empty} if its
interior is empty, i.e. contains no traces or other loops.
(Note that for a semiadequate diagram, no loops implies
no traces.) A loop is separating if both its interior
and exterior are non-empty (see \S\reference{S3c}).

\subsubsection{Cases A,B,C}

Let us first consider the case that the (edges of the)
triangle lie(s) on the same side of the outer cycle loops
as the outer cycle.

{\bf Case A.} The triangle and outer cycle have at most one
common loop. This case will be easily out. If there is no common 
loop, then we have, prior to attachments, something like
the figure (A1) in \eqref{csA}.
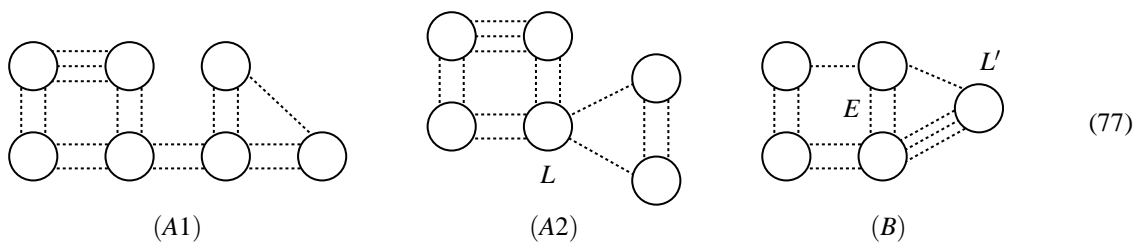
\begin{eqn}\label{csA}
\begin{array}{c@{\quad\qquad}c@{\quad\qquad}c}
\diag{8mm}{5.6}{2.4}{
  {\piclinedash{0.05 0.05}{0.01}
   \picline{0.2 0.2}{5.2 0.2}
   \picline{0.2 0.6}{5.2 0.6}
   \picline{0.2 2.15}{2.2 2.15}
   \picline{0.2 1.9}{2.2 1.9}
   \picline{0.2 1.65}{2.2 1.65}
   \picline{5.2 0.6}{3.6 2.1}
  }
  \picmultigraphics{3}{1.6 0}{
    {\piclinedash{0.05 0.05}{0.01}
     \picline{0.2 0.2}{0.2 1.8}
     \picline{0.6 0.2}{0.6 1.8}
    }
    \picfilledcircle{0.4 1.9}{0.4}{}
  }
  \picmultigraphics{4}{1.6 0}{
    \picfilledcircle{0.4 0.4}{0.4}{}
  }
  % \picfilledcircle{2.0 1.9}{0.4}{}
} &
\diag{8mm}{4.2}{3.4}{
  \pictranslate{0 1}{
    {\piclinedash{0.05 0.05}{0.01}
     \picline{0.2 0.2}{2.0 0.2}
     \picline{0.2 0.6}{2.0 0.6}
     \picline{0.2 2.15}{2.2 2.15}
     \picline{0.2 1.9}{2.2 1.9}
     \picline{0.2 1.65}{2.2 1.65}
     % \picline{0.2 2.2}{2.2 2.2}
     % \picline{0.2 1.9}{2.2 1.9}
     % \picline{0.2 1.6}{2.2 1.6}
     \picline{3.6 -0.5}{3.6 1.2}
     \picline{4.0 -0.5}{4.0 1.2}
     \picline{3.6 -0.5}{2.0 0.4}
     \picline{3.6 1.2}{2.0 0.4}
    }
    \picmultigraphics{2}{1.6 0}{
      {\piclinedash{0.05 0.05}{0.01}
       \picline{0.2 0.2}{0.2 1.8}
       \picline{0.6 0.2}{0.6 1.8}
      }
      \picfilledcircle{0.4 1.9}{0.4}{}
    }
    \picmultigraphics{2}{1.6 0}{
      \picfilledcircle{0.4 0.4}{0.4}{}
    }
    \picfilledcircle{3.8 -0.5}{0.4}{}
    \picfilledcircle{3.8 1.2}{0.4}{}
    \picputtext{2 -0.4}{$L$}
  }
} &
\diag{8mm}{4.2}{2.4}{
  \pictranslate{0 0}{
    {\piclinedash{0.05 0.05}{0.01}
     \picline{0.2 0.2}{2.0 0.2}
     \picline{0.2 0.6}{2.0 0.6}
     \picline{0.2 1.9}{2.2 1.9}
     \picline{3.8 1.1}{2.0 0.1}
     \picline{3.8 1.3}{2.0 0.3}
     \picline{3.8 1.5}{2.0 0.5}
     \picline{3.4 1.5}{2.0 2.1}
    }
    \picmultigraphics{2}{1.6 0}{
      {\piclinedash{0.05 0.05}{0.01}
       \picline{0.2 0.2}{0.2 1.8}
       \picline{0.6 0.2}{0.6 1.8}
      }
      \picfilledcircle{0.4 1.9}{0.4}{}
    }
    \picmultigraphics{2}{1.6 0}{
      \picfilledcircle{0.4 0.4}{0.4}{}
    }
    \picfilledcircle{3.6 1.2}{0.4}{}
    \picputtext{3.8 2.0}{$L'$}
  }
  \picputtext{1.5 1.2}{$E$}
} 
\\
(A1) & (A2) & (B)
\end{array}
\end{eqn}
% \begin{eqn}\label{uu}
% \diag{1cm}{6}{2.4}{
%   {\piclinedash{0.05 0.05}{0.01}
%    \picline{0.2 0.2}{5.2 0.2}
%    \picline{0.2 0.6}{5.2 0.6}
%    \picline{0.2 2.2}{2.2 2.2}
%    \picline{0.2 1.9}{2.2 1.9}
%    \picline{0.2 1.6}{2.2 1.6}
%    \picline{5.2 0.6}{3.6 2.1}
%   }
%   \picmultigraphics{3}{1.6 0}{
%     {\piclinedash{0.05 0.05}{0.01}
%      \picline{0.2 0.2}{0.2 1.8}
%      \picline{0.6 0.2}{0.6 1.8}
%     }
%     \picfilledcircle{0.4 1.9}{0.4}{}
%   }
%   \picmultigraphics{4}{1.6 0}{
%     \picfilledcircle{0.4 0.4}{0.4}{}
%   }
%   % \picfilledcircle{2.0 1.9}{0.4}{}
% }\,,
% \end{eqn}
To these cycles we must attach loops as in \eqref{attach}. However,
we see that (A1) of \eqref{csA} has already 3 connected sum 
factors, and so atoms. Thus we have only one attachment left by
corollary \ref{cr3l}, but this cannot make the diagram $D$ prime.

If the outer cycle and triangle have one common loop, then
we have the figure (A2) in \eqref{csA},
and two attachments left. In order to create two cycles in
$IG(A(D))$, we must attach both loops inside the common loop
$L$, so that their legs become intertwined with at least three
edges outside $L$. This is handled in \S\reference{Sca}.

We are next left to treat the cases B and C.

{\bf Case B.} The triangle and outer cycle have two common loops.
This is case (B) in \eqref{csA}. Here we have again 3 attachments
\eqref{attach} to perform. There is a degenerate instance of this
case, which we call case B', in which the common connection $E$
of the outer cycle and triangle goes on both sides of the
triangle loop $L'$.
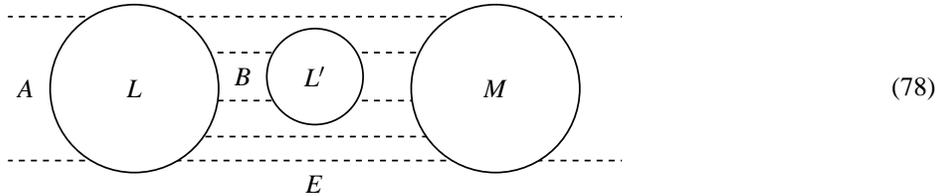
\begin{eqn}\label{B'}
\diag{1.6cm}{5.5}{2}{
  \piclinewidth{70}
  \pictranslate{1.25 0}{
  {\piclinedash{0.05 0.05}{0.01}
   \picline{-1.05 0.4}{4.05 0.4}
   \picline{0.2 0.6}{3.2 0.6}
   \picline{0.2 1.3}{3.2 1.3}
   \picline{0.2 0.9}{3.2 0.9}
   \picline{-1.05 1.6}{4.05 1.6}
  }
  \picfilledcircle{0 1}{0.7}{$L$}
  \picfilledcircle{1.5 1.1}{0.4}{$L'$}
  \picfilledcircle{3 1}{0.7}{$M$}
  \picputtext{1.5 0.2}{$E$}
  \picputtext{0.9 1.1}{$B$}
  \picputtext{-0.9 1.0}{$A$}
 }
}
\end{eqn}
This peculiarity will cost us some extra effort later.
At least we can assume that, if $M$ is empty, $L$ has
attachments inside that connect the regions $A$ and
$B$. (Otherwise we can return to case B by a mutation.)

Still we will see below that Cases B and B' are relatively
simple. The other case requires the most work.

{\bf Case C.} The triangle and outer cycle lie on opposite 
sides of a loop. In that case, we start with the outer cycle 
\eqref{oaz} (without insistence of odd length), and have three
attachments to perform, but among them there is one of the sort 
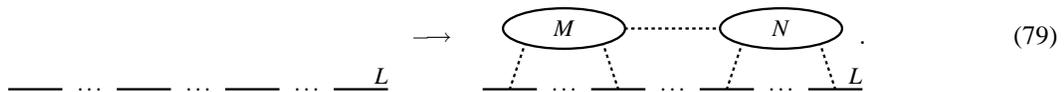
\begin{eqn}\label{tatt}
{\small
\diag{9mm}{5.5}{2.5}{
  \picmultigraphics{4}{1.6 0}{
    \picline{0.5 0 x}{0.5 0.8 x}
  }
  \picmultigraphics{3}{1.6 0}{
    \picputtext{0.5 1.2 x}{$\dots$}
  }
  % \picline{0.5 1.6 x}{0.5 2.4 x}
  % \picputtext{0.5 2.8 x}{$\dots$}
  % \picline{0.5 3.2 x}{0.5 4.0 x}
  \picputtext{5.5 0.7}{$L$}
}
\quad\lra\quad
\diag{9mm}{5.5}{2.5}{
  \picmultigraphics{4}{1.6 0}{
    \picline{0.5 0 x}{0.5 0.8 x}
  }
  \picmultigraphics{3}{1.6 0}{
    \picputtext{0.5 1.2 x}{$\dots$}
  }
  {\piclinedash{0.05 0.05}{0.01}
   \picline{1.2 1.4}{4.4 1.4}
  }
  \picmultigraphics{2}{3.2 0}{
   {\piclinedash{0.05 0.05}{0.01}
    \picline{0.5 0.4 x}{1.4 0.7 x}
    \picline{0.5 2.0 x}{1.4 1.7 x}
   }
   \picfilledellipse{1.2 1.4}{0.9 0.3}{}
  }
  \picputtext{1.2 1.4}{$M$}
  \picputtext{4.4 1.4}{$N$}
  \picputtext{5.5 0.7}{$L$}
}
}\,.
\end{eqn}
The loops $M,N$ will be often referred to as the \em{triangle loops}.
Again the same arguments concerning attachments and case distinction 
apply for the $B$-state of $D_k$. Often it will help to switch from 
$A(D_k)$ to $B(D_k)$, and to see that there some feature becomes 
violated. We will exhibit several features in the next lemmas.
These properties must apply to both the $A$- and $B$-state.

Let us in advance fix the following \em{complexity}. We regard
case A as simpler than case B, which in turn is simpler than case C.
We choose $D_k$ mirrored so that $A(D_k)$ is the simplest w.r.t.
this hierarchy. Thus, for example, if in the treatment of ($A(D_k)$
being of) case C we establish that $B(D_k)$ is of case A or B, then
we can consider this $D_k$ as dealt with. Note that to establish
case A or B, it is enough to identify two cycles whose edges do
not lie on opposite sides of some loop.

Next let us reintroduce the notion of \em{depth} of loops. This is 
relatively easily adapted from the previous section. 
% Let us leave out case A, which is already finished.

In our case A or B the additional loop(s) $L'$ on the
triangle but not on the outer cycle has/have depth 0.
Loops attached outside to $L'$ (i.e. on the side that contains the
edges of the triangle and outer cycle) have depth 0, loops attached
inside $L'$ have depth 1, loops attached inside those have depth 2
etc.

In case C the assignment of depth to the loops needs to be clarified
only for the triangle attachment \eqref{tatt}. This is done in the
obvious way, by saying that these 2 loops would have the depth of
a loop that would be attached instead of them by \eqref{attach}.

Let us again call an \em{edge} a parallel equivalence class of crossing
traces in $A(D)$. The multiplicity of an edge and the property to 
be simple are then obviously defined. Again we often regard a letter 
naming an edge as a variable that indicates its multiplicity. 

With this set up, we assume again $D=D_k$ is \em{flat}, i.e. of 
minimal loop depth sum in its mutation equivalence class. Also we
assume it \em{edge-reduced}, i.e. within the diagrams of the same
minimal total loop depth, it should have the smallest number of edges
(equivalence classes of parallel traces). In the following, we will
several times apply a mutation to bring a diagram into a more
favorable shape. While we do not say it each time explicitly,
it is to be understood that such a mutation can be, and is, chosen
so as to preserve the minimal loop depth and edge-reducedness.

\subsection{Preparatory considerations}

\subsubsection{Control on position of attachments}

In order to deal with cases B and C we will need several preparatory
lemmas. We postpone case B' to \S\reference{CBB'}.

\begin{lemma}
The attachment \eqref{tatt} in case C is made inside an
outer cycle loop (i.e. $L$ is an outer cycle loop, and $M,N$
are attached inside $L$, and of depth $1$).
\end{lemma}

\proof Let $L$ be the loop the attachment \eqref{tatt} is made at.
If $L$ is not an outer cycle loop, then $L$ is itself an attached 
loop. So after $L$ and \eqref{tatt} we still have no cycle in $IG(A)$,
but only (at most) one attachment \eqref{attach} left. It is 
easy to see that, in which way we even choose to make this attachment, 
we cannot create two cycles in $IG(A)$. So $L$ is an outer cycle 
loop. Then by definition of case C, \eqref{tatt} is made inside $L$.
(Otherwise we have case A.) \qed

An edge adjacent to an attached loop is called a \em{leg}. 
Edges between the two attached loops $M,N$ in \eqref{tatt} may 
not be regarded as legs. The valence of an attached loop using
\eqref{attach} is as before the number of its legs. This number
coincides with the number of regions touched by the legs, under
exclusion of peculiarities \eqref{pecu}. Because of primeness
and $\gm=4$, the number of regions touched is always 2 or 3.
We do not use a notion of valence for the loops $M,N$ of
\eqref{tatt}, but we will distinguish cases by the number of
regions (again 2 or 3) touched by legs of $M$ or $N$.

Under exclusion of peculiarities \eqref{pecu} (which we will soon
accomplish with lemma \reference{lmZ'}), the triangle
attachment \eqref{tatt} is of the shape shown on the left:
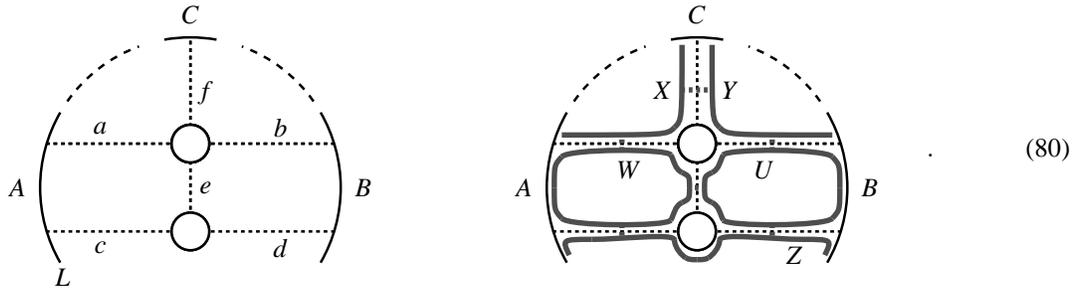
\begin{eqn}\label{tat2}
\diag{1cm}{6}{4}{
  \pictranslate{3 1.5}{
    \picmultigraphics[S]{2}{-1 1}{
      \piccirclearc{0 0}{2}{-30 30}
      {\piclinedash{0.05 0.05}{0.01}
       \picline{2 -17 polar x p 0 x}{0 2}
       \piclinedash{0.1 0.1}{0}
       \piccirclearc{0 0}{2}{35 70}
      }
      \picmultigraphics[S]{3}{1 -1}{
	{\piclinedash{0.05 0.05}{0.01}
	 \picline{2 17 polar}{2 17 polar x p 0 x}
	}
	\picfilledcircle{2 17 polar x p 0 x}{0.25}{}
      }
    }
    \piccirclearc{0 0}{2}{80 100}
    \picputtext{0.2   0.0}{$e$}
    \picputtext{0.2   1.3}{$f$}
    \picputtext{1.2   0.8}{$b$}
    \picputtext{1.2  -0.8}{$d$}
    \picputtext{-1.2  0.8}{$a$}
    \picputtext{-1.2 -0.8}{$c$}
    \picputtext{-2.3 -0.0}{$A$}
    \picputtext{ 2.3 -0.0}{$B$}
    \picputtext{0.0  2.3}{$C$}
  }
  \picputtext{ 1.3 0.3}{$L$}
}\qquad
\diag{1cm}{6}{4}{
  \pictranslate{3 1.5}{
    \picmultigraphics[S]{2}{-1 1}{
      \piccirclearc{0 0}{2}{-30 30}
      {\piclinedash{0.05 0.05}{0.01}
       \picline{2 -17 polar x p 0 x}{0 2}
       \piclinedash{0.1 0.1}{0}
       \piccirclearc{0 0}{2}{35 70}
      }
      \picmultigraphics[S]{3}{1 -1}{
	{\piclinedash{0.05 0.05}{0.01}
	 \picline{2 17 polar}{2 17 polar x p 0 x}
	}
	\picfilledcircle{2 17 polar x p 0 x}{0.25}{}
      }
    }
    \piccirclearc{0 0}{2}{80 100}
    \picputtext{-2.3 -0.0}{$A$}
    \picputtext{ 2.3 -0.0}{$B$}
    \picputtext{ 0.45 1.3}{$Y$}
    \picputtext{-0.45 1.3}{$X$}
    \picputtext{1.3 -0.9}{$Z$}
    \picputtext{-0.9 0.24}{$W$}
    \picputtext{0.9 0.24}{$U$}
    \picputtext{0.0  2.3}{$C$}
    { \piclinewidth{14}
      \picgraycol{0.3}
      \picPSgraphics{0 setlinecap}
      \picmultigraphics[S]{2}{-1 1}{
         \opencurvepath{0.2 1.9}{0.2 0.9}{0.4 0.7}{1.8 0.7}{}
	 \opencurvepath{1.7 -0.9}{1.8 -0.69}{1 -0.69}{0.3 -0.6}
	   {0.3 -0.8}{0.2 -0.91}{0.1 -0.96}{0 -0.96}{}
	 \picmultigraphics[S]{2}{1 -1}{
	   \opencurvepath{0.1 0}{0.1 0.15}{0.3 0.25}{0.3 0.35}{0.4 0.52}
	     {1.7 0.52}{1.9 0.3}{1.9 0}{}
	   {\piclinedash{0.05 0.05}{0.01}
	    \picline{1 0.5}{1 0.7}
	   }
	 }
      }
      {\piclinedash{0.05 0.05}{0.01}
       \picline{-0.2 1.3}{0.2 1.3}
       \picline{-0.12 0.0}{0.12 0.0}
      }
    }
  }
}\,.
\end{eqn}
Here we assume the regions $A,B,C$ outside $L$ along each segment
drawn in solid line to be the same (i.e. solid line segments
are not touched by traces from outside $L$). For future reference,
let us call the edges of the triangle $a$ to $f$. Hereby it is
possible that some of $a,\dots,f$ is (of multiplicity) zero. In
fact, this necessarily occurs, as shows the following lemma.

\begin{lemma}\label{lm1014}
We have in case $C$ no triangle attachment of the type \eqref{tat2}
with empty triangle loops and none of $a,b,c,d$ being $0$. In other
words, one of the loops of the triangle has a leg that touches
at most one of $A$ or $B$.
\end{lemma}

\proof First consider the case $f=0$ or $f=2$.
By connectivity again at least of $a,b,c,d$ is $\ge 2$, 
and then by direct observation $B(D)$, whose loops are shown
on the right of \eqref{tat2}, has two cycles not lying
on opposite sides of a loop. Thus $B(D)$ is of type $B$.

Next assume $f\ge 3$. If $X=Y$, we have a cycle in $B(D)$ made
up of $X$ and the loops between the traces of $f$. The triple
$X,V,W$ gives a second cycle, so $B(D)$ is of type $A$.

If $X\ne Y$, then we still have the cycles $U,W,Z$ and the
one containing $U,W,X,Y$. They are distinct even if $Z$ is
one of $X$ or $Y$. So $B(D)$ is of type $B$.
Finally, if $f=1$, then by adequacy of $D$
necessarily $X\ne Y$, and then the same two cycles exist. \qed

\begin{lemma}\label{lA1}
There is at most one non-empty outer cycle loop in case A.
\end{lemma}

\proof To create a cycle in $IG(A)$, we need both attachments 
to be within the same loop. \qed

\begin{lemma}\label{lC1}
We can assume, up to mutation, that there is at most one
non-empty outer cycle loop in case C.
\end{lemma}

\proof Compare the proof of lemma \reference{l1ne}. We assume
we have loops attached within two outer cycle loops and derive
a contradiction.

If we had an attachment inside some of the triangle loops,
then this would be the third attachment, and we cannot create
2 cycles in $IG(A)$. So the triangle loops are empty.

Then if the triangle attachment identifies only two regions $A$
and $B$, we have \eqref{tat2} with $f=0$, but $a,b,c,d$ non-zero
to have two cycles in $IG(A)$. We ruled this out by lemma \ref{lm1014}.
(Note that, with only $A$ and $B$ touched, we cannot have two legs
touching the same region on either side of an attached loop, as
in \eqref{pecu}, because then $D$ would not be prime.)
% $M,N$ will not touch any other region, and with the third
% attachement \eqref{attach}, we cannot create a loop in $IG(A)$.

Thus the triangle attachment identifies 3 regions, say $A,B,C$.
(Here we must include the situation similar to \eqref{pecu}, but
with $M,N$ connected by a trace.)
Now, an attachment in any other outer cycle loop, must identify
only 2 or all 3 of these regions. In former case, one can
achieve by a mutation (which does not spoil flatness or
edge-reducedness) that only one outer cycle loop is non-empty.
In latter case, we have a modification of \eqref{eqw}, in which
$X$ and $Y$ may contain several loops, but their legs still
touch only $A,B,C$. Then, however, again $IG(A)$ has no cycle.
\qed

\begin{lemma}\label{lB1}
There is at most one non-empty outer cycle or triangle loop in case B.
\end{lemma}

\proof If we had 3 attachments \eqref{attach} divided into
several outer cycle or triangle loops, then we must create
two cycles in $A(D)$ with two of them $M,N$, attached in
some (outer cycle or triangle) loop $L$. So each of the 
loops $M,N$ inside $L$ must be intertwined with $3$ edges
outside $L$. Now, since we have 3 attachments, and one is
inside another loop different from $L$, we cannot have outside
attachments to $L$. This means that $M,N$ must be intertwined
with 3 edges outside $L$ from the outer cycle or triangle.
While there are two loops $L$ with three outside connections,
one cannot (unlike in case B', which we will treat later)
install $M,N$ to intertwine with all of them,
without their legs intersecting, or $\gm>4$. \qed

Let (as before) $\tl D$ be obtained from $D=D_k$ by making all 
edges in $A(D)$ simple (i.e. a single trace). We call again 
an edge (parallel equivalence class of traces) in $A(D)$ to 
be \em{inadequate} if, when reduced to a single crossing, this
crossing is $B$-inadequate in $\tl D$, i.e. has a trace in 
$B(\tl D)$ that connects a loop to itself. 

% ]
% a loop of
% depth 0 attached outside a 

\begin{lemma}\label{x2}
If $x$ is an inadequate edge, then $x=2$.
\end{lemma}

\proof In order $D_k$ to be $B$-adequate, we must have 
$x\ge 2$. If $x>2$, then we have a cycle (triangle, or outer 
cycle) in $B(D)$ of simple edges, i.e. each pair of neighbored
loops on this cycle is connected only by one trace. Moreover, 
this cycle has only one loop with non-empty interior, and 
all its other loops have only two adjacent traces. If $B(D)$ is
of type $A$, then it is easy to see that with two attachments
\eqref{attach} one cannot create two cycles in $IG(B(D))$. This
means that $B(D)$ must be of type $C$. However, since one of 
its cycles has only simple edges, one sees again that with the 
two other attachments \eqref{attach} we cannot create two
cycles in $IG(B)$. \qed

Again an attached loop (though not some in \eqref{tatt})
of valence $2$ is inadequate, if both its legs are inadequate. 

\begin{corr}\label{cfg}
There are no inadequate loops.
\end{corr}

\proof If both their legs have multiplicity 2, the diagram is
disconnected. \qed

\begin{lemma}
There is no valence three loop, attached by \eqref{attach},
which has non-empty interior.
\end{lemma}

\proof Case A is easy, so consider cases B, B' and C.
Assume $L$ were such a loop, and $M$, resp. $M,N$,
the loop(s) attached inside $L$ using \eqref{attach}, resp. 
\eqref{tatt} (in case C). The attachments of $L$, $M$ and 
possibly $N$ do not create any cycle in $IG(A)$. By assuming 
flatness of $D$, we see also that attaching $L$, $M$ and 
possibly $N$ augments $\gm$ by at least 4 (see \eqref{Fig17}).
This means that 
the third attached loop $P$ is either outside an outer cycle 
loop (or the extra triangle loop in case B), or inside $L$.
In either case it is easy to see that this attachment cannot
create any cycle in $IG(A)$. \qed

\begin{lemma}
There is no valence two loop, attached by \eqref{attach},
which has non-empty interior.
\end{lemma}

\proof Again case A is easily settled, and excluded.
Assume $L$ were such a loop, and $M$, resp. $M,N$,
the loop(s) attached inside $L$ using \eqref{attach}, resp. 
\eqref{tatt} (in case C). If \eqref{tatt} is used inside $L$,
we have only one instance of \eqref{attach} left, but cannot
create two cycles in $IG(A)$ with it. The same situation occurs
if we use \eqref{attach} twice inside $L$. So we can assume that
the loop $M$ is the only loop inside $L$. If now $M$ has valence
$2$, then by flatness of $D$, we must have the picture \eqref{Fig5},
and $M$ is inadequate, in contradiction to corollary
\ref{cfg}. Thus $M$ must have valence 3, and we have
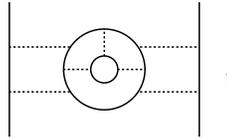
\[
\diag{6mm}{5}{3}{
  \picline{0.4 0}{0.4 3}
  \picline{4.6 0}{4.6 3}
  {\piclinedash{0.05 0.05}{0.01}
   \picline{1 0.4 x}{1 4.6 x}
   \picline{2 0.4 x}{2 4.6 x}
  }
  \picfilledcircle{2.5 1.5}{0.9}{}
  {\piclinedash{0.05 0.05}{0.01}
   \picline{1.6 1.5}{3.4 1.5}
   \picline{2.5 1.5}{2.5 2.4}
  }
  \picfilledcircle{2.5 1.5}{0.3}{}
}\ ,
\]
or a super-configuration thereof, but then $\gm\ge 6$. \qed

\begin{corr}\label{crdp}
There are no loops of depth $\ge 3$, and the only possible
ones of depth $2$ are (in case C and) attached inside triangle 
loops \eqref{tatt}. \qed
\end{corr}

% \begin{lemma}
% There is at most one attachment \eqref{attach} inside a
% triangle loop (i.e. such attached by \eqref{tatt}).
% \end{lemma}
% 
% \proof If there were two such, we would exhaust the
% three attachments we can perfrom, and have at most 1
% cycle in $IG(A)$. \qed

We will later rule out, during the detailed treatment of 
case C, depth-2 loops completely (lemma \reference{ldp2}).

The next step is to get disposed of things like \eqref{pecu}.

\begin{lemma}\label{lmZ'B}
In case A or B, there is no region touched in two different
segments by legs of (possibly different) attachments \eqref{attach}
on the opposite side.
\end{lemma}

\proof In case A, we have no outside attachment to $L$ (both
attachments must be inside to create cycles in $IG(A)$). In case B,
the proof is the same as for lemma \reference{lmZ}. \qed

\begin{lemma}\label{lmZ'}
In case C, there is no region touched in two different segments by
legs of attachments \eqref{attach} or \eqref{tatt} on the opposite side.
\end{lemma}

% NOTE THAT WE DID NOT RULE OUT THE CASE THAT attachments of
% different loops are made on different segments of the same region.
% now this is done!

\proof If $Z$ in \eqref{pecu} were attached inside an outer cycle
loop $L$, then we would have at most one loop outside
attached to $L$ (the second attachment \eqref{attach}),
and thus would reduce the number of legs by a mutation.
Thus $Z$ is attached outside $L$.
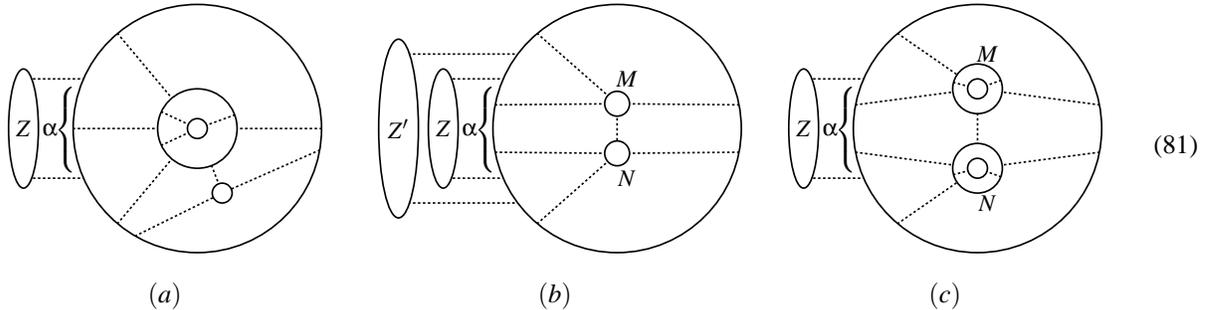
\begin{eqn}\label{pecu2}
\begin{array}{c@{\kern6mm}c@{\kern6mm}c}
\diag{6.6mm}{6.3}{5.9}{
  \small
  \pictranslate{3.8 3}{
   {\piclinedash{0.05 0.05}{0.01}
    \picline{-3.4 1.0}{0 1.0}
    \picline{-3.4 -1.0}{0 -1.0}
   }
   \picfilledellipse{0 -3.5 x}{0.3 1.2}{$Z$}
   \picfilledcircle{0 0}{2.5}{}
   \picputtext{0.39 -2.75 x}{$\ap\!\!\left\{\ry{0.7cm}\right.$}%}
   {\piclinedash{0.05 0.05}{0.01}
    \picline{2.5 0 polar}{0 0}
    \picline{2.5 180 polar}{0 0}
    \picline{2.5 130 polar}{0 0}
    \picline{2.5 -130 polar}{0 0}
    \picline{0 0}{0.5 -1.3}
    \picline{2.5 -120 polar}{0.5 -1.3}
    \picline{2.5 -10 polar}{0.5 -1.3}
   }
   \picfilledcircle{0 0}{0.8}{}
   \picfilledcircle{0.5 -1.3}{0.2}{}
   {\piclinedash{0.05 0.05}{0.01}
    \picline{0.8 20 polar}{0 0}
    \picline{0.8 155 polar}{0 0}
    \picline{0.8 -155 polar}{0 0}
   }
   \picfilledcircle{0 0}{0.2}{}
  }
}
&
\diag{6.6mm}{7.5}{5.9}{
  \small
  \pictranslate{5 3}{
   {\piclinedash{0.05 0.05}{0.01}
    \picline{-4.4 1.5}{0 1.5}
    \picline{-4.4 -1.5}{0 -1.5}
    \picline{-3.4 1.0}{0 1.0}
    \picline{-3.4 -1.0}{0 -1.0}
   }
   \picfilledellipse{0 -3.5 x}{0.3 1.2}{$Z$}
   \picfilledellipse{0 -4.4 x}{0.4 1.8}{$Z'$}
   \picfilledcircle{0 0}{2.5}{}
   \picputtext{0.39 -2.75 x}{$\ap\!\!\left\{\ry{0.7cm}\right.$}%}
   {\piclinedash{0.05 0.05}{0.01}
    \picline{2.5 11 polar}{0 0.5}
    \picline{2.5 169 polar}{0 0.5}
    \picline{2.5 130 polar}{0 0.5}
    \picline{2.5 -169 polar}{0 -0.5}
    \picline{2.5 -11 polar}{0 -0.5}
    \picline{2.5 -130 polar}{0 -0.5}
    \picline{0 0.5}{0 -0.5}
   }
   \picfilledcircle{0 0.5}{0.25}{}
   \picfilledcircle{0 -0.5}{0.25}{}
   \picputtext{0.2 1.0}{$M$}
   \picputtext{0.2 -1.0}{$N$}
  }
}
&
\diag{6.6mm}{6.3}{5.9}{
  \small
  \pictranslate{3.8 3}{
   {\piclinedash{0.05 0.05}{0.01}
    \picline{-3.4 1.0}{0 1.0}
    \picline{-3.4 -1.0}{0 -1.0}
   }
   \picfilledellipse{0 -3.5 x}{0.3 1.2}{$Z$}
   \picfilledcircle{0 0}{2.5}{}
   \picputtext{0.39 -2.75 x}{$\ap\!\!\left\{\ry{0.7cm}\right.$}%}
   {\piclinedash{0.05 0.05}{0.01}
    \picline{2.5 11 polar}{0 0.8}
    \picline{2.5 169 polar}{0 0.8}
    \picline{2.5 130 polar}{0 0.8}
    \picline{2.5 -169 polar}{0 -0.8}
    \picline{2.5 -11 polar}{0 -0.8}
    \picline{2.5 -130 polar}{0 -0.8}
    \picline{0 0.5}{0 -0.5}
   }
   \picfilledcircle{0 0.8}{0.5}{}
   \picfilledcircle{0 -0.8}{0.5}{}
   \picmultigraphics[S]{2}{1 -1}{
     \pictranslate{0 0.8}{
       {\piclinedash{0.05 0.05}{0.01}
	\picline{0.5 20 polar}{0 0}
	\picline{0.5 165 polar}{0 0}
	% \picline{0.5 -135 polar}{0 0}
       }
       \picfilledcircle{0 0}{0.2}{}
     }
   }
   \picputtext{0.2 1.5}{$M$}
   \picputtext{0.2 -1.5}{$N$}
  }
}
\\
(a) & (b) & (c)
\end{array}
\end{eqn}

Again $Z$ has valence $2$, else we would have $\gm>4$, or a
composite diagram $D$. If $Z$ is not inadequate, then either (i)
there is a depth-2 loop inside one of the triangle loops $M,N$, as
in \eqref{pecu2}(a), or (ii) there is a second attachment outside
$L$, as in \eqref{pecu2}(b). In either case it is impossible to have
the triangle \eqref{tatt} and the second attachment \eqref{attach}
inside $L$ of depth 1. In particular, the segments of $L$ on either
side of $\ap$ must be touched from the inside by legs of some of 
the two triangle loops, and only of those loops.

If (i) holds, then we need both $M,N$ to be intertwined with
$Z$, in order to have two cycles in $IG(A)$. In this case,
however, we need depth-2 loops within both $M$ and $N$ to
avoid that $Z$ is inadequate, as in \eqref{pecu2}(c).
However, then we need 4 attachments.

Thus we have a second attachment $Z'$ outside $L$. We claim that
$Z$ and $Z'$ are parallel. For the argument compare also the proof
of lemma \reference{lm2l} below. 

Let us say that $Z$ and $Z'$ \em{enclose} if they look like
in \eqref{pl} or the left diagram in \eqref{op2}, but not like
in \eqref{nenc} and the right diagram in \eqref{op2}.

If $Z$ and $Z'$ would not enclose, then there are three regions
outside $L$ (the two separated by the legs of $Z$, and the one
enclosed between $L$ and $Z'$) that must be touched by legs
inside $L$, and because $\gm=4$, there are no further ones. Now
because $D$ is prime, we must have $Z$ and $Z'$ being separated
by traces from $L$ to some other outer cycle loop (like in
\eqref{nenc}, and unlike the right diagram in \eqref{op2}).
But then again $Z$ would be inadequate.

Thus $Z$ and $Z'$ enclose. Then, since there is
no loop of depth $>1$, again it is easy to see that
if $Z$ and $Z'$ were not parallel, one would be inadequate.
(See case \ref{CS1} in the proof of lemma \ref{lm2l}.)

Now we have established that $Z$ and $Z'$ are parallel. Again if
some of $M$ or $N$ has no leg on the segment $\ap$ of $L$, then
it intertwines with at most two connections outside $L$,
and then we have at most one cycle in $IG(A)$. Therefore, 
we must have the picture \eqref{pecu2}(b).

We may assume w.l.o.g., using the two cycles in $IG(A)$ and
up to mutation, that $M,N$ intertwine the connection between
$L$ and an outer cycle loop on the left as drawn below.
Similarly, we can assume that between the two outer cycle
edges to the right at least the lower one is there.
Then we have (recall convention \reference{cthi}): 
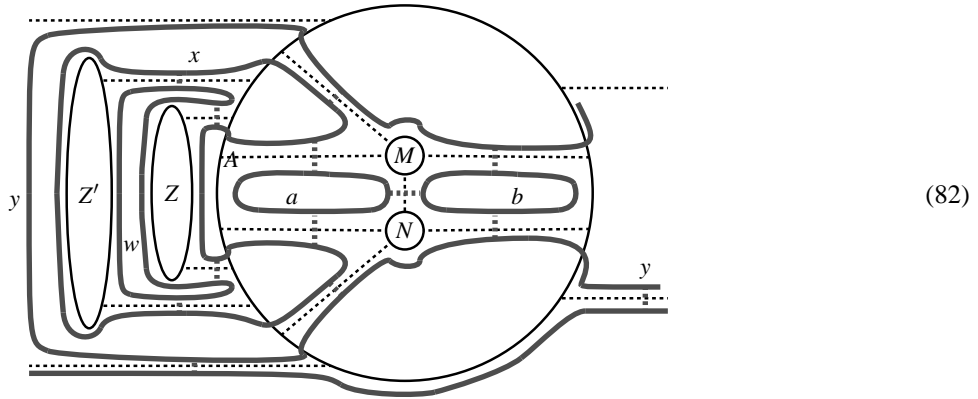
\begin{eqn}\label{mik}
\diag{10mm}{8.5}{5.9}{
  \small
  \pictranslate{5 3}{
   {\piclinedash{0.05 0.05}{0.01}
    \picline{-5 2.3}{0 2.3}
    \picline{-5 -2.3}{0 -2.3}
    \picline{0 -1.4}{3.5 -1.4}
    \picline{0 1.4}{3.5 1.4}
    \picline{-4.2 1.5}{0 1.5}
    \picline{-4.2 -1.5}{0 -1.5}
    \picline{-3.1 1.0}{0 1.0}
    \picline{-3.1 -1.0}{0 -1.0}
   }
   \picfilledellipse{0 -3.1 x}{0.27 1.16}{$Z$}
   \picfilledellipse{0 -4.2 x}{0.3 1.8}{$Z'$}
   \picfilledcircle{0 0}{2.5}{}
   {\piclinedash{0.05 0.05}{0.01}
    \picline{2.5 11 polar}{0 0.5}
    \picline{2.5 169 polar}{0 0.5}
    \picline{2.5 131 polar}{0 0.5}
    \picline{2.5 -169 polar}{0 -0.5}
    \picline{2.5 -11 polar}{0 -0.5}
    \picline{2.5 -131 polar}{0 -0.5}
    \picline{0 0.5}{0 -0.5}
   }
   \picfilledcircle{0 0.5}{0.25}{$M$}
   \picfilledcircle{0 -0.5}{0.25}{$N$}
  }
  \pictranslate{0 3}{
   {
    \piclinewidth{14}
    \picgraycol{0.3}
    \picPSgraphics{0 setlinecap}
    \picmultigraphics{2}{2.5 0}{
    \curvepath{4.0 0.25}{3.9 0.26}{3.4 0.27}{2.9 0.3}{2.7 0}
      {2.7 -0.2}{3.4 -0.25}{4.3 -0.25}{4.8 -0.2}{4.8 -0.0}{4.7 0.25}{}
    }
    \pictranslate{1.5 0}{
     \picmultigraphics[S]{2}{1 -1}{
      \opencurvepath{-0.3 0}{-0.3 1.3}{0.0 1.4}{1.2 1.45}{1.2 1.05}
	{0.1 1.5}{0.02 0.5}{0.0 0.35}{0.0 0.15}{0.0 0}{}
     }
    }
    \picscale{1 -1}{
      \opencurvepath{0 0}{0 2.0}{0.5 2.2}{3.9 2.3}{3.4 2.0}{4.7 0.8}
      {4.9 1}{5.2 1}{5.2 0.8}{5.5 0.6}{7.8 0.6}{7.3 1.2}{7.4 1.25}
      {8.4 1.25}{}
    }
    \opencurvepath{0 0}{0 2.0}{0.5 2.2}{3.9 2.3}{3.4 2.0}{4.7 0.8}
      {4.9 1}{5.2 1}{5.2 0.8}{5.5 0.6}{7.6 0.6}{7.3 1.2}{}
    \picmultigraphics[S]{2}{1 -1}{
     \opencurvepath{0.37 0}{0.42 1.2}{0.49 1.98}{0.87 1.98}{1.04 1.6}
     {1.4 1.6}{2.8 1.6}{2.9 1.7}{3.3 1.8}{3.5 1.6}{4.7 0.6}{2.6 0.6}
     {2.7 0.9}{2.38 0.9}{2.3 0.8}{2.3 0}{}
    }
    \opencurvepath{0 -2.4}{3.9 -2.4}{4.4 -2.7}{5.2 -2.7}{5.6 -2.6}
      {6.3 -2.39}{6.8 -2.0}{7.2 -1.57}{7.6 -1.57}{8.5 -1.57}{}
    {\piclinedash{0.05 0.05}{0.01}
     \picline{2.2 -2.4}{2.2 -2.2}
     \picline{8.2 -1.25}{8.2 -1.6}
     \picline{4.8 0}{5.2 0}
    }
    \picmultigraphics[S]{2}{1 -1}{
     \piclinedash{0.05 0.05}{0.01}
     \picline{6.2 0.6}{6.2 0.3}
     \picline{3.8 0.7}{3.8 0.3}
     \picline{2.5 0.9}{2.5 1.2}
     \picline{2.0 1.6}{2.0 1.4}
     \picline{4 1.2}{4.2 1.4}
    }
   }
   \picputtext[d]{6.5 -0.15}{$b$}
   \picputtext[d]{3.5 -0.15}{$a$}
   \picputtext{2.2 1.81}{$x$}
   \picputtext{8.2 -1.01}{$y$}
   \picputtext{-0.19 -0.1}{$y$}
   \picputtext{2.7 0.45}{$A$}
   \picputtext{1.36 -0.7}{$w$}
 }
}
\end{eqn}
Note that the loop $y$ (of $B(D)$) closes on the left as shown,
because of the neighboring loop to $L$ on the outer cycle (which
we omitted).

Now if we look at the $B$-state of \eqref{mik}, we see that
there is a cycle of length at least $4$ (containing $a,b,x,y$),
which then must be the outer cycle. (This cycle may get longer if
multiplicities of edges in the $A$-state are augmented.)
Now we observe also that this
cycle has two loops, $x$ and $y$, with non-empty interior. This
contradicts lemmas \ref{lA1}, \ref{lC1} and \ref{lB1}, though.

Note that in lemma \reference{lC1} we took the freedom to
apply a mutation. However, the situation in its proof where
this was necessary does not occur here. The valence-2 attachment
$w$ inside $x$ touches a region $A$ not touched by traces from
within $y$.
\qed

So in particular all valencies of attachments \eqref{attach}
are 2 or 3, and \eqref{tatt} has \em{always} the shape \eqref{tat2},
with some of $a$ to $f$ being $0$, except for possible depth-2 loops
(see corollary \reference{crdp}) attached inside the triangle loops.

\begin{lemma}\label{lm6}
We can assume, up to mutation, there is at most one
edge of multiplicity $\ge 3$.
\end{lemma}

\proof An edge $x$ of multiplicity $\ge 3$ gives rise to a cycle
of length $\ge 4$. To see this, let $a$ and $b$ be the loops
in $B(\tl D_k)$ connected by the trace of the simple crossing
$x$. As before, if $a\ne b$, then by primeness of $D_k$,
$a$ and $b$ must be connected outside $x$ (i.e. via traces
not dual to those in $x$) in $B(D_k)$.
If this cycle is of length $3$, then $x=3$, and $a=b$, i.e.
$x$ is inadequate. This contradicts lemma \ref{x2}, though. 

Now we have only one cycle of length at least 4. So
if $x,y\ge 3$, then they must give rise each to a single 
cycle in $B(D_k)$, and both cycles must coincide. Then again 
it is easy to see that by mutations (that preserve the total
loop depth, in order to not to collide with our assumption the
diagram to be flat), one of $x$ or $y$ can be made to $\le 2$.
\qed

\begin{lemma} \label{lmoz}
The outer cycle has length at least 5.
\end{lemma}

\proof 
That the outer cycle has length at least 5 follows directly from
\eqref{ABl}, with our assumption $k>0$, and the observation that
there are exactly 4 loops in $A(D)$ outside the outer cycle.
\qed

\begin{lemma}
There exist no two pairs of parallel loops, and no triple
of parallel loops.
\end{lemma}

\proof The two pairs are out immediately, since we have only
three attachments. The triple is out similarly, since, by the 
lack of further attachments, we would have $\gm=2$. \qed

% I MOVED THIS TO CASE C SINCE IN CASE B IT IS TRIV THAT
% WE HAVE ONLY ONE OUTSIDE ATTACHMENT
% \begin{lemma}
% If two loops are attached using \eqref{attach} outside
% to outer cycle loop, then they are parallel 
% \[
% \]
% or they are not within the same two outside loop edges.
% \[
% \]
% \end{lemma}
% 
% So what this lemma says is, for example, that the following 
% scenario is out:
% \[
% \]
%
% 
% \proof ...

\subsubsection{Going over to initial diagrams}

Let us now, again, temporarily eliminate the dependence of
$D_k$ on $k$ by considering only the simplest diagram that
comes in question. (In \S\reference{Scmp} below we explain
how we recover the whole sequence from that diagram.) Again we will see 

\begin{lemma}\label{lmcmp}
Assume there is at most one non-empty outer cycle loop and
at most one multiple edge on the outer cycle between two
empty outer cycle loops. Then connectivity (number of components)
of $D_k$ and all invariants we test so far, namely those in
\eqref{71}, \eqref{gm4}, and corollary \reference{cr3l},
are not changed if
\begin{enumerate}
\item we change, keeping the parity, an edge of multiplicity 
$>2$ to $4$ or $5$, or
\item change the length of the outer cycle to $5$ or $6$ by 
joining pairs of simple edges between outer cycle loops:
\begin{eqn}\label{2lp}
\diag{7mm}{3.5}{1}{
  {\piclinedash{0.05 0.05}{0.02}
   \picline{0 0.5}{3.5 0.5}
  }
  \picfilledcircle{1.05 0.5}{0.3}{}
  \picfilledcircle{2.45 0.5}{0.3}{}
}\quad\lra\quad
\diag{7mm}{1.5}{1}{
  {\piclinedash{0.05 0.05}{0.02}
   \picline{0 0.5}{1.5 0.5}
  }
}
\,.
\end{eqn}
\end{enumerate}
\end{lemma}

\proof To justify this we need to exclude one pathological 
situation in either points. 

The only troublesome situation in case 1 is when we may 
destroy a triangle in $B(D)$ by changing a leg of multiplicity 
3 to 5. This leg must be then inadequate. (Otherwise all 
cycles it belongs to have length at least 4.) In the proof 
of lemma \ref{lm6}, we argued this possibility out, though. 
% However, an inadequate edge has 
% multiplicity $2$ by lemma \ref{x2}. 

For case 2 we will need to see that by adjusting the 
length of the outer cycle, we do not make a multiple edge 
in $B(D)$ simple, which may affect $IG(B(D))$ (by deleting an 
isolated vertex). %For this purpose we prove lemma \ref{lmoz}.
But since at least two simple edges remain in the outer cycle
of $A(D)$ by assumption, the graph $IG(B(D))$ is not 
affected either. \qed

Observe that the premise of the lemma is needed only for
part 2. So again we can assume
% THIS SHOULD BE MOVED RIGHT AFTER RULING OUT EDGE MULT 3
% AND TWO EDGES OF MULT 4,5
\[
\mbox{\em all edges have multiplicities 1,2,4,5, and at 
most one has 4 or 5}.
\]
Let us call the edge of multiplicity $4$ or $5$ the \em{twist edge}.

With this assumption we have

\begin{lemma}\label{l2v}
The leg multiplicities of valence-2-loops are
$(1,1)$, $(1,2)$, $(1,4)$, $(1,5)$, and $(2,5)$.
\end{lemma}

\proof Similar to that of lemma \reference{lemmaW}. However, now
$(1,2)$ and $(1,5)$ come in, because we have a triangle resp. may
have an outer cycle of even length in $B(D)$. \qed 

It is easy to prove that the adjustment (also in point 2) of
lemma \ref{lmcmp} can be done in case B. Case C is more awkward 
and, after lemma \reference{lC1}, the completion will be postponed
to lemma \ref{lC2}. Case B' is dealt with in lemma \reference{ghp'},
and corollary \ref{cor10.4} (a direct consequence of lemma \ref{ghp}).

\begin{lemma}
In case B there is at most one outer cycle loop with 
non-empty interior, and there are no multiple edges between 
empty outer cycle loops.
\end{lemma}

\proof % It remains to estimate the number of simple edges.
That there is at most one non-empty outer cycle loop was
shown in lemma \reference{lB1}. It is easy to see that with
3 attachments we can create at most 2 cycles in $IG(A(D))$:
by the argument in the proof of lemma \reference{lB1}, there
cannot be two loops inside $L$ intertwined (each) with 4
connections outside $L$ in case of an outside attachment,
and intertwined with 3 connections outside $L$ in case of
3 attachments inside $L$. (Note again that \eqref{pecu}
is excluded, because of lemma \reference{lmZ'B}.)
% In case B, it is easy to see that with 3 attachments we can
% create at most 2 cycles in $IG(A(D))$. Since we need these 2 
% cycles, we see that all 3 attachments must be made at the same 
% loop of the outer cycle or triangle. So there is at most one 
% outer cycle loop with non-empty interior.
Therefore, there
cannot be a multiple edge between the other at least 4 outer 
cycle loops, for it would crate an isolated vertex in $IG(A)$.
\qed

\begin{lemma}
In case A there is at most one outer cycle loop with 
non-empty interior, and there is at most one multiple edge between 
empty outer cycle loops.
\end{lemma}

\proof We argued that case (A1) in \eqref{csA} is excluded. In (A2)
we need both attachments to be inside $L$, and with $L$ having
at most 4 multiple outside connections, we have at most 3 cycles
in $IG(A)$. So we can have at most one multiple edge between empty
outer cycle loops (giving an isolated vertex in $IG(A)$). \qed

\subsubsection{Computations\label{Scmp}}

It will be again inavoidable at some points of the proof 
to use a computer and check a certain number of explicit 
patterns. (Or at least, their manual case-by-case discussion
would lead to no reasonable argument and exposition.) While 
it is not helpful to enter into full (implementational) 
details, we clarify some features of this procedure in advance,
and concentrate on the mathematical part of the work later.

Since we have two cycles now, in turned in practice (due to the
patterns to check) too technically cumbersome to try to make 
the diagrams positive (by making both cycles to be of even length).
In other words, the convenience we would have from using pre-existing 
software to calculate $V$ would not longer compensate for the effort 
in adapting the input to and output from this computation. So it was 
better to write an entirely new program which tests all Jones and
Kauffman semiadequacy invariants directly from the $A$-state. 

The state was encoded by recording the basepoints of the traces 
in cyclic order along each loop. Hereby we chose the orientation
of loops to be \em{coherent}, i.e. so that each trace looks
locally like\ \ $
\diag{5mm}{1}{1.3}{
  \picvecline{0 0}{0 1.3}
  \picvecline{1 1.3}{1 0}
  {\piclinedash{0.09 0.09}{0.01}
   \picline{0 0.6}{1 0.6}
  }
}$~. In contrast to the positive orientation in the proof
of lemma \ref{lmps}, the coherent orientation always exists, but 
it does not necessarily extend to a link diagram orientation.
This explains the inconvenience we overcome when abandoning
the calculation from (oriented) link diagrams.

The coherent orientation also allows one to switch easily
between the $A$- and $B$-state (and test latter's invariants
too). We also programmed a component number test.

After we implemented these checks, we ran as a test the patterns (a)
and (b) of figure \reference{Fig19} in the previous section. We
obtained the 12 vectors (a)-(n) of \eqref{css} within about a
minute. The previous calculation of $V$ took about half an hour.
This speed-up should be related also to the feature of the 
Jones and Kauffman semiadequacy invariants to be of 
polynomial complexity, in contrast to the full polynomial 
invariants they are derived from, which are NP hard.

To make things even more convenient, we implemented also
the recovery of the sequence of $D_k$ from the pattern.
We used as input patterns in which the outer cycle
has length 5. First, if the total number of crossings
is odd, then we extended the outer cycle to length 6
(by putting an extra loop connected by two simple edges
to its neighbors).

If all legs have multiplicity at most 2, then we have a
sporadic case, in which we can only augment the length 
of the outer cycle by multiples of 2 (by adding two loops
connected by only simple edges, in the opposite direction
to \eqref{2lp}). Now we used that $A$- and
$B$-state must have the same number of loops \eqref{ABl}.
Thus there is at most one diagram we can obtain (``at most'' 
because the $B$-state may have already fewer loops than the 
$A$-state before we even extend the outer cycle in $A(D)$). 
We designed our program to report this diagram, but in 
practice it turned out that this case never occurred. 

If we have an edge of multiplicity 4 or 5, the twist edge, then
first we adjust it (within its parity) so that $A$- and $B$-state
have the same number of loops. (In practice it turned out that
in most cases we had to change here multiplicity 5 to 3.)

Then the sequence of diagrams $\hD_k$ is obtained by adding
2 to the multiplicity of that edge and simultaneously extending
the outer cycle by 2. Since this augments the crossing number
of the diagram by 4, we can discard diagrams which
(after equating the number of $A$- and $B$-state loops) are of 
crossing number not divisible by 4. Otherwise, we designed our 
program to output the DT code of the first four diagrams of each 
such series (which were then subjected to Vassiliev invariant
tests).

We have a total of 6 places where we use a computation.
These are designated so. The programs used are available
at \cite{program}.

We collect at the end of the proof all sequences coming out of
our computations. We refer at this place already to figure
\ref{Fig20'} on p.~\pageref{Fig20'}, where the diagram of each
sequence is shown for $k=1$. See \S\ref{SSI} for more explanation. 
We will in fact encounter some sequences repeatedly. For the
calculations it is too cumbersome and not worthwhile to implement
all conditions that lead to the particular case where the calculation
is performed, and that make this case to be mutually exclusive from
the others. (In fact, extraneous checks lend more safety to the
calculation.) Even during the calculation in the same case there are
various symmetries which are too hard to see and remove in advance.

\subsection{Restrictions on outside attachments\label{S10.5}}

In this subsection, before we start a detailed separate
discussion of cases A, B and C, we derive restrictions on how 
loops can be attached (using \eqref{attach}) to loops on 
the outer cycle (or in case B the extra triangle loop).
Case A is trivial from this point of view, so we exclude it
from our treatment.

\begin{lemma}\label{lm2l}
If two loops of valence 2 are attached using \eqref{attach} 
outside to outer cycle loop $L$, then \em{either} they are parallel 
\em{or} they are not within the same two outside loop edges.
\end{lemma}

Let us illustrate the two options by pictures. (We discussed, though,
the situation already partially in the proof of lemma \ref{lmZ'}.)
Parallel loops attached outside to $L$ are
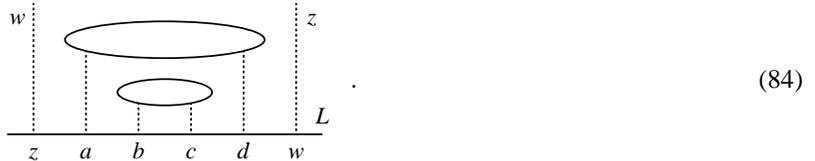
\begin{eqn}\label{pl}
{\small
\diag{7mm}{6}{3}{
  \picline{0 0.5}{6 0.5}
  {\piclinedash{0.05 0.05}{0.01}
   \picline{0.5 3}{0.5 0.5}
   \picline{5.5 3}{5.5 0.5}
   \picline{1.5 2.3}{1.5 0.5}
   \picline{4.5 2.3}{4.5 0.5}
   \picline{2.5 1.3}{2.5 0.5}
   \picline{3.5 1.3}{3.5 0.5}
  }
  \picfilledellipse{3 2.3}{1.9 0.35}{}
  \picfilledellipse{3 1.3}{0.9 0.25}{}
   \picputtext[d]{0.5 0.05}{$z$}
   \picputtext[d]{5.5 0.05}{$w$}
   \picputtext[d]{1.5 0.05}{$a$}
   \picputtext[d]{2.5 0.05}{$b$}
   \picputtext[d]{3.5 0.05}{$c$}
   \picputtext[d]{4.5 0.05}{$d$}
   \picputtext[d]{6 0.7}{$L$}
   \picputtext{5.8 2.7}{$z$}
   \picputtext{0.2 2.7}{$w$}
}}\quad.
\end{eqn}
Here the outside of $L$ is above in the diagram, and 
the edges $w,z$ are understood to connect $L$ to some
other outer cycle loop. Also (unlike below) we understand
that no further loops are attached on the part of $L$ drawn
as horizontal line, \em{except} between $b$ and $c$.

Loops not within the same two outside loop edges look like:
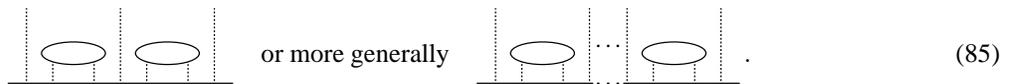
\begin{eqn}\label{nenc}
\diag{5mm}{6}{2.5}{
  \picline{0 0.5}{6 0.5}
  {\piclinedash{0.05 0.05}{0.01}
   \picline{0.5 2.5}{0.5 0.5}
   \picline{5.5 2.5}{5.5 0.5}
   \picline{3.0 2.5}{3.0 0.5}
  }
   \picmultigraphics{2}{2.5 0}{
     \picstroke
     {\piclinedash{0.05 0.05}{0.01}
      \picline{1.2 1.3}{1.2 0.5}
      \picline{2.3 1.3}{2.3 0.5}
     }
     \picfilledellipse{1.75 1.3}{0.85 0.3}{}
   }
}\quad\mbox{or more generally}\quad
\diag{5mm}{7}{2.5}{
  \picline{0 0.5}{3.1 0.5}
  \picline{7 0.5}{3.9 0.5}
  \picputtext{3.5 0.5}{\small .\,.\,.}
  \picputtext{3.5 1.5}{\small .\,.\,.}
  {\piclinedash{0.05 0.05}{0.01}
   \picline{0.5 2.5}{0.5 0.5}
   \picline{6.5 2.5}{6.5 0.5}
   \picline{3.0 2.5}{3.0 0.5}
   \picline{4.0 2.5}{4.0 0.5}
  }
   \picmultigraphics{2}{3.5 0}{
     \picstroke
     {\piclinedash{0.05 0.05}{0.01}
      \picline{1.2 1.3}{1.2 0.5}
      \picline{2.3 1.3}{2.3 0.5}
     }
     \picfilledellipse{1.75 1.3}{0.85 0.3}{}
   }
}\,.
\end{eqn}
(Again the outside is above $L$, but something may now be
attached from below to $L$.)

So what this lemma says is, for example, that the following 
scenarios are out:
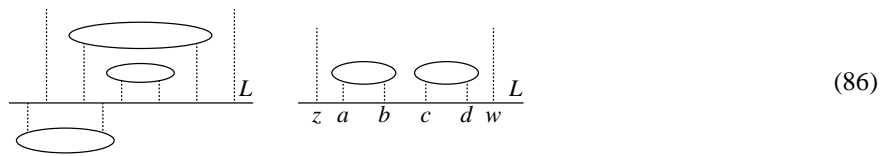
\begin{eqn}\label{op2}
{\small
\diag{5mm}{6.5}{4}{
  \pictranslate{0.5 1}{
    \picline{-0.5 0.5}{6 0.5}
    {\piclinedash{0.05 0.05}{0.01}
     \picline{0.5 3}{0.5 0.5}
     \picline{5.5 3}{5.5 0.5}
     \picline{1.5 2.3}{1.5 0.5}
     \picline{4.5 2.3}{4.5 0.5}
     \picline{2.5 1.3}{2.5 0.5}
     \picline{3.5 1.3}{3.5 0.5}
    }
    \picfilledellipse{3 2.3}{1.9 0.35}{}
    \picfilledellipse{3 1.3}{0.9 0.25}{}
    {\piclinedash{0.05 0.05}{0.01}
     \picline{2.0 0.5}{2.0 -0.5}
     \picline{-0.0 0.5}{-0.0 -0.5}
    }
    \picfilledellipse{1 -0.5}{1.3 0.33}{}
    \picputtext{5.8 0.9}{$L$}
  }
}\quad
\diag{5mm}{6}{4}{
  \pictranslate{0.5 1}{
  \picline{0 0.5}{6 0.5}
  {\piclinedash{0.05 0.05}{0.01}
   \picline{0.5 2.5}{0.5 0.5}
   \picline{5.2 2.5}{5.2 0.5}
  }
   \picmultigraphics{2}{2.2 0}{
     \picstroke
     {\piclinedash{0.05 0.05}{0.01}
      \picline{1.2 1.3}{1.2 0.5}
      \picline{2.3 1.3}{2.3 0.5}
     }
     \picfilledellipse{1.75 1.3}{0.85 0.3}{}
   }
   \picputtext{5.8 0.9}{$L$}
   \picputtext[d]{0.5 0.0}{$z$}
   \picputtext[d]{5.2 0.0}{$w$}
   \picputtext[d]{1.2 0.0}{$a$}
   \picputtext[d]{2.3 0.0}{$b$}
   \picputtext[d]{3.4 0.0}{$c$}
   \picputtext[d]{4.5 0.0}{$d$}
  }
}
}
\end{eqn}
\proof If we have two loops attached outside, in order 
$IG(A)$ to have a cycle, we must have case C, and the triangle 
\eqref{tatt} is the only attachment inside $L$.

Assume the two loops $M$, $N$ are attached within the same arc of $L$
separated by basepoints of edges to other outer cycle loops.

\begin{caselist}
\case\label{CS1}
$M$ and $N$ are enclosed. This is the situation \eqref{pl},
but allowing for attachments from below on the horizontal line.
Under exclusion of inadequate loops, $b$ cannot be alone on an
interval between trace basepoints from inside $L$. So it must be
grouped with $a$ or $c$ on the same interval. In latter case $D$
is composite, so $b$ is grouped with $a$. Similarly $c$ must be
grouped with $d$, and then $M$ and $N$ are parallel.

\case $M$ and $N$ are not enclosed. This is the situation in the
right picture of \eqref{op2}. The argument in the previous case
modifies to show that $b$ is grouped with $c$, $d$ with $w$, and
$a$ with $z$. So we have the picture (a) in \eqref{44}.
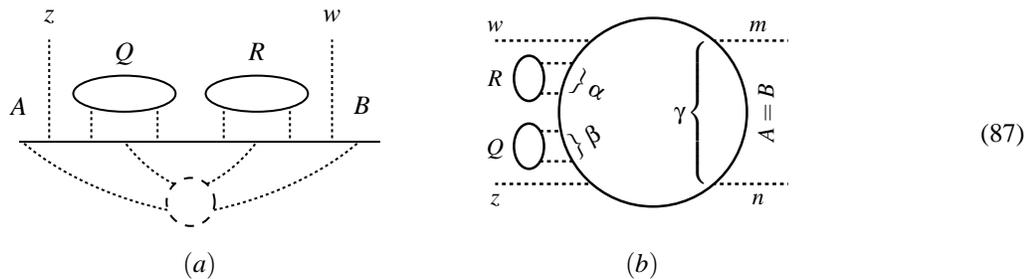
\begin{eqn}\label{44}
\begin{array}{c@{\kern1.5cm}c}
\diag{8mm}{6}{4}{
  \pictranslate{0 1}{
  \picline{0 0.5}{6 0.5}
  {\piclinedash{0.05 0.05}{0.01}
   \picline{0.5 2.2}{0.5 0.5}
   \picline{5.2 2.2}{5.2 0.5}
  }
   \picmultigraphics{2}{2.2 0}{
     \picstroke
     {\piclinedash{0.05 0.05}{0.01}
      \picline{1.2 1.3}{1.2 0.5}
      \picline{2.3 1.3}{2.3 0.5}
     }
     \picfilledellipse{1.75 1.3}{0.85 0.3}{}
   }
   \picputtext{1.75 2.0}{$Q$}
   \picputtext{3.95 2.0}{$R$}
  }
  {\piclinedash{0.05 0.05}{0.01}
  \pictranslate{2.85 0.7}{
    \picmultigraphics[S]{2}{-1 1}{
       \piccurve{0 0}{0.4 0.1}{0.8 0.4}{1.1 0.8}
       \piccurve{0 -0.4}{0.8 -0.3}{2.0 0.1}{2.8 0.8}
    }
    \piclinedash{0.15 0.15}{0.01}
    \picfilledcircle{0 -0.2}{0.4}{}
  }
  }
  \picputtext{0 2.1}{$A$}
  \picputtext{5.7 2.1}{$B$}
  \picputtext{0.5 3.6}{$z$}
  \picputtext{5.2 3.6}{$w$}
} &
\small
\diag{1cm}{3.9}{2.5}{
  {\piclinedash{0.05 0.05}{0.01}
   \picline{0 0.3}{3.9 0.3}
   \picline{0 2.2}{3.9 2.2}
  }
  \picmultigraphics{2}{0 0.9}{
    {\piclinedash{0.05 0.05}{0.01}
     \picline{0.5 0.6}{1.5 0.6}
     \picline{0.5 1.0}{1.5 1.0}
    }
    \picfilledellipse{0.45 0.8}{0.2 0.3}{}
  }
  \picfilledcircle{2.1 1.25}{1.25}{}
  \picputtext{0 0.1}{$z$}
  \picputtext{0 2.4}{$w$}
  \picputtext{3.5 0.1}{$n$}
  \picputtext{3.5 2.4}{$m$}
  \pictranslate{3.6 1.25}{
    \picrotate{90}{\picputtext{0 0}{$A=B$}}
  }
  \picputtext{0.0 0.8}{$Q$}
  \picputtext{0.0 1.7}{$R$}
  \pictranslate{2.63 1.7}{
    \picputtext{0 0}{$\gm\left\{\ry{3.4em}\right.$} %}
  }
  \pictranslate{1.25 1.65}{
    \picrotate{-30}{ 
      \picputtext{0 0}{$\left.%{
	\ry{0.9em}\right\}\,\ap$%
      }
    }
  }
  \pictranslate{1.15 0.92}{
    \picrotate{28}{ 
      \picputtext{0 0}{$\left.%{
	\ry{0.9em}\right\}\,\bt$%
      }
    }
  }
} \\
\ry{1.4em}(a) & (b)
\end{array}
\end{eqn}
The dashed circle in the lower part means that we do not specify
how traces are connected to the loops $M,N$ inside the circle.
A trace that touches regions $A$ or $B$ must be there
by primeness of $D$. However, since $\gm=4$, we see that
actually then $A=B$. So outside $L$ we have the picture (b)
in \eqref{44}, with $\ap$, $\bt$, $\gm$ being the segments
of $L$ on which legs of $M$ and $N$ are attached from inside $L$.

We assume that $z,w$ resp. $m,n$ connect $L$ to its two neighboring
loops on the outer cycle $O$ resp. $P$. We must allow that some
of $z,w,m,n$ be 0. This covers all possibilities how the edges
outside $L$ connect to $O,P$, with one exception. This exception
is the option that $A$ lies between edges that connect $L$ to
the same outer cycle loop $O$ as $z$ and $w$. This situation is,
however, recurred to the other cases by a mutation, as explained
after \eqref{exss}.

\begin{caselist}
\case Assume first that some of $z,w,m,n$ is $0$. Up to mutations,
we may assume w.l.o.g. that $n=0$. Then the connections of $L$ to
the two inside loops $M$ and $N$ are intertwined with at most one
outer cycle connection from $L$. To have two cycles in $IG(A)$, this
connection must be then indeed intertwined with both $M$ and $N$,
i.e. $z\ne0\ne w$. So we have the following picture, up to symmetry:
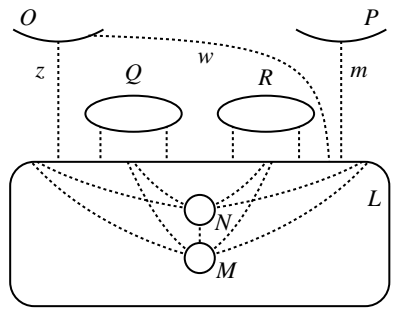
\[
{\small
\diag{8mm}{6}{6}{
  \pictranslate{0 2.7}{
  \picovalbox{2.85 -0.7}{6.3 2.4}{0.4}{}
  {\piclinedash{0.05 0.05}{0.01}
   \picline{0.5 2.5}{0.5 0.5}
   \picline{5.2 2.5}{5.2 0.5}
   \piccurve{1.0 2.6}{4.4 2.5}{4.9 2.2}{5.0 0.5}
  }
  \piccirclearc{0.5 4}{1.5}{-120 -60}
  \piccirclearc{5.2 4}{1.5}{-120 -60}
   \picmultigraphics{2}{2.2 0}{
     \picstroke
     {\piclinedash{0.05 0.05}{0.01}
      \picline{1.2 1.3}{1.2 0.5}
      \picline{2.3 1.3}{2.3 0.5}
     }
     \picfilledellipse{1.75 1.3}{0.8 0.3}{}
   }
   \picputtext{0 2.9}{$O$}
   \picputtext{5.7 2.9}{$P$}
  }
  \pictranslate{2.85 2.4}{
   {\piclinedash{0.05 0.05}{0.01}
    \picmultigraphics[S]{2}{-1 1}{
       \piccurve{0 0}{0.4 0.1}{0.8 0.4}{1.1 0.8}
       \piccurve{0 0}{0.8 0.1}{2.0 0.4}{2.7 0.8}
       \piccurve{0 -0.8}{0.8 -0.6}{2.0 -0.1}{2.8 0.8}
       \piccurve{0 -0.8}{0.4 -0.6}{0.8 -0.1}{1.2 0.8}
    }
    \picline{0 0}{0 -0.8}
   }
   \picfilledcircle{0 -0.0}{0.25}{}
   \picfilledcircle{0 -0.8}{0.25}{}
   \picputtext{0.4 -0.2}{$N$}
   \picputtext{0.45 -1.0}{$M$}
   \picputtext{-1.1 2.3}{$Q$}
   \picputtext{ 1.1 2.2}{$R$}
   \picputtext{ 2.9 0.2}{$L$}
   \picputtext{ 0.1 2.52}{$w$}
   \picputtext{-2.65 2.3}{$z$}
   \picputtext{ 2.65 2.3}{$m$}
  } 
}
}\unskip
\]
The understanding of legs whose traces cross is that for each such
intersecting pair at least one leg must be of multiplicity 0. All
edges outside $L$ are of non-zero multiplicity.

Now the connection between $L$ and $P$ is not intertwined with
any of the (connections to the) loops inside $L$. In order to have
then two cycles in $IG(A)$, we need that $M$ and $N$ are intertwined
with all of $Q$, $R$ and $O$. However, then legs from $M$ and $N$
must be installed inside $L$ on all intervals $\ap$, $\bt$, $\gm$
of picture (b) in \eqref{44}. That is easily seen to be impossible
without leg intersections.

\case So we have the case that in \eqref{44} (b) all of
$z,w,m,n$ are non-zero. % THERE WAS AN ERROR HERE FIXED

Since $\chi(IG(A))=-1$, we see that $M$ and $N$ are intertwined
with at least 3 connections outside $L$. Thus they must both have
legs on $\gm$ of picture (b) in \eqref{44}. A simple check shows
then that we have (up to symmetry) only two options left. These
cases are shown below, assuming all edges drawn have non-zero
multiplicity. In case (b) we have drawn simultaneously the $B$-state
(following convention \reference{cthi}).
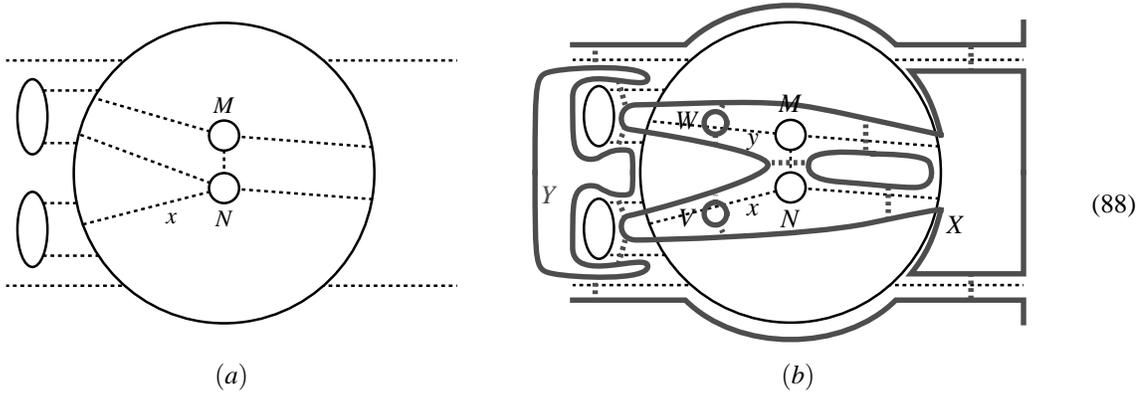
\begin{eqn}\label{45}
\begin{array}{c@{\kern1.5cm}c}
\small
\diag{1cm}{6}{4}{
  {\piclinedash{0.05 0.05}{0.01}
   \picline{0 0.5}{6 0.5}
   \picline{0 3.5}{6 3.5}
  }
   \pictranslate{0 1.25}{
  \picmultigraphics{2}{0 1.5}{
    {\piclinedash{0.05 0.05}{0.01}
     \picline{0.5 -0.35}{1.5 -0.35}
     \picline{0.5 0.35}{1.5 0.35}
    }
    \picfilledellipse{0.35 0.0}{0.2 0.5}{}
  }
  }
  \pictranslate{2.9 2}{
  \picfilledcircle{0 0}{2}{}
  {\piclinedash{0.05 0.05}{0.01}
   \picline{2 10 polar}{0 0.5}
   \picline{2 -10 polar}{0 -0.2}
   \picline{0 0.5}{0 -0.2}
   \picline{0 0.5}{2 150 polar}
   \picline{0 -0.2}{2 165 polar}
   \picline{0 -0.2}{2 200 polar}
  }
  \picfilledcircle{0 0.5}{0.2}{}
  \picfilledcircle{0 -0.2}{0.2}{}
  \picputtext[d]{0 0.8}{$M$}
  \picputtext[u]{0 -0.5}{$N$}
  \picputtext{-0.7 -0.6}{$x$}
  }
} &
\diag{1cm}{6}{4}{
  {\piclinedash{0.05 0.05}{0.01}
   \picline{0 0.5}{6 0.5}
   \picline{0 3.5}{6 3.5}
  }
   \pictranslate{0 1.25}{
  \picmultigraphics{2}{0 1.5}{
    {\piclinedash{0.05 0.05}{0.01}
     \picline{0.5 -0.35}{1.5 -0.35}
     \picline{0.5 0.35}{1.5 0.35}
    }
    \picfilledellipse{0.35 0.0}{0.2 0.4}{}
  }
  }
  \pictranslate{2.9 2}{
    \picfilledcircle{0 0}{2}{}
    {\piclinedash{0.05 0.05}{0.01}
     \picline{2 10 polar}{0 0.5}
     \picline{2 -10 polar}{0 -0.2}
     \picline{0 0.5}{0 -0.2}
     \picline{0 0.5}{2 160 polar}
     \picline{0 -0.2}{2 200 polar}
    }
    \picfilledcircle{0 0.5}{0.2}{}
    \picfilledcircle{0 -0.2}{0.2}{}
  }
  \piclinewidth{14}
  \picgraycol{0.3}
  \picPSgraphics{2 setlinecap}
  \pictranslate{0 2}{
   \picputtext{-0.3 -0.3}{$Y$}
   \picmultigraphics[S]{2}{1 -1}{
    \picstroke{
      \picline{6 2}{6 1.7}
      \piclineto{4.3 1.7}
      \piccurveto{3.6 2.4}{2.2 2.4}{1.5 1.7}
      \piclineto{0.0 1.7}
      \picline{6 0}{6 1.35}
      \piclineto{4.5 1.35}
      \piccurveto{4.7 1.1}{4.8 0.9}{4.9 0.5}%{5 0}
      \opencurvepath{-0.5 0}{-0.5 1.3}{0.0 1.4}{1 1.45}{1 1.05}{0.0 1.5}{0.0 0.5}
        {0.0 0.35}{0.15 0.15}{0.3 0.17}{0.4 0.23}{0.5 0.3}{0.8 0.3}
	{0.8 0}{}
    }
    {\piclinedash{0.05 0.05}{0.01}
     \picPSgraphics{0 setlinecap}
     \picline{0.7 0.6}{0.6 0.35}
     \picline{0.7 0.9}{0.6 1.25}
     \picline{0.3 1.38}{0.3 1.7}
     \picline{5.3 1.35}{5.3 1.7}
    }
   }
   \picstroke{
    \opencurvepath{4.9 0.5}{3.1 0.9}{2.1 0.95}{1.1 0.9}
    {0.6 0.9}{0.6 0.6}{0.9 0.55}{2.6 0.3}{2.6 -0.1}{0.9 -0.55}
    {0.6 -0.6}{0.6 -0.9}{1.1 -0.9}{3.1 -0.85}{4.8 -0.5}{4.9 -0.5}{}
    % \picPSgraphics{cp}
    \curvepath{4.0 0.25}{3.9 0.26}{3.7 0.27}{3.2 0.3}{3.1 0}
      {3.2 -0.1}{3.9 -0.15}{4.3 -0.2}{4.8 -0.2}{4.8 -0.0}{4.7 0.25}{}
   }
   \piccircle{1.9 -0.55}{0.15}{}
   \piccircle{1.9 0.66}{0.15}{}
   {\piclinedash{0.05 0.05}{0.01}
    \picPSgraphics{0 setlinecap}
    \picline{2.65 0.13}{3.15 0.13}
    \picline{2.0 0.9}{1.95 0.8}
    \picline{1.9 0.5}{1.9 0.4}
    \picline{2.0 -0.8}{1.95 -0.7}
    \picline{1.9 -0.4}{1.9 -0.3}
    \picline{3.9 0.3}{3.9 0.7}
    \picline{4.2 -.2}{4.2 -0.65}
   }
  }
  \picgraycol{0.0}
  \pictranslate{2.9 2}{
    \picputtext[d]{0 0.8}{$M$}
    \picputtext[u]{0 -0.5}{$N$}
    \picputtext{-0.5 -0.5}{$x$}
    \picputtext{-0.5 0.45}{$y$}
    \picputtext{2.2 -0.7}{$X$}
    \picputtext{-1.36 0.7}{$W$}
    \picputtext{-1.36 -0.6}{$V$}
  }
} \\
\ry{2.2em}(a) & (b)
\end{array}
\end{eqn}
Case (a) can be handled as in lemma \ref{lm1014}.
The edge $x$ is inadequate, so $x=2$. Then by connectivity,
one of the other legs of $M$ or $N$ must be multiple. Then one
exhibits in $B(D)$ two cycles not separated by a loop. So
$B(D)$ is of type $A$ or $B$, and we are done.

In case (b), we have $x=y=2$ by inadequacy. We have drawn above
again into the diagram the loops and traces in $B(D)$ by
thicker, and gray, lines. The property $x=y=2$ (in $A(D)$) explains
the small loops drawn near them in $B(D)$. Let us call these
small loops $V$ and $W$. The $B$-state loops $X$ and $Y$
close as indicated, because they touch a segment of the
$A$-state loops neighbored to $L$ on the outer cycle.
With that it is easy to see that $X$ is the only separating
loop in $B(D)$. Now $X$ has at most 3 multiple connections inside.
(When augmenting multiplicity of some of the edges in $A(D)$, further
small loops enter into $B(D)$, and multiple connections and their
intertwining may only be spoiled, but never newly
created.) Two of these connections, those from $V$ and $W$,
are intertwined with at most one connection outside $X$, the one to
$Y$. Then we see that $IG(B)$ cannot have any cycle.
\end{caselist}
\end{caselist}

With this lemma \reference{lm2l} is proved.
\qed
%til here anything was cleaned up

The following lemma makes a series of technical assumptions.
However, these are often satisfied, and the argument will
be needed repeatedly, so it is better to single it out.

\begin{lemma}\label{UV}
Consider a diagram $D$ in case B, B' or C with the following
properties:
\begin{itemize}
\item 
Assume $B(D)$ has a single separating loop $X$, and
$\le 3$ multiple connections inside (resp. outside), each 
intertwined with at most $2$ multiple connections outside
(resp. inside) $X$.

\item
Let $U,V$ be the regions of $A(D)$, 
%for the diagram $\tD$ in \eqref{ABC}
which contain traces of all but at most one (in case B, B')
connection of the outer cycle in their boundary.
(For case $C$ it is as in \eqref{oaz}.) Assume further no
attachment legs (of both \eqref{tatt} and \eqref{attach})
touch upon the boundary of $U$ or $V$.

\item $B(D')$ has at least 7 loops, where $D'$ is as in definition
\reference{cfiq}. (See also remark \ref{remUV} after the lemma.)

\item $IG(A(D))$ has at most $4$ cycles, and at most one non-empty
outer cycle loop.
\end{itemize}
Then such a diagram $D$ cannot occur among our
diagrams $D_k$.
\end{lemma}

\begin{rem}\label{remUV}
Let us make the following remarks concerning the 7-loop
condition on $B(D')$. (These points will be used implicitly when
invoking the lemma several times below.) More precisely, as the
proof will show, it is
enough that $B(D)$ has at least 7 loops, when we assume the
minimal admissible multiplicity for each edge of $A(D)$, which
is not an edge between empty outer cycle loops.

Moreover, to find such 7 loops, it is enough to have one loop in
$B(D)$ connected by a trace to $X$, whose connection to $X$ is not
intertwined with any connection on the opposite side of $X$.
This is because we already have $X$, and at least 5 other loops
are needed whose connections to create 2 cycles in $IG(B)$.
\end{rem}

\proof Now, by the second assumption, the boundaries of $U,V$
correspond to loops in $B(D)$, which we denote again by $U,V$.
The connection $(U,V)$ between loops $U$ and $V$ in $B(D)$ 
is not intertwined with any other connection. Since, by the
first condition, we can create at most two cycles in $IG(B)$, 
this means that $(U,V)$ must be simple (else it would yield an 
isolated vertex in $IG(B)$). 

Thus there exists at most one simple edge between empty outer
cycle loops in $A(D)$. On the other hand, $IG(A(D))$ has at most
$4$ cycles, so there are at most two multiple edges between
empty outer cycle loops in $A(D)$. Since the outer cycle has
length at least 5, we see that it must have length exactly 5,
and exactly two of the three edges between empty outer cycle
loops in $A(D)$ must be multiple.

With an outer cycle in $A(D)$ of length 5, we have in $A(D)$ exactly
9 loops. Thus $B(D)$ has also 9 loops, and $c(D)=20$ by \eqref{ABl}.
We are thus already down to a finite number of cases. Even
those are easily dealt with by hand, assuming $B(D)$ has at
least 7 loops, when all outer cycle edges between empty
loops in $A(D)$ are simple.

Now we are supposed to have two multiple outer cycle edges between
empty loops in $A(D)$. By connectivity reasons their multiplicities
cannot be $(2,2)$. So making outer cycle edges multiple, we must 
add at least 3 traces in $A(D)$, and so at least 3 loops to $B(D)$. 
Then, however, $B(D)$ has at least 10 loops, a contradiction.

The completes the proof of lemma \reference{UV}. \qed

\begin{lemma}\label{lm3v}
In cases B or C, there is no loop $E$ of valence 3 attached (using 
\eqref{attach}) outside to an outer cycle loop $L$.
\end{lemma}

\proof Assume there is such $E$. We will rule
gradually out all possibilities.

First observe that there must be three regions outside $L$
touched by legs of loops attached inside $L$. Two of these
regions are enclosed by $L$ and $E$. (The third one is needed 
to avoid $D$ becoming composite.)

Case $B$ is easy. We already know that there is only one loop
on the outer cycle or triangle with non-empty interior. The 
only option we have, avoiding inadequate loops, is
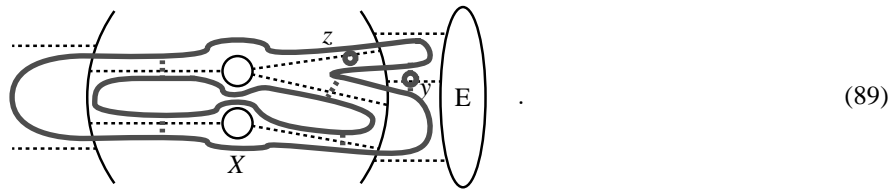
\begin{eqn}\label{..1}
\diag{1cm}{6}{3}{
  \pictranslate{3 1.5}{
    \picmultigraphics[S]{2}{-1 1}{
      \piccirclearc{0 0}{2}{-35 35}
    }
    {\piclinedash{0.05 0.05}{0.01}
     \picline{2 170 polar}{2 170 polar x p 0 x}
     \picline{2 190 polar}{2 190 polar x p 0 x}
     \picline{2 170 polar x p 0 x}{2 18 polar}
     \picline{2 170 polar x p 0 x}{2 -3 polar}
     \picline{2 190 polar x p 0 x}{2 -20 polar}
    }
    \picfilledcircle{2 170 polar x p 0 x}{0.2}{}
    \picfilledcircle{2 190 polar x p 0 x}{0.2}{}
    {\piclinedash{0.05 0.05}{0.01}
     \picline{2 160 polar}{2 160 polar x p -3 x}
     \picline{2 200 polar}{2 200 polar x p -3 x}
     \picline{2 25 polar}{2 25 polar x p 3 x}
     \picline{2 -25 polar}{2 -25 polar x p 3 x}
     \picline{2 6 polar}{2 6 polar x p 3 x}
    }
    \picfilledellipse{3 0}{0.3 1.2}{E}
    \picputtext{1.2  0.8}{$z$}
    \picputtext{2.5 0.1}{$y$}
    \picputtext{0.0 -0.9}{$X$}
    { \piclinewidth{14}
      \picgraycol{0.3}
      \picPSgraphics{2 setlinecap}
      \picstroke{
	\curvepath{2.4 0.75}{2.3 0.75}{2.0 0.71}{0.4 0.5}{0.4 0.8}
	  {-0.4 0.8}{-0.4 0.5}{-1.8 0.55}{-2.3 0.55}{-3.0 0.5}
	  {-3.0 -0.5}{-2.3 -0.55}{-1.8 -0.55}{-0.4 -0.5}
	  {-0.4 -0.7}{0.4 -0.7}{0.4 -0.6}{1.3 -0.71}{2.3 -0.77}
	  {2.5 -0.71}{2.6 -0.3}{2.5 -0.0}{2.0 0.1}{1.5 0.25}{1.2 0.29}
	  {1.2 0.34}{2.0 0.4}
	  {2.5 0.4}{2.6 0.6}{}
	\curvepath{-0.6 0.25}{-0.7 0.25}{-1.9 0.25}{-1.9 0.1}
	  {-1.9  -0.28}{-0.3 -0.28}{-0.25 -.04}{0.17 -0.04}
	  {0.3 -0.22}{0.4 -0.36}{1.5 -0.5}{1.8 -0.5}{1.8 -0.1}
	  {0.6 0.08}{0.3 0.15}{0.2 0.09}
	  {0.0 0.02}{-0.2 .02}{-0.2 0.18}{}{}
      }
      \piccircle{1.5 0.52}{0.07}{}
      \piccircle{2.3 0.24}{0.08}{}
      {\piclinedash{0.05 0.05}{0.01}
       \picPSgraphics{0 setlinecap}
       \picline{-1.0 0.25}{-1.0 0.55}
       \picline{-1.0 -0.23}{-1.0 -0.6}
       \picline{1.2 0.0}{1.4 0.33}
       \picline{1.4 -0.5}{1.4 -0.73}
       \picline{2.3 0.08}{2.3 0.16}
       \picline{2.3 0.32}{2.3 0.49}
      }
    }
  }
}
\qquad.
\end{eqn}
Here all (thin, black) edges in $A(D)$ are assumed to be of
non-zero multiplicity.
%...[1]
Again $y$ and $z$ are inadequate, so $y=z=2$. We look at $B(D)$. 
Let $Z$ be the loop between the two traces of $z$. We see 
that in $B(D)$ there is only one separating loop, call it $X$. 
It has at most two multiple connections inside, and one of them, 
to $Z$, is intertwined with at most two multiple connections 
outside $X$. Then $IG(B)$ has at most one cycle.
%end...[1]

Consider case C. If we have a second outer cycle loop with 
non-empty interior, then in this interior must be the third 
attachment, using \eqref{attach}. However, it is easy to see that 
if this attachment has valence 2, then it is inadequate, while 
if it has valence 3, then $\gm>6$. So we can assume there is at 
most one outer cycle loop with non-empty interior.
%
%...[2]
%
Let below $L$ be this non-empty outer cycle loop.

\begin{caselist} % case 1
\case Inside $L$ is the second attachment $M$ \eqref{attach}
and the triangle \eqref{tatt}.

\begin{caselist}
\case %1.1
The valence of $M$ is 3. If $M$ is attached within some of
the triangle loops, then using flatness and \eqref{Fig17} we
see that the other triangle loop must not intertwine anything
outside $L$ (i.e. the diagram on the right of \eqref{Fig17}
must become composite). Then, however, $IG(A)$ has no cycle.
So $M$ is of depth $1$ within $L$, and both triangle loops
are empty.

Then the legs of the triangle identify only two
regions $A,B$ outside $L$, and by lemma \ref{lm1014}
only one of the 2 triangle loops is intertwined with other
connections. This implies that region $A$ must be between two
edges of the same outside connection of $L$. So we have, up to
mutation,
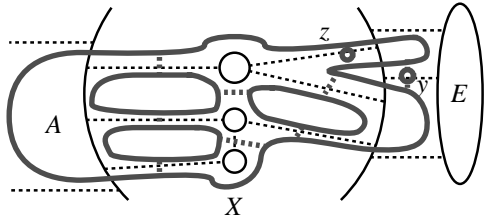
\[
% \fbox{
\diag{1cm}{6.4}{4}{
  \pictranslate{3 2.5}{
    \picmultigraphics[S]{2}{-1 1}{
      \piccirclearc{0 0}{2}{-45 35}
    }
    {\piclinedash{0.05 0.05}{0.01}
     \picline{2 170 polar}{2 170 polar x p 0 x}
     \picline{2 190 polar}{2 190 polar x p 0 x}
     \picline{2 170 polar x p 0 x}{2 18 polar}
     \picline{2 170 polar x p 0 x}{2 -3 polar}
     \picline{2 190 polar x p 0 x}{2 -20 polar}
     \picline{0 -0.9}{2 190 polar x p 0 x}
     \picline{0 -0.9}{2 210 polar}
    }
    \picfilledcircle{2 170 polar x p 0 x}{0.2}{}
    \picfilledcircle{2 190 polar x p 0 x}{0.15}{}
    {\piclinedash{0.05 0.05}{0.01}
     \picline{2 160 polar}{2 160 polar x p -3 x}
     \picline{2 220 polar}{2 220 polar x p -3 x}
     \picline{2 25 polar}{2 25 polar x p 3 x}
     \picline{2 -25 polar}{2 -25 polar x p 3 x}
     \picline{2 6 polar}{2 6 polar x p 3 x}
    }
    \picfilledcircle{0 -0.9}{0.15}{}
    \picfilledellipse{3 0}{0.3 1.2}{$E$}
    \picputtext{1.2  0.8}{$z$}
    \picputtext{2.5 0.1}{$y$}
    \picputtext{0.0 -1.5}{$X$}
    \picputtext{-2.4 -0.4}{$A$}
    { \piclinewidth{14}
      \picgraycol{0.3}
      \picPSgraphics{2 setlinecap}
      \picstroke{
	\curvepath{2.4 0.75}{2.3 0.75}{2.0 0.71}{0.4 0.5}{0.4 0.8}
	  {-0.4 0.8}{-0.4 0.5}{-1.8 0.55}{-2.3 0.55}{-3.0 0.5}
	  {-3.0 -1.1}{-2.3 -1.15}{-1.8 -1.15}{-0.3 -1.0}
	  {-0.3 -1.3}{0.3 -1.3}{0.4 -0.4}{1.3 -0.71}{2.3 -0.77}
	  {2.5 -0.71}{2.6 -0.3}{2.5 -0.0}{2.0 0.1}{1.5 0.25}{1.2 0.29}
	  {1.2 0.34}{2.0 0.4}
	  {2.5 0.4}{2.6 0.6}{}
	\curvepath{-0.6 0.25}{-0.7 0.25}{-1.9 0.25}{-1.9 0}
	  {-1.9  -0.28}{-0.3 -0.28}{-0.2 -.12}{-0.2 0.08}{-0.4 0.2}{}
	\curvepath{-0.5 -.45}{-0.7 -.45}{-1.9 -.45}{-1.7 -0.7}
	  {-1.6  -0.9}{-0.3 -0.9}{-0.2 -.72}{-0.2 -0.59}{-0.2 -0.52}{}
        \curvepath{0.5 0.08}{0.6 0.08}{1.8 -0.2}{1.8 -0.4}{1.5 -0.6}
	  {0.6 -0.36}{0.3 -0.30}{0.2 -0.12}{0.2 0.1}{}
      }
      \piccircle{1.5 0.52}{0.07}{}
      \piccircle{2.3 0.24}{0.08}{}
      {\piclinedash{0.05 0.05}{0.01}
       \picPSgraphics{0 setlinecap}
       \picline{-1.0 0.25}{-1.0 0.55}
       \picline{-1.0 -0.23}{-1.0 -0.45}
       \picline{-1.0 -0.9}{-1.0 -1.1}
       \picline{1.2 0.0}{1.4 0.33}
       % \picline{1.4 -0.5}{1.4 -0.73}
       \picline{0.2 0.02}{-0.2 0.02}
       \picline{2.3 0.08}{2.3 0.16}
       \picline{2.3 0.32}{2.3 0.49}
       \picline{-0.2 -.6}{0.4 -0.7}
       \picline{0.9 -0.45}{0.8 -0.61}
      }
    }
  }
%}
}\,.
\]
Then we apply the same argument as for \eqref{..1}.
% (LABEL EQN b4 ...[1] & refer to eqn)

\case %1.2
The valence of $M$ is 2.
We have up to mutations and symmetries two options:
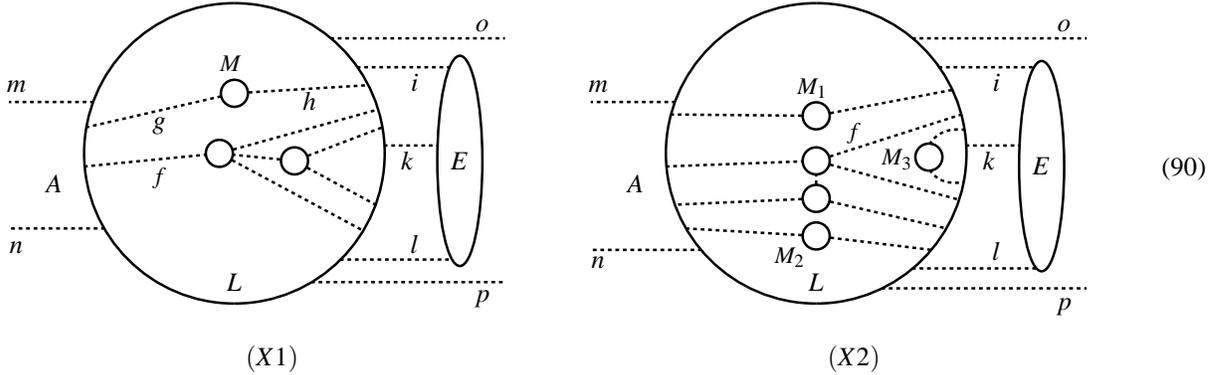
\begin{eqn}\label{77}
\begin{array}{c@{\qquad}c}
\diag{1cm}{7}{5}{
  \pictranslate{3 2.5}{
    \picmultigraphics[S]{2}{-1 1}{
      \piccirclearc{0 0}{2}{-90 90}
    }
    {\piclinedash{0.05 0.05}{0.01}
     \picline{2 170 polar}{0 0.8}
     \picline{2 185 polar}{-0.2 0}
     \picline{0 0.8}{2 27 polar}
     \picline{-0.2 0}{2 17 polar}
     \picline{-0.2 0}{2 -31 polar}
     \picline{0.8 -0.1}{-0.2 0}
     \picline{0.8 -0.1}{2 -20 polar}
     \picline{0.8 -0.1}{2 10 polar}
    }
    \picfilledcircle{0 0.8}{0.18}{}
    \picfilledcircle{-0.2 0}{0.18}{}
    \picfilledcircle{0.8 -0.1}{0.18}{}
    {\piclinedash{0.05 0.05}{0.01}
     \picline{2 160 polar}{2 160 polar x p -3 x}
     \picline{2 210 polar}{2 210 polar x p -3 x}
     \picline{2 35 polar}{2 35 polar x p 3 x}
     \picline{2 -45 polar}{2 -45 polar x p 3 x}
     \picline{2 -59 polar}{2 -59 polar x p 3.6 x}
     \picline{2 50 polar}{2 50 polar x p 3.6 x}
     \picline{2 3 polar}{2 3 polar x p 3 x}
    }
    \picfilledellipse{3 -0.1}{0.3 1.4}{$E$}
    % \picputtext{1.2  0.8}{$z$}
    \picputtext{2.3 -0.1}{$k$}
    \picputtext{2.4 0.95}{$i$}
    \picputtext{2.4 -1.2}{$l$}
    \picputtext{0.0 -1.72}{$L$}
    \picputtext{-2.4 -0.4}{$A$}
    \picputtext{-0.05 1.2}{\small$M$}
    \picputtext{-1.0 -0.3}{\small$f$}
    \picputtext{3.3 1.7}{$o$}
    \picputtext{3.3 -1.9}{$p$}
    \picputtext{-2.9 0.9}{$m$}
    \picputtext{-2.9 -1.25}{$n$}
    \picputtext{-1 0.4}{$g$}
    \picputtext{1 0.7}{$h$}
  }
}
& \diag{1cm}{7}{5}{
  \pictranslate{3 2.5}{
    \picmultigraphics[S]{2}{-1 1}{
      \piccirclearc{0 0}{2}{-90 90}
    }
    {\piclinedash{0.05 0.05}{0.01}
     \picline{2 165 polar}{0 0.5}
     \picline{2 185 polar}{0 -0.1}
     \picline{0 0.5}{2 25 polar}
     \picline{0 -0.1}{2 15 polar}
     \picline{0 -0.1}{2 -18 polar}
     \picline{0 -0.6}{2 -30 polar}
     \picline{0 -0.6}{0 -0.1}
     \picline{0 -0.6}{2 200 polar}
     \picline{0 -1.1}{2 210 polar}
     \picline{0 -1.1}{2 -40 polar}
     \piccurve{2 9 polar}{1.3 0.4}{1.3 -0.5}{2 -11 polar}
    }
    \picfilledcircle{0 0.5}{0.18}{}
    \picfilledcircle{0 -0.1}{0.18}{}
    \picfilledcircle{0 -0.6}{0.18}{}
    \picfilledcircle{0 -1.1}{0.18}{}
    \picfilledcircle{1.5 -0.05}{0.18}{}
    {\piclinedash{0.05 0.05}{0.01}
     \picline{2 160 polar}{2 160 polar x p -3 x}
     \picline{2 220 polar}{2 220 polar x p -3 x}
     \picline{2 35 polar}{2 35 polar x p 3 x}
     \picline{2 -50 polar}{2 -50 polar x p 3 x}
     \picline{2 -64 polar}{2 -64 polar x p 3.6 x}
     \picline{2 50 polar}{2 50 polar x p 3.6 x}
     \picline{2 3 polar}{2 3 polar x p 3 x}
    }
    \picfilledellipse{3 -0.17}{0.3 1.4}{$E$}
    % \picputtext{1.2  0.8}{$z$}
    \picputtext{2.3 -0.1}{$k$}
    \picputtext{2.4 0.95}{$i$}
    \picputtext{2.4 -1.3}{$l$}
    \picputtext{0.0 -1.72}{$L$}
    \picputtext{-2.4 -0.4}{$A$}
    \picputtext{0.5 0.3}{\small$f$}
    \picputtext{1.08 -0.05}{\small$M_3$}
    \picputtext{-0.05 0.9}{\small$M_1$}
    \picputtext{-0.35 -1.35}{\small$M_2$}
    \picputtext{3.3 1.7}{$o$}
    \picputtext{3.3 -2.0}{$p$}
    \picputtext{-2.9 0.9}{$m$}
    \picputtext{-2.9 -1.45}{$n$}
  }
} \\
(X1) & (X2)
\end{array}
\end{eqn}
In the second case we indicated the possible positions
for $M$ as $M_i$. Only one of these positions is used. An
attachment $M$ of depth 2 inside some of the triangle loops is
prohibited in (X1) because $IG(A)$ will have no two cycles, and
in (X2) because $\gm$ becomes $>4$ (and otherwise by corollary
\reference{crdp}).

The position
$M_3$ in $(X1)$ is excluded, because $IG(A)$ has no cycle.
For the same reason, adopting the nomenclature of \eqref{tat2}
for the triangle edges, the edge $f$ in (X1) has non-zero
multiplicity. The same is true in (X2), because otherwise
$M_2$ would give $\gm=2$, while $M_{1,3}$ would be inadequate.
(We indicated only the edge $f$ of \eqref{tat2} in the
above diagrams, to call to mind that the orientation 
of the pattern \eqref{tat2} differs slightly.) Moreover,
$M_{1,2}$ are symmetric in (X1).

Here is the first situation where it appeared more useful
to apply a computer to check the cases, rather than discussing 
them manually one by one. However, as a preparation,
we need still to argue out first some of them by hand.

Remember that to apply our program, and work with the 
simplest diagram according to lemma \ref{lmcmp}, we 
need to show that there is at most one multiple edge 
between empty outer cycle loops. (We already assured 
that there is at most one non-empty outer cycle loop.)

So we need to show that $IG(A)$ has at most 3 cycles,
i.e. rule out the cases that $IG(A)$ has $\ge 4$ cycles.
(One could likely deal with the other cases using similar
arguments, rather than using a computer, but this discussion 
would become too long and tedious.)

\begin{sublemma}\label{lmm3}
$IG(A)$ has at most 3 cycles.
\end{sublemma}

\proof
Now assume $IG(A)$ has $\ge 4$ cycles. This means that 
there are 3 loops inside $L$, two from the triangle and 
one from \eqref{attach}, and the connections from $L$ to 
them are intertwined with the one to $E$ and the two outer 
cycle loops neighbored to $L$. Then $IG(A)=K_{3,3}$,
the complete bipartite graph on $3+3$ vertices.

In both cases the assumed 4 cycles force $A$ to be between edges of 
the same outside connection of $L$. Similarly the outside attached 
loop $E$ is in a region $B$ between edges of the same connection. 
Thus the multiplicities of all edges in (X1), (X2) are non-zero, 
\em{except} for the choice between $M_i$ for $M$ in (X2), and for 
some of the edges of the triangle.
%end...[2]

%...[3]
\begin{figure}[htb]
\begin{eqn}\label{strs}
\begin{array}{c}
\begin{array}{c@{\qquad}c}
\diag{1cm}{7}{4.5}{
  \pictranslate{3 2.5}{
    \picmultigraphics[S]{2}{-1 1}{
      \piccirclearc{0 0}{2}{-90 90}
    }
    {\piclinedash{0.05 0.05}{0.01}
     \picline{2 170 polar}{2 170 polar x p 0 x}
     \picline{2 190 polar}{2 190 polar x p 0 x}
     \picline{2 170 polar x p 0 x}{2 18 polar}
     \picline{2 190 polar x p 0 x}{2 -20 polar}
     \picline{2 190 polar x p 0 x}{2 170 polar x p 0 x}
     \picline{2 207 polar}{0 -0.96}
     \picline{2 -35 polar}{0 -0.96}
    }
    \picfilledcircle{2 170 polar x p 0 x}{0.17}{}
    \picfilledcircle{2 190 polar x p 0 x}{0.17}{}
    \picfilledcircle{0 -0.96}{0.17}{}
    {\piclinedash{0.05 0.05}{0.01}
     \picline{2 65 polar}{2 65 polar x p 3 x}
     \picline{2 -70 polar}{2 -70 polar x p 3 x}
     \picline{2 160 polar}{2 160 polar x p -3 x}
     \picline{2 214 polar}{2 214 polar x p -3 x}
     \picline{2 25 polar}{2 25 polar x p 3 x}
     \picline{2 -45 polar}{2 -45 polar x p 3 x}
     \picline{2 6 polar}{2 6 polar x p 3 x}
    }
    \picfilledellipse{3 -0.2}{0.3 1.33}{$E$}
    \picputtext{1.2  0.8}{$z$}
    \picputtext{1.8 0.2}{$y$}
    \picputtext{0.0 -1.4}{$X$}
    \picputtext{0.0 3.8 x}{$Y$}
    \picputtext{-2.4 -1.45}{$V$}  
    \picputtext{-2.4 1.1}{$U$}  
    { \piclinewidth{14}
      \picgraycol{0.3}
      \picPSgraphics{2 setlinecap}
      \picstroke{
	\curvepath{2.4 0.75}{2.3 0.75}{2.0 0.71}{0.4 0.5}{0.4 0.8} 
	  {-0.4 0.8}{-0.4 0.5}{-1.8 0.55}{-2.3 0.55}{-3.0 0.5}
	  {-3.0 -1.0}{-2.3 -1.05}{-1.8 -1.05}{-0.4 -1.0}{-0.4 -1.2}
	  {0.4 -1.2}{0.4 -1.1}{1.3 -1.31}{2.3 -1.37}{2.5 -1.31}
	  {2.6 -0.9}{2.5 0.0}{2.0 0.0}{2.1 -0.3}{2.0 -0.65}{1.8 -0.6}
	  {1.5 -0.55}{0.6 -0.36}{0.3 -0.30}{0.2 -0.12}{0.2 0.1}
	  {0.4 0.08}{0.5 0.18}{0.6 0.28}{2.0 0.4}{2.5 0.4}{2.6 0.6}{}
	\curvepath{-0.6 0.25}{-0.7 0.25}{-1.9 0.25}{-1.9 0}
	  {-1.9  -0.28}{-0.3 -0.28}{-0.2 -.12}{-0.2 0.08}{-0.4 0.2}{}
	\curvepath{-0.5 -.45}{-0.7 -.45}{-1.9 -.45}{-1.7 -0.7}
	  {-1.6 -0.8}{-0.3 -0.8}{-0.2 -.76}{-0.2 -0.75}
          {0.2 -.75}{0.3 -0.85}{1.7 -1.1}{1.7 -1.0}{1.6 -0.65}
	  {0.6 -0.59}{0.3 -0.55}{0.2 -0.57}{0.15 -0.6}{-0.15 -0.6}
	  {-0.2 -0.5}{}
        \opencurvepath{-3 0.8}{-2 0.8}{-1 2.2}{0 2.2}{0.5 2.2}
	  {1.1 1.9}{2.7 1.9}{3.1 1.9}{}
        \curvepath{2.1 1.7}{1.3 1.7}{1.7 1.3}{1.9 0.9}
	  {2.7 0.9}{2.9 1.4}{3.3 1.4}{3.5 0}{3.3 -1.7}{3.0 -1.7}
	  {2.7 -1.6}{1.8 -1.5}{1.5 -1.65}{1.5 -1.8}{2.5 -1.8} 
	  {3.5 -1.8}{4 -1.5}{4.0 0}{4.0 1.2}{3.6 1.5}{3.35 1.65}{} 
         \opencurvepath{-3 -1.2}{-1.8 -1.2}{-1 -2.1}{0 -2.2}{1.3 -2.0}
	   {3.3 -2.0}{}
      }
      \piccircle{1.5 0.52}{0.07}{}
      \piccircle{2.3 0.24}{0.08}{}
      {\piclinedash{0.05 0.05}{0.01}
       \picPSgraphics{0 setlinecap}
       \picline{-1.0 0.25}{-1.0 0.55}
       \picline{-1.0 -0.23}{-1.0 -0.5}
       \picline{-1.0 -0.75}{-1.0 -1.1}
       \picline{1.2 -0.5}{1.2 -0.68}
       \picline{1.2 -1.2}{1.2 -1.08}
       \picline{0.2 0.02}{-0.2 0.02}
       \picline{2.3 0.01}{2.3 0.16}
       \picline{2.3 0.32}{2.3 0.49}
       \picline{2.3 0.7}{2.3 0.9}
       \picline{2.3 1.7}{2.3 1.9}
       \picline{2.2 -1.4}{2.2 -1.6}
       \picline{2.2 -1.9}{2.2 -2.0}
       \picline{-2.2 -1.25}{-2.2 -1.0}
       \picline{-2.2 0.6}{-2.2 0.85}
      }
    }
  }
}
&
\diag{1cm}{7}{4.5}{  % (**)
  \pictranslate{3 2.5}{  % (**)
    \picmultigraphics[S]{2}{-1 1}{  % (**)
      \piccirclearc{0 0}{2}{-90 90}  % (**)
    }  % (**)
    {\piclinedash{0.05 0.05}{0.01}  % (**)
     \picline{2 170 polar}{2 170 polar x p 0 x}  % (**)
     \picline{2 190 polar}{2 190 polar x p 0 x}  % (**)
     \picline{2 170 polar x p 0 x}{2 18 polar}  % (**)
     \picline{2 190 polar x p 0 x}{2 -20 polar}  % (**)
     \picline{2 190 polar x p 0 x}{2 170 polar x p 0 x}  % (**)
     \picline{2 170 polar x p 0 x}{2 -3 polar}  % (**)
     \picline{2 65 polar}{2 65 polar x p 3 x}  
     \picline{2 -70 polar}{2 -70 polar x p 3 x} 
     \picline{2 207 polar}{0 -0.96}  % (**)
     \picline{2 -35 polar}{0 -0.96}  % (**)
    }  % (**)
    \picfilledcircle{2 170 polar x p 0 x}{0.17}{}  % (**)
    \picfilledcircle{2 190 polar x p 0 x}{0.17}{}  % (**)
    \picfilledcircle{0 -0.96}{0.17}{}  % (**)
    {\piclinedash{0.05 0.05}{0.01}  % (**)
     \picline{2 160 polar}{2 160 polar x p -3 x}  % (**)
     \picline{2 214 polar}{2 214 polar x p -3 x}  % (**)
     \picline{2 25 polar}{2 25 polar x p 3 x}  % (**)
     \picline{2 -45 polar}{2 -45 polar x p 3 x}  % (**)
     \picline{2 6 polar}{2 6 polar x p 3 x}  % (**)
    }  % (**)
    \picfilledellipse{3 -0.2}{0.3 1.33}{$E$}  
    \picputtext{1.2  0.8}{$z$}  % (**)
    \picputtext{2.5 0.1}{$y$}  % (**)
    \picputtext{0.0 -1.4}{$X$}  % (**)
    \picputtext{0.0 3.8 x}{$Y$}  
    \picputtext{-2.4 -1.45}{$V$}  
    \picputtext{-2.4 1.1}{$U$}  
    { \piclinewidth{14}  % (**)
      \picgraycol{0.3}  % (**)
      \picPSgraphics{2 setlinecap}  % (**)
      \picstroke{  % (**)
	\curvepath{2.4 0.75}{2.3 0.75}{2.0 0.71}{0.4 0.5}{0.4 0.8}  % (**)
	  {-0.4 0.8}{-0.4 0.5}{-1.8 0.55}{-2.3 0.55}{-3.0 0.5}  % (**)
	  {-3.0 -1.0}{-2.3 -1.05}{-1.8 -1.05}{-0.4 -1.0}  % (**)
	  {-0.4 -1.2}{0.4 -1.2}{0.4 -1.1}{1.3 -1.31}{2.3 -1.37}  % (**)
	  {2.5 -1.31}{2.6 -0.9}{2.5 -0.0}{2.0 0.1}{1.5 0.25}{1.2 0.29}  % (**)
	  {1.2 0.34}{2.0 0.4}{2.5 0.4}{2.6 0.6}{}  % (**)
	\curvepath{-0.6 0.25}{-0.7 0.25}{-1.9 0.25}{-1.9 0}  % (**)
	  {-1.9  -0.28}{-0.3 -0.28}{-0.2 -.12}{-0.2 0.08}{-0.4 0.2}{}  % (**)
        \curvepath{0.5 0.08}{0.6 0.08}{1.8 -0.2}{1.8 -0.4}{1.5 -0.6}  % (**)
	  {0.6 -0.36}{0.3 -0.30}{0.2 -0.12}{0.2 0.1}{}  % (**)
	\curvepath{-0.5 -.45}{-0.7 -.45}{-1.9 -.45}{-1.7 -0.7}  % (**)
	  {-1.6 -0.8}{-0.3 -0.8}{-0.2 -.76}{-0.2 -0.75}  % (**)
          {0.2 -.75}{0.3 -0.85}{1.7 -1.1}{1.7 -1.0}{1.6 -0.65}  % (**)
	  {0.6 -0.59}{0.3 -0.55}{0.2 -0.57}{0.15 -0.6}{-0.15 -0.6}  % (**)
	  {-0.2 -0.5}{}  % (**)
        \curvepath{2.1 1.7}{1.3 1.7}{1.7 1.3}{1.9 0.9}
	  {2.7 0.9}{2.9 1.4}{3.3 1.4}{3.5 0}{3.3 -1.7}{3.0 -1.7}
	  {2.7 -1.6}{1.8 -1.5}{1.5 -1.65}{1.5 -1.8}{2.5 -1.8} 
	  {3.5 -1.8}{4 -1.5}{4.0 0}{4.0 1.2}{3.6 1.5}{3.35 1.65}{} 
        \opencurvepath{-3 0.8}{-2 0.8}{-1 2.2}{0 2.2}{0.5 2.2}
	  {1.1 1.9}{2.7 1.9}{3.1 1.9}{}
        \opencurvepath{-3 -1.2}{-1.8 -1.2}{-1 -2.1}{0 -2.2}{1.3 -2.0}
	   {3.3 -2.0}{}
      }  % (**)
      \piccircle{1.5 0.52}{0.07}{}  % (**)
      \piccircle{2.3 0.24}{0.08}{}  % (**)
      {\piclinedash{0.05 0.05}{0.01}  % (**)
       \picPSgraphics{0 setlinecap}  % (**)
       \picline{-1.0 0.25}{-1.0 0.55}  % (**)
       \picline{-1.0 -0.23}{-1.0 -0.5}  % (**)
       \picline{-1.0 -0.75}{-1.0 -1.1}  % (**)
       \picline{1.2 0.0}{1.4 0.33}  % (**)
       \picline{1.2 -0.5}{1.2 -0.68}  % (**)
       \picline{1.2 -1.2}{1.2 -1.08}  % (**)
       \picline{0.2 0.02}{-0.2 0.02}  % (**)
       \picline{2.3 0.08}{2.3 0.16}  % (**)
       \picline{2.3 0.32}{2.3 0.49}  % (**)
       \picline{2.3 0.7}{2.3 0.9}
       \picline{2.3 1.7}{2.3 1.9}
       \picline{2.2 -1.4}{2.2 -1.6}
       \picline{2.2 -1.9}{2.2 -2.0}
       \picline{-2.2 -1.25}{-2.2 -1.0}
       \picline{-2.2 0.6}{-2.2 0.85}
      }  % (**)
    }  % (**)
  }  % (**)
}  % (**)
\\
(*) & (**) 
\end{array}\\
\begin{array}{c}
\diag{1cm}{7}{5}{  % (***)
  \pictranslate{3 2.5}{  % (***)
    \picmultigraphics[S]{2}{-1 1}{  % (***)
      \piccirclearc{0 0}{2}{-90 90}  % (***)
    }  % (***)
    {\piclinedash{0.05 0.05}{0.01}  % (***)
     \picline{2 170 polar}{2 170 polar x p 0 x}  % (***)
     \picline{2 190 polar}{2 190 polar x p 0 x}  % (***)
     \picline{2 170 polar x p 0 x}{2 18 polar}  % (***)
     \picline{2 190 polar x p 0 x}{2 3 polar}  % (***)
     \picline{2 190 polar x p 0 x}{2 -12 polar}  % (***)
     \picline{0 -0.9}{2 -30 polar}  % (***)
     \picline{0 -0.9}{2 190 polar x p 0 x}  % (***)
     \picline{0 -0.9}{2 210 polar}  % (***)
    }  % (***)
    \picfilledcircle{2 170 polar x p 0 x}{0.2}{}  % (***)
    \picfilledcircle{2 190 polar x p 0 x}{0.15}{}  % (***)
    {\piclinedash{0.05 0.05}{0.01}  % (***)
     \picline{2 160 polar}{2 160 polar x p -3 x}  % (***)
     \picline{2 220 polar}{2 220 polar x p -3 x}  % (***)
     \picline{2 25 polar}{2 25 polar x p 3 x}  % (***)
     \picline{2 -40 polar}{2 -40 polar x p 3 x}  % (***)
     \picline{2 -60 polar}{2 -60 polar x p 3.6 x}  % (***)
     \picline{2 50 polar}{2 50 polar x p 3.6 x}  % (***)
     \picline{2 -3 polar}{2 -3 polar x p 3 x}  % (***)
    }  % (***)
    \picfilledcircle{0 -0.9}{0.15}{}  % (***)
    \picfilledellipse{3 -0.2}{0.3 1.2}{$E$}  % (***)
    % \picputtext{1.2  0.8}{$z$}  % (***)
    \picputtext{2.3 -0.0}{$y$}  % (***)
    \picputtext{0.0 -1.5}{$X$}  % (***)
    \picputtext{0.0 3.8 x}{$Y$}
    \picputtext{-2.4 -1.65}{$V$}  % (***)
    \picputtext{-2.4 1.1}{$U$}  % (***)
    { \piclinewidth{14}  % (***)
      \picgraycol{0.3}  % (***)
      \picPSgraphics{2 setlinecap}  % (***)
      \picstroke{  % (***)
	\curvepath{2.4 0.75}{2.3 0.75}{2.0 0.71}{0.4 0.5}{0.4 0.8}  % (***)
	  {-0.4 0.8}{-0.4 0.5}{-1.8 0.55}{-2.3 0.55}{-3.0 0.5}  % (***)
	  {-3.0 -1.1}{-2.3 -1.15}{-1.8 -1.15}{-0.3 -1.0}  % (***)
	  {-0.3 -1.3}{0.3 -1.3}{0.4 -1.1}{1.3 -1.11}{2.3 -1.17}  % (***)
	  {2.5 -1.11}{2.6 -0.7}{2.5 -0.6}{2.0 -0.3}{1.5 -0.28}  % (***)
	  {1.2 -0.25}{1.2 -0.16}{2.0 -0.04}{2.5 0.1}{2.6 0.4}{}  % (***)
	\curvepath{-0.5 -.45}{-0.7 -.45}{-1.9 -.45}{-1.7 -0.7}  % (***)
	  {-1.6  -0.9}{-0.3 -0.9}{-0.2 -.72}{-0.2 -0.59}{-0.2 -0.52}{}  % (***)
        \curvepath{0.5 -.45}{0.6 -.45}{1.8 -0.5}{1.8 -0.6}{1.5 -0.9}  % (***)
	  {0.6 -0.86}{0.3 -0.80}{0.2 -0.72}{0.2 -0.5}{}  % (***)
        \curvepath{0.5 -.15}{0.6 -.14}{1.8 0.1}{1.8 0.5}{1.5 0.5}  % (***)
	  {0.6 0.4}{0.4 0.23}{-0.0 -0.0}{-0.3 0.2}{-0.6 0.25}  % (***)
	  {-0.7 0.25}{-1.9 0.25}{-1.9 0}{-1.9  -0.28}{-0.2 -.12}  % (***)
	  {-0.2 -0.18}{0.0 0.02}{0.2 -0.18}{0.3 -0.17}{}  % (***)
        \curvepath{4.0 0}{4.0 1.2}{3 1.5}{1.6 1.5}{1.85 0.95}{2.6 0.85}  % (***)
	  {3.05 1.2}{3.25 1.06}{3.45 0}{3.35 -1.1}{3.2 -1.3}{3.1 -1.55}  % (***)
	  {2.7 -1.4}{1.8 -1.35}{1.5 -1.5}{1.5 -1.65}{2.5 -1.65}  % (***)
	  {3.5 -1.65}{4 -1.3}{}  % (***)
        \opencurvepath{-3 0.8}{-2 0.8}{-1 2.2}{0 2.2}{0.5 2.1}
	  {1.3 1.67}{2.7 1.67}{3.1 1.67}{}
        \opencurvepath{-3 -1.4}{-1.8 -1.4}{-1 -2.1}{0 -2.2}{1.3 -1.9}
	   {3.3 -1.9}{}
      }  % (***)
      \piccircle{2.3 -0.2}{0.11}{}  % (***)
      {\piclinedash{0.05 0.05}{0.01}  % (***)
       \picPSgraphics{0 setlinecap}  % (***)
       \picline{-1.0 0.25}{-1.0 0.55}  % (***)
       \picline{-1.0 -0.23}{-1.0 -0.45}  % (***)
       \picline{-1.0 -0.9}{-1.0 -1.1}  % (***)
       \picline{1.2 0.0}{1.3 -0.23}  % (***)
       \picline{1.2 0.5}{1.2 0.66}  % (***)
       \picline{1.3 -0.5}{1.4 -0.23}  % (***)
       % \picline{1.4 -0.5}{1.4 -0.73}  % (***)
       \picline{2.3 -.48}{2.3 -.31}  % (***)
       \picline{2.3 -.1}{2.3 0.1}  % (***)
       \picline{-0.2 -.6}{0.2 -0.65}  % (***)
       \picline{1.3 -1.15}{1.3 -0.81}  % (***)
       \picline{1.9 -1.15}{1.9 -1.41}  % (***)
       \picline{2.1 -1.65}{2.1 -1.91}  % (***)
       \picline{2.2 0.75}{2.2 0.91}  % (***)
       \picline{2.5 1.65}{2.5 1.45}  % (***)
       \picline{-2.2 -1.35}{-2.2 -1.1}
       \picline{-2.2 0.6}{-2.2 0.85}
      }  % (***)
    }  % (***)
  }  % (***)
} \\  % (***)
(***)
\end{array}
\end{array}
\end{eqn}
\caption{\label{ffg}}
\end{figure}

There are, up to mutations, 3 cases, see figure \reference{ffg}.
(Here we assume all edge multiplicities to be non-zero.)
The $B$-state loops $X$ and $Y$ are drawn closed since they
are bounded by the outer cycle ($A$-state) loops neighboring
to $L$. Observe again that all edges $y$ and $z$ in $(*)$ and
$(**)$, are inadequate, so of multiplicity two.

These cases are ruled out by the following argument. Again 
there is a unique separating loop $X$ in $B(D)$. It has at 
most 3 multiple connections inside, each intertwined with
at most two connections outside. Also, $IG(A(D))$ has 
at most $4$ cycles. Taking into account the ($B$-state
loops coming from the multiplicity 2) inadequate edges
$y$ and $z$, we see that in all cases $B(D)$ has at least
7 loops, even when all outer cycle edges between empty 
loops in $A(D)$ are simple. Thus we can apply lemma \ref{UV}.
This completes the proof of sublemma \reference{lmm3}. \qed

%end...[3]

With these cases ruled out, we justified the use of our
computer algorithm on the patterns of \eqref{77}.
% can use the computer. We have the following patterns...

The edges are named as in \eqref{77}, with the triangle
edges following \eqref{tat2}. To fix the remaining names for
the edges, let also in (X2) $g,h$ be the edges of $M$, and
call the 3 edges between other empty outer cycle loops $x,y,z$
(we use, as explained in \S\ref{Scmp}, an outer cycle
of length 5).

Let us specify what conditions we have on
the multiplicities thus named.
% ...[4]
% (multiplicities of the edges)
In case X1:
\begin{itemize}
\item $f>0$ (as argued above), $k=2$ (inadequate)
\item $e=1$, otherwise it gives an isolated vertex in $IG$, and
  $\chi(IG(A))>-1$. 
\item $m,n,o,p>0$ (otherwise $IG(A)$ has at most 1 cycle)
  and not all are even (by connectivity)
\item $x=y=z=1$ (because of $\chi(IG(A))=-1$, as for $e=1$).
\item $i,l>0$ and not both even (by connectivity)
\item not all of $a,b,c,d,f$ even (by connectivity)
\item the multiplicities for the valence-2 loop legs follow 
  lemma \ref{l2v}
\end{itemize}
In case X2:
\begin{itemize}
\item $f>0$ (as argued above), $k=2$ (inadequate)
\item $o,p>0$ and one of $m,n>0$ (w.l.o.g. up to mutation);
  $m,n,o,p$ are not all even (by connectivity)
\item $y=z=1$; $x,e$ are not both $>1$, and if one is $>1$, 
  then both $m,n>0$ (because of $\chi(IG(A))=-1$).
\item $i,l>0$ and not both even (by connectivity)
\item not all of $a,b,c,d,f$ even (by connectivity)
% \item if $f=0$, then $a\ge b,c,d$ (w.l.o.g. up to mutation) and
  % one of $a,b,c,d$ is 0 (by lemma \ref{lm1014})
\item the multiplicities for the valence-2 loop legs 
  follow lemma \ref{l2v}
\end{itemize}

%(quote calc.)
\begin{clc}
% aqv23tst8a.C trdg_X[12].in
Without strong effort at speed-up, the computer calculation
took about 27 minutes, and ruled out all cases.
\end{clc}

\end{caselist}

\case %2
Inside $L$ is only the triangle, and outside $L$ are the
two attachements \eqref{attach}. 

First note that then always $IG(A)$ has $\le 3$ cycles.
To have at least 2 cycles, it is necessary that at least
one outer cycle connection of $L$ is intertwined by
the triangle. There are 3 cases:
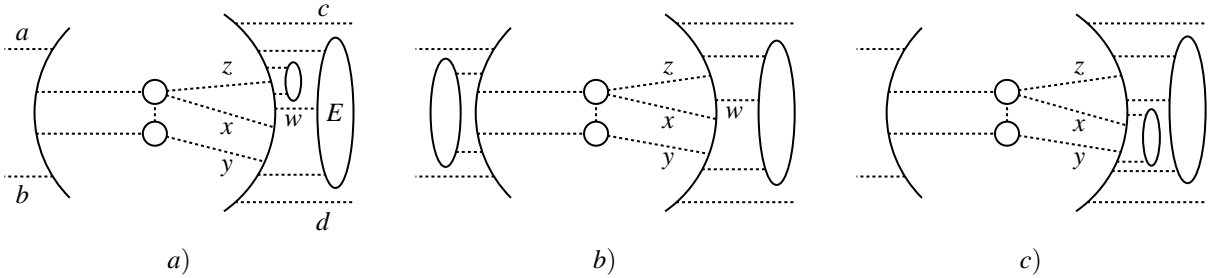
\[
\begin{array}{c@{\qquad\ }c@{\qquad\ }c}
\diag{8mm}{5.8}{3.5}{
  \pictranslate{2.5 1.8}{
    \piccirclearc{0 0}{2}{135 225}
    \piccirclearc{0 0}{2}{-55 55}
    {\piclinedash{0.05 0.05}{0.01}
     \picline{2 170 polar}{2 170 polar x p 0 x}
     \picline{2 190 polar}{2 190 polar x p 0 x}
     \picline{2 170 polar x p 0 x}{2 15 polar}
     \picline{2 170 polar x p 0 x}{2 -7 polar}
     \picline{2 190 polar x p 0 x}{2 -24 polar}
     \picline{2 190 polar x p 0 x}{2 170 polar x p 0 x}
     \picline{2 48 polar}{2 48 polar x p 3.3 x}
     \picline{2 -48 polar}{2 -48 polar x p 3.3 x}
    }
    \picfilledcircle{2 170 polar x p 0 x}{0.2}{}
    \picfilledcircle{2 190 polar x p 0 x}{0.2}{}
    {\piclinedash{0.05 0.05}{0.01}
     \picline{2 148 polar}{2 148 polar x p -2.5 x}
     \picline{2 212 polar}{2 212 polar x p -2.5 x}
     \picline{2 22 polar}{2 22 polar x p 2.3 x}
     \picline{2  9 polar}{2  9 polar x p 2.3 x}
     \picline{2 31 polar}{2 31 polar x p 3 x}
     \picline{2 -31 polar}{2 -31 polar x p 3 x}
     \picline{2 2 polar}{2 2 polar x p 3 x}
    }
    \picfilledellipse{3 0.0}{0.3 1.25}{$E$}
    \picfilledellipse{2.3 0.52}{0.13 0.32}{}
    \picputtext{1.2 -0.25}{$x$}
    \picputtext{-2.2 1.3}{$a$}
    \picputtext{-2.2 -1.35}{$b$}
    \picputtext{2.8 1.7}{$c$}
    \picputtext{2.8 -1.8}{$d$}
    \picputtext{1.2 0.75}{$z$}
    \picputtext{2.3 -0.13}{$w$}
    \picputtext{1.2 -0.87}{$y$}
  }
} &
\diag{8mm}{6.3}{3.5}{
  \pictranslate{3 1.8}{
    \piccirclearc{0 0}{2}{135 225}
    \piccirclearc{0 0}{2}{-55 55}
    {\piclinedash{0.05 0.05}{0.01}
     \picline{2 170 polar}{2 170 polar x p 0 x}
     \picline{2 190 polar}{2 190 polar x p 0 x}
     \picline{2 170 polar x p 0 x}{2 18 polar}
     \picline{2 170 polar x p 0 x}{2 -3 polar}
     \picline{2 190 polar x p 0 x}{2 -20 polar}
     \picline{2 190 polar x p 0 x}{2 170 polar x p 0 x}
     \picline{2 48 polar}{2 48 polar x p 3.3 x}
     \picline{2 -48 polar}{2 -48 polar x p 3.3 x}
    }
    \picfilledcircle{2 170 polar x p 0 x}{0.2}{}
    \picfilledcircle{2 190 polar x p 0 x}{0.2}{}
    {\piclinedash{0.05 0.05}{0.01}
     \picline{2 148 polar}{2 148 polar x p -3 x}
     \picline{2 212 polar}{2 212 polar x p -3 x}
     \picline{2 161 polar}{2 161 polar x p -2.5 x}
     \picline{2 199 polar}{2 199 polar x p -2.5 x}
     \picline{2 28 polar}{2 28 polar x p 3 x}
     \picline{2 -28 polar}{2 -28 polar x p 3 x}
     \picline{2 6 polar}{2 6 polar x p 3 x}
    }
    \picfilledellipse{3 0}{0.3 1.2}{}
    \picfilledellipse{-2.5 0}{0.25 0.9}{}
    \picputtext{1.2 -0.15}{$x$}
    \picputtext{1.2 0.75}{$z$}
    \picputtext{2.3 0.01}{$w$}
    \picputtext{1.2 -0.75}{$y$}
  }
} &
\diag{8mm}{5.8}{3.5}{
  \pictranslate{2.5 1.8}{
    \piccirclearc{0 0}{2}{135 225}
    \piccirclearc{0 0}{2}{-55 55}
    {\piclinedash{0.05 0.05}{0.01}
     \picline{2 170 polar}{2 170 polar x p 0 x}
     \picline{2 190 polar}{2 190 polar x p 0 x}
     \picline{2 170 polar x p 0 x}{2 18 polar}
     \picline{2 170 polar x p 0 x}{2 -6 polar}
     \picline{2 190 polar x p 0 x}{2 -19 polar}
     \picline{2 190 polar x p 0 x}{2 170 polar x p 0 x}
     \picline{2 48 polar}{2 48 polar x p 3.3 x}
     \picline{2 -48 polar}{2 -48 polar x p 3.3 x}
    }
    \picfilledcircle{2 170 polar x p 0 x}{0.2}{}
    \picfilledcircle{2 190 polar x p 0 x}{0.2}{}
    {\piclinedash{0.05 0.05}{0.01}
     \picline{2 148 polar}{2 148 polar x p -2.5 x}
     \picline{2 212 polar}{2 212 polar x p -2.5 x}
     \picline{2 28 polar}{2 28 polar x p 3 x}
     \picline{2 -29 polar}{2 -29 polar x p 3 x}
     \picline{2 -24 polar}{2 -24 polar x p 2.3 x}
     \picline{2 -1 polar}{2 -1 polar x p 2.3 x}
     \picline{2 6 polar}{2 6 polar x p 3 x}
    }
    \picfilledellipse{3 0.05}{0.3 1.22}{}
    \picfilledellipse{2.4 -0.41}{0.14 0.47}{}
    \picputtext{1.2 -0.75}{$y$}
    \picputtext{1.2 -0.21}{$x$}
    \picputtext{1.2 0.75}{$z$}
  }
} \\
\ry{1.7em}a) & b) & c)
\end{array}
\]
All edges are assumed to be of non-zero multiplicity,
except possibly $x$. and/or at most one of $a,b,c,d$.
(We labeled latter for simplicity only in case a).) 
If $x$ is not there, $y$ is inadequate, so $y=2$. 

Up to mutations, we may assume the intertwined outer cycle 
connection of $L$ to be on either side. In other words, a 
region bounded by two edges of this connection is either 
where $E$ is attached (so $c,d\ne 0$), or where the legs of 
the triangle loops inside $L$ going to the left in the 
diagrams are touching $L$ (here $a,b\ne 0$).

In case a), let us assume that the intertwined outer cycle
connection is on the left of $L$ in the diagram, i.e. $a\ne 
0\ne b$. Now $z$ is inadequate, so $z=2$. In $B(D)$ again the 
only separating loop $X$ has at most 3 multiple connections
inside. The one, coming from $z=2$ in $A(D)$, is intertwined
with at most one outside connection. Another connection comes
from the trace in $B(D)$ dual to $y$, if $x=1$ or $(x,y)=(0,2)$.
This connection is intertwined with at most two connections outside
$X$. This already suffices to conclude that $\chi(IG(B))>-1$.

In case b), assume w.l.o.g. $c,d\ne 0$. Now $z,w$ are inadequate, 
and with the $B$-state loops coming from $z=w=2$ we count 
already 7 loops. Again the loops $U$, $V$ have a connection 
not intertwined with anything. 
% (If one
% of $a$ or $b$ is zero, one of $U$ and $V$ has a connection to the 
% separating $B$-state loop $X$. However, even in this event this
% connection is not intertwined with any of the connections 
% inside $X$.) 
Thus we can use the same argument of lemma \ref{UV} with which we 
ruled out \eqref{strs}.

Finally, case c) is handled similarly to case a). We assume
again $a\ne 0\ne b$. The only difference to the argument for
case a) is that the connection inside $X$ coming from $y$ (again
if $x=1$ or $(x,y)=(0,2)$) is intertwined with at most one
connection outside $X$, while the one coming from the inadequate
edge $z=2$ is intertwined with at most two connections outside $X$.
\end{caselist}

With this lemma \reference{lm3v} is proved.
\qed
%end...[4]

Now we can go back to sublemma \ref{lmm3}, which operated
under the setting of lemma \reference{lm3v}, and extend it
fully using lemma \reference{lm3v}.

\begin{lemma}\label{lmm3'}
In case C, there are at most three cycles in $IG(A(D))$.
\end{lemma}

\proof
% If $IG(A)$ has $\ge 4$ cycles, then as before $IF(A)=K_{3,3}$.
% By exclusion of a valence-3 outeside attachement to $L$ with
% lemma \reference{lm3v}, we have
Assume $IG(A(D))$ has $\ge 4$ cycles. As in the proof 
of lemma \ref{lm3v}, we have then that $IG(A)=K_{3,3}$,
with the triangle \eqref{tatt} and one loop \eqref{attach} 
inside $L$, and the other loop \eqref{attach} outside $L$. 
Now we know, by lemma \ref{lm3v}, that the attachment 
outside $L$ has valence 2. The only option we have then is
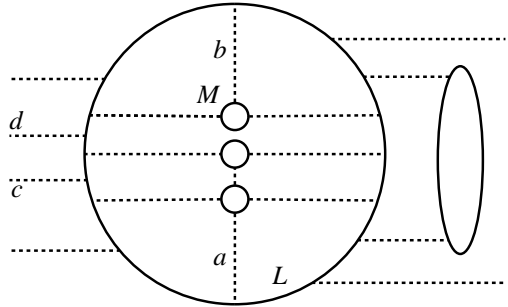
\[
\diag{1cm}{7}{4}{
  \pictranslate{3 2}{
    \picmultigraphics[S]{2}{-1 1}{
      \piccirclearc{0 0}{2}{-90 90}
    }
    {\piclinedash{0.05 0.05}{0.01}
     \picline{2 165 polar}{0 0.5}
     \picline{2 165 polar}{0 0.5}
     \picline{2 180 polar}{0 -0.0}
     \picline{0 0.5}{2 15 polar}
     \picline{0 0.5}{2 90 polar}
     \picline{0 -0.6}{2 -90 polar}
     \picline{0 -0.0}{2 -0 polar}
     \picline{0 -0.6}{2 -18 polar}
     \picline{0 -0.6}{0 -0.1}
     \picline{0 -0.6}{2 198 polar}
     % \picline{0 -1.1}{2 210 polar}
     % \picline{0 -1.1}{2 -28 polar}
    }
    \picfilledcircle{0 0.5}{0.18}{}
    \picfilledcircle{0 -0.0}{0.18}{}
    \picfilledcircle{0 -0.6}{0.18}{}
    % \picfilledcircle{0 -1.1}{0.18}{}
    {\piclinedash{0.05 0.05}{0.01}
     \picline{2 173 polar}{2 173 polar x p -3 x}
     \picline{2 190 polar}{2 190 polar x p -3 x}
     \picline{2 150 polar}{2 150 polar x p -3 x}
     \picline{2 220 polar}{2 220 polar x p -3 x}
     \picline{2 31 polar}{2 31 polar x p 3 x}
     \picline{2 -35 polar}{2 -35 polar x p 3 x}
     \picline{2 -59 polar}{2 -59 polar x p 3.6 x}
     \picline{2 50 polar}{2 50 polar x p 3.6 x}
    }
    \picfilledellipse{3 -0.08}{0.3 1.25}{}
    % \picputtext{1.2  0.8}{$z$}
    \picputtext{0.6 -1.62}{$L$}
    \picputtext{-0.2 1.4}{$b$}
    \picputtext{-0.2 -1.4}{$a$}
    \picputtext{-0.35 0.79}{$M$}
    \picputtext{-2.9 0.45}{$d$}
    \picputtext{-2.9 -.5}{$c$}
  }
}
\]
The edges $a,b,c,d$ are optional, but using $\gm=4$ we see
that exactly one of these four must be there (i.e. $\ne 0$).
Again we assume all other edges have non-zero multiplicity,
and the two edges that leave from $L$ to the left resp.
right connect it to the same outer cycle loop. 

In fact, $a$ or $b$ may connect to a region between two outer cycle
traces from $L$ to the same neighbored outer cycle loop. Then, however,
we can apply lemma \ref{UV}. (Note that, because by inadequacy $a$
or $b$ is $2$, and $c=d=0$, by connectivity at least one of the other
6 legs of the 3 loops inside $L$ is multiple, so we find 7 loops in
$B(D)$ by remark \reference{remUV}.)
% 
% 
% If $a\ne 0$, then by inadequacy, $a=2$, and can we can apply
% lemma \reference{lm1014} (where $f$ is our $a$). Similarly,
% if $b\ne 0$, then lemma \reference{lm1014} applies with $f=0$.
% If $d\ne 0$, then $M$ becomes inadequate.
% 

Now if $c=0$, then we can use lemma \reference{lm1014}
(where $f$ is our $a$). So it remains to consider $c\ne 0$.
This option is dealt with by the argument that was applied
to the diagrams in figure \reference{ffg}. (Again there
are two distinct $B$-state loops for the regions $U$
and $V$ in \eqref{oaz}.) \qed

\subsection{Case A\label{Sca}}

We argued that from the options (A1) and (A2) in \eqref{csA},
only (A2) is relevant. To have two cycles in $IG(A)$, both
attachments we perform must be inside $L$, the common loop
of the outer cycle and triangle. Both loops inside $L$ must
intertwine at least (the same) three connections outside
$L$. In fact, if (A2) is drawn as in \eqref{csA}, we easily
see that this is impossible without $\gm\ge 6$. The minorly
less trivial case results again from the situation in
which we have the phenomenon \eqref{B'} for the triangle,
and/or the triangle lies between edges of the same outer
cycle connection.
% and/or the analogous phenomenon
% for the outer cycle give a non-trivial case.
Even with this complication, the case is relatively easy to settle,
albeit with some computation.

Using $\gm=4$ and the exclusion of inadequate loops, we see
that from both loops inside $L$, one, $N$, must be of valence 3,
and the other $M$ of valence 2. The inadequate leg $e$ of $N$ is
of valence $2$. So we have (see \eqref{tat2'}) in $B(D)$ a unique
separating loop $X$, which has at most two connections inside
intertwined with something outside, and both %of them
are intertwined with at most 3 connections outside $X$.

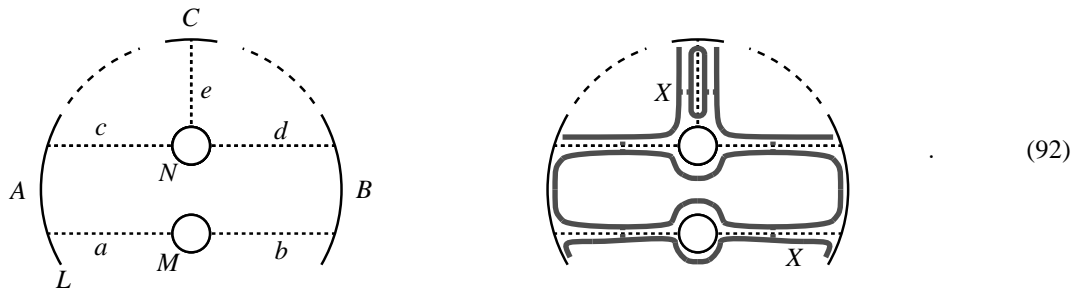
\begin{eqn}\label{tat2'}
\diag{1cm}{6}{4}{
  \pictranslate{3 1.5}{
    \picmultigraphics[S]{2}{-1 1}{
      \piccirclearc{0 0}{2}{-30 30}
      {\piclinedash{0.05 0.05}{0.01}
       \picline{2 17 polar x p 0 x}{0 2}
       \piclinedash{0.1 0.1}{0}
       \piccirclearc{0 0}{2}{35 70}
      }
      \picmultigraphics[S]{3}{1 -1}{
	{\piclinedash{0.05 0.05}{0.01}
	 \picline{2 17 polar}{2 17 polar x p 0 x}
	}
	\picfilledcircle{2 17 polar x p 0 x}{0.25}{}
      }
    }
    \piccirclearc{0 0}{2}{80 100}
    \picputtext{0.2   1.3}{$e$}
    \picputtext{1.2   0.8}{$d$}
    \picputtext{1.2  -0.8}{$b$}
    \picputtext{-1.2  0.8}{$c$}
    \picputtext{-1.2 -0.8}{$a$}
    \picputtext{-2.3 -0.0}{$A$}
    \picputtext{ 2.3 -0.0}{$B$}
    \picputtext{0.0  2.3}{$C$}
    \picputtext{-0.3  0.23}{$N$}
    \picputtext{-0.3  -0.97}{$M$}
  }
  \picputtext{ 1.3 0.3}{$L$}
}\qquad
\diag{1cm}{6}{4}{
  \pictranslate{3 1.5}{
    \picmultigraphics[S]{2}{-1 1}{
      \piccirclearc{0 0}{2}{-30 30}
      {\piclinedash{0.05 0.05}{0.01}
       \picline{2 17 polar x p 0 x}{0 2}
       \piclinedash{0.1 0.1}{0}
       \piccirclearc{0 0}{2}{35 70}
      }
      \picmultigraphics[S]{3}{1 -1}{
	{\piclinedash{0.05 0.05}{0.01}
	 \picline{2 17 polar}{2 17 polar x p 0 x}
	}
	\picfilledcircle{2 17 polar x p 0 x}{0.25}{}
      }
    }
    \piccirclearc{0 0}{2}{80 100}
    \picputtext{-0.45 1.3}{$X$}
    \picputtext{1.3 -0.9}{$X$}
    { \piclinewidth{14}
      \picgraycol{0.3}
      \picPSgraphics{0 setlinecap}
      \picmultigraphics[S]{2}{-1 1}{
         \opencurvepath{0.25 1.9}{0.25 0.9}{0.4 0.7}{1.8 0.7}{}
	 \opencurvepath{1.7 -0.9}{1.8 -0.69}{1 -0.69}{0.3 -0.6}
	   {0.3 -0.8}{0.2 -0.91}{0.1 -0.96}{0 -0.96}{}
	 \picmultigraphics[S]{2}{1 -1}{
	   \opencurvepath{0.0 0.15}{0.1 0.15}{0.3 0.25}{0.3 0.35}{0.4 0.52}
	     {1.7 0.52}{1.9 0.3}{1.9 0}{}
	   {\piclinedash{0.05 0.05}{0.01}
	    \picline{1 0.5}{1 0.7}
	   }
	 }
      }
      \picovalbox{0 1.43}{0.18 0.9}{0.08}{}
      {\piclinedash{0.05 0.05}{0.01}
       \picline{-0.24 1.3}{-0.12 1.3}
       \picline{0.24 1.3}{0.12 1.3}
      }
    }
  }
}\,.
\end{eqn}

Let $a,b,c,d$ be the 4 other legs of $M,N$. Then since $e=2$,
one of $a,b,c,d$ is $>1$. So by lemma \reference{UV} (and
remark \reference{remUV}), the region $C$ outside $L$ touched
by $e$ is among one of the two regions $U,V$ in \eqref{oaz},
and by symmetry we may assume $C=U$.

We have then the following picture. (The interior of $L$ is as
on the left of \eqref{tat2'}.)
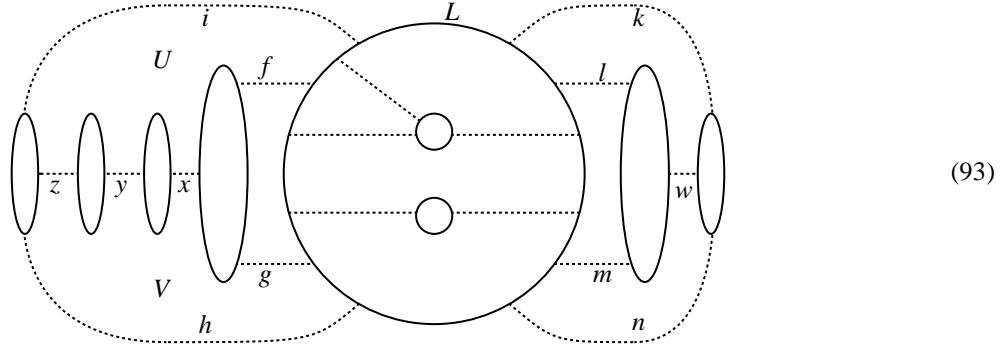
\begin{eqn}\label{dgA}
\diag{8mm}{12}{6}{
  \pictranslate{7 3}{
    {\piclinedash{0.05 0.05}{0.01}
     \picline{-3.7 1.5}{3.7 1.5}
     \picline{-3.7 -1.5}{3.7 -1.5}
     \picline{-3.7 0}{-6.7 0}
     \picline{3.7 0}{4.7 0}
     \picmultigraphics[S]{2}{1 -1}{
       \opencurvepath{2.5 60 polar}{2 2.8}{4.65 2.8}{4.65 0.8}{4.65 0}{4.5 0}
       {}
       \picscale{-1 1}{
         \opencurvepath{2.5 60 polar}{2 2.8}{6.9 2.8}{6.9 0.8}{6.7 0}{6.5 0}
         {}
       }
     }
    }
    \picmultigraphics[S]{2}{-1 1}{
      \picfilledellipse{3.5 0}{0.4 1.8}{}
      \picfilledellipse{4.6 0}{0.22 1.0}{}
    }
    \picfilledellipse{-5.7 0}{0.22 1.0}{}
    \picfilledellipse{-6.8 0}{0.22 1.0}{}
    \picfilledcircle{0 0}{2.5}{}
    \picputtext{0.3 2.7}{$L$}
    {\piclinedash{0.05 0.05}{0.01}
     \picline{0 0.7}{2.5 130 polar}
     \picline{2.5 165 polar}{2.5 15 polar}
     \picline{2.5 -165 polar}{2.5 -15 polar}
    }
    \picfilledcircle{0 0.7}{0.3}{}
    \picfilledcircle{0 -0.7}{0.3}{}
    \picputtext{-2.8 1.8}{$f$}
    \picputtext{-2.8 -1.7}{$g$}
    \picputtext{-3.8 -2.5}{$h$}
    \picputtext{-3.8 2.6}{$i$}
    \picputtext{2.8 1.7}{$l$}
    \picputtext{2.8 -1.7}{$m$}
    \picputtext{3.4 -2.5}{$n$}
    \picputtext{3.4 2.6}{$k$}
    \picputtext{4.15 -0.3}{$w$}
    \picputtext{-4.15 -0.2}{$x$}
    \picputtext{-5.2 -0.2}{$y$}
    \picputtext{-6.3 -0.2}{$z$}
    \picputtext{-4.5 1.9}{$U$}
    \picputtext{-4.5 -1.9}{$V$}
  }
}
\end{eqn}

Now we said that draw the outer cycle like \eqref{oaz}.
In the above diagram we made a unique exception to this rule.
We choose, in deviation from elsewhere, the position of the
outer cycle so that the curve giving the connected sum in
(A2) goes through infinity. (This curve is supposed to intersect
$L$ in exactly two points, but should not give a connected sum
decomposition of $D$, since it intersects traces of loops to be
attached inside $L$.) This exception is made here merely for
better display.

In \eqref{dgA}, we have the following conditions on edge multiplicities:
\begin{itemize}
\item All edges are non-zero, except at most one of $g,h,i,k,l,m,n$.
\item $f=2$ (inadequate); w.l.o.g. $y=z=1$, and $x,w$ are not both $>1$
  (because of $\chi(IG(A))=-1$)
\item If one of $g,h,i,k,l,m,n$ is zero, then $x=w=1$ (because of
  $\chi(IG(A))=-1$).
\item $e=2$ (inadequate), $c,d$ not both even (by connectivity)
\item the multiplicities for the valence-2 loop legs $a,b$
  follow lemma \ref{l2v}
\end{itemize}

\begin{clc}
% aqv23tst8d.C trdg_A.in
The computer checked these cases within 40 seconds, and ruled all out.
This finishes case A.
\end{clc}

% 
% To find the 7 loops in $B(D)$ note that $B(D)$ has a unique separating
% loop $X$, which has a connection inside not intertwined with
% anything (coming from the one of $a,b,c,d>1$), and at least 5
% other loops are needed to have two cycles in $IG(B(D))$.

\subsection{Cases B and B'\label{CBB'}}

Now we are prepared enough to finish off case B completely.
Later we move on to case C, which requires more work.
However, let us first advance with the degenerate case B'
as far as we did with case B.

We use the diagram \eqref{B'}.

\begin{lemma}\label{ghp'}
There is only one non-empty outer cycle or triangle loop
in case B', which is w.l.o.g. $L$.
\end{lemma}

\proof If both of $L$ and $M$ contain at most one attachment,
then $IG(A)$ will have no cycles. So we can assume up to
symmetry that w.l.o.g. $L$ contains at least 2 of the 3
attachments, inside or outside. 

Assume that $L$ has 2 attachments, so the third attachment
is inside some loop $N\ne L$. Note that then we have no outside
attachments to $L,N$, so we have no different legs connecting
to the same region. We distinguish the following 4 situations:
\def\labelenumi{(\theenumi}
\def\theenumi{\roman{enumi})}
\begin{enumerate}
\item\label{(i)} The third attachment has valence $3$, and
  the attachments inside $L$ identify $2$ of these $3$ regions.
\item\label{(ii)}
  The attachments inside $L$ identify $3$ regions, and the
one in $N$ identifies 2 of the three regions.
\item\label{(iii)}
  The attachments both inside $L$ and $N$ identify the same
  2 regions each, or 
\item\label{(iv)} The attachments inside $L,N$ identify
  the same triple of regions.
\end{enumerate}
In case \ref{(iii)} the two pairs of identified regions must
coincide, for otherwise the attachment within $N$ would be inadequate.
Cases \reference{(i)},  \reference{(ii)}, and \reference{(iii)}
can then be reduced to one non-empty loop by mutations.

In case \ref{(iv)}, if $N=M$, then one has a modification of
\eqref{eqw}. (Now $X$ in \eqref{eqw} is $L$ here, $Y$ is $M$,
and $A,B$ may not be those of \eqref{B'}.)
Now one sees that the attachments inside $L$ intertwine at most
two connections (to $M$ and $N$), and the attachments inside
$L,M$ intertwine at most one common connection (the one between
$L$ and $M$). Thus $IG(A)$ cannot have 2 cycles.

If $N\ne M$, then by the argument after \eqref{B'}, we may assume
that $A,B$ of \eqref{B'} are among the $A,B,C$ of \eqref{eqw}.
However, $A,B$ are not simultaneously neighbored to any other loop
than $L$, a contradiction. 

Thus we have only one non-empty loop on the outer cycle/triangle.
If this loop is different from $L$, we can reduce case B' to
case B by mutations. So we can assume that it is $L$, and the
lemma is proved. \qed

% In both cases one sees also that the above potential traces $e''$
% can be eliminated by a mutation (that moves $e''$ into the region $Q$
% resp. $U$ of \eqref{uipo}). So we did not draw them in figure \ref{B12}.

Since now $M$ is empty in \eqref{B'}, we can assume that attachments
inside $L$ touch, or attachments outside $L$ are made within, both
regions $A$ and $B$ in \eqref{B'}.

\begin{lemma}\label{ghp}
In case B', $IG(A)$ has at most 3 cycles.
\end{lemma}

\proof Assume we have $\ge 4$ cycles. Then all 3 attachments are
inside $L$, and $IG(A)=K_{3,3}$, with exactly 4 cycles. 
We have the following three options, under exclusion of
inadequate loops. (We draw for space reasons the refinement
of \eqref{B'} only until the left half of $L'$.) In case (a),
we allow the edge $y$ to connect to any region outside $L$
different from $A,B$.
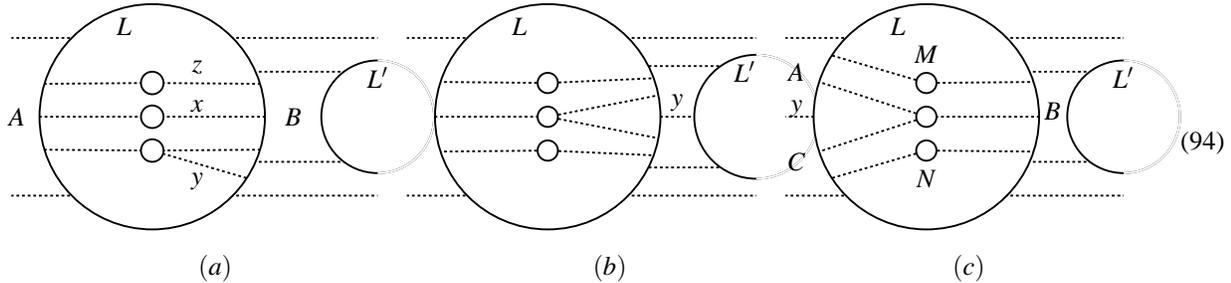
\begin{eqn}\label{B'1}
\begin{array}{ccc}
% \fbox{
\diag{7.5mm}{6.5}{4}{
 \pictranslate{0.0 0}{
  {\piclinedash{0.05 0.05}{0.01}
   \picline{-0.5 3.4}{6 3.4}
   \picline{-0.5 0.6}{6 0.6}
   \picline{2 2.8}{6 2.8}
   \picline{2 1.2}{6 1.2}
  }
  \picfilledcircle{6 2}{1.0}{}
  {\picgraycol{1}\piccirclearc{6 2}{1.0}{-90 90}}
  % \piccirclearc{6 2}{1.0}{90 270}
  \picfilledcircle{2 2}{2}{}
  \picputtext{1.5 3.6}{$L$} 
  \picputtext{6 2.7}{$L'$} 
  \picputtext{4.5 2}{$B$}
  \picputtext{-.4 2}{$A$}
  \picputtext{2.8 0.9}{$y$}
  \picputtext{2.8 2.2}{$x$}
  \picputtext{2.8 2.9}{$z$}
  \pictranslate{2 2}{
   {\piclinedash{0.05 0.05}{0.01}
    \picline{2 180 polar}{0 0}
    \picline{2 0 polar}{0 0}
    \picline{2 -163 polar}{0 -0.6}
    \picline{2 -17 polar}{0 -0.6}
    \picline{2 163 polar}{0 0.6}
    \picline{2 17 polar}{0 0.6}
    \picline{0 -0.6}{2 -33 polar}
   }
   \picfilledcircle{0 0}{0.2}{}
   \picfilledcircle{0 0.6}{0.2}{}
   \picfilledcircle{0 -0.6}{0.2}{}
  }
 }
}
% } 
&
\diag{7.5mm}{6.2}{4}{
  {\piclinedash{0.05 0.05}{0.01}
   \picline{-0.5 3.4}{5.7 3.4}
   \picline{-0.5 0.6}{5.7 0.6}
   \picline{2 2.0}{5.7 2.0}
   \picline{2 2.9}{5.7 2.9}
   \picline{2 1.1}{5.7 1.1}
  }
  \picfilledcircle{5.7 2}{1.1}{}
  {\picgraycol{1}\piccirclearc{5.7 2}{1.1}{-90 90}}
  % \piccirclearc{6 2}{1.0}{90 270}
  \picfilledcircle{2 2}{2}{}
  \picputtext{1.5 3.6}{$L$} 
  \picputtext{5.5 2.8}{$L'$} 
  \picputtext{2.3 4.3 x}{$y$} 
  \pictranslate{2 2}{
    {\piclinedash{0.05 0.05}{0.01}
     \picline{0 0.6}{2 20 polar}
     \picline{0 0.6}{2 162 polar}
     \picline{0 -0.6}{2 -20 polar}
     \picline{0 -0.6}{2 198 polar}
     \picline{-2 0}{0 0}
     \picline{0 0}{2 11 polar}
     \picline{0 0}{2 -11 polar}
    }
    \picfilledcircle{0 0.0}{0.18}{}
    \picfilledcircle{0 0.6}{0.18}{}
    \picfilledcircle{0 -0.6}{0.18}{}
  }
}
&
\diag{7.5mm}{5.5}{4}{
  {\piclinedash{0.05 0.05}{0.01}
   \picline{-0.5 3.4}{5.5 3.4}
   \picline{-0.5 0.6}{5.5 0.6}
   \picline{2 2.8}{5.5 2.8}
   \picline{2 1.2}{5.5 1.2}
   \picline{2 2.0}{-0.5 2.0}
  }
  \picfilledcircle{5.5 2}{1.0}{}
  {\picgraycol{1}\piccirclearc{5.5 2}{1.0}{-90 90}}
  \picfilledcircle{2 2}{2}{}
  \picputtext{2 0.9}{$N$}
  \picputtext{2 3.1}{$M$}
  \picputtext{1.5 3.6}{$L$} 
  \picputtext{5.5 2.7}{$L'$} 
  \picputtext{-0.3 2.2}{$y$} 
  \picputtext{-0.3 2.8}{$A$}
  \picputtext{-0.3 1.2}{$C$}
  \picputtext{4.24 2.1}{$B$}
  \pictranslate{2 2}{
   \picscale{-1 1}{
    {\piclinedash{0.05 0.05}{0.01}
     \picline{0 0.6}{2 33 polar}
     \picline{0 0.6}{2 162 polar}
     \picline{0 -0.6}{2 -33 polar}
     \picline{0 -0.6}{2 198 polar}
     \picline{-2 0}{0 0}
     \picline{0 0}{2 18 polar}
     \picline{0 0}{2 -18 polar}
    }
    \picfilledcircle{0 0.0}{0.18}{}
    \picfilledcircle{0 0.6}{0.18}{}
    \picfilledcircle{0 -0.6}{0.18}{}
   }
  }
}
\\
\ry{1.6em}\kern-2mm(a) & (b) & (c)
\end{array}
\end{eqn}

Let us first deal with (a). Since $y$ is inadequate, $y=2$,
and then w.l.o.g. we have $x,z\ge 2$ by connectivity.
Then, however, it is easy to see that $B(D)$ has two
cycles, which share a single loop, but are not on opposite
sides of it. Thus $B(D)$ is of type A, and we are done.

Case (c) is similarly simple. Again $y=2$ by inadequacy. Then
we observe that $B(D)$ has a single separating loop $X$. It
has again at most two multiple connections inside, each
intertwined with at most two multiple connections outside
$X$. This implies that $IG(B)$ has at most one cycle. There
is also a modification of this case, in which, say, $M$
connects $B,C$ rather than $A,B$, and is parallel to $N$.
This is handled by the argument for case (a).

Case (b) is handled as case (c).
\qed

\begin{corr}\label{cor10.4}
In case B', there exists at most one multiple edge
between empty outer cycle loops.
\end{corr}

\proof As before, using lemma \reference{lmoz}, lemma
\reference{ghp}, and $\chi(IG(A))=-1$. \qed

\begin{lemma}
In case B', there is no pair of distinct attachment
legs touching the same region.
\end{lemma}

\proof Let $P$ be an attachment to $L$, on either side of which
there are legs touching from the opposite side of $L$.

First, again $P$ cannot be inside $L$. If so, and $L$ had only
one attachment outside, $P$ would intertwine at most one
connection, and $IG(A)$ has no cycle. If $L$ had both other
attachments outside, $D$ will be composite.

So $P$ is outside $L$. Again because $\gm=4$ and $D$ is
prime, $P$ has valence 2. Recall that legs inside $L$
should touch $A$ and $B$ in \eqref{B'}. Since they
can touch at most 3 regions, this means that $P$ must be
attached within $A$ or $B$, and legs inside $L$ touch only
$A$, $B$ and the region $C$ between $L$ and $P$. So by a
mutation we can assume w.l.o.g. that $P$ is attached within $B$.
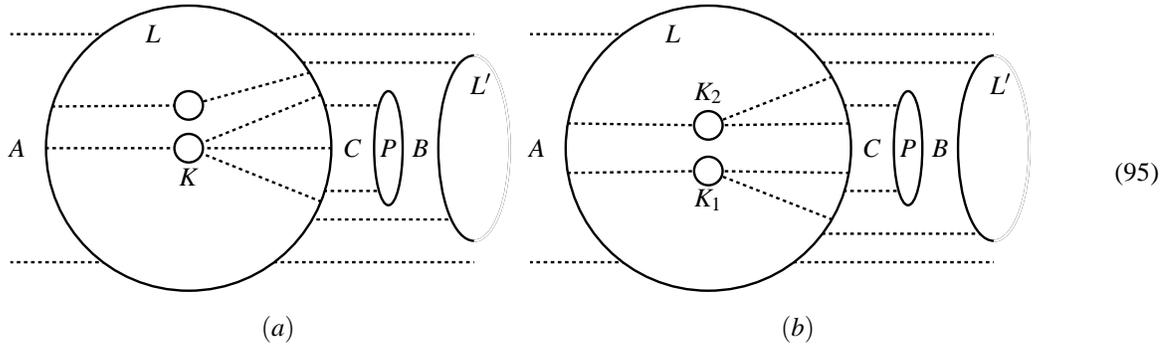
\begin{eqn}\label{B''1}
\begin{array}{c@{\qquad}c}
\diag{9.5mm}{6.5}{4}{
  {\piclinedash{0.05 0.05}{0.01}
   \picline{-0.5 3.6}{6 3.6}
   \picline{-0.5 0.4}{6 0.4}
   \picline{2 2.6}{4.8 2.6}
   \picline{2 1.4}{4.8 1.4}
   \picline{2 3.2}{6 3.2}
   \picline{2 1.0}{6 1.0}
  }
  \picfilledellipse{6 2}{0.5 1.3}{}
  \picfilledellipse{4.8 2}{0.2 0.8}{$P$}
  {\picgraycol{1}\picellipsearc{6 2}{0.5 1.3}{-90 90}}
  % \piccirclearc{6 2}{1.0}{90 270}
  \picfilledcircle{2 2}{2}{}
  \picputtext{1.5 3.6}{$L$} 
  \picputtext{6.1 2.9}{$L'$} 
  \picputtext{5.25 2}{$B$}
  \picputtext{4.3 2}{$C$}
  \picputtext{-.4 2}{$A$}
  \picputtext{2 1.6}{$K$}
  \pictranslate{2 2}{
   {\piclinedash{0.05 0.05}{0.01}
    \picline{2 180 polar}{0 0}
    \picline{2 0 polar}{0 0}
    \picline{2 22 polar}{0 0}
    \picline{2 -22 polar}{0 0}
    \picline{2 163 polar}{0 0.6}
    \picline{2 32 polar}{0 0.6}
   }
   \picfilledcircle{0 0}{0.2}{}
   \picfilledcircle{0 0.6}{0.2}{}
  }
} &
\diag{9.5mm}{6.5}{4}{
  {\piclinedash{0.05 0.05}{0.01}
   \picline{-0.5 3.6}{6 3.6}
   \picline{-0.5 0.4}{6 0.4}
   \picline{2 2.6}{4.8 2.6}
   \picline{2 1.4}{4.8 1.4}
   \picline{2 3.2}{6 3.2}
   \picline{2 0.8}{6 0.8}
  }
  \picfilledellipse{6 2}{0.5 1.3}{}
  \picfilledellipse{4.8 2}{0.2 0.8}{$P$}
  {\picgraycol{1}\picellipsearc{6 2}{0.5 1.3}{-90 90}}
  % \piccirclearc{6 2}{1.0}{90 270}
  \picfilledcircle{2 2}{2}{}
  \picputtext{1.5 3.6}{$L$} 
  \picputtext{6.1 2.9}{$L'$} 
  \picputtext{5.25 2}{$B$}
  \picputtext{4.3 2}{$C$}
  \picputtext{-.4 2}{$A$}
  \picputtext{2 1.3}{$K_1$}
  \picputtext{2 2.8}{$K_2$}
  \pictranslate{2 2}{
   {\piclinedash{0.05 0.05}{0.01}
    \picline{2 190 polar}{0 -0.32}
    \picline{2 -10 polar}{0 -0.32}
    \picline{2 -30 polar}{0 -0.32}
    \picline{2 170 polar}{0 0.32}
    \picline{2 10 polar}{0 0.32}
    \picline{2 30 polar}{0 0.32}
   }
   \picfilledcircle{0 -0.32}{0.2}{}
   \picfilledcircle{0 0.32}{0.2}{}
  }
} \\
\ry{1.6em}(a) & (b)
\end{array}
\end{eqn}
Clearly, given $P$ is outside, it is the only outside attached
loop. Now $P$ is inadequate, unless there is an attachment
(of depth 2) inside $K$ in case (a) of \eqref{B''1}, or both
$K_{1,2}$ in case (b). But then after 3 attachments $IG(A)$ will
have no cycle. \qed

In particular, we assured now that in case B' too, all
attachments have valence 2 or 3.

\begin{lemma}
In case B', there is no attachment outside $L$ of valence 3.
\end{lemma}

\proof Assume there were such an attachment $P$. The three regions
identified by legs inside $L$ must be then, w.l.o.g. up to mutation,
$A$ of \eqref{B'}, and the two regions between $P$ and $L$, whereby
$P$ is attached inside $B$ of \eqref{B'}. Up to symmetries, the
only option we have, avoiding inadequate loops, is
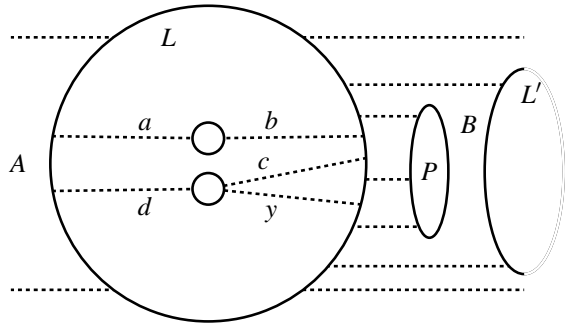
\[
\diag{10.5mm}{6.5}{4}{
  {\piclinedash{0.05 0.05}{0.01}
   \picline{-0.5 3.6}{6 3.6}
   \picline{-0.5 0.4}{6 0.4}
   \picline{2 2.6}{4.8 2.6}
   \picline{2 1.8}{4.8 1.8}
   \picline{2 1.2}{4.8 1.2}
   \picline{2 3.0}{6 3.0}
   \picline{2 0.7}{6 0.7}
  }
  \picfilledellipse{6 1.9}{0.5 1.3}{}
  \picfilledellipse{4.8 1.9}{0.24 0.84}{$P$}
  {\picgraycol{1}\picellipsearc{6 1.9}{0.5 1.3}{-90 90}}
  % \piccirclearc{6 2}{1.0}{90 270}
  \picfilledcircle{2 2}{2}{}
  \picputtext{1.5 3.6}{$L$} 
  \picputtext{6.1 2.9}{$L'$} 
  \picputtext{5.3 2.5}{$B$}
  \picputtext{-.4 2}{$A$}
  \picputtext{2.8 1.4}{$y$}
  \picputtext{1.2 2.5}{$a$}
  \picputtext{2.8 2.55}{$b$}
  \picputtext{2.7 2.0}{$c$}
  \picputtext{1.2 1.46}{$d$}
  \pictranslate{2 2}{
   {\piclinedash{0.05 0.05}{0.01}
    \picline{2 190 polar}{0 -0.32}
    \picline{2 2 polar}{0 -0.32}
    \picline{2 -15 polar}{0 -0.32}
    \picline{2 170 polar}{0 0.32}
    \picline{2 10 polar}{0 0.32}
   }
   \picfilledcircle{0 -0.32}{0.2}{}
   \picfilledcircle{0 0.32}{0.2}{}
  }
} \\
\]
Again $y$ is inadequate, so $y=2$. Then $B(D)$ has a separating loop
$X$, with at most two multiple connections inside. One of those,
the one dual to $y$, is intertwined with at most two connections
outside $X$. So $IG(B)$ has no two cycles. \qed

With this all work until the end of \S\reference{S10.5} is extended
to case $B'$. We continue treating $B$ and $B'$ simultaneously.

Using lemma \ref{lB1}, \ref{lmZ'B} and \ref{lm3v}, we are left
with the following choices to perform 3 attachments \eqref{attach}
in case B, so as to create 2 cycles in $IG(A)$, have $\gm=4$, and
avoid inadequate loops. The region outside $L$ touched by $y$ in (a)
and (b) of \eqref{CsB} may \em{a priori} be arbitrary.

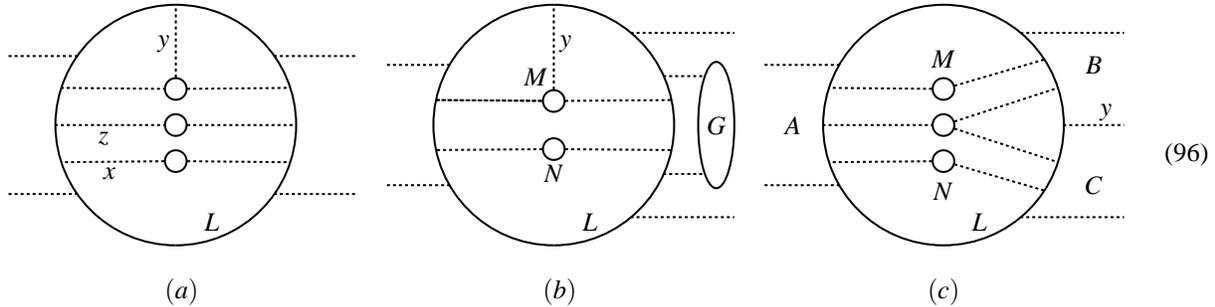
\begin{eqn}\label{CsB}
\begin{array}{ccc}
\diag{8mm}{5.8}{4}{
  \pictranslate{2.8 2}{
    \picmultigraphics[S]{2}{-1 1}{
      \piccirclearc{0 0}{2}{-90 90}
    }
    {\piclinedash{0.05 0.05}{0.01}
     \picline{2 162 polar}{0 0.6}
     \picline{2 180 polar}{0 -0.0}
     \picline{0 0.6}{2 18 polar}
     \picline{0 0.6}{2 90 polar}
     % \picline{0 -0.6}{2 -90 polar}
     \picline{0 -0.0}{2 -0 polar}
     \picline{0 -0.6}{2 -18 polar}
     % \picline{0 -0.6}{0 -0.1}
     \picline{0 -0.6}{2 198 polar}
     % \picline{0 -1.1}{2 210 polar}
     % \picline{0 -1.1}{2 -28 polar}
    }
    \picfilledcircle{0 0.6}{0.18}{}
    \picfilledcircle{0 -0.0}{0.18}{}
    \picfilledcircle{0 -0.6}{0.18}{}
    % \picfilledcircle{0 -1.1}{0.18}{}
    {\piclinedash{0.05 0.05}{0.01}
     % \picline{2 173 polar}{2 173 polar x p -3 x}
     % \picline{2 190 polar}{2 190 polar x p -3 x}
     \picline{2 145 polar}{2 145 polar x p -2.8 x}
     \picline{2 215 polar}{2 215 polar x p -2.8 x}
     % \picline{2 31 polar}{2 31 polar x p 3 x}
     % \picline{2 -35 polar}{2 -35 polar x p 3 x}
     \picline{2 -35 polar}{2 -35 polar x p 3 x}
     \picline{2 35 polar}{2 35 polar x p 3 x}
    }
    % \picfilledellipse{3 -0.08}{0.3 1.25}{}
    \picputtext{-1.2 -0.2}{$z$}
    \picputtext{-1.1 -0.8}{$x$}
    \picputtext{0.6 -1.62}{$L$}
    \picputtext{-0.2 1.4}{$y$}
    % \picputtext{-0.2 -1.4}{$a$}
    % \picputtext{-0.35 0.79}{$M$}
    % \picputtext{-2.9 0.45}{$d$}
    % \picputtext{-2.9 -.5}{$c$}
  }
}
&
\diag{8mm}{5.8}{4}{
  \pictranslate{2.8 2}{
    \picmultigraphics[S]{2}{-1 1}{
      \piccirclearc{0 0}{2}{-90 90}
    }
    {\piclinedash{0.05 0.05}{0.01}
     \picline{2 168 polar}{0 0.4}
     \picline{0 0.4}{2 12 polar}
     \picline{0 0.4}{2 90 polar}
     \picline{0 0.4}{2 168 polar}
     \picline{0 -0.4}{2 -12 polar}
     \picline{0 -0.4}{2 192 polar}
    }
    \picfilledcircle{0 0.4}{0.18}{}
    % \picfilledcircle{0 -0.0}{0.18}{}
    \picfilledcircle{0 -0.4}{0.18}{}
    % \picfilledcircle{0 -1.1}{0.18}{}
    {\piclinedash{0.05 0.05}{0.01}
     % \picline{2 173 polar}{2 173 polar x p -3 x}
     % \picline{2 190 polar}{2 190 polar x p -3 x}
     \picline{2 150 polar}{2 150 polar x p -2.8 x}
     \picline{2 210 polar}{2 210 polar x p -2.8 x}
     \picline{2 24 polar}{2 24 polar x p 2.7 x}
     \picline{2 -24 polar}{2 -24 polar x p 2.7 x}
     \picline{2 -50 polar}{2 -50 polar x p 3 x}
     \picline{2 50 polar}{2 50 polar x p 3 x}
    }
    \picfilledellipse{2.7 0.}{0.3 1.05}{$G$}
    % \picputtext{-1.2 -0.2}{$z$}
    % \picputtext{-1.1 -0.7}{$x$}
    \picputtext{0.6 -1.62}{$L$}
    \picputtext{0.2 1.4}{$y$}
    % \picputtext{-0.2 -1.4}{$a$}
    % \picputtext{-0.35 0.79}{$M$}
    % \picputtext{-2.9 0.45}{$d$}
    % \picputtext{-2.9 -.5}{$c$}
    \picputtext{-0.3 0.8}{$M$}
    \picputtext{0 -0.8}{$N$}
  }
} & 
\diag{8mm}{6}{4}{
  \pictranslate{3 2}{
    \picmultigraphics[S]{2}{-1 1}{
      \piccirclearc{0 0}{2}{-90 90}
    }
    {\piclinedash{0.05 0.05}{0.01}
     \picline{0 0.6}{2 33 polar}
     \picline{0 0.6}{2 162 polar}
     \picline{0 -0.6}{2 -33 polar}
     \picline{0 -0.6}{2 198 polar}
     \picline{-2 0}{0 0}
     \picline{0 0}{2 18 polar}
     \picline{0 0}{2 -18 polar}
    }
    \picfilledcircle{0 0.0}{0.18}{}
    \picfilledcircle{0 0.6}{0.18}{}
    \picfilledcircle{0 -0.6}{0.18}{}
    {\piclinedash{0.05 0.05}{0.01}
     % \picline{2 173 polar}{2 173 polar x p -3 x}
     % \picline{2 190 polar}{2 190 polar x p -3 x}
     \picline{2 150 polar}{2 150 polar x p -3 x}
     \picline{2 210 polar}{2 210 polar x p -3 x}
     \picline{2 0 polar}{2 0 polar x p 3 x}
     % \picline{2 -24 polar}{2 -24 polar x p 2.7 x}
     \picline{2 -50 polar}{2 -50 polar x p 3 x}
     \picline{2 50 polar}{2 50 polar x p 3 x}
    }
    % \picfilledellipse{2.7 0.}{0.3 1.05}{}
    % \picputtext{-1.2 -0.2}{$z$}
    % \picputtext{-1.1 -0.7}{$x$}
    \picputtext{0.6 -1.62}{$L$}
    \picputtext{2.7 0.25}{$y$}
    % \picputtext{-0.2 -1.4}{$a$}
    % \picputtext{-0.35 0.79}{$M$}
    % \picputtext{-2.9 0.45}{$d$}
    % \picputtext{-2.9 -.5}{$c$}
     \picputtext{0 1.1}{$M$}
     \picputtext{0 -1.1}{$N$}
     \picputtext{-2.5 0}{$A$}
     \picputtext{2.5 1.0}{$B$}
     \picputtext{2.5 -1.0}{$C$}
  }
} \\
\ry{1.9em}(a) & (b) & (c)
\end{array}
\end{eqn}

Now we can treat case B' uniformly, since we know that we must
connect the regions $A,B$ in \eqref{B'}. The modification to
\eqref{CsB} is the that traces going out from $L$ on the right
are supposed to connect to both $L'$ and $M$ in \eqref{B'},
with the innermost two (excluding $y$ in \eqref{CsB} (c))
connecting to $L'$. (So regions $B,C$ in \eqref{CsB} (c)
lie between $L$ and $L'$. Only in case (b) we could place the
right legs of $M,N$ so that $G$ is intertwined with only one
of them, but then $G$ would be inadequate.)

Cases (a) and (c) are ruled out as in the proof of lemma \ref{ghp}.
% Let us first deal with (a). Since $y$ is inadequate, $y=2$,
% and then w.l.o.g. we have $x,z\ge 2$ by connectivity.
% Then, however, it is easy to see that $B(D)$ has two
% cycles, which share a single loop, but are not on opposite
% sides of it. Thus $B(D)$ is of type A, and we are done.
% 
% Case (c) is similarly simple. Again $y=2$ by inadequacy. Then
% we observe that $B(D)$ has a single separating loop $x$. It
% has again at most two multiple connections inside, each
% intertwined with at most two multiple connections outside
% $x$. This implies that $IG(B)$ has at most one cycle. There
% is also a modification of this case, in which, say, $M$
% connects $A,C$ rather than $A,B$, and is parallel to $N$.
% This is handled by the argument for case (a).
% 
Thus, for the rest of this subsection, we consider case (b) in
\eqref{CsB}. We draw in \eqref{Bb} the $A$- and $B$-state (we
waived, though, on drawing the traces in $B(D)$, which are obvious).
For the $A$-state, our understanding is that $e.e'$, as well
as $x,x'$, connect $L$ to the same outer cycle or triangle loop.
All drawn edges have non-zero multiplicity. To keep the case
general, we allow then that loops end from outside on the
segments $C,D$ of $L$, as well as on the lower one $E$ between
the endpoints of $x'$ and $e'$.

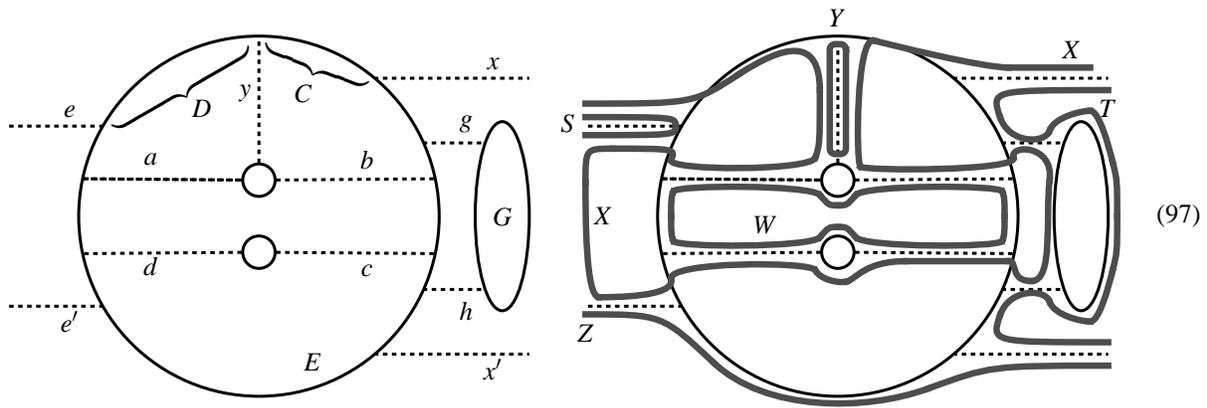
\begin{eqn}\label{Bb}
\diag{12mm}{5.8}{4}{
  \pictranslate{2.8 2}{
    \picmultigraphics[S]{2}{-1 1}{
      \piccirclearc{0 0}{2}{-90 90}
    }
    {\piclinedash{0.05 0.05}{0.01}
     \picline{2 168 polar}{0 0.4}
     \picline{0 0.4}{2 12 polar}
     \picline{0 0.4}{2 90 polar}
     \picline{0 0.4}{2 168 polar}
     \picline{0 -0.4}{2 -12 polar}
     \picline{0 -0.4}{2 192 polar}
    }
    \picfilledcircle{0 0.4}{0.18}{}
    % \picfilledcircle{0 -0.0}{0.18}{}
    \picfilledcircle{0 -0.4}{0.18}{}
    % \picfilledcircle{0 -1.1}{0.18}{}
    {\piclinedash{0.05 0.05}{0.01}
     % \picline{2 173 polar}{2 173 polar x p -3 x}
     % \picline{2 190 polar}{2 190 polar x p -3 x}
     \picline{2 150 polar}{2 150 polar x p -2.8 x}
     \picline{2 210 polar}{2 210 polar x p -2.8 x}
     \picline{2 24 polar}{2 24 polar x p 2.7 x}
     \picline{2 -24 polar}{2 -24 polar x p 2.7 x}
     \picline{2 -50 polar}{2 -50 polar x p 3 x}
     \picline{2 50 polar}{2 50 polar x p 3 x}
    }
    \picfilledellipse{2.7 0.}{0.3 1.05}{$G$}
    % \picputtext{-1.2 -0.2}{$z$}
    \picputtext{2.6 1.7}{$x$}
    \picputtext{2.6 -1.7}{$x'$}
    \picputtext{-2.1 1.17}{$e$}
    \picputtext{-2.1 -1.17}{$e'$}
    \picputtext{0.6 -1.62}{$E$}
    \pictranslate{1.9 108 polar}{
      \picrotate{-60}{%{
	\picputtext[l]{0 0}{$\ \left.\ry{3.4em}\right\}$}
      }
    }
    \pictranslate{2.01 61.5 polar}{
      \picrotate{-111}{%{
	\picputtext[l]{0 0}{$\ \left.\ry{2.3em}\right\}$}
      }
    }
    \picputtext{-0.15 1.4}{$y$}
    \picputtext{1.2 0.63}{$b$}
    \picputtext{-1.2 0.63}{$a$}
    \picputtext{-1.2 -0.55}{$d$}
    \picputtext{1.2 -.6}{$c$}
    \picputtext{2.3 1.05}{$g$}
    \picputtext{2.3 -1.05}{$h$}
    \picputtext{1.35 118 polar}{$D$}
    \picputtext{1.45 70 polar}{$C$}
  }
} \qquad
\diag{12mm}{5.8}{4}{
  \pictranslate{2.8 2}{
    \picmultigraphics[S]{2}{-1 1}{
      \piccirclearc{0 0}{2}{-90 90}
    }
    {\piclinedash{0.05 0.05}{0.01}
     \picline{2 168 polar}{0 0.4}
     \picline{0 0.4}{2 12 polar}
     \picline{0 0.4}{2 90 polar}
     \picline{0 0.4}{2 168 polar}
     \picline{0 -0.4}{2 -12 polar}
     \picline{0 -0.4}{2 192 polar}
    }
    \picfilledcircle{0 0.4}{0.18}{}
    % \picfilledcircle{0 -0.0}{0.18}{}
    \picfilledcircle{0 -0.4}{0.18}{}
    % \picfilledcircle{0 -1.1}{0.18}{}
    {\piclinedash{0.05 0.05}{0.01}
     % \picline{2 173 polar}{2 173 polar x p -3 x}
     % \picline{2 190 polar}{2 190 polar x p -3 x}
     \picline{2 150 polar}{2 150 polar x p -2.8 x}
     \picline{2 210 polar}{2 210 polar x p -2.8 x}
     \picline{2 24 polar}{2 24 polar x p 2.7 x}
     \picline{2 -24 polar}{2 -24 polar x p 2.7 x}
     \picline{2 -50 polar}{2 -50 polar x p 3 x}
     \picline{2 50 polar}{2 50 polar x p 3 x}
    }
    \picfilledellipse{2.7 0.}{0.3 1.05}{}
    % \picputtext{-1.2 -0.2}{$z$}
    % \picputtext{-1.1 -0.7}{$x$}
    % \picputtext{0.6 -1.62}{$L$}
    % \picputtext{-0.2 1.4}{$y$}
    % \picputtext{-0.2 -1.4}{$a$}
    % \picputtext{-0.35 0.79}{$M$}
    % \picputtext{-2.9 0.45}{$d$}
    % \picputtext{-2.9 -.5}{$c$}
    { \piclinewidth{14}
      \picgraycol{0.3}
      \picPSgraphics{2 setlinecap}
      \opencurvepath{-2.8 1.1}{-1.8 1.1}{-1.8 0.9}{-2.8 0.9}{}
      \opencurvepath{-2.5 0.75}{-1.8 0.75}{-1.9 0.50}{-0.3 0.50}
	{-0.2 0.8}{-0.2 0.8}{-0.2 2}{-0.6 1.9}{-1.6 1.25}{-2.8 1.25}
	{}
      % \opencurvepath{-2.8 -1.1}{-1.8 -1.1}{-1.9 -0.5}{-0.3 -.5}
      \opencurvepath{-2.5 0.75}{-2.8 0.75}{-2.8 0}{-2.8 -0.9}
        {-2.5 -0.9}{-1.8 -0.9}{-1.9 -0.5}{-0.3 -.5}
	{-0.2 -0.7}{0 -0.75}{0.2 -0.7}{0.3 -.5}{1.2 -0.5}{1.8 -0.5}
	{2 -0.5}{1.9 -0.72}{2.4 -0.72}{2.3 0}{2.4 0.8}{1.9 0.8}
	{2 0.5}{0.3 0.5}{0.2 0.6}{0.2 2}{0.6 1.9}{1.6 1.65}{2.2 1.65}
	{2.8 1.65}{}
      \picovalbox{0 1.3}{0.18 1.22}{0.06}{}
      \picmultigraphics[S]{2}{1 -1}{
        \picmultigraphics[S]{2}{-1 1}{
	  \opencurvepath{0 0.12}{0.1 0.12}{0.3 .35}{1.8 .35}{1.85 0.2}
	  {1.85 0}{}
	}
	\opencurvepath{3 1.4}{1.7 1.4}{1.9 0.8}{2.3 0.8}{2.4 1.1}
	  {2.7 1.17}{2.9 1.17}{3.0 0.8}{3.1 0.4}{3.1 0}{}
      }
      \opencurvepath{-2.8 -1.09}{-1.95  -1.09}{-0.8 -2.17}{0.3 -2.17}
	{1.7 -1.66}{3 -1.66}{}
    }
    \picputtext{-2.6 0}{$X$}
    \picputtext{2.6 1.85}{$X$}
    \picputtext{0 2.2}{$Y$}
    \picputtext{-2.8 -1.3}{$Z$}
    \picputtext{-0.8 -0.1}{$W$}
    \picputtext{3 1.2}{$T$}
    \picputtext{-3 1.03}{$S$}
  }
} 
\end{eqn}\\[-1mm]

% Again there are not multiple connections outside $L$ except the two
% drawn in \eqref{Bb}. 

Note that we drew the $B$-state assuming there are no additional
edges on $C,D,E$ in the $A$-state. For $E$ adding edges will not
alter the following argument. For $C,D$ we will just see that
in fact there must not be any further edges.

Namely, we have again only one separating loop $X$ in $B(D)$.
Again, $X$ closes on the left as shown because $e,e'$ connect
to the same loop in $A(D)$. $X$ contains in its interior only
two multiple connections, to $W$ and $Y$, the loop obtained from
the inadequate edge $y=2$ in $A(D)$. 

Moreover, $X$ has at most 3 multiple connections in its exterior,
which can be intertwined with $Y$ or $W$, the connections to $Z$,
$T$ and $S$. In particular, $y$ does not end on a segment of $L$
between edges that connect $L$ from outside to the same loop;
otherwise $Y$ intertwines at most two connections outside
$X$, and $IG(B)$ has no two cycles. Moreover, in order $T$
to intertwine with both $Y,W$, we see that $C$ must be
empty, $x=1$, and at least one of $g,h$, say $g$, is $1$.
We call this condition \em{$y$ has distance 2 from $G$}.
Also, in order $S$ to intertwine, we must have either that
$D$ is empty and $e=2$, or $e=1$, and $D$ contains exactly
one endpoint of a trace that connects $L$ to a loop different
from the one $e$ does.

Note that we could also have something like
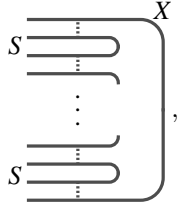
\[
\diag{6mm}{3}{4}{
  { \piclinewidth{14}
    \picgraycol{0.3}
    \picPSgraphics{2 setlinecap}
    \picarcangle{0 4}{3 4}{3 2}{0.5}
    \picarcangleto{3 0}{0 0}{0.5}
    \pictranslate{0 2}{ \picmultigraphics[S]{2}{1 -1}{
      \picarcangle{0 1.6}{2 1.6}{2 1.4}{0.2}
      \picarcangleto{2 1.2}{0 1.2}{0.2}
      \picarcangle{0 0.8}{2 0.8}{2 0.6}{0.2}
      {\piclinedash{0.05 0.05}{0.01}
       \picPSgraphics{0 setlinecap}
       \picline{1.1 1.6}{1.1 2}
       \picline{1.1 1.2}{1.1 0.8}
      }
    }}
  }
  \pictranslate{0 2}{
    \picmultigraphics{3}{0 -0.3}{
      \picputtext{1.1 0.3}{.}
    }
    \picputtext{3 2.2}{$X$}
    \picputtext{-0.3 1.44}{$S$}
    \picputtext{-0.3 -1.44}{$S$}
  }
}\ ,
\]
but this is excluded because $D$ is prime, and $X$ is the only
separating loop in $B(D)$, so $S$ is empty on the side opposite to
$X$.

Note that the loop $T$ closes on the right, as $X$ closes on the
left. We have then a cycle in $B(D)$ that contains $X,T,Z$.

Consider now the $B$ state of the diagram $\tl D$ of (B) in
\eqref{csA} before the attachments \eqref{attach} are made to $L$.
We draw  the connections of $x,x'$ and $e,e'$ as simple edges, i.e.
the diagram $(\tl D)'$ (see definition \reference{cfiq}).
\begin{eqn}\label{uipo}
\def\mycv{
  \opencurvepath{-0.75 -0.1}{-0.4 -0.1}{-0.2 -0.37}{0 -0.37}{}
}
\def\mycvv{\picmultigraphics[S]{2}{-1 1}{\mycv}}
\diag{1cm}{6}{2.2}{
 \pictranslate{0.25 0}{
  {\piclinedash{0.05 0.05}{0.01}
   \picline{-0.25 1.5}{5.75 1.5}
   \picline{2 1.5}{2.7 0.2}
   \picline{3.5 1.5}{2.7 0.2}
  }
  \picfilledcircle{2.7 0.2}{0.3}{$C$}
  \picfilledcircle{0.5 1.5}{0.3}{}
  \picfilledcircle{5.0 1.5}{0.3}{}
  \picfilledcircle{2.0 1.5}{0.3}{$A$}
  \picfilledcircle{3.5 1.5}{0.3}{$B$}
  \picputtext{2.7 1.1}{$Q$}
  \picputtext{2.7 1.9}{$E$}
  \picputtext{1.2 2.1}{$U$}
  \picputtext{1.2 1.0}{$V$}
  { \piclinewidth{14}
    \picgraycol{0.3}
    \picPSgraphics{0 setlinecap}
    \pictranslate{0.5 1.5}{
      \picscale{1 -1}{
        \picmultigraphics{4}{1.5 0}{\mycvv}
      }
      \picmultigraphics{2}{4.5 0}{\mycvv}
      \pictranslate{1.5 0}{\mycv}
      \pictranslate{3 0}{\picscale{-1 1}{\mycv}}
    }
    \opencurvepath{2 1.13}{2.1 1.13}{2.4 0.6}{2.3 0.3}{2.36 -0.25}
      {3.04 -0.25}{3.1 0.3}{3.0 0.6}{3.4 1.13}{3.5 1.13}{}
    \pictranslate{2.72 1.07}{
      \picmultigraphics[rt]{3}{120}{
	 \picmultigraphics[S]{2}{-1 1}{
           \opencurvepath{0.42 0.20}{0.4 0.25}{0.31 0.35}{0 0.35}{}
	 }
      }
    }
    {\piclinedash{0.05 0.05}{0.01}
     \picmultigraphics{3}{1.5 0}{\picline{1.25 1.4}{1.25 1.6}}
     \picline{2.2 0.9}{2.4 1.0}
     \picline{3.1 1.0}{3.25 0.9}
    }
  }
 }
}
\end{eqn}
This diagram $\tl D'$  already has a cycle in the $B$-state, call
it $R=(Q,V,U)$. Now the effect of the attachments \eqref{attach}
on $B(D)$ is that $T$ and one of $X$ or $S$ is added,
and $X$ obtains 2 new loops in its
interior, $Y,W$. Again by connectivity, since $y=2$, we have some
of $a,b,c,d>1$, so that $B(D)$ contains a cycle in the interior
of $X$. 

There is also the option of case B' that $A,B$ are connected by
traces $e''$ below $C$ in \eqref{uipo}. We will soon rule out
this option, but for the moment just observe that the arguments
below apply to case B' in the same way as for case B.

If now $y$ connects to a region outside $L$ different from $U,V$, then
we can use lemmas \reference{ghp'} and  \ref{ghp} and apply lemma
\ref{UV}: we have at least $7$ loops in $B(D)$ (namely, $X$, and at
least 3 loops on its either side, taking into account
$\max(a,b,c,d)>1$),
and the loops $U,V$ of \eqref{uipo} remain non-separating, with
a multiple edge. Thus the attachments inside $L$ are done so as
$y$ to touch $U$ or $V$.

Next note that under the attachments inside $L$, the cycle $R$ in
\eqref{uipo} also persists in $D$, unless legs inside $L$ touch
between $x$ and $x'$ in
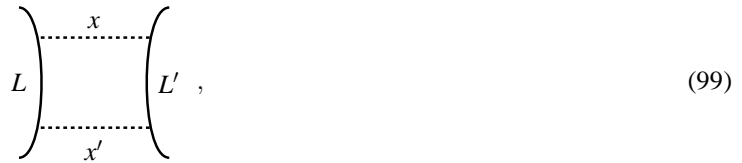
\begin{eqn}\label{LL'}
\diag{1cm}{2}{2}{
  \piccurve{2 0}{1.6 0}{1.6 2}{2 2}
  \piccurve{0 2}{0.4 2}{0.4 0}{0 0}
  \picclip{
    \picline{0 0}{2 0}
    \piccurveto{1.6 0}{1.6 2}{2 2}
    \piclineto{0 2}
    \piccurveto{0.4 2}{0.4 0}{0 0}
  }{
    \piclinedash{0.05 0.05}{0.01}
    \picline{0 0.4}{2 0.4}
    \picline{0 1.6}{2 1.6}
  }
  \picputtext{1 1.8}{$x$}
  \picputtext{0 1}{$L$}
  \picputtext{2 1}{$L'$}
  \picputtext{1 0.1}{$x'$}
}\quad,
\end{eqn}
and $x$, $x'$ lie in the common connection $E$ of the outer cycle and
triangle in $A(D)$ in \eqref{uipo}. Otherwise the cycle $X,T,Z$ in
\eqref{Bb} does not contain edges dual to those in $E$, and so it
is not part of $R$.
Then, counting the cycle inside $X$, the addition of $T$ in $B(D)$
under the attachments in $A(D)$ to $L$ creates a third cycle in $B(D)$.

\begin{figure}
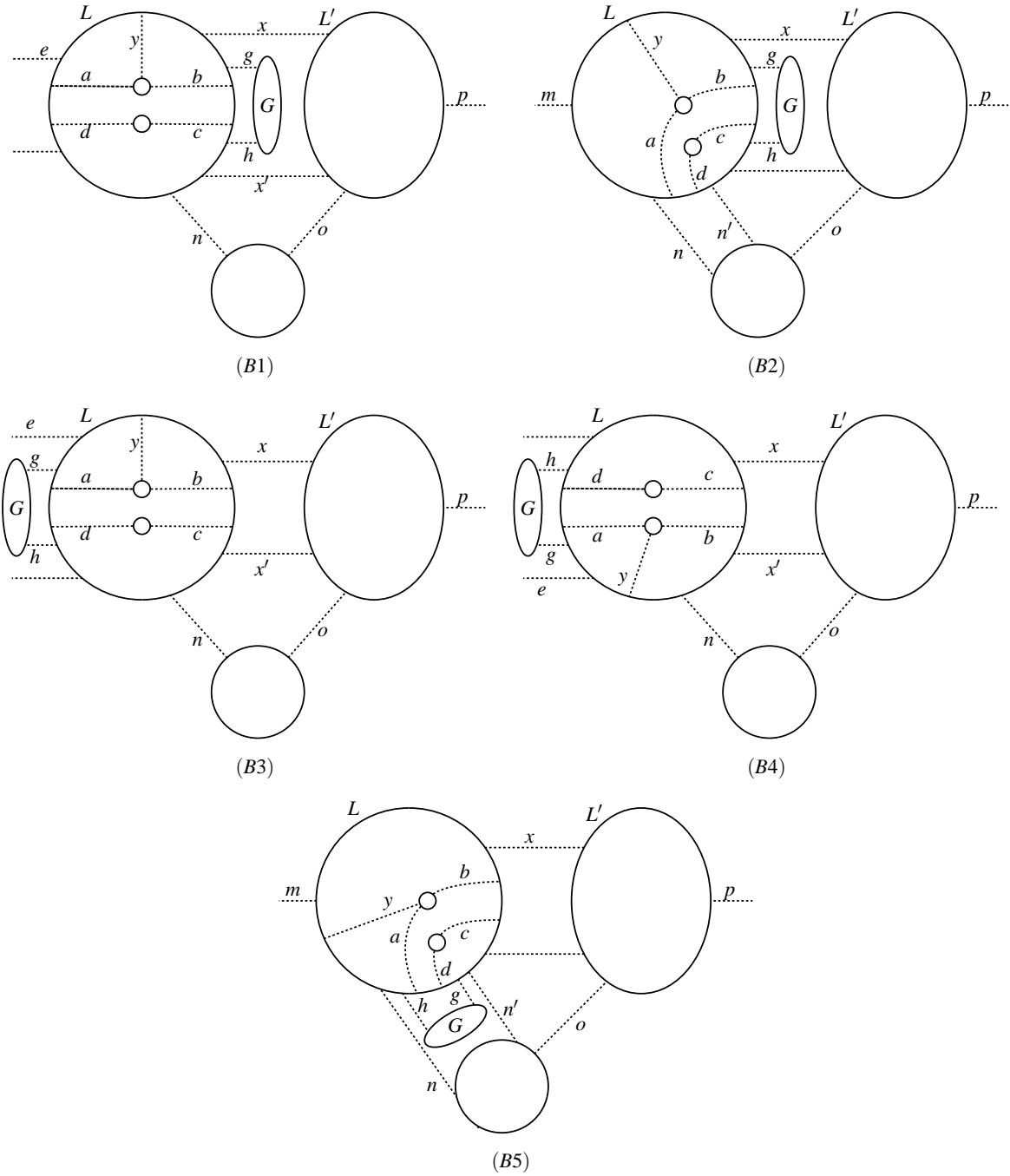

\[
{\small
\begin{array}{cc} 
\diag{7.1mm}{10.5}{7}{
  \pictranslate{2.8 5}{
    \picmultigraphics[S]{2}{-1 1}{
      \piccirclearc{0 0}{2}{-90 90}
    }
    {\piclinedash{0.05 0.05}{0.01}
     \picline{2 168 polar}{0 0.4}
     \picline{0 0.4}{2 12 polar}
     \picline{0 0.4}{2 168 polar}
     \picline{0 0.4}{2 90 polar}
     \picline{0 -0.4}{2 -12 polar}
     \picline{0 -0.4}{2 192 polar}
     \picline{0 7.4 x}{0 6 x}
     \picline{2.5 -4}{0 6 x}
     \picline{2.5 -4}{2 -72 polar}
    }
    \picfilledcircle{0 0.4}{0.18}{}
    % \picfilledcircle{0 -0.0}{0.18}{}
    \picfilledcircle{0 -0.4}{0.18}{}
    % \picfilledcircle{0 -1.1}{0.18}{}
    {\piclinedash{0.05 0.05}{0.01}
     % \picline{2 173 polar}{2 173 polar x p -3 x}
     % \picline{2 190 polar}{2 190 polar x p -3 x}
     \picline{2 150 polar}{2 150 polar x p -2.8 x}
     \picline{2 210 polar}{2 210 polar x p -2.8 x}
     \picline{2 24 polar}{2 24 polar x p 2.7 x}
     \picline{2 -24 polar}{2 -24 polar x p 2.7 x}
     \picline{2 -50 polar}{2 -50 polar x p 5 x}
     \picline{2 50 polar}{2 50 polar x p 5 x}
    }
    \picfilledellipse{2.7 0.}{0.3 1.05}{$G$}
    \picfilledellipse{5 0}{2 1.5 x}{}
    \picfilledcircle{2.5 -4}{1.0}{}
    % \picputtext{-1.2 -0.2}{$z$}
    \picputtext{2.6 1.7}{$x$}
    \picputtext{2.6 -1.7}{$x'$}
    \picputtext{-2.1 1.17}{$e$}
    \picputtext{-0.15 1.4}{$y$}
    \picputtext{1.2 0.63}{$b$}
    \picputtext{-1.2 0.63}{$a$}
    \picputtext{-1.2 -0.55}{$d$}
    \picputtext{1.2 -.6}{$c$}
    \picputtext{2.3 1.05}{$g$}
    \picputtext{2.3 -1.05}{$h$}
    \picputtext{6.9 0.2}{$p$}
    \picputtext{3.9 -2.7}{$o$}
    \picputtext{1.2 -2.9}{$n$}
    \picputtext{-1.2 2.0}{$L$}
    \picputtext{4 1.9}{$L'$}
  }
} & \diag{7.1mm}{10.0}{7}{
  \pictranslate{2.8 5}{
    \picmultigraphics[S]{2}{-1 1}{
      \piccirclearc{0 0}{2}{-90 90}
    }
    {\piclinedash{0.05 0.05}{0.01}
     \picline{2 114 polar}{0.4 0.0}
     \piccurve{0.4 0.0}{-0.05 -0.3}{-0.3 -1.1}{2 -85 polar}
     \piccurve{0.4 0.0}{0.5 0.3}{1.3 0.4}{2 12 polar}
     \piccurve{0.6 -0.9}{0.7 -0.5}{1.3 -0.4}{2 -12 polar}
     \piccurve{0.6 -0.9}{0.5 -1.1}{0.5 -1.4}{2 -69 polar}
     \picline{0 7.4 x}{0 6 x}
     \picline{2.0 -4}{0 6 x}
     \picline{2.0 -3.1}{2 -60 polar}
     \picline{2.0 -4.9}{2 -98 polar}
    }
    \picfilledcircle{0.4 0.0}{0.18}{}
    % \picfilledcircle{0 -0.0}{0.18}{}
    \picfilledcircle{0.6 -0.9}{0.18}{}
    % \picfilledcircle{0 -1.1}{0.18}{}
    {\piclinedash{0.05 0.05}{0.01}
     % \picline{2 173 polar}{2 173 polar x p -3 x}
     % \picline{2 190 polar}{2 190 polar x p -3 x}
     \picline{2 180 polar}{2 180 polar x p -2.8 x}
     \picline{2 24 polar}{2 24 polar x p 2.7 x}
     \picline{2 -24 polar}{2 -24 polar x p 2.7 x}
     \picline{2 -45 polar}{2 -45 polar x p 5 x}
     \picline{2 45 polar}{2 45 polar x p 5 x}
    }
    \picfilledellipse{2.7 0.}{0.3 1.05}{$G$}
    \picfilledellipse{5 0}{2 1.5 x}{}
    \picfilledcircle{2.0 -4}{1.0}{}
    % \picputtext{-1.2 -0.2}{$z$}
    \picputtext{2.6 1.6}{$x$}
    \picputtext{-0.15 1.4}{$y$}
    \picputtext{1.2 0.63}{$b$}
    \picputtext{1.2 -0.73}{$c$}
    \picputtext{0.8 -1.45}{$d$}
    \picputtext{-0.3 -.8}{$a$}
    \picputtext{2.3 1.05}{$g$}
    \picputtext{2.3 -1.05}{$h$}
    % \picputtext{-0.6 0.2}{$y$}
    \picputtext{-2.5 0.2}{$m$}
    \picputtext{6.9 0.2}{$p$}
    \picputtext{0.3 -3.2}{$n$}
    \picputtext{1.3 -2.7}{$n'$}
    \picputtext{3.7 -2.7}{$o$}
    \picputtext{-1.2 2.0}{$L$}
    \picputtext{4 1.9}{$L'$}
  }
}\\
\ry{1.6em}(B1) & (B2) \\[5mm]
\diag{7.1mm}{10.5}{7}{
  \pictranslate{2.8 5}{
    \picmultigraphics[S]{2}{-1 1}{
      \piccirclearc{0 0}{2}{-90 90}
    }
    {\piclinedash{0.05 0.05}{0.01}
     \picline{2 168 polar}{0 0.4}
     \picline{0 0.4}{2 12 polar}
     \picline{0 0.4}{2 168 polar}
     \picline{0 0.4}{2 90 polar}
     \picline{0 -0.4}{2 -12 polar}
     \picline{0 -0.4}{2 192 polar}
     \picline{0 7.4 x}{0 6 x}
     \picline{2.5 -4}{0 6 x}
     \picline{2.5 -4}{2 -72 polar}
    }
    \picfilledcircle{0 0.4}{0.18}{}
    % \picfilledcircle{0 -0.0}{0.18}{}
    \picfilledcircle{0 -0.4}{0.18}{}
    % \picfilledcircle{0 -1.1}{0.18}{}
    {\piclinedash{0.05 0.05}{0.01}
     % \picline{2 173 polar}{2 173 polar x p -3 x}
     % \picline{2 190 polar}{2 190 polar x p -3 x}
     \picline{2 130 polar}{2 130 polar x p -2.8 x}
     \picline{2 230 polar}{2 230 polar x p -2.8 x}
     \picline{2 24 180 + polar}{2 24 180 + polar x p -2.7 x}
     \picline{2 -24 180 + polar}{2 -24 180 + polar x p -2.7 x}
     \picline{2 -30 polar}{2 -30 polar x p 5 x}
     \picline{2 30 polar}{2 30 polar x p 5 x}
    }
    \picfilledellipse{-2.7 0.}{0.3 1.05}{$G$}
    \picfilledellipse{5 0}{2 1.5 x}{}
    \picfilledcircle{2.5 -4}{1.0}{}
    % \picputtext{-1.2 -0.2}{$z$}
    \picputtext{2.6 1.3}{$x$}
    \picputtext{2.6 -1.3}{$x'$}
    \picputtext{-2.4 1.8}{$e$}
    \picputtext{-0.15 1.4}{$y$}
    \picputtext{1.2 0.63}{$b$}
    \picputtext{-1.2 0.63}{$a$}
    \picputtext{-1.2 -0.55}{$d$}
    \picputtext{1.2 -.6}{$c$}
    \picputtext{-2.3 1.05}{$g$}
    \picputtext{-2.3 -1.05}{$h$}
    \picputtext{6.9 0.2}{$p$}
    \picputtext{3.9 -2.7}{$o$}
    \picputtext{1.2 -2.9}{$n$}
    \picputtext{-1.2 2.0}{$L$}
    \picputtext{4 1.9}{$L'$}
  }
} & \diag{7.1mm}{10.5}{7}{
  \pictranslate{2.8 5}{
    \picmultigraphics[S]{2}{-1 1}{
      \piccirclearc{0 0}{2}{-90 90}
    }
    {\piclinedash{0.05 0.05}{0.01}
     \picline{2 168 polar}{0 0.4}
     \picline{0 0.4}{2 12 polar}
     \picline{0 0.4}{2 168 polar}
     \picline{0 -0.4}{2 255 polar}
     \picline{0 -0.4}{2 -12 polar}
     \picline{0 -0.4}{2 192 polar}
     \picline{0 7.4 x}{0 6 x}
     \picline{2.5 -4}{0 6 x}
     \picline{2.5 -4}{2 -72 polar}
    }
    \picfilledcircle{0 0.4}{0.18}{}
    % \picfilledcircle{0 -0.0}{0.18}{}
    \picfilledcircle{0 -0.4}{0.18}{}
    % \picfilledcircle{0 -1.1}{0.18}{}
    {\piclinedash{0.05 0.05}{0.01}
     % \picline{2 173 polar}{2 173 polar x p -3 x}
     % \picline{2 190 polar}{2 190 polar x p -3 x}
     \picline{2 130 polar}{2 130 polar x p -2.8 x}
     \picline{2 230 polar}{2 230 polar x p -2.8 x}
     \picline{2 24 180 + polar}{2 24 180 + polar x p -2.7 x}
     \picline{2 -24 180 + polar}{2 -24 180 + polar x p -2.7 x}
     \picline{2 -30 polar}{2 -30 polar x p 5 x}
     \picline{2 30 polar}{2 30 polar x p 5 x}
    }
    \picfilledellipse{-2.7 0.}{0.3 1.05}{$G$}
    \picfilledellipse{5 0}{2 1.5 x}{}
    \picfilledcircle{2.5 -4}{1.0}{}
    % \picputtext{-1.2 -0.2}{$z$}
    \picputtext{2.6 1.3}{$x$}
    \picputtext{2.6 -1.3}{$x'$}
    \picputtext{-2.4 -1.8}{$e$}
    \picputtext{-0.65 -1.5}{$y$}
    \picputtext{1.2 0.69}{$c$}
    \picputtext{-1.2 0.69}{$d$}
    \picputtext{-1.2 -0.65}{$a$}
    \picputtext{1.2 -.65}{$b$}
    \picputtext{-2.2 1.05}{$h$}
    \picputtext{-2.2 -1.05}{$g$}
    \picputtext{6.9 0.2}{$p$}
    \picputtext{3.9 -2.7}{$o$}
    \picputtext{1.2 -2.9}{$n$}
    \picputtext{-1.2 2.0}{$L$}
    \picputtext{4 1.9}{$L'$}
  }
}\\
\ry{1.6em}(B3) & (B4) \\[4mm]
\multicolumn{2}{c}{
\diag{7.1mm}{10.0}{7}{
  \pictranslate{2.8 5}{
    \picmultigraphics[S]{2}{-1 1}{
      \piccirclearc{0 0}{2}{-90 90}
    }
    {\piclinedash{0.05 0.05}{0.01}
     \picline{2 204 polar}{0.4 0.0}
     \piccurve{0.4 0.0}{-0.05 -0.3}{-0.3 -1.1}{2 -85 polar}
     \piccurve{0.4 0.0}{0.5 0.3}{1.3 0.4}{2 12 polar}
     \piccurve{0.6 -0.9}{0.7 -0.5}{1.3 -0.4}{2 -12 polar}
     \piccurve{0.6 -0.9}{0.5 -1.1}{0.5 -1.4}{2 -69 polar}
     \picline{0 7.4 x}{0 6 x}
     \picline{2.0 -4}{0 6 x}
     \picline{2.5 -3.3}{2 -50 polar}
     \picline{1.5 -4.9}{2 -108 polar}
    }
    \picfilledcircle{0.4 0.0}{0.18}{}
    % \picfilledcircle{0 -0.0}{0.18}{}
    \picfilledcircle{0.6 -0.9}{0.18}{}
    % \picfilledcircle{0 -1.1}{0.18}{}
    {\piclinedash{0.05 0.05}{0.01}
     % \picline{2 173 polar}{2 173 polar x p -3 x}
     % \picline{2 190 polar}{2 190 polar x p -3 x}
     \picline{2 180 polar}{2 180 polar x p -2.8 x}
     \picline{2 -58 polar}{1.5 -2.4}
     \picline{2 -94 polar}{0.6 -3.1}
     \picline{2 -35 polar}{2 -35 polar x p 5 x}
     \picline{2 35 polar}{2 35 polar x p 5 x}
    }
    \pictranslate{1.0 -2.7}{\picrotate{120}{
    \picfilledellipse{0 0.}{0.3 .75}{}
    }}
    \picfilledellipse{5 0}{2 1.5 x}{}
    \picfilledcircle{2.0 -4}{1.0}{}
    % \picputtext{-1.2 -0.2}{$z$}
    \picputtext{2.6 1.35}{$x$}
    \picputtext{-0.45 -.02}{$y$}
    \picputtext{1.2 0.63}{$b$}
    \picputtext{1.2 -0.73}{$c$}
    \picputtext{0.8 -1.45}{$d$}
    \picputtext{-0.3 -.8}{$a$}
    \picputtext{1.0 -2.05}{$g$}
    \picputtext{.3 -2.25}{$h$}
    \picputtext{1.0 -2.67}{$G$}
    % \picputtext{-0.6 0.2}{$y$}
    \picputtext{-2.5 0.2}{$m$}
    \picputtext{6.9 0.2}{$p$}
    \picputtext{0.5 -4.0}{$n$}
    \picputtext{2.2 -2.3}{$n'$}
    \picputtext{3.7 -2.7}{$o$}
    \picputtext{-1.2 2.0}{$L$}
    \picputtext{4 1.9}{$L'$}
  }
}
}
\\
\multicolumn{2}{c}{\ry{1.6em}(B5)}
\end{array}
}
\]
\caption{Patterns in cases B1 to B5.\label{B12}}
\end{figure}

Thus we can assume that $x$ in \eqref{LL'} belongs to the
connection $E$, that is, (w.l.o.g. up to symmetry) $L=A$ and
$L'=B$ for $A,B$ in \eqref{uipo}. With the previously observed
properties on the edges $e$, $x$ and $y$, we are left with 5
options for the attachments \eqref{attach}, see figure \ref{B12}.

The patterns of figure \ref{B12} we have now are very explicit.
Additionally, we have strong restrictions on edge multiplicities:
% The conditions on edge multiplicities in figure \ref{B12} are 
\begin{itemize}
\item $m=n=o=p=1$ (where $m=1$ applies only if $m$ is drawn),
  and all not drawn outer cycle connections are simple (because
  $\chi(IG(A))=-1$),
\item w.l.o.g. $g=1$ (by the arguments on $C,D$ in \eqref{Bb});
  for (B1), $e=2$, $x=1$, for (B2) $x=1$, for (B3) $e=1$,
  $x=2$, for (B4) $e=x'=1$, for (B5) $x=1$ (for the same reasons),
\item $y=2$ (by inadequacy),
\item one of $a,b,c,d$ is even, but not both $(a,b)$ or both
  $(c,d)$ are even (by connectivity).
\end{itemize}

\begin{clc}
% aqv23tst8b.C trdg_B[12345].in
The computer checked these cases (though still to many to
write one by one here) within a few seconds. With the explanation
of \S\reference{Scmp}, the sequences our computation finds are those
of figure \reference{Fig20'} for $i=1,2,8,9$. This concludes case B.
\end{clc}

In case B', i.e. if the above named potential traces $e''$ exist,
we must have a trace ending from inside $L$ between $n$ and $n'$
in (B2) or (B5), following the remark after \eqref{B'}. If we attach
$G$ somewhere else than between $n$ and $n'$, the three regions
identified by traces inside $L$ are determined, and it is impossible
to install $y$ so that it connects $U$ or $V$ in \eqref{uipo}. Thus
we should try to attach $G$ between $n$ and $n'$. However, with
$e''$ there, it is then impossible to install $y$ to have distance
$2$ from $G$ and to connect $U$ or $V$. With this case B' is finished,
too.

\subsection{Case C}

First we try again to clear the way for applying computations, by
establishing the premise of lemma \ref{lmcmp}. We need the below
lemma \ref{lC2}.

\begin{lemma}\label{lC2}
There is at most one multiple edge between empty outer cycle loops.
\end{lemma}

\proof 
This lemma is immediately clear from lemmas \ref{lmm3'} and \ref{lC1}.
\qed

Next we eliminate loops of depth $\ge 2$ completely.

\begin{lemma}\label{ldp2}
There are no depth-2 loops.
\end{lemma}

\proof With corollary \reference{crdp}, we are left only with
depth-2 loops inside triangle loops.

If there is a depth-2 loop of valence 3, then we can argue
using \eqref{Fig17} and $\gm=4$.

Next assume one of the triangle loops $M$ has exactly one
valence-2 loop $N$ inside. Then we have by flatness the
case (a) of \eqref{dp2q}.
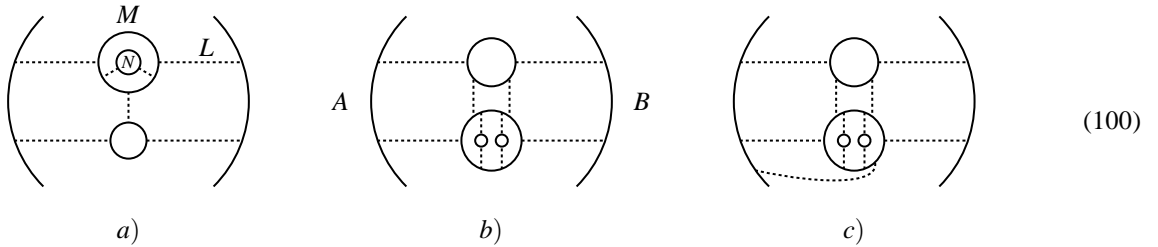
\begin{eqn}\label{dp2q}
\begin{array}{c@{\qquad\ }c@{\qquad\ }c}
\diag{8mm}{5}{3.5}{
  \pictranslate{2.5 1.8}{
    \picputtext{1.3 0.9}{$L$}
    \piccirclearc{0 0}{2}{135 225}
    \piccirclearc{0 0}{2}{-45 45}
    {\piclinedash{0.05 0.05}{0.01}
     \picline{2 180 19 - polar}{2 19 polar}
     \picline{2 180 19 + polar}{2 -19 polar}
     \picline{2 180 19 - polar x p 0 x}{2 180 19 + polar x p 0 x}
    }
    \picfilledcircle{2 180 19 - polar x p 0 x}{0.5}{}
    \picfilledcircle{2 180 19 + polar x p 0 x}{0.3}{}
    \pictranslate{2 180 19 - polar x p 0 x}{
      {\piclinedash{0.05 0.05}{0.01}
       \picline{0 0}{0.5 -30 polar}
       \picline{0 0}{0.5 -150 polar}
      }
      {\scriptsize
      \picfilledcircle{0 0}{0.2}{$N$}
      }
      % {\tiny\picputtext{0.3 50 polar}{$N$}}
      \picputtext{0.8 90 polar}{$M$}
    }
  }
} &
\diag{8mm}{5}{3.5}{
  \pictranslate{2.5 1.8}{
    \piccirclearc{0 0}{2}{135 225}
    \piccirclearc{0 0}{2}{-45 45}
    {\piclinedash{0.05 0.05}{0.01}
     \picline{2 180 19 - polar}{2 19 polar}
     \picline{2 180 19 + polar}{2 -19 polar}
     \picline{2 180 19 - polar x p 0.3 x}{2 180 19 + polar x p 0.3 x}
     \picline{2 180 19 - polar x p -0.3 x}{2 180 19 + polar x p -0.3 x}
    }
    \picfilledcircle{2 180 19 - polar x p 0 x}{0.4}{}
    \picfilledcircle{2 180 19 + polar x p 0 x}{0.5}{}
    \pictranslate{2 180 19 + polar x p 0 x}{
      {\piclinedash{0.05 0.05}{0.01}
       \picline{0.5 -70 polar}{0.5 70 polar}
       \picline{0.5 -110 polar}{0.5 110 polar}
      }
      \picfilledcircle{0.5 -70 polar p 0}{0.1}{}
      \picfilledcircle{0.5 -70 polar p N 0}{0.1}{}
    }
    \picputtext{-2.5 0}{$A$}
    \picputtext{2.5 0}{$B$}
  }
} &
\diag{8mm}{5}{3.5}{
  \pictranslate{2.5 1.8}{
    \piccirclearc{0 0}{2}{135 225}
    \piccirclearc{0 0}{2}{-45 45}
    {\piclinedash{0.05 0.05}{0.01}
     \picline{2 180 19 - polar}{2 19 polar}
     \picline{2 180 19 + polar}{2 -19 polar}
     \picline{2 180 19 - polar x p 0.3 x}{2 180 19 + polar x p 0.3 x}
     \picline{2 180 19 - polar x p -0.3 x}{2 180 19 + polar x p -0.3 x}
     \piccurve{2 215 polar}{-0.5 -0.5 0.9 -}{0.6 -0.5 0.9 -}{0.3 0 0.9 -}
    }
    \picfilledcircle{2 180 19 - polar x p 0 x}{0.4}{}
    \picfilledcircle{2 180 19 + polar x p 0 x}{0.5}{}
    \pictranslate{2 180 19 + polar x p 0 x}{
      {\piclinedash{0.05 0.05}{0.01}
       \picline{0.5 -70 polar}{0.5 70 polar}
       \picline{0.5 -110 polar}{0.5 110 polar}
      }
      \picfilledcircle{0.5 -70 polar p 0}{0.1}{}
      \picfilledcircle{0.5 -70 polar p N 0}{0.1}{}
    }
  }
} \\
\ry{1.0em}a) & b) & c)
\end{array}
\end{eqn}
Note that the other triangle loop must be empty, otherwise
we cannot create 2 cycles in $IG(A)$. 

Now if $N$ is not supposed to be inadequate, we need traces
from outside $L$ like on the left below\\[5mm]
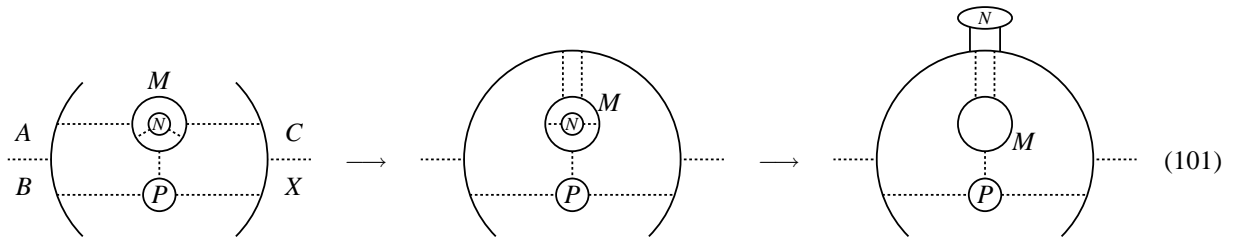
\begin{eqn}\label{89}
\diag{7.2mm}{5.6}{3.5}{
  \pictranslate{2.8 1.8}{
    \piccirclearc{0 0}{2}{135 225}
    \piccirclearc{0 0}{2}{-45 45}
    {\piclinedash{0.05 0.05}{0.01}
     \picline{2 180 19 - polar}{2 19 polar}
     \picline{2 180 19 + polar}{2 -19 polar}
     \picline{2 180 19 - polar x p 0 x}{2 180 19 + polar x p 0 x}
     \picline{-2.8 0}{-2. 0}
     \picline{2.8 0}{2. 0}
    }
    \picputtext{-2.5 0.5}{$A$}
    \picputtext{-2.5 -0.5}{$B$}
    \picputtext{2.5 0.5}{$C$}
    \picputtext{2.5 -0.5}{$X$}
    \picfilledcircle{2 180 19 - polar x p 0 x}{0.5}{}
    \picfilledcircle{2 180 19 + polar x p 0 x}{0.3}{$P$}
    \pictranslate{2 180 19 - polar x p 0 x}{
      {\piclinedash{0.05 0.05}{0.01}
       \picline{0 0}{0.5 -30 polar}
       \picline{0 0}{0.5 -150 polar}
      }
      {\scriptsize
      \picfilledcircle{0 0}{0.2}{$N$}
      }
      % {\tiny\picputtext{0.3 50 polar}{$N$}}
      \picputtext{0.8 90 polar}{$M$}
    }
  }
}
\quad\lra\quad
\diag{7.2mm}{5.6}{3.5}{
  \pictranslate{2.8 1.8}{
    \piccirclearc{0 0}{2}{-45 225}
    {\piclinedash{0.05 0.05}{0.01}
     \picline{2 85 polar}{2 85 polar p 0.6}
     \picline{2 95 polar}{2 95 polar p 0.6}
     \picline{2 180 19 + polar}{2 -19 polar}
     \picline{2 180 19 - polar x p 0 x}{2 180 19 + polar x p 0 x}
     \picline{-2.8 0}{-2. 0}
     \picline{2.8 0}{2. 0}
    }
    \picfilledcircle{2 180 19 - polar x p 0 x}{0.5}{}
    \picfilledcircle{2 180 19 + polar x p 0 x}{0.3}{$P$}
    \pictranslate{2 180 19 - polar x p 0 x}{
      {\piclinedash{0.05 0.05}{0.01}
       \picline{0 0}{0.5 -0 polar}
       \picline{0 0}{0.5 -180 polar}
      }
      {\scriptsize
      \picfilledcircle{0 0}{0.2}{$N$}
      }
      % {\tiny\picputtext{0.3 50 polar}{$N$}}
      \picputtext{0.8 30 polar}{$M$}
    }
  }
}
\quad\lra\quad
\diag{7.2mm}{5.6}{3.5}{
  \pictranslate{2.8 1.8}{
    \piccirclearc{0 0}{2}{-45 225}
    {\piclinedash{0.05 0.05}{0.01}
     \picline{2 85 polar}{2 85 polar p 0.6}
     \picline{2 95 polar}{2 95 polar p 0.6}
     \picline{2 180 19 + polar}{2 -19 polar}
     \picline{2 180 19 - polar x p 0 x}{2 180 19 + polar x p 0 x}
     \picline{-2.8 0}{-2. 0}
     \picline{2.8 0}{2. 0}
    }
    \picfilledcircle{2 180 19 - polar x p 0 x}{0.5}{}
    \picfilledcircle{2 180 19 + polar x p 0 x}{0.3}{$P$}
    \picline{2 82 polar}{2 82 polar p 2.6}
    \picline{2 98 polar}{2 98 polar p 2.6}
    \picfilledellipse{2.6 0 x}{0.5 0.2}{\scriptsize $N$}
    \picputtext{0.8 25 polar}{$M$}
  }
}
\end{eqn}
Now if there are other traces from outside $L$ both between $A$ and $C$,
and between $B$ and $X$, then $\gm\ge 6$. Thus w.l.o.g. (up to flyping
$N$ into $P$), $A=C$. Then the legs of $M$ become parallel,
and we can mutate $N$, as shown above, into a loop of depth $0$,
in contradiction to flatness.

The last option is to have two valence-2 loops attached inside
a triangle loop $M$. Diagrams (b) and (c) and \eqref{dp2q} display
the only options we have, if one is to create two cycles in
$IG(A)$ and avoid inadequate loops. Because of $\gm=4$, the regions
$A,B$ outside $L$ cannot be separated by traces as on the left of
\eqref{89}. On the other hand, the regions $A$ and $B$ must be
between edges that connect $L$ to the same neighboring loop on
the outer cycle, in order to have 2 cycles in $IG(A)$. 

The $B$-state for diagram (b) of \eqref{dp2q}
looks then as diagram (a) in \eqref{dp2q'}.
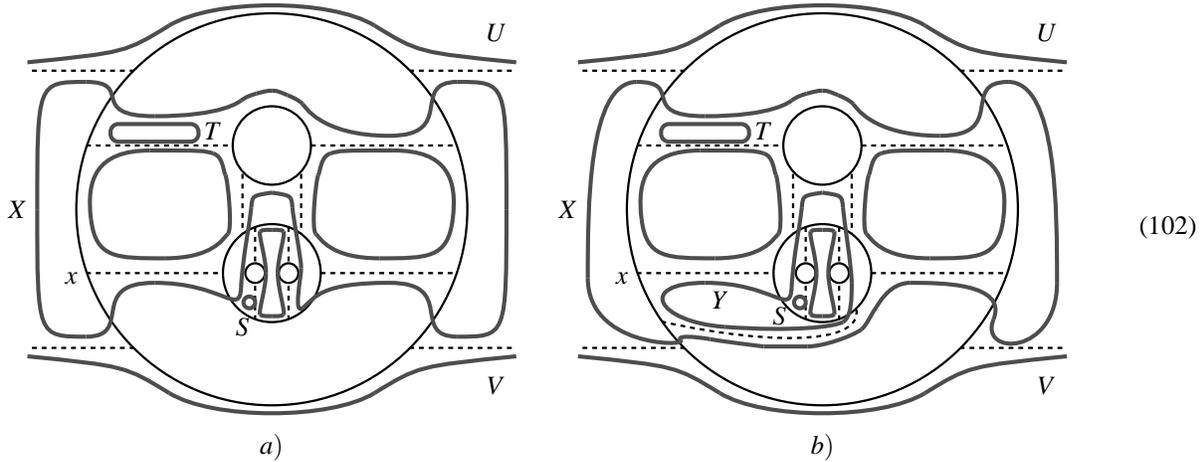
\begin{eqn}\label{dp2q'}
\begin{array}{c@{\qquad\ }c}
\diag{1.3cm}{5}{4.4}{
  \piclinewidth{44}
  \pictranslate{2.5 2.2}{
    \piccircle{0 0}{2}{}
    {\piclinedash{0.05 0.05}{0.01}
     \picline{2 45 polar}{2 45 polar x p 2.5 x}
     \picline{2 -45 polar}{2 -45 polar x p 2.5 x}
     \picline{2 135 polar}{2 135 polar x p -2.5 x}
     \picline{2 -135 polar}{2 -135 polar x p -2.5 x}
     \picline{2 180 19 - polar}{2 19 polar}
     \picline{2 180 19 + polar}{2 -19 polar}
     \picline{2 180 19 - polar x p 0.3 x}{2 180 19 + polar x p 0.3 x}
     \picline{2 180 19 - polar x p -.3 x}{2 180 19 + polar x p -0.3 x}
    }
    \picputtext{-2.05 -0.7}{$x$}
    \picfilledcircle{2 180 19 - polar x p 0 x}{0.4}{}
    \picfilledcircle{2 180 19 + polar x p 0 x}{0.5}{}
    \pictranslate{2 180 19 + polar x p 0 x}{
      {\piclinedash{0.05 0.05}{0.01}
       \picline{0.5 -70 polar}{0.5 70 polar}
       \picline{0.5 -110 polar}{0.5 110 polar}
      }
      \picfilledcircle{0.5 -70 polar p 0}{0.1}{}
      \picfilledcircle{0.5 -70 polar p N 0}{0.1}{}
    }
    \picputtext{-2.6 0}{$X$}
    \picputtext{2.3 1.8}{$U$}
    \picputtext{2.3 -1.8}{$V$}
    \picputtext{-.3 -1.2}{$S$}
    \picputtext{-.6 0.8}{$T$}
    { \piclinewidth{24}
      \picgraycol{0.3}
      \picPSgraphics{0 setlinecap}
      \picmultigraphics[S]{2}{1 -1}{
        \opencurvepath{-2.5 1.5}{-1.5 1.6}{-0.6 2.15}{0.6 2.15}
	  {1.5 1.6}{2.5 1.5}{}
      }
      \piccircle{-0.23 -0.95}{0.06}{}
      \picmultigraphics[S]{2}{-1 1}{
        \opencurvepath{0 -0.21}{0.1 -0.21}{0.14 -0.25}{0.06 -0.5}
	  {0.06 -0.8}{0.14 -1.05}{0.1 -1.09}{0 -1.09}{}
	\curvepath{0.6 0.57}{0.7 0.6}{1.8 0.6}{1.9 0.3}{1.8 -0.5}
	  {1.6 -0.5}{0.6 -0.5}{0.5 -0.3}{0.4 -0.2}{0.45 0.3}{0.5 0.4}{}
      }
      \opencurvepath{0 0.17}{0.2 0.14}{0.22 -0.1}{0.3 -0.5}{0.3 -0.7}
        {0.21 -1.1}{0.44 -0.9}
	{0.5 -0.75}{1.9 -0.75}{1.4 -1.3}{2.4 -1.3}{2.4 0}{}
      \picscale{-1 1}{
	\opencurvepath{0 0.17}{0.2 0.14}{0.22 -0.1}{0.3 -0.5}{0.33 -1.1}
	  {0.5 -0.75}{1.9 -0.75}{1.4 -1.3}{2.4 -1.3}{2.4 0}{}
      }
      \opencurvepath{0 1.21}{0.05 1.21}{0.27 1.12}{0.3 1.09}{0.33 1.06}
	{0.5 0.75}{1.9 0.75}{1.4 1.3}{2.4 1.3}{2.4 0}{}
      \picscale{-1 1}{
	\opencurvepath{0 1.21}{0.05 1.21}{0.27 1.15}{0.3 1.12}{0.33 1.1}
	  {0.5 0.95}{1.9 0.95}{1.4 1.3}{2.4 1.3}{2.4 0}{}
      }
      \picovalbox{-1.2 0.78}{0.9 0.18}{0.08}{}
    }
  }
} &
\diag{1.3cm}{5}{4.4}{
  \piclinewidth{44}
  \pictranslate{2.5 2.2}{
    \piccircle{0 0}{2}{}
    {\piclinedash{0.05 0.05}{0.01}
     \picline{2 45 polar}{2 45 polar x p 2.5 x}
     \picline{2 -45 polar}{2 -45 polar x p 2.5 x}
     \picline{2 135 polar}{2 135 polar x p -2.5 x}
     \picline{2 -135 polar}{2 -135 polar x p -2.5 x}
     \picline{2 180 19 - polar}{2 19 polar}
     \picline{2 180 19 + polar}{2 -19 polar}
     \picline{2 180 19 - polar x p 0.3 x}{2 180 19 + polar x p 0.3 x}
     \picline{2 180 19 - polar x p -0.3 x}{2 180 19 + polar x p -0.3 x}
     \piccurve{2 215 polar}{-0.5 -0.5 0.9 -}{0.6 -0.5 0.9 -}{0.3 0 0.9 -}
    }
    \picputtext{-2.05 -0.7}{$x$}
    \picfilledcircle{2 180 19 - polar x p 0 x}{0.4}{}
    \picfilledcircle{2 180 19 + polar x p 0 x}{0.5}{}
    \pictranslate{2 180 19 + polar x p 0 x}{
      {\piclinedash{0.05 0.05}{0.01}
       \picline{0.5 -70 polar}{0.5 70 polar}
       \picline{0.5 -110 polar}{0.5 110 polar}
      }
      \picfilledcircle{0.5 -70 polar p 0}{0.1}{}
      \picfilledcircle{0.5 -70 polar p N 0}{0.1}{}
    }
    \picputtext{-2.6 0}{$X$}
    \picputtext{2.3 1.8}{$U$}
    \picputtext{2.3 -1.8}{$V$}
    \picputtext{-.44 -1.08}{$S$}
    \picputtext{-1.04 -0.93}{$Y$}
    \picputtext{-.6 0.8}{$T$}
    { \piclinewidth{24}
      \picgraycol{0.3}
      \picPSgraphics{0 setlinecap}
      \picmultigraphics[S]{2}{1 -1}{
        \opencurvepath{-2.5 1.5}{-1.5 1.6}{-0.6 2.15}{0.6 2.15}
	  {1.5 1.6}{2.5 1.5}{}
      }
      \piccircle{-0.23 -0.95}{0.06}{}
      \picmultigraphics[S]{2}{-1 1}{
        \opencurvepath{0 -0.21}{0.1 -0.21}{0.14 -0.25}{0.06 -0.5}
	  {0.06 -0.8}{0.14 -1.05}{0.1 -1.09}{0 -1.09}{}
	\curvepath{0.6 0.57}{0.7 0.6}{1.8 0.6}{1.9 0.3}{1.8 -0.5}
	  {1.6 -0.5}{0.6 -0.5}{0.5 -0.3}{0.4 -0.2}{0.45 0.3}{0.5 0.4}{}
      }
      \picscale{-1 1}{
	\opencurvepath{0 0.17}{0.2 0.14}{0.22 -0.1}{0.3 -0.5}{0.33 -1.1}
	  {0.5 -0.75}{1.9 -0.75}{1.7 -1.15}{0.3 -1.25}{-.12 -1.20}
	  {-.3 -1.1}{-0.3 -0.9}{-0.3 -0.6}{-0.22 -0.1}{-0.2 0.14}
	  {0 0.17}{}
      }
      \opencurvepath{0 1.21}{0.05 1.21}{0.27 1.12}{0.3 1.09}{0.33 1.06}
	{0.5 0.75}{1.9 0.75}{1.4 1.3}{2.4 1.3}{2.4 0}
	{2.4 -1.5}{1.6 -1.5}{1.9 -0.75}{1.2 -0.75}{0.6 -0.75}
	{0.5 -1.15}{0.3 -1.3}{0 -1.4}{-0.3 -1.4}{-0.6 -1.4}
	{}
      \picscale{-1 1}{
	\opencurvepath{0 1.21}{0.05 1.21}{0.27 1.15}{0.3 1.12}{0.33 1.1}
	  {0.5 0.95}{1.9 0.95}{1.4 1.3}{2.4 1.3}{2.4 0}
	  {2.4 -1.4}{1.4 -1.4}{1.5 -1.28}{1.3 -1.3}
	  {1.0 -1.35}{0.6 -1.4}
	  {}
      }
      \picovalbox{-1.2 0.78}{0.9 0.18}{0.08}{}
    }
  }
}
\\
\ry{1.0em}a) & b) 
\end{array}
\end{eqn}
Here again the extra loops $S,T$ result from edges in $A(D)$
which must be multiple for connectivity reasons. ($T$ may be,
instead of its displayed position, also connected by the dual 2
edges to $x=2$.) Then this is ruled out by lemma \reference{UV}.

The case c) of \eqref{dp2q} leads to the diagram b) of \eqref{dp2q'}
(with the same remarks on $S,T$ as \eqref{dp2q'} a)). For this diagram,
$B(D)$ has two separating loops, $X$ and $Y$, but each one has at most
one multiple connection on one side. Thus $IG(B)$ has no cycle. \qed

% \begin{lemma}\label{lmZ'}
% There is no region touched by two different legs of the
% same attachment \eqref{attach} or \eqref{tatt}.
% \end{lemma}
% 
% \begin{lemma}\label{lC3}
% \end{lemma}
% 
% \[
% \]
% but then

\begin{lemma}\label{l45}
Assume $M$ is a loop of an attachment, and $M$ has two edges.
Then, with the exception of the 2 sequences of figure \ref{Fig20'}
for $i=4,10$, neither of them has multiplicity $4$ or $5$.
In particular, the list of options of lemma \reference{l2v}
reduces to $(1,1)$ and $(1,2)$.
\end{lemma}

\proof If there were such an edge, then we would have in $B(D)$
an outer cycle with at most one multiple connection. Thus the
claim then follows from the next lemma. \qed

\begin{lemma}\label{l46}
The outer cycle has at least 2 multiple connections, except the 
sequences of figure \ref{Fig20'} for $i=4,10$.
\end{lemma}

\proof Assume it has at most one multiple connection.
Using lemma \reference{lm3v}, avoiding inadequate loops,
we are left with only few possibilities to create
two cycles in $IG(A)$. Clearly the outer cycle loop
$L$ with non-empty interior must be among those of
the multiple connection. Moreover, we need at least one
attachment outside $L$.

\begin{caselist}
\case
Assume first there is one attachment \eqref{attach}
outside $L$, and one inside.

Then we have two options, up to mutation, which are shown below
as a) and b).
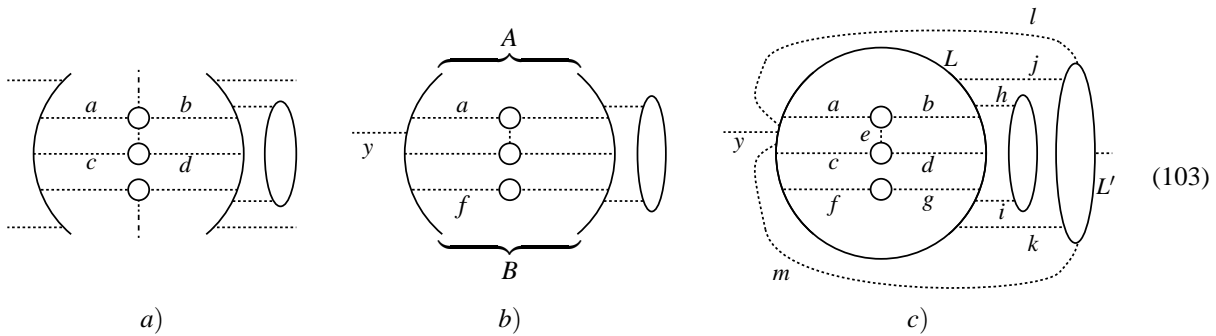
\begin{eqn}\label{91}
\begin{array}{c@{\qquad}c@{\qquad}c}
\diag{7mm}{5.5}{4}{
  \pictranslate{2.5 2}{
  {\piclinedash{0.05 0.05}{0.01}
   \picline{-2.5 1.4}{3 1.4}
   \picline{-2.5 -1.4}{3 -1.4}
   \picline{0 -0.9}{2.7 -0.9}
   \picline{0 0.9}{2.7 0.9}
  } 
  \picfillnostroke{
   \piccircle{0 0}{2}{}
  }
  {\piclinedash{0.05 0.05}{0.01}
   \picline{2 20 polar}{2 180 20 - polar}
   \picline{2 -20 polar}{2 180 -20 - polar}
   \picline{0 1.1}{0 0}
   \picline{0 -2 x}{0 2 x}
   \picline{0 -1.1}{0 2 -20 polar x p 0 x}
   {\small
   \picputtext[d]{0.9 0.8}{$b$}
   \picputtext{0.9 -.2}{$d$}
   \picputtext[d]{-.9 0.8}{$a$}
   \picputtext{-.9 -.2}{$c$}
   }
   {\piclinedash{0.1 0.1}{0.14}
    \picline{0 1.1}{0 1.6}
    \picline{0 -1.1}{0 -1.6}
   }
  }
  \piccirclearc{0 0}{2}{-50 50}
  \piccirclearc{0 0}{2}{130 230}
  \picfilledellipse{2.7 0}{0.3 1}{}
  % \picfilledellipse{-2.7 0}{0.3 1}{}
  \picfilledcircle{2 20 polar x p 0 x}{0.2}{}
  \picfilledcircle{2 -20 polar x p 0 x}{0.2}{}
  \picfilledcircle{0 0}{0.2}{}
 }
}
& 
\diag{7mm}{6}{4.5}{
 \pictranslate{3 2.25}{
  {\piclinedash{0.05 0.05}{0.01}
   \picline{-3 0.4}{0 0.4}
   \picline{0 -0.9}{2.7 -0.9}
   \picline{0 0.9}{2.7 0.9}
  } 
  \picfillnostroke{
   \piccircle{0 0}{2}{}
  }
  {\small
  \picputtext{-.9 0.9}{$a$}
  \picputtext{-2.7 0.1}{$y$}
  }
  {\piclinedash{0.05 0.05}{0.01}
   \picline{2 20 polar}{2 180 20 - polar}
   \picline{2 -20 polar}{2 180 -20 - polar}
   \picline{0 20 polar x p 0 x}{0 0}
   \picline{0 -2 x}{0 2 x}
   \picline{0 0}{2 20 polar x p 0 x}
  }
  \piccirclearc{0 0}{2}{-50 50}
  \piccirclearc{0 0}{2}{130 230}
  \picfilledellipse{2.7 0}{0.27 1.1}{}
  % \picfilledellipse{-2.7 0}{0.3 1}{}
  \picfilledcircle{2 20 polar x p 0 x}{0.2}{}
  \picfilledcircle{2 -20 polar x p 0 x}{0.2}{}
  \picfilledcircle{0 0}{0.2}{}
  \picmultigraphics[S]{2}{1 -1}{
    \picputtext{1.6 -90 polar}{
      $\underbrace{\kern2.7\unitlength}_{}$
    }
  }
  \picputtext{0 -2.22}{$B$}
  \picputtext{0 2.22}{$A$}
  \picputtext[u]{-.9 -0.8}{$f$}
 }
}
&
\diag{7mm}{7.4}{5}{
  \pictranslate{3 2.5}{
  {\piclinedash{0.05 0.05}{0.01}
   \picline{-3 0.4}{0 0.4}
   \picline{0 -0.9}{2.7 -0.9}
   \picline{0 0.9}{2.7 0.9}
   \picline{0 -1.4}{3.7 -1.4}
   \picline{0 1.4}{3.7 1.4}
   \picline{3.7 0}{4.4 0}
   \opencurvepath{2 165 polar}{-2.5 1}{-1.9 2.5}
     {3 2.4}{3.8 1.8}{3.9 0}{}
   \opencurvepath{2 175 polar}{-2.5 -0.1}{-1.9 -2.8}
     {3 -2.7}{3.8 -1.8}{3.9 0}{}
  } 
  \picfill{
   \piccircle{0 0}{2}{}
  }
  {\small
  % \picputtext{-.9 0.9}{$a$}
  \picputtext[d]{-.9 0.8}{$a$}
  \picputtext[d]{0.9 0.8}{$b$}
  \picputtext{-.9 -.2}{$c$}
  \picputtext{0.9 -.2}{$d$}
  \picputtext{-0.3 0.33}{$e$}
  \picputtext[u]{-.9 -0.8}{$f$}
  \picputtext[u]{0.9 -0.8}{$g$}
  \picputtext{2.3 1.13}{$h$}
  \picputtext{2.3 -1.1}{$i$}
  \picputtext{2.9 1.71}{$j$}
  \picputtext{2.9 -1.7}{$k$}
  \picputtext{2.9 2.6}{$l$}
  \picputtext[l]{1.2 1.8}{$L$}
  \picputtext[l]{4.1 -.6}{$L'$}
  \picputtext{-1.9 -2.3}{$m$}
  \picputtext{-2.7 0.2}{$y$}
  }
  {\piclinedash{0.05 0.05}{0.01}
   \picline{2 20 polar}{2 180 20 - polar}
   \picline{2 -20 polar}{2 180 -20 - polar}
   \picline{0 20 polar x p 0 x}{0 0}
   \picline{0 -2 x}{0 2 x}
   \picline{0 0}{2 20 polar x p 0 x}
  }
  \piccirclearc{0 0}{2}{-50 50}
  \piccirclearc{0 0}{2}{130 230}
  \picfilledellipse{2.7 0}{0.27 1.1}{}
  \picfilledellipse{3.7 0}{0.37 1.7}{}
  % \picfilledellipse{-2.7 0}{0.3 1}{}
  \picfilledcircle{2 20 polar x p 0 x}{0.2}{}
  \picfilledcircle{2 -20 polar x p 0 x}{0.2}{}
  \picfilledcircle{0 0}{0.2}{}
 }
}
\\
\ry{1.5em}a) & b) & c)
\end{array}
\end{eqn}
All legs are assumed non-empty, except a few in a). There, among
the two vertical legs, exactly one of the two is supposed to be
there (because of $\gm=4$). Also, some outer cycle edges may be zero.
(It is not implied here that edges going out to the left/right
from $L$ connect to the outer cycle loop on the same side of $L$.
Also, it must be allowed that an outgoing trace stands in fact
for several edges, whose traces have neighbored basepoints on $L$.)

\begin{caselist}
\case In case \eqref{91} a), we must have all of $a,b,c,d$ non-zero
to have 2 cycles in $IG(A)$. But then we can apply lemma \ref{lm1014}.

\case Now consider case b). Let $L'$ resp. $L''$ be the
neighbors of $L$ in the outer cycle. % on the left resp. right. 
We assume that the connection $(L,L')$ is the multiple one, and
$(L,L'')$ is the simple one.
% one (simple) edge that connects $L$ to $L'$, while if there are
% both edges, one connects $L$ to $L''$.
% 
% In case a) one immediately observes that $B(D)$ is of type
% A or B. (In fact, using that one of $a,b,c,d$ is $>1$ by
% connectivity, one finds in $B(D)$ 3 cycles.)
In order to have two cycles in $IG(A)$, on both segments $A$ and $B$
of \eqref{91} b) there must be a trace connecting $L$ to $L'$.

\begin{caselist}
\case If $y$ consists of only one edge, one has first $y=2$
by inadequacy, and also $a=2$ the same way. Now
it is easy to see that $B(D)$ has a single separating loop
$X$, with at most 3 multiple connections inside.

\begin{caselist} 
\case Assume first the leftmost
(i.e. neighbored to $a$ resp. $f$ in \eqref{91} b)) traces outside
$L$ on the segments $A$ and $B$ are both part of the multiple outer
cycle connection $(L,L')$. Then we have diagram a) or b) of \eqref{91'},
with the convention that all edges drawn are non-zero, except possibly
$y$. All 3
multiple connections inside $X$ in $B(D)$ are intertwined each with
at most 2 connections outside $X$, and we can apply lemma \ref{UV}. 

\def\mydg#1{
\diag{7mm}{9.3}{4.5}{
 \pictranslate{4.3 2.25}{
  {\piclinedash{0.05 0.05}{0.01}
   \picline{-3 0.42}{0 0.42}
   \picline{-3 0.27}{0 0.27}
   \picline{0 -0.9}{2.7 -0.9}
   \picline{0 0.9}{2.7 0.9}
  } 
  \picfill{
   \piccircle{0 0}{2}{}
  }
  {\small
  \picputtext{-.9 1.0}{$a$}
  % \picputtext{-.9 -1}{$f$}
  }
  {\piclinedash{0.05 0.05}{0.01}
   \picline{2 158 polar}{2 158 polar x p 0 x}
   \picline{2 162 polar}{2 162 polar x p 0 x}
   \picline{2 20 polar}{2 20 polar x p 0 x}
   \picline{2 -20 polar}{2 180 -20 - polar}
   \picline{0 20 polar x p 0 x}{0 0}
   \picline{0 -2 x}{0 2 x}
   \picline{0 0}{2 20 polar x p 0 x}
  }
  % \piccirclearc{0 0}{2}{-50 50}
  % \piccirclearc{0 0}{2}{130 230}
  \picfilledellipse{2.7 0}{0.27 1.1}{}
  % \picfilledellipse{-2.7 0}{0.3 1}{}
  \picfilledcircle{2 20 polar x p 0 x}{0.2}{}
  \picfilledcircle{2 -20 polar x p 0 x}{0.2}{}
  \picfilledcircle{0 0}{0.2}{}
  % \picmultigraphics[S]{2}{1 -1}{
    % \picputtext{1.6 -90 polar}{
      % $\underbrace{\kern2.7\unitlength}_{}$
    % }
  % }
  % \picputtext{0 -2.22}{$B$}
  \picputtext{0 2.28}{$L$}
  % \picputtext[u]{-.9 -0.8}{$f$}
  #1
 }
}}
\begin{eqn}\label{91'}
\begin{array}{c@{\qquad}c}
% \mydg{
%   {\piclinedash{0.05 0.05}{0.01}
%    \picline{2 135 polar}{2 135 polar x p -3 x}
%    \picline{2 -135 polar}{2 -135 polar x p -3 x}
%    \picline{-3 0}{-4.3 0}
%    % \picmultigraphics[S]{2}{1 -1}{
%      \piccurve{2 55 polar}{2 2}{3.5 2}{4 0}
%    % }
%    \picline{4 0}{5 0}
%   }
%   \picfilledellipse{-3 0}{0.5 1.8}{$L'$}
%   \picfilledcircle{4 0}{0.5}{$L''$}
% } 
% &
% \mydg{
%   {\piclinedash{0.05 0.05}{0.01}
%    \picline{2 135 polar}{2 135 polar x p -3 x}
%    \picline{2 -135 polar}{2 -135 polar x p -3 x}
%    \picline{-3 0}{-4.3 0}
%    \picscale{1 -1}{
%      \piccurve{2 55 polar}{2 2}{3.5 2}{4 0}
%    }
%    \picline{4 0}{5 0}
%   }
%   \picfilledellipse{-3 0}{0.5 1.8}{$L'$}
%   \picfilledcircle{4 0}{0.5}{$L''$}
% } 
% \\
\mydg{
  {\piclinedash{0.05 0.05}{0.01}
   \picline{2 -135 polar}{2 -135 polar x p -3 x}
   \picline{2 135 polar}{2 135 polar x p -3 x}
   \picline{-3 0}{-4.3 0}
   \opencurvepath{2 55 polar}{2 1.8}{3.3 1.6}{3.4 0}{3.2 -2}{2 -2.4}
     {-3.2 -2.4}{-3.2 0}{}
   \picline{4 0}{5 0}
  % \picmultigraphics[S]{2}{1 -1}{
    \piccurve{2 70 polar}{2 2.5}{3.5 2.5}{4.3 0}
  % }
  }
  \picfilledellipse{-3 0}{0.5 1.8}{$L'$}
  \picfilledcircle{4.2 0}{0.5}{$L''$}
  \picputtext{-1 -0.4}{\small $S$}
  \picputtext{1 0.3}{\small $T$}
  \picputtext{2.1 -1.5}{$R$}
  \picputtext{-2.2 1.65}{$x$}
  \picputtext{2.4 -2.29}{$y$}
} 
&
\mydg{
  {\piclinedash{0.05 0.05}{0.01}
   \picline{2 -135 polar}{2 -135 polar x p -3 x}
   \picline{2 135 polar}{2 135 polar x p -3 x}
   \picline{-3 0}{-4.3 0}
   \picscale{1 -1}{
     \piccurve{2 70 polar}{2 2.5}{3.5 2.5}{4.3 0}
     \opencurvepath{2 55 polar}{2 1.8}{3.3 1.6}{3.4 0}{3.2 -2}{2 -2.6}
       {-2 -2.6}{-2.8 -2.2}{-3.0 -1.8}{-3.0 0}{}
   }
   \picline{4 0}{5 0}
  }
  \picfilledellipse{-3 0}{0.5 1.8}{$L'$}
  \picfilledcircle{4.2 0}{0.5}{$L''$}
  \picputtext{-1 -0.4}{\small $S$}
  \picputtext{1 0.3}{\small $T$}
  \picputtext{2.1 1.5}{$R$}
  \picputtext{-2.2 -1.65}{$x$}
  \picputtext{2.4 2.49}{$y$}
} 
\\
\ry{1.5em}a) & b)
% \ry{1.5em}c) & d) \\[5mm]
\end{array}
\end{eqn}

\case
So one of the two leftmost traces outside $L$ on segments $A$
or $B$ of \eqref{91} b) is the simple outer cycle connection
$(L,L'')$, as in diagram a) or b) of \eqref{91'} with $x=0$
(and all other edges non-zero). In either case again a look
at the $B$ state shows that the unique separating loop $X$
has at most 3 multiple connections inside, coming from regions
$S,T$ of $A(D)$, and one, $W$, dual to $a=2$. Each of $S,T,W$
is intertwined with at most 2 connections outside $X$. However,
the pair of intertwined outside connections is not the same
for all of $S,T,W$: the connection of $X$ to the loop (in
$B(D)$) of the region $R$ (in $A(D)$) intertwines one of $S$
and $W$, but not the other. So again $IG(B)$ has no two cycles.
\end{caselist}

\case 
However, if $y$ consists of several edges, the situation is less easy.
Then we have c) of \eqref{91}. By assumption, a single trace from $y$
connects $L$ to the outer cycle loop $L''$ (different from $L'$ on the
right).

All edges drawn in \eqref{91} c) are assumed non-zero, except possibly
$l,m$. It is easy to see that one of $l,m$ must be 0 by lemma 
\reference{UV}, but of course (by assumption) not both $l=m=0$. Also,
$(h,i)$ follow lemma \ref{l2v}, and $a=2$ by inadequacy (independently
on whether $l=0$ or $l\ne 0$). By looking at $B(D)$, we can assume
w.l.o.g. that $b=c=f=g=1$, otherwise $B(D)$ is of type $B$.
Next $e=1$, because otherwise we would need 3 cycles in $IG(A(D))$,
which we cannot create. Furthermore, all outer cycle edges,
except possibly $j,k$ (and $l,m$), have multiplicity one. 

\begin{clc}
% aqv23tst8b.C trdg_prep2.in
This input was again processed by the computer. The result, after
a few seconds, are the sequences of figure \reference{Fig20'} for
$i=4,10$.
\end{clc}

\end{caselist}
\end{caselist}

\case 
Assume next we have two attachments \eqref{attach} outside $L$
(and none inside). By lemma \reference{lm2l}, we have two
options.

\begin{caselist}
\case The two outside loops are parallel.

\begin{caselist}
\case The two outside loops are attached in a region
  between edges of the multiple outer cycle connection of $L$.
  See diagram (a) in \eqref{zx1}.
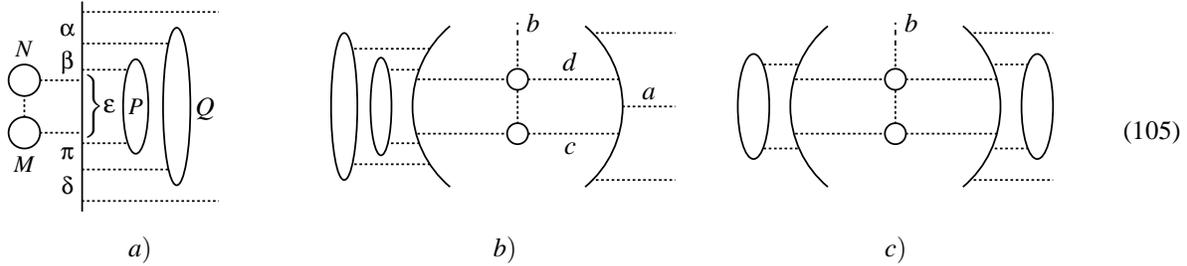
\begin{eqn}\label{zx1}
\begin{array}{c@{\qquad\ }c@{\qquad\ }c}
\diag{7mm}{5}{4}{
  \picline{1.4 0}{1.4 4}
  {\piclinedash{0.05 0.05}{0.01}
   \picline{0.3 1.5}{1.4 1.5}
   \picline{0.3 1.5}{0.3 2.5}
   \picline{0.3 2.5}{1.4 2.5}
   \picline{1.4 1.3}{2.4 1.3}
   \picline{1.4 2.7}{2.4 2.7}
   \picline{1.4 0.8}{3.2 0.8}
   \picline{1.4 3.2}{3.2 3.2}
   \picline{1.4 0.2}{4 0.2}
   \picline{1.4 3.8}{4 3.8}
  }
  \picfilledcircle{0.3 1.5}{0.3}{}
  \picfilledcircle{0.3 2.5}{0.3}{}
  \picfilledellipse{2.42 2}{0.24 0.9}{\small $P$}
  \picfilledellipse{3.2 2}{0.26 1.5}{}
  \picputtext[l]{1.37 2.26}{%{
    $\left.\ry{1.5em}\right\}\!\eps$
  }
  {\small
   \picputtext{1.12 3.47}{$\ap$}
   \picputtext{1.12 2.89}{$\bt$}
   \picputtext{1.12 1.11}{$\pi$}
   \picputtext{1.12 0.51}{$\dl$}
   \picputtext{0.3 3.1}{$N$}
   \picputtext{0.3 0.9}{$M$}
  }
  \picputtext{3.74 2}{$Q$}
}
& 
\diag{7mm}{6.5}{4}{
  \pictranslate{3.5 2}{
  {\piclinedash{0.05 0.05}{0.01}
   \picline{0 1.4}{3 1.4}
   \picline{0 -1.4}{3 -1.4}
   \picline{-2.6 -0.7}{0 -0.7}
   \picline{-2.6 0.7}{0 0.7}
   \picline{-3.3 -1.1}{0 -1.1}
   \picline{-3.3 1.1}{0 1.1}
   \picline{0 0}{3 0}
  } 
  \picfillnostroke{
   \piccircle{0 0}{2}{}
  }
  {\piclinedash{0.05 0.05}{0.01}
   \picline{2 15 polar}{2 180 15 - polar}
   \picline{2 -15 polar}{2 180 -15 - polar}
   \picline{0 1.1}{2 -15 polar x p 0 x}
   \picputtext{0.3 1.6}{$b$}
   \picputtext{1.0 -0.8}{$c$}
   \picputtext{1.0 0.86}{$d$}
   \picputtext{2.5 0.25}{$a$}
   {\piclinedash{0.1 0.1}{0.14}
    \picline{0 1.1}{0 1.6}
   }
  }
  \piccirclearc{0 0}{2}{-50 50}
  \piccirclearc{0 0}{2}{130 230}
  \picfilledellipse{-3.3 0}{0.25 1.4}{}
  \picfilledellipse{-2.6 0}{0.2 0.93}{}
  \picfilledcircle{2 15 polar x p 0 x}{0.2}{}
  \picfilledcircle{2 -15 polar x p 0 x}{0.2}{}
 }
} &
\diag{7mm}{6}{4}{
  \pictranslate{3 2}{
  {\piclinedash{0.05 0.05}{0.01}
   \picline{0 1.4}{3 1.4}
   \picline{0 -1.4}{3 -1.4}
   \picline{-2.7 -0.8}{2.7 -0.8}
   \picline{-2.7 0.8}{2.7 0.8}
  } 
  \picfillnostroke{
   \piccircle{0 0}{2}{}
  }
  {\piclinedash{0.05 0.05}{0.01}
   \picline{2 15 polar}{2 180 15 - polar}
   \picline{2 -15 polar}{2 180 -15 - polar}
   \picline{0 1.1}{2 -15 polar x p 0 x}
   \picputtext{0.3 1.6}{$b$}
   {\piclinedash{0.1 0.1}{0.14}
    \picline{0 1.1}{0 1.6}
   }
  }
  \piccirclearc{0 0}{2}{-50 50}
  \piccirclearc{0 0}{2}{130 230}
  \picfilledellipse{2.7 0}{0.3 1}{}
  \picfilledellipse{-2.7 0}{0.3 1}{}
  \picfilledcircle{2 15 polar x p 0 x}{0.2}{}
  \picfilledcircle{2 -15 polar x p 0 x}{0.2}{}
 }
}
\\
\ry{1.6em}a) & b) & c)
\end{array}
\end{eqn}
Clearly for two cycles in $IG(A)$, both $M,N$ must be intertwined
with $P,Q$, so they have legs on the segment $\eps$ of $L$. 
If none of $M,N$ had a leg on one of $\ap,\bt,\pi,\dl$, then
using that some of the legs of $P,Q$ must be multiple by
connectivity reasons, we see that $B(D)$ is of type A. Let
thus $\nu\in\{\ap,\bt,\pi,\dl\}$ be this segment touched. (The
following argument immediately clarifies that $\nu$ must be
unique.) Now, since $M,N$ must be intertwined also with the
(only multiple) outer cycle connection of $L$, there must
be a unique (because of $\gm=4$) third region $\mu\nin\{\ap,
\bt,\pi,\dl,\eps\}$ outside $L$ touched by a leg of $M,N$, and
because of $\chi(IG)=-1$, it is touched by legs of both $M,N$.
Then lemma \reference{lm1014} applies with $f=2$. 

\case The two outside loops are attached in a region between
an edge of the multiple outer cycle connection of $L$ and the
edge of the simple connection.  See diagram (b) in \eqref{zx1}.
Because of $\chi(IG)=-1$, all edges of $M,N$ are there, except
possibly $b$. By lemma \ref{lm1014}, there must be an edge $a$
at the indicated position. Then by $\gm=4$, we see $b=0$.
Now, by inadequacy we have $a=c=d=2$. Then we can apply
lemma \reference{UV} to exclude this case.
\end{caselist}
% here gamma4 -sheet is to enter

\case The two outside loops are in different regions outside $L$.
% bootom of gamma1-sheet
By $\gm=4$ and $\chi(IG(A))=-1$ we see that we have the diagram
c) in \eqref{zx1},
and all 6 edges of the triangle attachment must be there.
Then we are done by lemma \reference{lm1014}.
\end{caselist}
\end{caselist}
\qed

\subsubsection{Case C1\label{SSC1}}

We distinguish two subcases of case C, depending on
whether the triangle loops have legs that (jointly) touch
2 or 3 regions out side $L$. We consider first the case C1 of
2 regions $A,B$.

We have then a triangle attachment of the type \eqref{tat2}
with $f=0$. By lemma \reference{lm1014}, one of $a,b,c,d$
in \eqref{tat2} must be 0. The triangle attachment looks
then like diagram (a) in \eqref{zx3}. This implies in particular
that for a cycle in $IG(A)$ we need at least one more attachment
inside $L$, so at most one attachment is outside $L$.

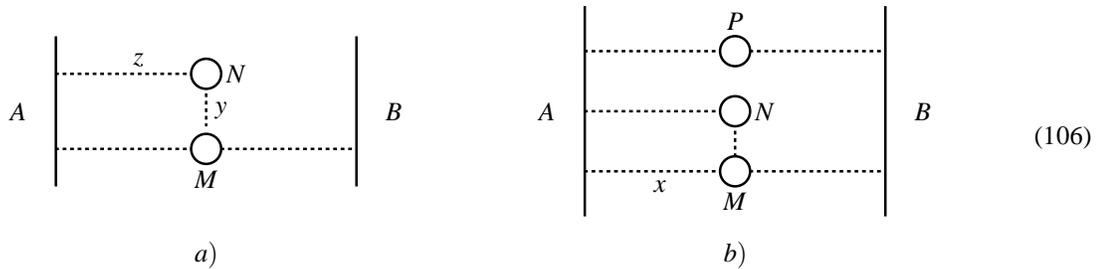
\begin{eqn}\label{zx3}
\begin{array}{c@{\kern3cm}c}
\diag{1cm}{4}{2.0}{
  \picline{0 0}{0 2.0}
  \picline{4 0}{4 2.0}
  {\piclinedash{0.05 0.05}{0.01}
   \picline{0 0.5}{4 0.5}
   \picline{2 0.8}{2 1.5}
   \picline{0 1.5}{2 1.5}
  }
  \picfilledcircle{2 1.5}{0.2}{}
  \picfilledcircle{2 0.5}{0.2}{}
  \picputtext{4.5 1.0}{$B$}
  \picputtext{-.5 1.0}{$A$}
  \picputtext{2 0.1}{$M$}
  \picputtext{2.4 1.5}{$N$}
  \picputtext{2.2 1.05}{$y$}
  \picputtext{1.1 1.7}{$z$}
} &
\diag{1cm}{4}{2.8}{
  \picline{0 0}{0 2.8}
  \picline{4 0}{4 2.8}
  {\piclinedash{0.05 0.05}{0.01}
   \picline{0 0.6}{4 0.6}
   \picline{0 2.2}{4 2.2}
   \picline{2 0.8}{2 1.4}
   \picline{0 1.4}{2 1.4}
  }
  \picfilledcircle{2 1.4}{0.2}{}
  \picfilledcircle{2 0.6}{0.2}{}
  \picfilledcircle{2 2.2}{0.2}{}
  \picputtext{4.5 1.4}{$B$}
  \picputtext{-.5 1.4}{$A$}
  \picputtext{2 0.2}{$M$}
  \picputtext{2.4 1.4}{$N$}
  \picputtext{2 2.6}{$P$}
  \picputtext{1 0.4}{$x$}
}
\\
\ry{1.6em}a) & b)
\end{array}
\end{eqn}

It will be useful to restrict first the valence-2 attachments
\eqref{attach}, touching upon $A,B$.

\begin{lemma}\label{opz}
In diagram (b) of \eqref{zx3},
there is no parallel loop to $P$, and all legs of $M,N,P$ are
simple, except possibly $x$.
\end{lemma}

\proof Otherwise $B(D)$ is of type A or B. \qed

Next we can eliminate parallel loops completely.

\begin{lemma}
There is no pair of parallel (valence-2) loops.
\end{lemma}

\proof Assume there were two such loops $P,Q$. By the previous
lemma, at least one of the two regions outside $L$ they touch
is different from $A,B$. We have, up to mutation, diagram a) in
\eqref{zx2}, with the option that $m=0$ if $P,Q$ touch one of $A,B$,
and $m\ne 0$ otherwise.
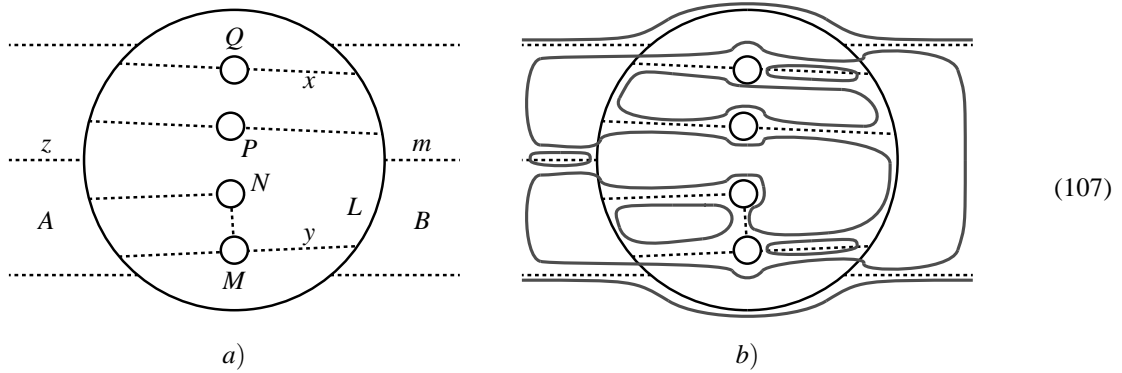
\begin{eqn}\label{zx2}
\begin{array}{c@{\qquad\ }c}
\diag{1cm}{6}{4}{
  \pictranslate{3 2}{
    {\piclinedash{0.05 0.05}{0.01}
     \picline{0 -3 x}{0 3 x}
    }
    \picfilledcircle{0 0}{2}{}
    {\piclinedash{0.05 0.05}{0.01}
     \picline{2 140 polar}{2 35 polar}
     \picline{2 -140 polar}{2 -35 polar}
     \picline{2 165 polar}{2 10 polar}
     \picline{2 -165 polar}{-0.05 -0.45}
     \picline{-0.05 -0.45}{-0.0 -1.2}
    }
    \picfilledcircle{-0.0 1.2}{0.18}{}
    \picfilledcircle{-0.05 0.45}{0.18}{}
    \picfilledcircle{-0.0 -1.2}{0.18}{}
    \picfilledcircle{-0.05 -0.45}{0.18}{}
    {\piclinedash{0.05 0.05}{0.01}
     \picscale{-1 1}{
     \picline{2 130 polar}{2 130 polar x p -3 x}
     \picline{2 230 polar}{2 230 polar x p -3 x}
     \picline{2 -50 polar}{2 -50 polar x p 3 x}
     \picline{2 50 polar}{2 50 polar x p 3 x}
     }
    }
    % \picfilledellipse{2.7 0.}{0.3 1.05}{}
    % \picputtext{-1.2 -0.2}{$z$}
    % \picputtext{-1.1 -0.7}{$x$}
    \picputtext{1.6 -0.62}{$L$}
    \picputtext{1.0 -0.99}{$y$}
    \picputtext{1.0 1.04}{$x$}
    % \picputtext{-0.2 -1.4}{$a$}
    % \picputtext{-0.35 0.79}{$M$}
    % \picputtext{-2.9 0.45}{$d$}
    % \picputtext{-2.9 -.5}{$c$}
     \picputtext{0 1.65}{$Q$}
     \picputtext{0.2 0.15}{$P$}
     \picputtext{0.35 -0.3}{$N$}
     \picputtext{0.0 -1.6}{$M$}
     \picputtext{-2.5 -0.8}{$A$}
     \picputtext{2.5 -0.8}{$B$}
     \picputtext{2.5 0.2}{$m$}
     \picputtext{-2.5 0.2}{$z$}
  }
} & 
\diag{1cm}{6}{4}{
  \pictranslate{3 2}{
    \picmultigraphics[S]{2}{-1 1}{
      \piccirclearc{0 0}{2}{-90 90}
    }
    {\piclinedash{0.05 0.05}{0.01}
     \picline{2 140 polar}{2 35 polar}
     \picline{2 -140 polar}{2 -35 polar}
     \picline{2 165 polar}{2 10 polar}
     \picline{2 -165 polar}{-0.05 -0.45}
     \picline{-0.05 -0.45}{-0.0 -1.2}
    }
    \picfilledcircle{-0.0 1.2}{0.18}{}
    \picfilledcircle{-0.05 0.45}{0.18}{}
    \picfilledcircle{-0.0 -1.2}{0.18}{}
    \picfilledcircle{-0.05 -0.45}{0.18}{}
    {\piclinedash{0.05 0.05}{0.01}
     \picscale{-1 1}{
     \picline{2 130 polar}{2 130 polar x p -3 x}
     \picline{2 230 polar}{2 230 polar x p -3 x}
     \picline{2 0 polar}{2 0 polar x p 3 x}
     \picline{2 -50 polar}{2 -50 polar x p 3 x}
     \picline{2 50 polar}{2 50 polar x p 3 x}
     }
    }
    % \picputtext{-0.2 -1.4}{$a$}
    % \picputtext{-0.35 0.79}{$M$}
    % \picputtext{-2.9 0.45}{$d$}
    % \picputtext{-2.9 -.5}{$c$}
    % \picputtext{0 1.65}{$Q$}
    % \picputtext{0.2 0.15}{$P$}
    % \picputtext{0.35 -0.3}{$N$}
    % \picputtext{0.0 -1.6}{$M$}
    { \piclinewidth{24}
      \picgraycol{0.3}
      \picPSgraphics{0 setlinecap}
      \curvepath{-1.0 -0.6}{-1.8  -0.6}{-1.6 -1.2}{-1.1 -1.1}
	{-0.1 -1.0}{-0.1 -0.6}{-0.8 -0.6}{}
      \picmultigraphics[S]{2}{1 -1}{
        \opencurvepath{-3.0 1.6}{-1.5 1.6}{-0.6 2.15}{0.6 2.15}
	  {1.5 1.6}{3.0 1.6}{}
	\opencurvepath{0 0.2}{-0.2 0.2}{-0.35 0.35}{-1.9 0.45}
	  {-2.0 0.1}{-2.5 0.2}{-3 0.2}{-3 1.3}{-2.5 1.4}
	  {-0.3 1.35}{-0.1 1.6}{0.1 1.6}
	  {0.3 1.35}{1.6 1.25}{1.5 1.5}{2.7 1.5}{2.9 0.9}{2.9 0}{}
        \pictranslate{3.26 -1.05}{
	  \picrotate{3}{
            \curvepath{-2.4 0.1}{-3 0.1}{-3 -0.1}
	      {-1.8 -0.1}{-1.8 0.1}{}
	  }
        }
      }
      \opencurvepath{0 0.2}{0.2 0.2}{0.25 0.35}{1.9 0.3}{1.9 -0.4}
	{1.9 -0.9}{0.25 -1.04}{0.2 -0.95}{0 -0.9}{0 -0.65}{0.2 -0.6}
	{0.25 -0.5}{0.18 -0.2}{0 -0.2}{}
      \curvepath{-1 1.15}{-1.5 1.2}{-1.9 0.45}{-0.3 0.5}{-0.25 0.75}
       {0.21 0.75}{0.25 0.48}{1.2 0.46}{1.9 0.44}{1.7 0.93}{0.8 0.95}
       {0.0 0.98}{-0.2 1}{-0.25 1.02}{-0.3 1.1}{}
      \pictranslate{0.11 0}{
        \curvepath{-2.6 0.1}{-3 0.1}{-3 -0.1}{-2.2 -0.1}{-2.2 0.1}{}
      }
    }
  }
} 
\\
\ry{1.9em}a) & b) 
\end{array}
\end{eqn}
Note that, by assumption, and in order to have 2 cycles in
$IG(A)$, all 5 other edges drawn outside $L$, different from $m$,
must be there. Thus again the regions $U,V$ of \eqref{oaz} remain
untouched. Now by inadequacy, $y=2$, and by connectivity,
w.l.o.g. up to mutation, $x\ge 2$. 

If $P,Q$ touch one of $A,B$, i.e. $m=0$, we see that again
we can apply lemma \reference{UV}. The diagram (b) in \eqref{zx2}
shows the $B$-state loops, taking account of $x\ge 2$ and $z=y=2$,
since $z$ is also inadequate if $m=0$.

In case $P,Q$ identify a pair of regions disjoint from $\{A,B\}$,
we have $m\ne 0$. Then we see that $B(D)$
has a pair of separating loops, each with at most one multiple
connection on one side. Thus $IG(B)$ cannot contain any cycle. \qed

\begin{lemma}\label{rew}
In diagram (a) in \eqref{zx3}, we can assume that $y=z=1$.
\end{lemma}

\proof With lemma \reference{l45}, it remains to rule
out the option $(y,z)=(1,2)$. (The case $y=z=2$ is out by
connectivity.) For this refine the
complexity of diagrams in case C, by saying that among
type C states $A(D)$ and $B(D)$, one is simpler, if
it has a triangle loop with only two edges connected to
it, but the triangle in the other state has no such loop.
With this set up, for $(y,z)=(1,2)$ in $A(D)$, go over
to $B(D)$. \qed

\begin{figure}[htb]
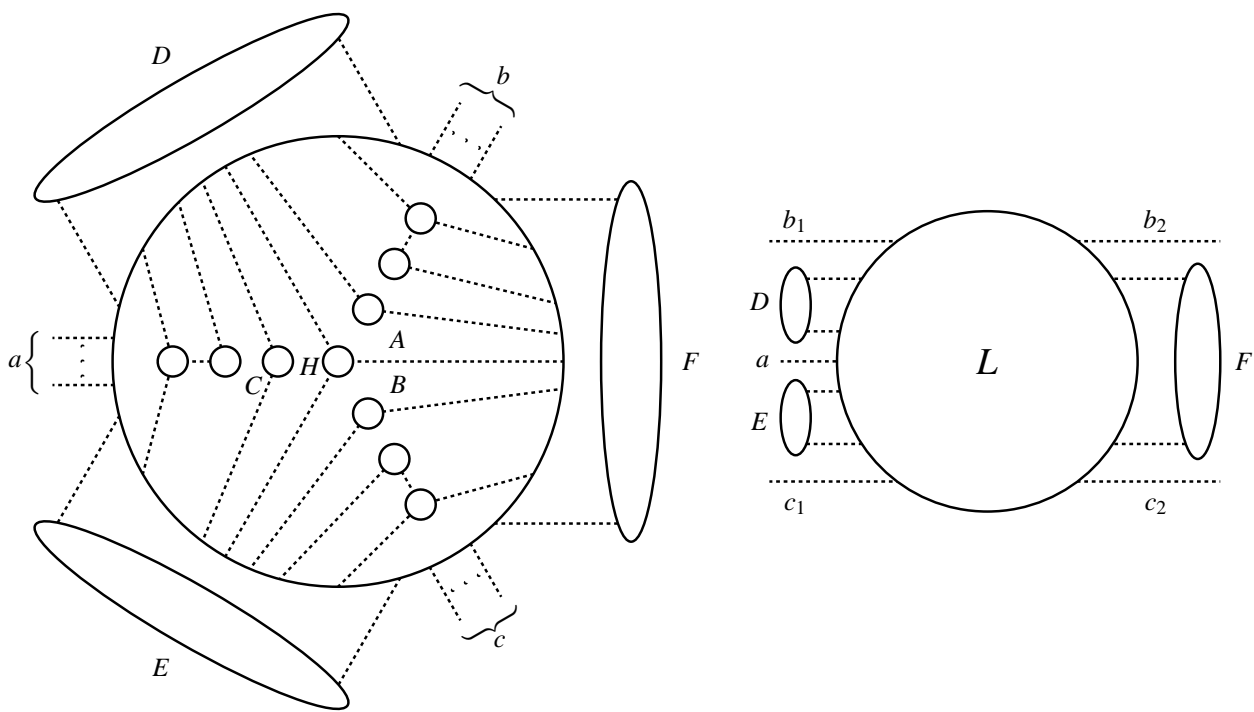

\begin{center}
\[
\diag{1cm}{8.1}{9.2}{
  \pictranslate{3.9 4.6}{
    \piccircle{0 0}{3}{}
    \picmultigraphics[rt]{3}{120}{
      \picscale{-1 1}{
       {\piclinedash{0.05 0.05}{0.01}
        \picline{3 30 polar}{0 2.2 x}
        \picline{3 -30 polar}{0 2.2 x}
        \picline{0 1.5 x}{0 2.2 x}
        \picline{0 1.5 x}{3 45 polar}
	\picline{0 0.8 x}{3 53 polar}
	\picline{0 0.8 x}{3 -53 polar}
       }
       \picfilledcircle{0 2.2 x}{0.2}{}
       \picfilledcircle{0 1.5 x}{0.2}{}
       \picfilledcircle{0 0.8 x}{0.2}{}
      }
      {\piclinedash{0.05 0.05}{0.01}
       \picline{0 3 x}{0 0}
       \picline{3 46 polar}{3 46 polar x p 3.9 x}
       \picline{3 -46 polar}{3 -46 polar x p 3.9 x}
       \picline{3 186 polar}{3 186 polar x p -3.8 x}
       \picline{3 -186 polar}{3 -186 polar x p -3.8 x}
       \picmultigraphics{3}{0 0.2}{
	 \piccircle{-3.4 -0.2}{0.01}{}
       }
       \picputtext[r]{-3.9 0.15}{$\left\{\ry{1.4em}\right.$}%}
      }
      \picfilledellipse{3.9 0}{0.4 2.4}{}
    }
    \picfilledcircle{0 0 x}{0.2}{}
    \picputtext{4.7 120 polar}{$D$}
    \picputtext{4.7 240 polar}{$E$}
    \picputtext{4.7   0 polar}{$F$}
    \picputtext{4.3 120 60 + polar}{$a$}
    \picputtext{4.3 240 60 + polar}{$c$}
    \picputtext{4.4   0 60 + polar}{$b$}
    \picputtext{0.8 0.3}{$A$}
    \picputtext{0.8 -0.3}{$B$}
    \picputtext{-0.37 -0.07}{$H$}
    \picputtext{-1.12 -0.31}{$C$}
  }
}\kern1.5cm
\diag{1cm}{6}{4}{
  {\piclinedash{0.05 0.05}{0.01}
   \picline{0 0.4}{6 0.4}
   \picline{0 3.6}{6 3.6}
   \picline{0.15 2}{2.9 2}
  }
  \pictranslate{0 1.25}{
   \picmultigraphics{2}{0 1.5}{
     {\piclinedash{0.05 0.05}{0.01}
      \picline{0.5 -0.35}{1.5 -0.35}
      \picline{0.5 0.35}{1.5 0.35}
     }
     \picfilledellipse{0.35 0.0}{0.2 0.5}{}
   }
  }
  \pictranslate{2.9 2}{
    {\piclinedash{0.05 0.05}{0.01}
     \picline{0 1.1}{2.8 1.1}
     \picline{0 -1.1}{2.8 -1.1}
    }
    \picfilledellipse{2.8 0}{0.3 1.3}{}
  }
  \pictranslate{2.9 2}{
   \picfilledcircle{0 0}{2}{\Large$L$}
   \picputtext[r]{-2.9 0.8}{$D$}
   \picputtext[r]{-2.9 -0.8}{$E$}
   \picputtext[r]{-2.9 -0.0}{$a$}
   \picputtext[r]{-2.4 -1.9}{$c_1$}
   \picputtext[r]{2.4 -1.9}{$c_2$}
   \picputtext[r]{-2.4 1.9}{$b_1$}
   \picputtext[r]{2.4 1.9}{$b_2$}
   \picputtext[l]{3.3 0}{$F$}
  }
} 
\]
\end{center}
\caption{\label{figC1}The pattern in case C1.}
\end{figure}

With all these preparations, we have now again, as in figure 
\reference{Fig19}, a pattern that unifies all remaining
possibilities for attachments, see figure \reference{figC1}.

The left part shows $L$ and its attachments. There are
now 3 triangles, so obviously our understanding is that
exactly one of them is there. Also, we argued in the beginning
of \S\reference{SSC1} that
among the outside attachments $D,E,F$ at most one is
there. If one is, we have to choose one, otherwise two, among the
3 valence-2 attachments $A,B,C$ and the valence-3 attachment
$H$ inside $L$.

The legs in $a,b,c$ go to the outer cycle loops neighbored to
$L$. This is shown on the right diagram. Again we can assume
that traces going from $L$ on the left resp. right connect $L$
to $L'$ resp. $L''$. Thus the traces that appear as $b$ on
the left may belong to two edges $b_{1,2}$, connecting $L$ to
$L'$ and $L''$. The same remark applies to $c$ and $c_{1,2}$.
The option that the regions outside $L$ into which $D,E,F$
are attached are between edges from $L$ to $L'$ (say) only,
is again recurred to the others by a flype.

We can assume again an outer cycle of length 5, with at most
one multiple edge between empty loops (lemma \reference{lC2}).

The multiplicities of edges conform to lemmas \ref{opz}, \ref{l45},
and \ref{rew}. Those of the outer cycle edges are as on the left of
part (a) of figure \reference{Fig19} (page \pageref{Fig19}), 
except that we change, in the notation of that figure, 
\[
\diag{1cm}{1.5}{1}{
  {\piclinedash{0.05 0.05}{0.01}
   \picline{0 0.5}{1.5 0.5}
  }
  \picputtext{0.75 0.73}{$b$}
}\quad\lra\quad
\diag{1cm}{2.5}{1}{
  {\piclinedash{0.05 0.05}{0.01}
   \picline{0 0.5}{2.5 0.5}
  }
  \picfilledcircle{1.25 0.5}{0.3}{}
  \picputtext{0.55 0.73}{$b$}
  \picputtext{2.05 0.73}{$e$}
}\,,
\]
and set $b=e=1$, w.l.o.g. up to mutation. (Otherwise still
$a,c,d\ge 1$ and $x,y,z\ge 0$, with the same further restrictions
as in \S\ref{SPAT}.)

\begin{clc}
% aqv23tst8a.C trdg_C1.in
Now we felt the diagrams are explicit enough to be checked
by computer. We also implemented a direct calculation of
$\chi(IG(A))$ during the successive attachments (rather than
calculating it from the diagram), in order to save the verification
of unnecessary cases. With this we could reduce the calculation
to 3 minutes, and it ruled out all possibilities.
\end{clc}

\subsubsection{Case C2}

The final case, C2, is when the triangle inside $L$
has legs touching on 3 regions outside $L$. The triangle
is then of the type \eqref{tat2}, with $f\ne 0$.

This case must again be dealt with by computer, but
first we simplify it.

\begin{lemma}\label{yuq}
There is at most one outside attachment to $L$.
\end{lemma}

\proof Assume we had two outside attached loops $R_{1,2}$.

\begin{caselist}
% See gm9=gm11
\case The two outside attachments are parallel. The two triangle
loops $M,N$ are the only loops inside $L$, and $f=2$ by inadequacy.
Then by lemma \ref{lm1014}
with $f=2$, one of $a,b,c,d$ must be 0. Still we need that both
$M,N$ are intertwined with at least 3 connections outside $L$.
Up to symmetry assume $R_{1,2}$ are attached in region $B$ of
\eqref{tat2}. This easily rules out the case that $b$, $c$ or
$d$ is 0, so assume w.l.o.g. $a=0$. We have then two situations,
shown by a) and b) of \eqref{zx5}:
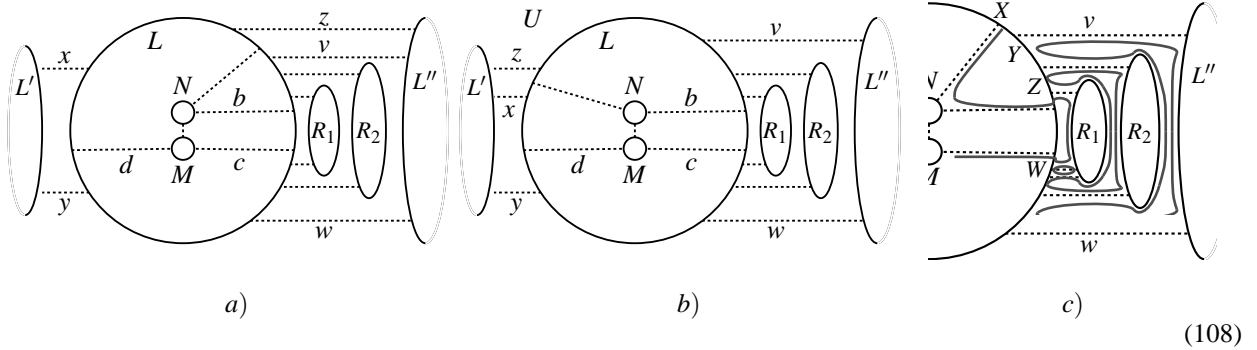
\begin{eqn}\label{zx5}
\begin{array}{c@{\quad}c@{\quad}c}
\diag{7.5mm}{7.5}{4.2}{
  \pictranslate{0.8 0}{
  {\piclinedash{0.05 0.05}{0.01}
   \picline{2.5 3.3}{6.3 3.3}
   \picline{2.5 3.8}{6.3 3.8}
   \picline{2.5 0.4}{6.3 0.4}
   \picline{-0.5 0.9}{2.5 0.9}
   \picline{-0.5 3.1}{2.5 3.1}
   \picline{2 2.6}{4.5 2.6}
   \picline{2 1.4}{4.5 1.4}
   \picline{2 3.0}{5.2 3.0}
   \picline{2 1.0}{5.2 1.0}
  }
  \picfilledellipse{4.5 2}{0.27 0.8}{\small$R_1$}
  \picfilledellipse{5.3 2}{0.3 1.2}{\small$R_2$}
  \picfilledellipse{6.3 2}{0.4 2.0}{}
  \picfilledellipse{-.8 2}{0.3 1.5}{}
  {\picgraycol{1}\picellipsearc{6.3 2}{0.4 2.0}{-90 90}
   \picellipsearc{-.8 2}{0.3 1.5}{90 -90}
  }
  \picfilledcircle{2 2}{2}{}
  \picputtext{1.5 3.6}{$L$} 
  \picputtext{6.3 2.9}{\small$L''$} 
  \picputtext{-.8 2.8}{\small$L'$}
  \picputtext{2 1.22}{$M$}
  \picputtext{2 2.8}{$N$}
  \picputtext{-0.1 3.3}{$x$}
  \picputtext{-0.1 0.7}{$y$}
  \picputtext{4.5 4.0}{$z$}
  \picputtext{4.5 3.5}{$v$}
  \picputtext{4.5 0.2}{$w$}
  \picputtext{3 2.6}{$b$}
  \picputtext{3 1.4}{$c$}
  \picputtext{1 1.4}{$d$}
  \pictranslate{2 2}{
   {\piclinedash{0.05 0.05}{0.01}
    \picline{2 190 polar}{0 -0.32}
    \picline{2 -10 polar}{0 -0.32}
    % \picline{2 170 polar}{0 0.32}
    \picline{2 10 polar}{0 0.32}
    \picline{2 47 polar}{0 0.32}
    \picline{0 -0.32}{0 0.32}
   }
   \picfilledcircle{0 -0.32}{0.2}{}
   \picfilledcircle{0 0.32}{0.2}{}
  }
  }
} & 
\diag{7.5mm}{7.5}{4.2}{
  \pictranslate{0.8 0}{
  {\piclinedash{0.05 0.05}{0.01}
   \picline{2.5 3.6}{6.3 3.6}
   \picline{2.5 2.6}{-0.5 2.6}
   \picline{2.5 0.4}{6.3 0.4}
   \picline{-0.5 0.9}{2.5 0.9}
   \picline{-0.5 3.1}{2.5 3.1}
   \picline{2 2.6}{4.5 2.6}
   \picline{2 1.4}{4.5 1.4}
   \picline{2 3.0}{5.2 3.0}
   \picline{2 1.0}{5.2 1.0}
  }
  \picfilledellipse{4.5 2}{0.27 0.8}{\small$R_1$}
  \picfilledellipse{5.3 2}{0.3 1.2}{\small$R_2$}
  \picfilledellipse{6.3 2}{0.4 2.0}{}
  \picfilledellipse{-.8 2}{0.3 1.5}{}
  {\picgraycol{1}\picellipsearc{6.3 2}{0.4 2.0}{-90 90}
   \picellipsearc{-.8 2}{0.3 1.5}{90 -90}
  }
  \picfilledcircle{2 2}{2}{}
  \picputtext{1.5 3.6}{$L$} 
  \picputtext{6.3 2.9}{\small$L''$} 
  \picputtext{-.8 2.8}{\small$L'$}
  \picputtext{2 1.22}{$M$}
  \picputtext{2 2.8}{$N$}
  \picputtext{-0.1 3.35}{$z$}
  \picputtext{-0.23 2.35}{$x$}
  \picputtext{-0.1 0.7}{$y$}
  \picputtext{4.5 3.8}{$v$}
  \picputtext{4.5 0.2}{$w$}
  \picputtext{0.2 4.0}{$U$}
  \picputtext{3 2.6}{$b$}
  \picputtext{3 1.4}{$c$}
  \picputtext{1 1.4}{$d$}
  \pictranslate{2 2}{
   {\piclinedash{0.05 0.05}{0.01}
    \picline{2 190 polar}{0 -0.32}
    \picline{2 -10 polar}{0 -0.32}
    % \picline{2 170 polar}{0 0.32}
    \picline{2 10 polar}{0 0.32}
    \picline{2 155 polar}{0 0.32}
    \picline{0 -0.32}{0 0.32}
   }
   \picfilledcircle{0 -0.32}{0.2}{}
   \picfilledcircle{0 0.32}{0.2}{}
  }
  }
} & 
\diag{8.5mm}{4.5}{4.2}{
  \pictranslate{-2 0}{
    \picclip{\picPSgraphics{2 0 m 2 4.3 l 6.5 4.3 l 6.5 0 l cp}}
    {
      {\piclinedash{0.05 0.05}{0.01}
       \picline{2.5 3.5}{6.3 3.5}
       \picline{2.5 0.4}{6.3 0.4}
       \picline{-0.5 0.9}{2.5 0.9}
       \picline{-0.5 3.15}{2.5 3.15}
       \picline{2 2.6}{4.5 2.6}
       \picline{2 1.4}{4.5 1.4}
       \picline{2 1.28}{4.5 1.28}
       \picline{2 3.0}{5.2 3.0}
       \picline{2 1.0}{5.2 1.0}
      }
      \picfilledellipse{4.5 2}{0.27 0.8}{\small$R_1$}
      \picfilledellipse{5.3 2}{0.3 1.2}{\small$R_2$}
      \picfilledellipse{6.3 2}{0.4 2.0}{}
      \picfilledellipse{-.8 2}{0.3 1.5}{}
      {\picgraycol{1}\picellipsearc{6.3 2}{0.4 2.0}{-90 90}
       \picellipsearc{-.8 2}{0.3 1.5}{90 -90}
      }
      \picfilledcircle{2 2}{2}{}
      \picputtext{1.5 3.6}{$L$} 
      \picputtext{6.3 2.9}{\small$L''$} 
      \picputtext{-.8 2.8}{\small$L'$}
      \picputtext{2 1.3}{$M$}
      \picputtext{2 2.8}{$N$}
      \picputtext{-0.1 3.3}{$x$}
      \picputtext{-0.1 0.7}{$y$}
      \picputtext{4.5 3.7}{$v$}
      \picputtext{4.5 0.2}{$w$}
      \pictranslate{2 2}{
       {\piclinedash{0.05 0.05}{0.01}
        \picline{2 190 polar}{0 -0.32}
        \picline{2 -10 polar}{0 -0.32}
        % \picline{2 170 polar}{0 0.32}
        \picline{2 10 polar}{0 0.32}
        \picline{2 57 polar}{0 0.32}
        \picline{0 -0.32}{0 0.32}
       }
       \picfilledcircle{0 -0.32}{0.2}{}
       \picfilledcircle{0 0.32}{0.2}{}
      }%
      \pictranslate{2 2}
      { \small
	\picputtext{1.35 1.2}{$Y$}
	\picputtext{1.15 1.9}{$X$}
	\picputtext{1.65 0.7}{$Z$}
	\picputtext{1.7 -0.6}{$W$}
	\piclinewidth{24}
        \picgraycol{0.3}
        \picPSgraphics{0 setlinecap}
	\picclip{\picPSgraphics{1.5 -1.3 m 1.5 2 l 4.5 2 l 4.5 -1.3 l cp}}{
          \picmultigraphics[S]{2}{1 -1}{
	    \opencurvepath{3.8 0}{3.8 1.0}{3.95 1.45}{1.5 1.45}{1.9 1.1}
	      {3.1 1.1}{3.3 1.35}{3.5 1.2}{3.7 0.5}{3.7 0}{}
	    \opencurvepath{2.92 0}{2.92 0.8}{3.09 0.96}{1.8 0.96}{1.9 0.7}
	      {2.3 0.7}{2.4 0.9}{2.6 0.9}{2.84 0.5}{2.84 0}{}
	  }
	}
	\picellipse{2.11 -0.61}{0.16 0.06}{}
	\opencurvepath{1.15 1.6}{0.4 0.7}{0.4 0.35}{2 0.35}{1.9 0.55}
	  {2.25 0.55}{2.15 0.2}{2.15 -0.3}{2.25 -0.5}{1.95 -0.5}
	  {2 -0.4}{0.4 -0.4}{}
      }
    }
  }
} 
\\
\ry{1.9em}a) & b) & c)
\end{array}
\end{eqn}
It is \em{a priori} to allow some of $v,w,x,y,z$ to be $0$.

Consider first case a) in \eqref{zx5}. We must have $x+z>0$
(otherwise $\gm=2$) and $y+w>0$ (otherwise $M$ intertwines only
two connections) and $v>0$ (otherwise $N$ intertwines only
two connections).

The multiplicities of the legs of $R_{1,2}$ are three times 1 and once
2. For this use lemma \ref{l45}; four simple legs are excluded by
connectivity, two times (2,1) because $B(D)$ becomes of type B.
So assume w.l.o.g. up to mutation that the lower leg of $R_1$
has multiplicity two.

Now consider a fragment of the $B$-state, as shown in part (c) of
\eqref{zx5}. (The lower part of the loop $Y$ is not drawn to remind
that $w$ may be not there.) It still shows that there exist two
cycles in $B(D)$, given by $XZY$ and $XZW$. Here we use that by
$x+z>0$ and $y+w>0$, we have $X\ne Y$. Thus $B(D)$ is of type B.

Consider then case b) in \eqref{zx5}. Then again $y+w>0$ (else $M$
intertwines only with 2) and now $z+v>0$ (else $N$ intertwines only
with 2); also $x>0$ (because else $\gm=2$), and so $x=2$ by inadequacy.
Then the same argument as for case a) applies. (The cycle in $B(D)$
formed by $X,Z,Y$ may contain an extra loop from the region $U$
in $A(D)$, if both $z$ and $v$ are non-zero.)

\case The two outside attachments $R_{1,2}$ are not parallel.
Then, by lemma \reference{lm2l}, they do not mutually enclose.

We use diagram \eqref{tat2}. Since the triangle is the only
attachment inside $L$, by lemma \reference{lm1014} for $f=2$,
we must have some of $a,b,c,d$ vanishing.

If now $c$ or $d$ is 0, or both of $a,b$ are 0, then one of 
the triangle loops $M,N$ in intertwined with at most one connection
outside $L$. Then we cannot have 2 cycles in $IG(A)$. So w.l.o.g.,
$a=0$, and then $b,c,d$ are non-zero. Then $c$ is inadequate, so
$c=2$ (and also $f=2$). 

Now the triangle loops $M,N$ have each two legs, but must intertwine
with at least 3 connections outside $L$. This easily implies
that none of the regions $A,B,C$ in \eqref{lDD} touched by legs
of $M,N$ can be any of the regions $U,V$ in \eqref{oaz}. For
the same reason, at least one of $R_{1,2}$, say $R_1$, must be
attached outside $L$ in a region different from $U,V$, i.e.
between edges $x,x'$ that connect $L$ to the same loop $L'$
on the outer cycle. 

Then we have at least 7 loops in $B(D)$: a separating loop $X$, 3
loops inside $X$ (two of which from $c=f=2$), and 3 outside $X$
(those along $U,V$, and one formed by $R_1$, its legs, $x$, $x'$
and arcs of $L,L'$). So with lemma \reference{UV} we are done.
\end{caselist}
\qed

\begin{lemma}\label{nopa}
There are no parallel loops.
\end{lemma}

\proof Parallel loops outside $L$ were ruled out in lemma
\reference{yuq}. For parallel loops inside $L$, note that
one could obtain \eqref{lDD} by deleting edges and valence-2
loops like
$
\diag{4mm}{2.5}{1}{
  {\piclinedash{0.15 0.15}{0.05}
   \picline{0 0.5}{2.5 0.5}
  }
  \picfilledcircle{1.25 0.5}{0.3}{}
}\,\lra\,
\diag{4mm}{1.5}{1}{
  {\piclinedash{0.15 0.15}{0.05}
   \picline{0 0.5}{1.5 0.5}
  }
}
$. Then the argument in the proof of lemma \reference{lempar}
modifies to show that $B(D)$ is of type A or B. \qed

With this we can already specify a common superconfiguration
of $A(D)$. Before we do so, though, let us make one final
simplification.

\begin{lemma}
In the attachment \eqref{tat2}, we have $e=1$.
\end{lemma}

\proof In \eqref{tat2}, assume $e\ge 2$. Since this then
gives an isolated vertex in $IG(A)$, we need at least 3 cycles therein.
Moreover, by lemma \reference{yuq}, we know that we must do an
attachment \eqref{attach} inside $L$, call it $P$, which must be then
of valence 2. It can be specified by choosing two of the three
regions $A,B,C$ touched by legs of the triangle loops $M,N$. 

If all of $a,b,c,d$ are non-zero in \eqref{tat2},
we use lemma \reference{lm1014}. % it is easy to see that
% any of the three choices of $P$ makes $B(D)$ to be of type A or B.
% (If we add a second attachment \eqref{attach} inside $L$, the
% situation does not change.) 

If now $c$ or $d$ is 0, or both of $a,b$ are 0, then one of 
the triangle loops $M,N$ in intertwined with at most one connection
outside $L$. Then we cannot have 3 cycles in $IG(A)$. So w.l.o.g.,
$a=0$, and then $b,c,d$ are non-zero. 

When $a=0$, among the three choices for $P$, only the one
identifying $A$ and $C$ avoids that $B(D)$ becomes of type A or B.
This, together with lemma \reference{nopa}, shows in particular that
$P$ must be the only attachment \eqref{attach} inside $L$, so that
the other one is outside $L$.

Next, by lemma \reference{UV} we can assume that one of $A,B,C$
coincides with one of the regions $U,V$ in \eqref{oaz}, say w.l.o.g $U$.
(Again $B(D)$ has only one separating loop, $X$, we need at least 5
others to have 2 cycles in $IG(B)$, and there is an extra loop
inside $X$, between the two edges of $e$, which does not
intertwine with anything. Thus again we have at least 7 loops.)
Then we have the following 3 choices:
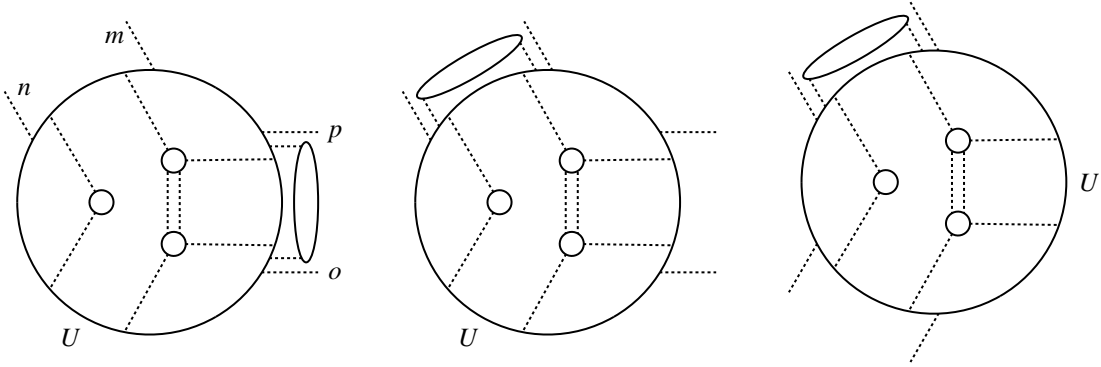
\[
\diag{8mm}{5.7}{5.3}{
  \pictranslate{2.5 2.5}{
    \picmultigraphics[rt]{3}{120}{
      {\piclinedash{0.05 0.05}{0.01}
       \picline{0.8 60 polar}{2.2 19 polar}
       \picline{0.8 -60 polar}{2.2 -19 polar}
      }
      \piccirclearc{0 0}{2.2}{-60 60}
    }
    {\piclinedash{0.05 0.05}{0.01}
     \picline{0.8 60 polar x 0.1 + x}{0.8 -60 polar x 0.1 + x}
     \picline{0.8 60 polar x 0.1 - x}{0.8 -60 polar x 0.1 - x}
    }
    \picmultigraphics[rt]{3}{120}{
      \picfilledcircle{0.8 60 polar}{0.2}{}
    }
    \picmultigraphics[rt]{2}{120}{
      {\piclinedash{0.05 0.05}{0.01}
       \picline{2.2 32 polar}{2.2 32 polar x p 2.8 x}
       \picline{2.2 -32 polar}{2.2 -32 polar x p 2.8 x}
      }
    }
    \picputtext{2.6 -120 polar}{$U$}
    \picputtext{2.8 102 polar}{$m$}
    \picputtext{2.8 138 polar}{$n$}
    \picputtext{3.3 -21 polar}{$o$}
    \picputtext{3.3 21 polar}{$p$}
    {\piclinedash{0.05 0.05}{0.01}
     \picline{2.2 25 polar}{2.2 25 polar x p 2.6 x}
     \picline{2.2 -25 polar}{2.2 -25 polar x p 2.6 x}
    }
    \picfilledellipse{2.6 0}{0.2 1}{}
  }
}\qquad
\diag{8mm}{5.5}{5.3}{
  \pictranslate{2.5 2.5}{
    \picmultigraphics[rt]{3}{120}{
      {\piclinedash{0.05 0.05}{0.01}
       \picline{0.8 60 polar}{2.2 19 polar}
       \picline{0.8 -60 polar}{2.2 -19 polar}
      }
      \piccirclearc{0 0}{2.2}{-60 60}
    }
    {\piclinedash{0.05 0.05}{0.01}
     \picline{0.8 60 polar x 0.1 + x}{0.8 -60 polar x 0.1 + x}
     \picline{0.8 60 polar x 0.1 - x}{0.8 -60 polar x 0.1 - x}
    }
    \picmultigraphics[rt]{3}{120}{
      \picfilledcircle{0.8 60 polar}{0.2}{}
    }
    \picmultigraphics[rt]{2}{120}{
      {\piclinedash{0.05 0.05}{0.01}
       \picline{2.2 32 polar}{2.2 32 polar x p 2.8 x}
       \picline{2.2 -32 polar}{2.2 -32 polar x p 2.8 x}
      }
    }
    \picputtext{2.6 -120 polar}{$U$}
    \picrotate{120}{
      {\piclinedash{0.05 0.05}{0.01}
       \picline{2.2 25 polar}{2.2 25 polar x p 2.6 x}
       \picline{2.2 -25 polar}{2.2 -25 polar x p 2.6 x}
      }
      \picfilledellipse{2.6 0}{0.2 1}{}
    }
  }
}\qquad
\diag{8mm}{5.0}{5.8}{
  \pictranslate{2.5 3.1}{
    \picmultigraphics[rt]{3}{120}{
      {\piclinedash{0.05 0.05}{0.01}
       \picline{0.8 60 polar}{2.2 19 polar}
       \picline{0.8 -60 polar}{2.2 -19 polar}
      }
      \piccirclearc{0 0}{2.2}{-60 60}
    }
    {\piclinedash{0.05 0.05}{0.01}
     \picline{0.8 60 polar x 0.1 + x}{0.8 -60 polar x 0.1 + x}
     \picline{0.8 60 polar x 0.1 - x}{0.8 -60 polar x 0.1 - x}
    }
    \picmultigraphics[rt]{3}{120}{
      \picfilledcircle{0.8 60 polar}{0.2}{}
    }
    \picrotate{120}{
      \picmultigraphics[rt]{2}{120}{
        {\piclinedash{0.05 0.05}{0.01}
         \picline{2.2 32 polar}{2.2 32 polar x p 2.8 x}
         \picline{2.2 -32 polar}{2.2 -32 polar x p 2.8 x}
        }
      }
    }
    \picputtext{2.6 0 polar}{$U$}
    \picrotate{120}{
      {\piclinedash{0.05 0.05}{0.01}
       \picline{2.2 25 polar}{2.2 25 polar x p 2.6 x}
       \picline{2.2 -25 polar}{2.2 -25 polar x p 2.6 x}
      }
      \picfilledellipse{2.6 0}{0.2 1}{}
    }
  }
}
\]
(Here we can assume that $o,p$ connect $L$ to $L''$, and $m,n$
connect $L$ to $L'$.)
In all cases, though, there is a loop inside $L$ intertwined with
at most one connection outside $L$, and so we cannot have $3$ cycles
in $IG(A(D))$.
\qed

\begin{figure}[htb]
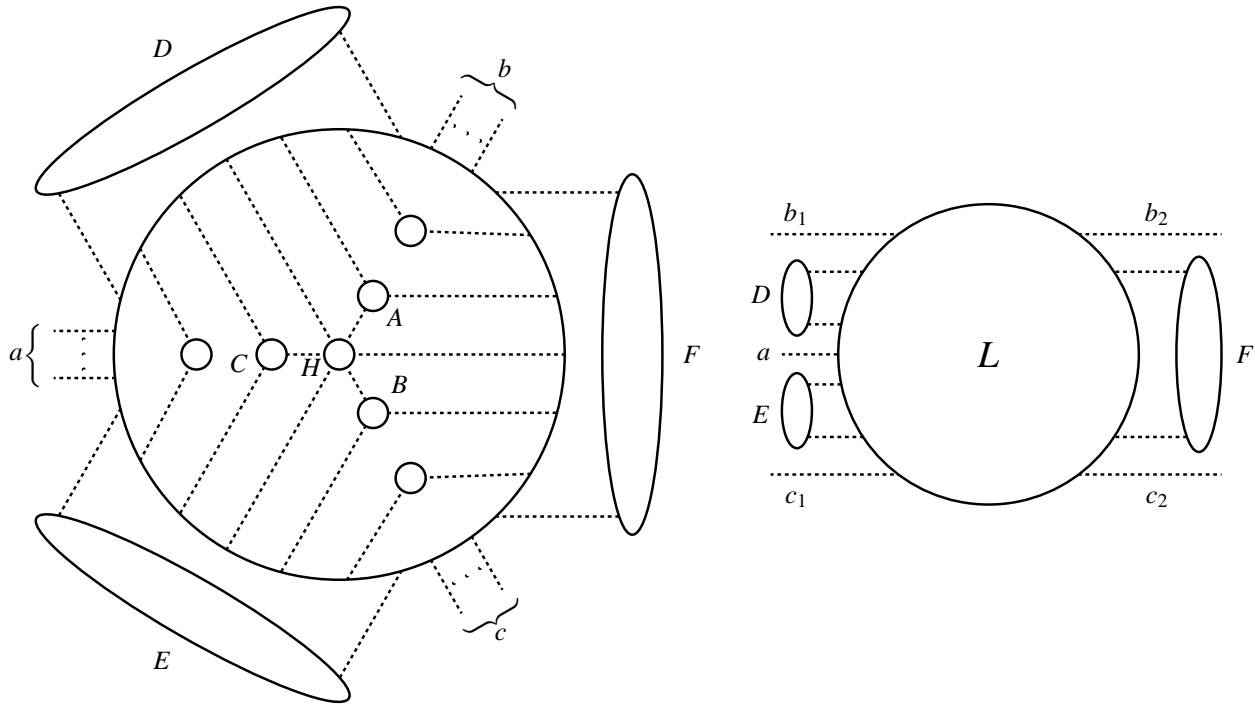

\begin{center}
\[
\diag{1cm}{8.1}{9.2}{
  \pictranslate{3.9 4.6}{
    \piccircle{0 0}{3}{}
    \picmultigraphics[rt]{3}{120}{
      \picscale{-1 1}{
       {\piclinedash{0.05 0.05}{0.01}
        \picline{3 28 polar}{0 1.9 x}
        \picline{3 -28 polar}{0 1.9 x}
        \picline{0 0.9 x}{0 0.0 x}
	\picline{0 0.9 x}{3 45 polar}
	\picline{0 0.9 x}{3 -45 polar}
       }
       \picfilledcircle{0 1.9 x}{0.2}{}
       \picfilledcircle{0 0.9 x}{0.2}{}
      }
      {\piclinedash{0.05 0.05}{0.01}
       \picline{0 3 x}{0 0}
       \picline{3 46 polar}{3 46 polar x p 3.9 x}
       \picline{3 -46 polar}{3 -46 polar x p 3.9 x}
       \picline{3 186 polar}{3 186 polar x p -3.8 x}
       \picline{3 -186 polar}{3 -186 polar x p -3.8 x}
       \picmultigraphics{3}{0 0.2}{
	 \piccircle{-3.4 -0.2}{0.01}{}
       }
       \picputtext[r]{-3.9 0.15}{$\left\{\ry{1.4em}\right.$}%}
      }
      \picfilledellipse{3.9 0}{0.4 2.4}{}
    }
    \picfilledcircle{0 0 x}{0.2}{}
    \picputtext{4.7 120 polar}{$D$}
    \picputtext{4.7 240 polar}{$E$}
    \picputtext{4.7   0 polar}{$F$}
    \picputtext{4.3 120 60 + polar}{$a$}
    \picputtext{4.3 240 60 + polar}{$c$}
    \picputtext{4.4   0 60 + polar}{$b$}
    \picputtext{0.75 0.5}{$A$}
    \picputtext{0.8 -0.4}{$B$}
    \picputtext{-0.37 -0.17}{$H$}
    \picputtext{-1.33 -0.11}{$C$}
  }
}
\kern1.5cm
\diag{1cm}{6}{4}{
  {\piclinedash{0.05 0.05}{0.01}
   \picline{0 0.4}{6 0.4}
   \picline{0 3.6}{6 3.6}
   \picline{0.15 2}{2.9 2}
  }
  \pictranslate{0 1.25}{
   \picmultigraphics{2}{0 1.5}{
     {\piclinedash{0.05 0.05}{0.01}
      \picline{0.5 -0.35}{1.5 -0.35}
      \picline{0.5 0.35}{1.5 0.35}
     }
     \picfilledellipse{0.35 0.0}{0.2 0.5}{}
   }
  }
  \pictranslate{2.9 2}{
    {\piclinedash{0.05 0.05}{0.01}
     \picline{0 1.1}{2.8 1.1}
     \picline{0 -1.1}{2.8 -1.1}
    }
    \picfilledellipse{2.8 0}{0.3 1.3}{}
  }
  \pictranslate{2.9 2}{
   \picfilledcircle{0 0}{2}{\Large$L$}
   \picputtext[r]{-2.9 0.8}{$D$}
   \picputtext[r]{-2.9 -0.8}{$E$}
   \picputtext[r]{-2.9 -0.0}{$a$}
   \picputtext[r]{-2.4 -1.9}{$c_1$}
   \picputtext[r]{2.4 -1.9}{$c_2$}
   \picputtext[r]{-2.4 1.9}{$b_1$}
   \picputtext[r]{2.4 1.9}{$b_2$}
   \picputtext[l]{3.3 0}{$F$}
  }
} 
\]
\end{center}
\caption{\label{figC2}The pattern in case C2.}
\end{figure}

The promised common superconfiguration of all options we have 
for $A(D)$ is obtained now by replacing in figure \reference{figC1}
the left diagram by the one on the left of figure \reference{figC2}.

The reading of figure \reference{figC2} is as follows. Let
$e,e',e''$ resp. be the edges between $H$ and $A,B,C$ resp.
Then we demand that exactly one of $e,e',e''$ is there,
and then it is $1$. It determines the orientation of the
attachment in \eqref{tat2}. For example, if $e=1$, then
$B,C$ and the traces touching upon them are not there.

With the choice between $e,e',e''$, which gives the edge $e$ of
\eqref{lDD}, also the edges $a,c,b,d,f$ are determined, and the
following conditions apply: 
\begin{itemize}
\item $f\ne 0$ (else case C1), but one of $a,b,c,d$ is $0$ (by
  lemma \ref{lm1014}),
\item $b+d>0$, $a+c>0$ (else case C1), and $c+d>0$ (else
  $e$ is isolated and $D$ composite),
\item if $a=b=0$, then $f=1$; if one of $c$ or $d$ is zero, the
  other is 1 (otherwise the connection between $L$ and one of $M,N$
  is an isolated vertex in $IG(A)$, and one cannot create 3 cycles
  therein with the other attachments), and
\item not all of $a,b,c,d,f$ are even (by connectivity).
\end{itemize}

Moreover, at most one of $D,E,F$ is attached (lemma \reference{yuq}).
If an outside attachment is there, we have 1, otherwise 2, of the
three valence-2 loops inside $L$.

The same remark as in case C1 on the outer cycle edges
and multiplicities apply.

\begin{clc}
% aqv23tst8a.C trdg_C2.in
At this stage we applied a computer to check the pattern. For
speed-up, again we kept track of $\chi(IG(A))$ during the successive
attachments, rather than calculating it from the (completely)
constructed diagram. With all theoretical input and optimization,
the calculation time could be reduced to 18 minutes. 
% Still this was the most difficult case of the calculation.
%and as mentioned in \S\ref{Scmp}, yields a non-empty result.

After adjusting the multiplicity of the twist edge to 3, so
that \eqref{ABl} holds, we obtain a list of 56 diagrams $\hD_0$,
from which we construct 56 sequences of diagrams $\hD_k$.
However, it turns out that only 7 sequences are distinct up
to flypes. (Again diagrams appear repeatedly, since it is too
hard to remove all symmetries in advance.) These indeed then
represent distinct sequences of knots, of which 3 consist
of amphicheiral ones. The knots for $k=1$ are shown in
figure \ref{Fig20'} as series $\hD^{(1)}$ to $\hD^{(7)}$.
\end{clc}

\subsection{Concluding part: writhe and Vassiliev invariant
  test\label{SSI}}

% MUST BE REDRAWN NICELY
\def\rrs#1#2{\begin{array}{c}#1\\ \ry{1.5em}\hD^{(#2)}_1\end{array}}
\def\rs#1#2{\unitxsize#1\relax\epsffile{t1-#2.eps}}
\def\rS{\rs{3.5cm}}
\begin{figure}
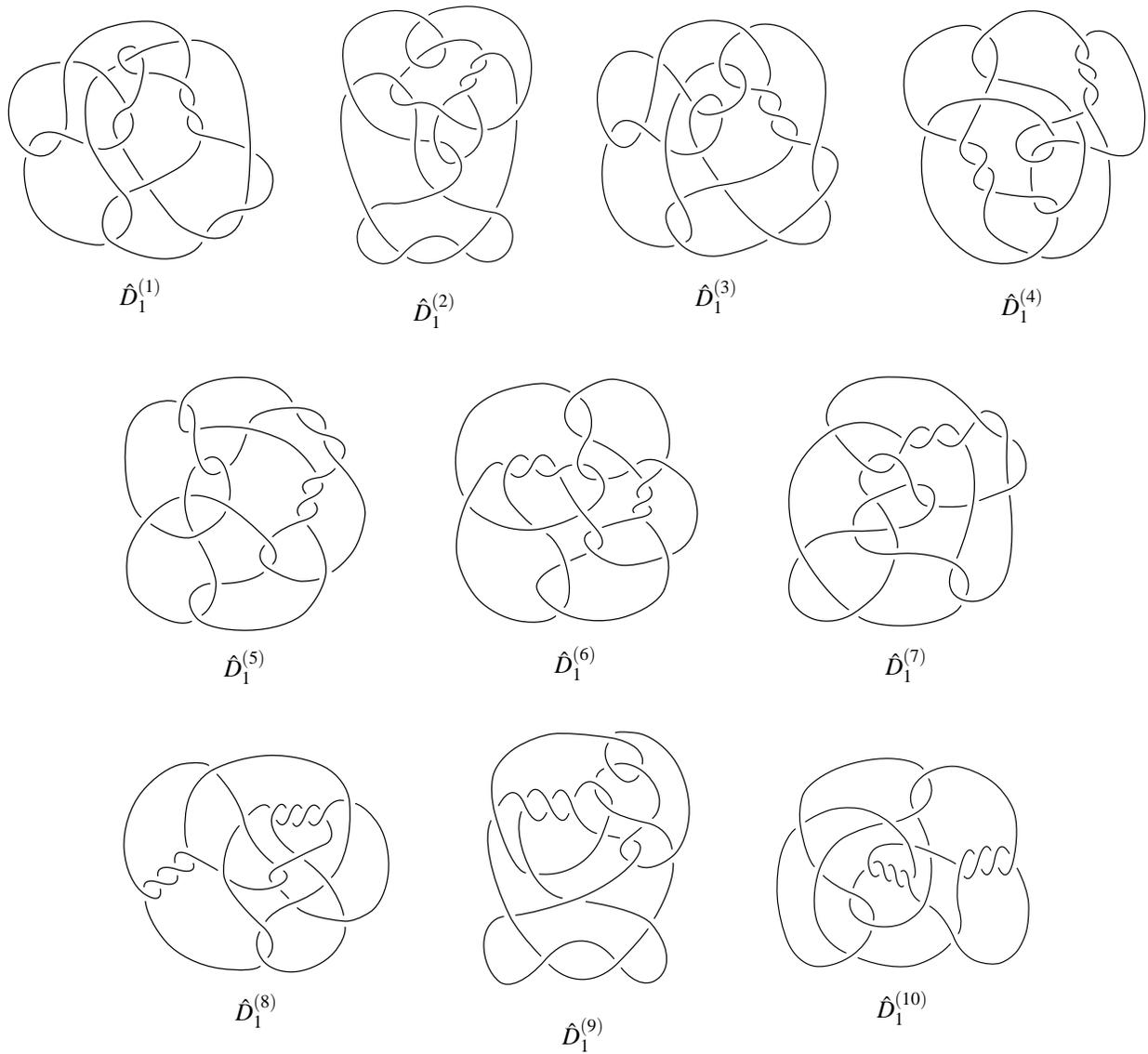
%[htb]
\begin{center}
\begin{eqnarray*}%{*4c}
& \rrs{\rS{17+4k_ach1}}{1} \quad\es \rrs{\rS{17+4k_ach2}}{2} \quad\es
  \rrs{\rS{17+4k_ach3}}{3} \quad\es \rrs{\rS{17+4k_ach4}}{4} & \\ [5mm]
& \rrs{\rS{17+4k_ach5}}{5} \qquad\es \rrs{\rs{3.5cm}{17+4k_ach6}}{6}
  \qquad\es \rrs{\rS{17+4k_ach7}}{7} &  \\[5mm]
& \rrs{\rS{17+4k_ach8}}{8} \qquad\es \rrs{\rs{3.8cm}{17+4k_ach9}}{9} 
  \qquad\es \rrs{\rs{3.3cm}{17+4k_ach10}}{10} &
\end{eqnarray*}
\end{center}
\caption{\label{Fig20'}The diagrams $\hD^{(i)}_k$ for $k=1$. Here
the parenthesized superscript refers to the number $i$ of the series
(and is omitted if $i$ is fixed). The subscript stands for the
number of the diagram in the series, so that $\hD_k$ has $16+4k$
crossings.}
\end{figure}

To sum up the proof so far, after check of semiadequacy
invariants, we are left with only 10 sequences $\{\hD^{(i)}_k\}_{k
\in\bN}$ of diagrams that could match a hypothetical $16+4k$-crossing
diagram $D_k$ of our amphicheiral knots $K_k$.

Figure \reference{Fig20'} shows the diagrams of these sequences for
$k=1$. The diagrams of higher $k$ are obtained by adding a full twist
(2 crossings) at each twist class of 3 or 4 crossings.

In order to distinguish $D_k$ from $\hD^{(i)}_k$, we use first the
writhe. The diagrams in the series for $i=1,2,8,9$ all have non-zero
writhe, and are excluded.

The the other 6 series, we use again Vassiliev invariants
to distinguish $\hD_k$ from our hypothetic diagrams $D_k$.
Considering $v_2$ and $v_3$ of \eqref{v23}, we have 5
distinct sequences of value pairs, up to negating $v_3$
(which is the result of taking mirror images). The values
for $k=1,\dots,4$ are given in the following table, in
comparison to those for $D_k$ in the last row:
\begin{eqn}\label{v23t}
\begin{array}{c||r|r|r|r||r|r|r|r||c}
k  &   1 &   2 &   3 &   4 &   1 &   2 &   3 &   4 & \\[1mm]
\hline\ry{1.3em}
v_2&  -4 &  -7 & -12 & -19 &   0 &   0 &   0 &   0 &v_3 \\
   &  -5 &  -7 &  -9 & -11 &   0 &   0 &   0 &   0 & \\
   &  -5 &  -9 & -13 & -17 &   1 &   2 &   3 &   4 & \\
   &  -6 & -13 & -22 & -33 &   0 &   0 &   0 &   0 & \\
   &  -8 & -15 & -24 & -35 &   0 &   0 &   0 &   0 & \\[1mm]
\hline\ry{1.3em}
D_k&  -7 & -11 & -15 & -19 &  0 &   0 &   0 &   0 & \\
\end{array}
\end{eqn}
By determining the braiding polynomial of $v_3$, one
easily finds that for the sequence with $v_3(\hD_k)\ne 0$
when $k\le 4$, the same holds also if $k\ge 5$. 
% (For two of them one can use also the non-zero writhe of $\hD_k$.)

A similar check of the braiding polynomial of $v_2$ of the
remaining 4 sequences $\hD_k$ shows that the only coincidence
with $v_2(D_k)$ occurs for $k=4$ in the first sequence in
\eqref{v23t}, where $v_2=-19$. In this case, one can use
$v_4$ of \cite{pv4}. It can be again evaluated on $D_4$ and
$\hD_4$ either directly, or (since these diagrams have 32
crossings) from its braiding polynomial obtained by examining
simpler diagrams for smaller $k$. We find that whenever
$v_2(\hD_4)=-19$, we have $v_4(\hD_4)=432$. On the other hand,
$v_4(D_4)=346$.

This completes the proof of theorem \ref{th10.1}, and then also
of theorem \ref{thm}. \qed

\noindent{\bf Acknowledgement.} I would wish to thank to
Sebastian Baader, Mikami Hirasawa, Efstratia Kalfagianni, Akio
Kawauchi, Hitoshi Murakami, Takuji Nakamura and Dan Silver for
stimulating me with their remarks into some of the investigations
described here. I also wish to thank to my JSPS host Prof. T.~Kohno
at U.~Tokyo for his support.

\end{document}

%% file: myeqn.tex
% typing `equation' breaks my fingers!
\newenvironment{eqn}{\begin{equation}}{\end{equation}\@ignoretrue}

% eqlabel=(#1)
\newenvironment{myeqn*}[1]{\begingroup\def\@eqnnum{\reset@font\rm#1}%
\xdef\@tempk{\arabic{equation}}\begin{equation}\edef\@currentlabel{#1}}
{\end{equation}\endgroup\setcounter{equation}{\@tempk}\ignorespaces}

% eqlabel=#1
\newenvironment{myeqn}[1]{\begingroup\let\eq@num\@eqnnum
\def\@eqnnum{\bgroup\let\r@fn\normalcolor % an extremely UGLY hack !!!
\def\normalcolor####1(####2){\r@fn####1#1}%
%\show\reset@font
\eq@num\egroup}%
\xdef\@tempk{\arabic{equation}}\begin{equation}\edef\@currentlabel{#1}}
{\end{equation}\endgroup\setcounter{equation}{\@tempk}\ignorespaces}

% eqlabel=(eqnnr) \qed
\newenvironment{myeqn**}{\begin{myeqn}{(%\arabic{equation}
%\show\theequation
\theequation)\es\es\mbox{\qed}}\edef\@currentlabel{\theequation}}
{\end{myeqn}\stepcounter{equation}}